\documentclass[12pt]{amsart}
\usepackage{latexsym, amssymb, amsfonts, amscd, multicol}

\theoremstyle{plain}
\newtheorem{Thm}{Theorem}[subsection]
\newtheorem{Cor}[Thm]{Corollary}

\newtheorem{Prop}[Thm]{Proposition}
\newtheorem{Lem}[Thm]{Lemma}
\newtheorem{Cl}[Thm]{Claim}
\newtheorem{Thm'}{Theorem}[section]

\theoremstyle{definition}
\newtheorem{Def}[Thm]{Definition}
\newtheorem{Rem}[Thm]{Remark}

\newtheorem{Emp}[Thm]{}

\newtheorem{Ex}[Thm]{Example}

\newtheorem{Con}[Thm]{Construction}
\newtheorem{Not}[Thm]{Notation}

\newtheorem{Not'}[Thm']{Notation}
\newtheorem{Emp'}[Thm']{}

\numberwithin{equation}{subsection}
\newcommand{\qlbar}{\overline{\B{Q}_l}}
\newcommand{\fqbar}{\overline{\fq}}
\newcommand{\om}{\omega}
\newcommand{\La}{\Lambda}

\newcommand{\ov}{\overline}

\newcommand{\fq}{\B{F}_q}

\newcommand{\B}[1]{\mathbb#1}
\newcommand{\cal}[1]{\mathcal{#1}}
\newcommand{\C}[1]{\cal#1}

\newcommand{\sr}{\operatorname{sr}}
\newcommand{\ssc}{\operatorname{sc}}
\newcommand{\der}{\operatorname{der}}
\newcommand{\nr}{\operatorname{nr}}
\newcommand{\GL}{\operatorname{GL}}
\newcommand{\SL}{\operatorname{SL}}
\newcommand{\tu}{\operatorname{tu}}
\newcommand{\tn}{\operatorname{tn}}

\newcommand{\isom}{\overset {\thicksim}{\to}}

\newcommand{\Om}{\Omega}

\newcommand{\si}{\sigma}

\newcommand{\lra}{\longrightarrow}
\newcommand{\lla}{\longleftarrow}
\newcommand{\hra}{\hookrightarrow}
\newcommand{\wt}{\widetilde}
\newcommand{\wh}{\widehat}
\newcommand{\Gm}{\Gamma}

\newcommand{\gm}{\gamma}
\newcommand{\ka}{\kappa}
\newcommand{\dt}{\delta}
\newcommand{\Dt}{\Delta}
\newcommand{\bs}{\backslash}

\newcommand{\m}{^{\times}}

\newcommand{\un}{\underline}

\newcommand{\tor}{\operatorname{tor}}

\newcommand{\al}{\alpha}

\newcommand{\la}{\lambda}

\newcommand{\rss}[1]{Subsection \ref{SS:#1}}
\newcommand{\rl}[1]{Lemma \ref{L:#1}}
\newcommand{\rn}[1]{Notation \ref{N:#1}}
\newcommand{\rcl}[1]{Claim \ref{C:#1}}
\newcommand{\rp}[1]{Proposition \ref{P:#1}}
\newcommand{\pg}[1]{\pageref{#1}}
\newcommand{\rr}[1]{Remark \ref{R:#1}}
\newcommand{\rc}[1]{Construction \ref{C:#1}}
\newcommand{\re}[1]{\ref{E:#1}}
\newcommand{\rco}[1]{Corollary \ref{C:#1}}

\newcommand{\rt}[1] {Theorem \ref{T:#1}}
\newcommand{\rd}[1]{Definition \ref{D:#1}}
\newcommand{\sm}{\smallsetminus}
\newcommand{\be}{\infty}

\newcommand{\Mat}{\operatorname{Mat}}

\newcommand{\IC}{\operatorname{IC}}
\newcommand{\inv}{\operatorname{inv}}
\newcommand{\pr}{\operatorname{pr}}
\newcommand{\Ker}{\operatorname{Ker}}
\newcommand{\Out}{\operatorname{Out}}
\newcommand{\Coker}{\operatorname{Coker}}
\newcommand{\Norm}{\operatorname{Norm}}
\newcommand{\Cent}{\operatorname{Cent}}
\newcommand{\val}{\operatorname{val}}

\newcommand{\im}{\operatorname{Im}}

\newcommand{\Aut}{\operatorname{Aut}}
\newcommand{\Gr}{\operatorname{Gr}}
\newcommand{\ab}{\operatorname{ab}}

\newcommand{\Inn}{\operatorname{Int}}
\newcommand{\sep}{\operatorname{sep}}
\newcommand{\Ind}{\operatorname{Ind}}  

\newcommand{\Ad}{\operatorname{Ad}}

\newcommand{\Isom}{\operatorname{Isom}}
\newcommand{\Gal}{\operatorname{Gal}}
\newcommand{\Tr}{\operatorname{Tr}}

\newcommand{\red}{\operatorname{red}}
\newcommand{\reg}{\operatorname{reg}}
\newcommand{\rgss}{\operatorname{rss}}
\newcommand{\Fr}{\operatorname{Fr}}

\newcommand{\Lie}{\operatorname{Lie}}
\newcommand{\ad}{\operatorname{ad}}
\newcommand{\Int}{\operatorname{Int}}
\newcommand{\Id}{\operatorname{Id}}

\newcommand{\ord}{\operatorname{ord}}

\newcommand{\End}{\operatorname{End}}
\newcommand{\Emb}{\operatorname{Emb}}

\newcommand{\rk}{\operatorname{rk}}
\newcommand{\Hom}{\operatorname{Hom}}

\newcommand{\lan}{\left\langle}
\newcommand{\ran}{\right\rangle}

\newcommand{\e}{\par \noindent}
\begin{document}
\title[Endoscopic decomposition]%
{On endoscopic decomposition of certain depth zero representations}

\author{David Kazhdan}
\author{Yakov Varshavsky}
\address{Institute of Mathematics, 
The Hebrew University of Jerusalem, Givat-Ram\\
Jerusalem  91904 Israel}
\email{kazhdan@math.huji.ac.il, vyakov@math.huji.ac.il }
\dedicatory{Dedicated to A.~Joseph on his 60th birthday}

\thanks{Both authors were supported by 
THE ISRAEL SCIENCE FOUNDATION (Grants No. 38/01 and 241/03)}
\date{November 2005}

\begin{abstract}
We construct an endoscopic decomposition for local $L$-packets 
associated to irreducible cuspidal Deligne--Lusztig 
representations. Moreover, the obtained decomposition is compatible 
with inner twistings.
\end{abstract}
\maketitle

\vskip 6truept
\centerline{\bf Introduction}
\vskip 6truept

Let $E$ be a local non-archimedean field, $\Gm\supset W\supset I$ the absolute Galois, 
the Weil and the inertia groups of $E$, respectively. Let $G$ be a reductive group over 
$E$, and let ${}^LG=\wh{G}\rtimes W$ be the complex Langlands dual group of $G$. 
Denote by $\C{D}(G(E))$ the 
space of invariant generalized functions on $G(E)$, that is, the space of 
$\Inn G$-invariant linear functionals on the space of locally constant 
compactly supported measures on $G(E)$.

Every admissible homomorphism $\la:W\to {}^LG$ (see \cite[$\S$ 10]{Ko1})  
gives rise to a finite group $S_{\la}:=\pi_0(Z_{\wh{G}}(\la)/Z(\wh{G})^{\Gm})$, 
where $Z_{\wh{G}}(\la)$ is the centralizer of $\la(W)$ in $\wh{G}$. 
To every conjugacy class $\ka$ of $S_{\la}$, Langlands \cite{La} associated an endoscopic 
subspace $\C{D}_{\ka,\la}(G(E))\subset\C{D}(G(E))$. For simplicity, we will restrict 
ourselves to the elliptic case, where $\la(W)$ does not lie in any proper Levi subgroup of ${}^LG$.

Langlands conjectured that every elliptic $\la$ corresponds to a finite set $\Pi_{\la}$,  
called an $L$-packet, of cuspidal irreducible representations of $G(E)$. 
Moreover, the subspace $\C{D}_{\la}(G(E))\subset\C{D}(G(E))$, generated by characters 
$\{\chi(\pi)\}_{\pi\in\Pi_{\la}}$, should have an endoscopic decomposition. 
More precisely, it is expected (\cite[IV, 2]{La})
that there exists a basis $\{a_{\pi}\}_{\pi\in \Pi_{\la}}$ of the space of central functions on  
$S_{\la}$ such that 
$\chi_{\ka,\la}:=\sum_{\pi\in \Pi_{\la}}a_{\pi}(\ka) \chi(\pi)$ belongs to 
$\C{D}_{\ka,\la}(G)$ for every conjugacy class $\ka$ of $S_{\la}$. 

The goal of this paper is to construct the endoscopic decomposition of 
$\C{D}_{\la}(G(E))$ for tamely ramified $\la$'s such that $Z_{\wh{G}}(\la(I))$ is 
a maximal torus. In this case, $G$ splits over
an unramified extension of $E$, and $\la$ factors through ${}^LT\hra{}^LG$ for an elliptic 
unramified maximal torus $T$ of $G$.

Each $\ka\in S_{\la}=\wh{T}^{\Gm}/Z(\wh{G})^{\Gm}$ gives rise to an elliptic endoscopic 
triple $\C{E}_{\ka,\la}$ for $G$, while characters of $S_{\la}$ are in bijection with 
conjugacy classes of embeddings $T\hra G$, stably conjugate to the inclusion. 
By the local Langlands correspondence for tori (\cite{La2}), a homomorphism 
$\la:W\to {}^LT$ defines a tamely ramified homomorphism $\theta:T(E)\to\B{C}\m$.
Therefore each character $a$ of $S_{\la}$ gives rise to an irreducible cuspidal 
representation $\pi_{a,\la}$ of $G(E)$ (denoted by $\pi_{a,\theta}$ in \rn{repr}). 

Our main theorem asserts that if the residual characteristic of $E$
is sufficiently large, 
then each $\chi_{\ka,\la}:=\sum_{a}a(\ka) \chi(\pi_{a,\la})$ is 
$\C{E}_{\ka,\la}$-stable (see \rd{stab}). 
More generally (see \rco{end} (b)), for each inner form $G'$ of $G$, we denote by 
$\chi'_{\ka,\la}$ the corresponding generalized function on $G'(E)$, and our main theorem
asserts that $\chi_{\ka,\la}$ and $\chi'_{\ka,\la}$ are ``$\C{E}_{\ka,\la}$-equivalent''
(see \rd{main} for a more precise term).
Though in this work we show this result only for local fields of characteristic zero, the  
case of local fields of positive characteristic follows by approximation (see \cite{KV2}).

Our argument goes as follows. First we prove the equivalence  
of the restrictions of $\chi_{\ka,\la}$ and $\chi'_{\ka,\la}$ to the subsets of  
topologically unipotent elements of $G(E)$ and $G'(E)$. 
If the residual characteristic of $E$ is sufficiently large, 
topologically unipotent elements of  $G(E)$ and $G'(E)$ can be identified with 
topologically nilpotent elements of the Lie algebras $\C{G}(E)$ and $\C{G}'(E)$, 
respectively. Thus we are reduced to an analogous assertion about generalized functions
on Lie algebras. Now the equivalence follows from a combination of a Springer hypothesis
(\rt{Spr}), which describes the trace of a Deligne--Lusztig representation in terms of 
Fourier transform of an orbit, 
and a generalization a theorem of Waldspurger \cite{Wa2} to inner forms, which asserts 
that  up to a sign, Fourier transform preserves the equivalence.

To prove the result in general, we use the topological Jordan decomposition (\cite{Ka2}).
We would like to stress that in order to prove just the stability of 
$\chi_{\ka,\la}$ 
one still needs a generalization of \cite{Wa2} to inner forms.

This paper is organized as follows. 

In the first section we give basic definitions and constructions of a rather general 
nature. In particular, most of the section is essentially a theory of endoscopy, which was 
developed by Langlands, Shelstad and Kottwitz. In order to incorporate 
both the case of algebraic groups and of Lie algebras we work in a more general context of 
algebraic varieties equipped with an action of $G^{\ad}$. 

More precisely, in Subsection 1.1 we recall basic properties, results and 
constructions concerning inner twistings and stable conjugacy. 
In Subsections 1.2 and 1.3 we give basic definitions and properties of dual groups 
and of endoscopic triples. Then in Subsection 1.4 we prove that 
certain subsets of the group $Z(\wh{G^{\ad}})^{\Gm}$ are actually subgroups. 
Unfortunately, this result is proven case-by-case.
In Subsection 1.5 we specialize previous results to the case of endoscopic triples over 
local non-archimedean fields.

In Subsection 1.6 we define the notions of stability and equivalence of generalized functions, while
in Subsection 1.7 we write down explicitly the condition for stability and equivalence 
for generalized functions coming from invariant locally $L^1$ functions.
Note that the notion of equivalence is much more subtle than that of stability.
In particular, it depends not just on an endoscopic triple but also on a triple
$(a,a';[b])$, consisting of compatible embeddings of maximal tori into $G$, $G'$ and 
the endoscopic group. 

We finish the section by Subsection 1.8 in which we study basic properties of 
certain equivariant maps from reductive groups to their Lie algebras, which we call 
quasi-logarithms. We use these maps to identify topologically unipotent elements 
of the group with topologically nilpotent elements of the Lie algebra. 

The second section is devoted to the formulation and the proof of the  main theorem. 
More precisely, in Subsection 2.1 we give two equivalent formulations of our main result.
In Subsection 2.2 we prove the equivalence of the restrictions of $\chi_{\ka,\la}$
and $\chi'_{\ka,\la}$ to topologically unipotent elements.

In Subsection 2.3 we rewrite character $\chi(\pi_{a,\theta})$ of $G(E)$
in terms of restrictions to  topologically unipotent elements of corresponding characters of 
centralizers $G_{\dt}(E)$. For this we use the topological Jordan decomposition. 
In the next Subsection 2.4 we compare endoscopic triples for the group $G$ and 
for its centralizers $G_{\dt}(E)$. Finally, in Subsections 2.5--2.7 we carry out the 
proof itself.

We finish the paper by two appendices of independent interest, 
crucially used in Subsection 2.2. In Appendix A we prove Springer hypothesis. 
In the case of large characteristic this result was proved by the first author in 
\cite{Ka1}. For the proof in general, we use Lusztig's 
interpretation of a trace of Deligne--Lusztig representation in terms of character sheaves 
\cite{Lu} and results of Springer \cite{Sp} on Fourier transform. 

In Appendix B we prove a generalization of both the theorem of Waldspurger \cite{Wa2} 
and that of Kazhdan--Polishchuk \cite[Thm. 2.7.1]{KP} (see also \rr{ch}).
Our strategy is very similar to those of \cite{Wa2} and \cite{KP}.
More precisely, using stationary phase principle and the results of Weil \cite{We},
we construct in Subsection B.2 certain measures whose Fourier transform can be explicitly 
calculated. Then in Subsection B.3 we extend our data over a local field to a 
corresponding data over a number field. Finally, in Section 1.4 we deduce our result 
from a simple form of the trace formula.

For the convenience of the reader, we also include a list of main terms and symbols, indexed by page number 
they first appear.


This work is an expanded version of the announcement \cite{KV}. In the process of writing, 
we have learned that DeBacker and  Reeder obtained similar results (see \cite{DBR}).
After the work was completed, it was pointed out to us that our scheme of the argument is similar 
to the one used by M\oe glin-Waldspurger in \cite{MW}.

We would like to thank the referee for his numerous valuable remarks.

\vskip 6truept
\centerline{\bf Notation and Conventions}
\vskip 6truept

For a finite abelian group $A$, we denote by $A^D$ the group of complex characters of $A$.

  For an algebraic group $G$, we denote by $G^0$,\label{g0} $Z(G)$, \label{zg} $G^{\ad}$, \label{gad} $G^{\der}$ and $W(G)$ \label{wg}
the connected component of the identity of $G$, the center of $G$, the adjoint group of $G$, 
the derived group of $G$, and the Weyl group of $G$, respectively. Starting from \re{toremb}, 
$G$ will be always assumed
to be reductive and connected, in which case we  denote by $G^{\ssc}$ \label{gss}
the simply connected covering of $G^{\der}$.

We denote by $\C{G}$, \label{lieg} $\C{H}$, \label{lieh} $\C{T}$ \label{liet} and $\C{L}$ \label{liel} Lie algebras
of algebraic groups $G$, $H$, $T$ and $L$, respectively.

Let an algebraic group $G$ acts on an algebraic variety $X$. For each $x\in X$, 
we denote by $G_x$ and $\C{G}_x$ \label{stab} the stabilizers of $x$ in $G$ and $\C{G}$, 
respectively. Explicitly, $\C{G}_x$ is the kernel of the differential at $g=1$ of the morphism
$G\to X\; (g\mapsto g(x))$. (Note that that $x$ and therefore also $G_x$ and $\C{G}_x$ 
will have different meaning starting from \re{BT}.) 

Each algebraic group acts on itself by inner automorphisms and by adjoint action on its Lie algebra.
For each $g\in G$ we denote by $\Inn g$ \label{intg} and $\Ad g$ \label{adg} the corresponding elements in 
$\Inn G\subset\Aut(G)$ and $\Ad G\subset \Aut(\C{G})$, respectively.

For a field $E$, we denote by $\ov{E}$ \label{ove} a fixed algebraic closure of $E$, 
and by $E^{\sep}$ \label{esep} the maximal separable extension of $E$ in $\ov{E}$.
$\Gm$ \label{Gm} will be always the absolute Galois group of $E$.
When $\Gm$ acts on a set $X$ we will write ${}^{\si}x$ instead of $\si(x)$.

For a reductive algebraic group $G$ over $E$, we denote by $\rk_E(G)$ the rank of $G$ over $E$, 
and put $e(G):=(-1)^{\rk_E(G^{\ad})}$. \label{eg}
We also set $e'(G):=e(G)e(G^*)$, where $G^*$ is the quasi-split inner form of $G$.
Then $e'(G)$ coincides with the sign defined by Kottwitz (\cite{Ko5}). 

Starting from \rss{local},  $E$ will be a local non-archimedean field with ring 
of integers $\C{O}$, \label{O} maximal ideal $\frak{m}$, \label{m} and residue field $\fq$ \label{fq} of characteristic
$p$. \label{p} We denote by $E^{\nr}$ \label{enr} the maximal unramified extension of $E$ in $\ov{E}$.
Starting from \rss{unip}, we will assume that the characteristic of $E$ is zero.

\tableofcontents

\section{Basic definitions and constructions}

\subsection{Stable conjugacy}

\begin*
\vskip 8truept

In this subsection we recall basic definitions and constructions concerning inner forms 
and stable conjugacy. 
 
Let $G$ be an algebraic group over a field $E$. Starting from \ref{E:toremb}, 
we will assume that $G$ is reductive and connected. 

Let $X$ be an algebraic variety  over $E$ 
(that is, a reduced scheme locally of finite type over $E$)
equipped with an action of $G^{\ad}$ (that is, with an action of $G$, trivial on $Z(G)$).
 
Our basic examples will be $X=G$ and $X=\C{G}$ with the natural action of $G^{\ad}$.
\end*

\begin{Emp} \label{E:inner}
{\bf Inner twistings.}
(a) Let $G'$ be an algebraic group over $E$.  Recall that 
an {\em inner twisting} \label{it} $\varphi:G\to G'$ is an isomorphism $\varphi:G_{E^{\sep}}\isom G'_{E^{\sep}}$ 
such that  for each $\si\in\Gm$ the automorphism 
$c_{\si}:=\varphi^{-1}{}^{\si}\varphi\in\Aut(G)$ is inner.
In this case, $\{c_{\si}\}_{\si}$ form a cocycle of $\Gm$ in 
$\Inn G=G^{\ad}$, and we denote by  
$\inv(\varphi)=\inv(G,G')\in H^1(E,G^{\ad})$ \label{invg'g} the corresponding cohomology 
class. 

(b) Two inner twistings are called {\em isomorphic} if they differ by an 
inner automorphism.
Then the map $(\varphi:G\to G')\mapsto\inv(G,G')$ gives a bijection between 
the set of isomorphism classes
of inner twistings of $G$ and $H^1(E,G^{\ad})$. 

(c) Each inner twisting $\varphi:G\to G'$ gives rise to a twisting 
$\varphi_X:X\to X'$, \label{varphiX}  where $X'$ is an algebraic variety over $E$ equipped 
with an action of $G'$, and
$\varphi_X$ is a $G_{E^{\sep}}\cong G'_{E^{\sep}}$-equivariant isomorphism  
$X_{E^{\sep}}\isom X'_{E^{\sep}}$. Explicitly, $X'$ is a twist of $X$ 
by the image of the cocycle 
$\{c_{\si}\}_{\si}\subset G^{\ad}$ in $\Aut(X)$. In particular, for each $\si\in\Gm$ 
we have ${}^{\si}\varphi_X=\varphi\circ c_{\si}$.

By construction, for each $x\in X$ and $g\in G$, we have $\varphi_X(g(x))=\varphi(g)(\varphi_X(x))$.

(d) An inner twisting $\varphi$ is called {\em trivial}, if  $\inv(\varphi)=1$.
Explicitly, $\varphi$ is trivial if and only if there exists $g\in G(E^{\sep})$
such that $\varphi\circ\Inn g$ induces an isomorphism $G\isom G'$ over $E$.
In particular, the identity map $\Id_G:G\to G$ is a trivial inner twisting. 
\end{Emp}

\begin{Def} \label{D:stconj} 
(a) Two points $x,x'\in X(E)$ are called 
{\em conjugate}, \label{conj} if there exists $g\in G(E)$ such that 
$x'=g(x)$.

(b) Let $\varphi:G\to G'$ be an inner twisting, and $\varphi_X:X\to X'$ the corresponding
twisting. Elements $x\in X(E)$ and $x'\in X'(E)$ are called {\em $E^{\sep}$-conjugate}, 
\label{esepconj}
if there exists $g\in G(E^{\sep})$ such that $x'=\varphi_X(g(x))$.

(c) When $G$ and $G_x$ (and hence also $G'$ and $G'_{x'}$) are connected reductive groups,
$E^{\sep}$-conjugate $x$ and $x'$ are also called {\em stably conjugate}.\label{stcon} 
\end{Def}

\begin{Rem} \label{R:stconj}
All of our examples will satisfy assumption (c). In this case, our notion of stable conjugacy 
generalizes the standard one (see \cite{Ko3}).
\end{Rem}


\begin{Emp} \label{E:inv}
{\bf Cohomological invariants.} 
Let $x\in X(E)$ and $x'\in X'(E)$ be $E^{\sep}$-conjugate elements.  
Denote by $G_{x,x'}$
the set of $g\in G(E^{\sep})$ such that $x'=\varphi_X(g(x))$. 

(a) Assume that $\varphi=\Id_G$. Then for each 
$g\in G_{x,x'}$, the map $\si\mapsto  g^{-1}{}^{\sigma}g$ defines a cocycle
of $\Gm$ in $G_x$. Moreover, the corresponding cohomology class  
$\inv(x,x')\in H^1(E,G_x)$ \label{invxx'} is independent of $g$. 
Furthermore, the correspondence $x'\mapsto \inv(x,x')$ gives a bijection between the set of 
conjugacy classes of $x'\in X(E)$ stably conjugate to $x$ and
$\Ker\,[H^1(E,G_x)\to H^1(E,G)]$ (compare \cite[4.1]{Ko2}).

(b) Let $\varphi$ be general. Then for each $g\in G_{x,x'}$,
the map $\si\mapsto  g^{-1}(\varphi^{-1}{}^{\si}\varphi){}^{\sigma}g$ defines a 
cocycle of $\Gm$ in $G_x/Z(G)=(G^{\ad})_x\subset\Inn G$. 
Moreover, the corresponding cohomology class 
$\ov{\inv}(x,x')\in H^1(E,(G^{\ad})_x)$ \label{ovinvxx'} is independent of $g$. 
Furthermore, the correspondence $x'\mapsto \ov{\inv}(x,x')$ gives a surjection 
from the set of 
conjugacy classes of $x'\in X'(E)$ stably conjugate to $x$ to the preimage of 
$\inv(G,G')\in H^1(E,G^{\ad})$ in $H^1(E,(G^{\ad})_x)$. 

When $\varphi=\Id_G$, then 
$\ov{\inv}(x,x')$ is the image of $\inv(x,x')$ under the natural projection
$H^1(E,G_x)\to H^1(E,(G^{\ad})_x)$.

(c) For each $g\in G_{x,x'}$, the map $h\mapsto \varphi(ghg^{-1})$ defines an inner twisting
$G_x\to G'_{x'}$. Moreover, the corresponding invariant 
$\inv(G_{x},G'_{x'})\in  H^1(E,(G_x)^{\ad})$ is just the image of $\ov{\inv}(x,x')$.
In particular, $G'_{x'}$ is canonically identified with $G_x$, if $G_x$ is abelian.

(d) Assume that  $G_x$ is abelian, and let $y\in X(E)$ and $y'\in X'(E)$ be $E^{\sep}$-conjugates 
of $x$ and $x'$. Then the identification $G'_{y'}=G'_{x'}=G_x=G_y$ identifies
$\ov{\inv}(y,y')$ with the product of $\ov{\inv}(x,x')$ and the images of 
 ${\inv}(y,x)$ and ${\inv}(x',y')$. Moreover, if $\varphi=\Id_G$, the same identification
identifies ${\inv}(y,y')$ with the product ${\inv}(y,x)\inv(x,x'){\inv}(x',y')$.
\end{Emp}

\begin{Emp} \label{E:coh}
{\bf Generalization }(compare \cite[(3.4)]{LS}). 
Let $\varphi:G\to G'$ be an inner twisting, $X_1,\ldots, X_k$  a $k$-tuple of 
algebraic varieties over $E$ equipped with an action of $G^{\ad}$, and  
$X'_1,\ldots, X'_k$ the corresponding inner twistings. 
Let $x_i\in X_i(E)$ and $x'_i\in X'_i(E)$ be $E^{\sep}$-conjugate
for each $i=1,\ldots,k$. 

Choose representatives $\wt{c}_{\si}\in G(E^{\sep})$ 
of $c_{\si}=\varphi^{-1}{}^{\si}\varphi\in G^{\ad}(E^{\sep})$ for all $\si\in\Gm$ 
and choose elements $g_i\in G(E^{\sep})$ such that $x'_i=\varphi_{X_i}(g_i(x_i))$
for $i=1,\ldots,k$. Then the map 
$\si\mapsto [(g_i^{-1}\wt{c}_{\si}{}^{\si}g_i)_i]\in G^k/Z(G)$ gives a cocycle of $\Gm$ in
$(\prod_i G_{x_i})/Z(G)$, independent of the choice of $\wt{c}_{\si}$'s, and 
the corresponding cohomology class
\[
\ov{\inv}((x_1,x'_1);\ldots;(x_k,x'_k))\in H^1(E,(\prod_i G_{x_i})/Z(G)) 
\]\label{ovinvkxx'}
of $[(g_i^{-1}\wt{c}_{\si}{}^{\si}g_i)_i]$ is independent of the $g_i$'s.

Note that $\ov{\inv}((x_1,x'_1);\ldots;(x_k,x'_k))$
lifts  both $(\ov{\inv}(x_i,x'_i)_i)\in \prod_i H^1(E,G_{x_i}/Z(G))$ and 
 $\Dt(\inv(G,G'))\in H^1(E,(G^k)/Z(G))$. (Here  
$\Dt:G^{\ad}\to (G^k)/Z(G)$ is the diagonal embedding.)
\end{Emp}

The following result follows immediately from definitions.

\begin{Lem} \label{L:propinv}
(a) $\ov{\inv}((x_1,x'_1);\ldots;(x_k,x'_k))$ depends only on the conjugacy classes of  
$x_i$'s and $x'_i$'s. 

(b) If $\varphi=\Id_G$, then 
$\ov{\inv}((x_1,x'_1);\ldots;(x_k,x'_k))$ is the 
image of  $((\inv(x_i,x'_i)_i)\in H^1(E,\prod_i G_{x_i})$.

(c) The canonical projection
$(\prod_{i=1}^k  G_{x_i})/Z(G)\to (\prod_{i=1}^{k-1}  G_{x_i})/Z(G)$ maps \\
$\ov{\inv}((x_1,x'_1);\ldots;(x_k,x'_k))$ to 
$\ov{\inv}((x_1,x'_1);\ldots;(x_{k-1},x'_{k-1}))$.

(d) The diagonal map $G_x/Z(G)\hra (G_x)^2/Z(G)$ maps 
$\ov{\inv}(x,x')\in  H^1(E,G_x/Z(G))$ to  $\ov{\inv}((x,x');(x,x'))\in H^1(E,(G_x)^2/Z(G))$.

(e) Assume that each $G_{x_i}$ is abelian, 
$y_i\in X_i(E)$ is an $E^{\sep}$-conjugate of $x_i$, and $y'_i\in X'_i(E)$ is 
an $E^{\sep}$-conjugate of $x'_i$. Then identifications 
$G_{x_i}=G_{y_i}=G_{x'_i}$ identify $\ov{\inv}((y_1,y'_1);\ldots;(y_k,y'_k))$ 
with the product of $\ov{\inv}((x_1,x'_1);\ldots;(x_k,x'_k))$ and the images
of $((\inv(y_i,x_i)_i)\in H^1(E,\prod_i G_{y_i})$ and 
$((\inv(x'_i,y'_i)_i)\in H^1(E,\prod_i G_{x'_i})$.

(f) Let $\varphi':G'\to G''$ be another inner twisting, $\varphi'_{X'}:X'\to X''$ 
the corresponding twisting of $X'$, and $x''_i\in X''(E)$ a stable conjugate of
$x_i$ and $x'_i$ for each $i=1,\ldots,k$. If each $G_{x_i}$ is abelian,  
then identifications $G_{x_i}=G_{x'_i}$ identify
$\ov{\inv}((x_1,x''_1);\ldots;(x_k,x''_k))\ov{\inv}((x_1,x'_1);\ldots;(x_k,x'_k))^{-1}$ with
$\ov{\inv}((x'_1,x''_1);\ldots;(x'_k,x''_k))$.
\end{Lem}

\begin{Not} \label{N:sconj}
(a) For each $x\in X(E)$, we denote by 
$[x]\subset X(E)$ \label{[x]} the $E^{\sep}$-conjugate class of $x$ and by 
$a_x:G_x\hra G$ \label{ax} the corresponding inclusion map.

(b) When $G$ is reductive and connected, we denote by 
$X^{\sr}$ \label{sr} the set of $x\in X$ such that $G_x\subset G$ is a maximal torus
and $\C{G}_x\subset\C{G}$ is a Lie algebra of $G_x$. 
We will call elements of $X^{\sr}$ {\em strongly regular}. \label{str}
\end{Not}

\begin{Rem}
The condition on $\C{G}_x$ holds automatically if the characteristic of $E$ is zero.
\end{Rem}

From now on we will assume that $G$ is reductive and connected, and
$T$ is a torus over $E$ of the same absolute rank as $G$.

\begin{Emp} \label{E:toremb}
{\bf Embedding of tori.}
(a)  There exists an affine variety $\un{\Emb}(T,G)$  over $E$
equipped with an action of $G^{\ad}$ such that for every extension $E'/E$ the set
$\un{\Emb}(T,G)(E')$ classifies embeddings $T_{E'}\hra G_{E'}$, and $G$ acts by conjugation.

To show the assertion, note that  both $G$ and $T$ split over $E^{\sep}$, 
therefore there exists an embedding
$\iota:T_{E^{\sep}}\hra G_{E^{\sep}}$. Consider the affine variety $\un{\Emb}_{\iota}(T,G):=
G\times_{\Norm_{G}(\iota(T))}\Aut(T)$ over $E^{\sep}$. Then 
the map 
$[g,\al]\mapsto (\Inn g)\circ\iota\circ\al$ defines a $G(E')$-equivariant bijection  
$\psi_{\iota}:\un{\Emb}_{\iota}(T,G)(E')\isom \un{\Emb}(T,G)(E')$ for every extension $E'$ of
$E^{\sep}$ (compare (the proof of) \rl{steinberg} (a) below). 

For every two embeddings $\iota_1,\iota_2:T_{E^{\sep}}\hra G_{E^{\sep}}$, 
there exists a unique isomorphism
$\un{\Emb}_{\iota_1}(T,G)\isom\un{\Emb}_{\iota_2}(T,G)$, compatible with the 
$\psi_{\iota_j}$'s, 
and we define $\un{\Emb}(T,G)$ be the inverse limit of the $\un{\Emb}_{\iota}(T,G)$'s.
Finally, since $\un{\Emb}_{\iota}(T,G)$ is a disjoint union of affine varieties, 
it descends to $E$.

(b) For each $a\in\un{\Emb}(T,G)$, the stabilizer $G_a=a(T)$ is a maximal torus
of $G$, which we will identify with $T$. It follows that $\un{\Emb}(T,G)^{\sr}=\un{\Emb}(T,G)$.
Also if $\varphi:G\to G'$ is an inner twisting, then the corresponding inner twisting  
$\un{\Emb}(T,G)'$ of $\un{\Emb}(T,G)$ is naturally isomorphic to $\un{\Emb}(T,G')$.
In particular, we can speak about stable conjugacy of embeddings $a:T\hra G$ and $a':T\hra G'$.

(c) If $x\in X^{\sr}(E)$ and $x'\in X'^{\sr}(E)$ are stably conjugate, then 
$a_x:G_x\hra G$ and 
$a_{x'}:G_x\cong G'_{x'}\hra G'$
are stably conjugate. Moreover, $\ov{\inv}(a_{x},a_{x'})=\ov{\inv}(x,x')$, (and  
$\inv(a_x,a_{x'})=\inv(x,x')$ when $\varphi=\Id_G$.)  

 Conversely, for every stably conjugate embeddings $a:T\hra G$ and
$a':T\hra G'$ and each $t\in T(E)$, elements $a(t)\in G(E)$ and $a'(t)\in G'(E)$ are 
$E^{\sep}$-conjugate.
\end{Emp}

\begin{Lem} \label{L:steinberg}

(a) Every conjugacy class $[a]$ of embeddings 
$T_{\ov{E}}\hra G_{\ov{E}}$ contains an  $E^{\sep}$-rational embedding $a$.

(b) If $G$ is quasi-split, then every $\Gm$-invariant conjugacy class $[a]$ of embeddings 
$T_{\ov{E}}\hra G_{\ov{E}}$ contains an  $E$-rational embedding $a:T\hra G$. 

(c) Let $\varphi:G\to G'$ be an inner twisting such that $G'$ is
quasi-split. Then for every 
embedding $a:T\hra G$  there exists an embedding 
$a':T\hra G'$ stably conjugate to $a$.
\end{Lem}
\begin{proof}

(a) Let $S\subset G$ be a maximal torus over $E$. Then there exists $a\in [a]$ 
such that $a(T_{\ov{E}})=S_{\ov{E}}$. Since both $T$ and $S$ split over $E^{\sep}$, 
we get that $a$ is $E^{\sep}$-rational.

(b) When $E$ is perfect, the assertion was shown in \cite[Cor. 2.2]{Ko3}. In general, the proof 
is similar: Since there is a $\Gm$-equivariant bijection between maximal tori of $G$ and 
those of $G^{\ssc}$, we can assume that $G$ is semisimple and simply connected. 
Next fix an $E^{\sep}$-rational $a'\in [a]$, which exists by (a). Since every 
homogeneous space for a connected group over a finite field has a rational point, 
the assertion holds in this case. Thus we may assume that $E$ is infinite, 
therefore there exists $t\in T(E)$ such that $a'(t)\in G(E^{\sep})$ is 
strongly regular.
 
The conjugacy class of $a'(t)\in G(E^{\sep})$ is 
$E^{\sep}$-rational and $\Gm$-invariant. Thus it is $E$-rational. 
By the theorem of Steinberg \cite[Thm 1.7]{St} 
(when $E$ is perfect) and Borel and Springer \cite[8.6]{BS} 
(in the general case) there exists $g\in G(E^{\sep})$ such that $g a'(t) g^{-1}\in G(E)$. 
Then $g a'(t) g^{-1}\in G^{\sr}(E)$, hence $a:=ga' g^{-1}:T\hra G$ is $E$-rational.

(c) Since $\varphi$ is an inner twisting, the conjugacy class of 
$\varphi\circ a:T_{\ov{E}}\hra G'_{\ov{E}}$ is $\Gm$-invariant.
Hence by (b), $[\varphi\circ a]$ 
contains an $E$-rational element $a'$, which by definition is stably conjugate to $a$.
\end{proof}

\begin{Cor} \label{C:steinberg}
Let  $\varphi:G\to G'$ be an inner twisting, and $\varphi_X:X\to X'$ the corresponding 
twisting.
If $G'$ is quasi-split, then  for every $x\in X^{\sr}(E)$ there exists $x'\in X'^{\sr}(E)$
stably conjugate to $x'$.
\end{Cor}
\begin{proof}
Denote by $\iota:G/G_x\hra X$ the canonical $G$-equivariant
embedding $[g]\mapsto g(x)$, and by $\iota':(G/G_x)'\hra X'$ the twisted map.
By \rl{steinberg} (c), there exists an embedding  $a'_x:G_x\hra G'$ stably conjugate to 
$a_x:G_x\hra G$. Moreover, $(G/G_x)'$ is $G'$-equivariantly isomorphic to 
$G'/a'_x(G_x)$. It follows that the image of $[1]\in G'/a'_x(G_x)(E)$ under 
$\iota'$ is stably conjugate to $x$.
\end{proof}

\begin{Def} \label{D:qisog}
By a {\em quasi-isogeny} \label{qi} we call a homomorphism $\pi:\wt{G}\to G$ such that
$\pi(Z(\wt{G}))\subset Z(G)$ and the induced homomorphism 
$\pi^{\ad}:\wt{G}^{\ad}\to G^{\ad}$ is an isomorphism.
\end{Def}

\begin{Emp} \label{E:qisog}
{\bf Quasi-isogenies.}
Let $\pi:\wt{G}\to G$ be a quasi-isogeny.

(a) Each inner twisting $\varphi:G\to G'$ gives rise to an inner twisting 
$\wt{\varphi}:\wt{G}\to\wt{G}'$ such that 
$\inv(\wt{G},\wt{G}')=\inv(G,G')$ (in $H^1(E,\wt{G}^{\ad})=H^1(E,G^{\ad})$).

(b) $X$ is equipped with an action of $\wt{G}$ trivial on $Z(\wt{G})$. 

(c) There is a $\pi$- and $\Gm$-equivariant bijection between embeddings of maximal tori
$a:T\hra G$ and the corresponding embeddings $\wt{a}:\wt{T}\hra \wt{G}$. Indeed, given $a$, 
put $\wt{T}=T\times_G\wt{G}:=\{t\in T, \wt{g}\in \wt{G}\,|\,a(t)=\pi(\wt{g})\}$ and
$\wt{a}(t,\wt{g}):=\wt{g}$. Conversely, given $\wt{a}$, define $a$ be the embedding
$\pi(\wt{a}(\wt{T}))\hra G$. In particular, $a_1$ and $a_2$ are stably conjugate if and only if 
$\wt{a_1}$ and $\wt{a_2}$ are such. 

We will call $\wt{a}$ {\em the lift of $a$} and $[\wt{a}]$ {\em the lift of $[a]$}.

(d) For each $i=1,\ldots,k$, let ${a}_i:{T}_i\hra {G}$ and ${a}'_i:{T}_i\hra {G}'$ 
be stable conjugate embeddings of maximal tori, and let  
 $\wt{a}_i:\wt{T}_i\hra \wt{G}$ and $\wt{a}'_i:\wt{T}_i\hra \wt{G}'$ 
be the lifts of  the $a_i$'s and the $a'_i$'s, respectively.
Then $\ov{\inv}((a_1,a'_1);\ldots;(a_k,a'_k))\in H^1(E,(\prod_i{T}_i)/Z({G})) $
is the image of  $\ov{\inv}((\wt{a}_1,\wt{a}'_1);\ldots;(\wt{a}_k,\wt{a}'_k))
\in H^1(E,(\prod_i\wt{T}_i)/Z(\wt{G})) $ 
under the canonical map $(\prod_i\wt{T}_i)/Z(\wt{G}))\to(\prod_i{T}_i)/Z({G})) $. 
\end{Emp}

\subsection{Preliminaries on dual groups} \label{SS:dual}
\begin*
\vskip 8truept
\end*
In this subsection, we will recall basic properties of Langlands dual groups. 
More specifically, we will study properties of triples $(G,H,[\eta])$ from \re{tr}. Constructions from this subsection will be later used in the case when $H$ is an endoscopic group for $G$.

\begin{Not} \label{N:dual}
For each  connected reductive group $G$ over a field $E$, we denote by $\wh{G}$ \label{dualg}
(or $\{G\}\:\wh{}$ ) the complex connected Langlands dual group, and by
$\rho_G:\Gm\to\Out(\wh{G})$ \label{rhog} the corresponding Galois action.
\end{Not}

\begin{Emp} \label{E:dual}
{\bf Basic properties of dual groups.}
(a) The map $G\mapsto (\wh{G},\rho_G)$ defines a surjection from the 
set of isomorphism classes of connected reductive groups over $E$ to 
that of pairs consisting of a connected
complex reductive group $\wh{G}$, and a continuous homomorphism 
$\rho:\Gm\to\Out(\wh{G})$. Moreover, 
$(\wh{G_1},\rho_{G_1})\cong(\wh{G_2},\rho_{G_2})$ if and only if $G_2$
is an inner twist of $G_1$. In particular, each pair $(\wh{G},\rho_G)$ comes from a 
unique quasi-split group $G$ over $E$.

(b) Let $T$ be a torus over $E$ of the same absolute rank as $G$.
Then there exists a canonical (hence $\Gm$-equivariant) bijection 
$[a]\mapsto \wh{[a]}$ \label{wha} between conjugacy classes of embeddings 
$T_{\ov{E}}\hra G_{\ov{E}}$ and conjugacy classes of embeddings $\wh{T}\hra\wh{G}$. 
In particular, $[a]$ is $\Gm$-invariant if and only if $\wh{[a]}$ is such.

(c) For each embeddings of maximal tori $T_{\ov{E}}\hra G_{\ov{E}}$ and $\wh{T}\hra\wh{G}$
related as in (b), the set of roots (resp. coroots) of $(G_{\ov{E}},T_{\ov{E}})$
is canonically identified with  the set of coroots (resp. roots) of $(\wh{G},\wh{T})$.
In particular, the Weyl group $W(\wh{G})$ is canonically identified with $W(G)$.


(d) Each quasi-isogeny $\pi:G_1\to G_2$ gives rise to a conjugacy class 
$[\wh{\pi}]$ \label{whpi} of quasi-isogenies $\wh{G_2}\to\wh{G_1}$. In particular, 
it induces a homomorphism $Z_{[\wh{\pi}]}:Z(\wh{G_2})\to Z(\wh{G_1})$.
\end{Emp}

\begin{Emp} \label{E:tr}
{\bf Triple.} For the rest of  this subsection, we fix a triple $(G,H,[\eta])$, consisting of
a connected reductive group  $G$ over a field $E$, a quasi-split reductive group $H$ over 
$E$ of the same absolute rank as $G$, and  a $\Gm$-invariant $\wh{G}$-conjugacy class 
$[\eta]$ of embeddings $\wh{H}\hra\wh{G}$. 
\end{Emp}

\begin{Emp} \label{E:pair}
{\bf Properties of the triple $(G,H,[\eta])$.}  
(a) Every stable conjugacy class $[b]$ of embeddings of maximal tori $T\hra H$ 
defines a $\Gm$-invariant conjugacy class  $\wh{[b]}$ of embeddings $\wh{T}\hra\wh{H}$ 
(by \re{dual}), hence  a $\Gm$-invariant conjugacy class  
$\wh{[b]}_G:=[\eta]\circ [\wh{b}]$ of embeddings $\wh{T}\hra\wh{G}$,
thus a $\Gm$-invariant conjugacy class $[b]_G$ \label{bg} of embeddings $T_{\ov{E}}\hra G_{\ov{E}}$.

(b) There are exist canonical ($\Gm$-equivariant) embeddings 
$Z(\wh{G})\hra Z(\wh{H})$, $Z(G)\hra Z(H)$ and $W(H)\hra W(G)$.

To see it, fix a maximal torus $\wh{T}\subset\wh{H}$ and 
an embedding $\eta:\wh{H}\hra\wh{G}$ from $[\eta]$. Then $\wh{T}$ is a maximal torus of 
$\wh{G}$, hence the set of roots (therefore also of coroots) of 
$(\wh{H},\wh{T})$ is naturally a subset of that of $(\wh{G},\wh{T})$. 
It follows that $W(H)=W(\wh{H})$ is naturally a subgroup
of $W(\wh{G})=W(G)$. Also by \re{dual} (c), the set of roots of 
$(H,T)$ is naturally a subset of that of $(G,T)$. 
Since $Z(\wh{G})\subset\wh{T}$ (resp. $Z(G)\subset T$) is the intersection of
kernels of all roots of $(\wh{G},\wh{T})$ (resp. of $(G,T)$), and similarly for $Z(\wh{H})$
(resp. $Z(H)$), we get an embedding $Z(\wh{G})\hra Z(\wh{H})$ (resp. $Z(G)\hra Z(H)$).

(c) $[\eta]$ naturally gives rise to a conjugacy class  $[\ov{\eta}]$ of embeddings 
$\wh{H/Z(G)}\hra\wh{G^{\ad}}$. Namely, each $\eta:\wh{H}\hra\wh{G}$ from $[\eta]$,
has a unique lift  $\ov{\eta}:\wh{H/Z(G)}\hra\wh{G^{\ad}}$, and we denote by 
 $[\ov{\eta}]$ the corresponding conjugacy class. 

By (b), $[\eta]$ thus induces a homomorphism
$Z_{[\ov{\eta}]}:Z(\wh{G^{\ad}})^{\Gm}\to\pi_0(Z(\wh{H/Z(G)})^{\Gm})$. \label{zbareta}
\end{Emp}


\begin{Emp} \label{E:example}
{\bf Example.}
Any embedding $a:T\hra G$ of a maximal torus gives rise to a triple  \re{tr} with $H=T$ and 
$[\eta]=[\wh{a}]$. In particular, it gives rise to a $\Gm$-equivariant embedding 
$Z(\wh{G})\hra\wh{T}$, hence to a homomorphism $\pi_0(Z(\wh{G})^{\Gm})\to\pi_0(\wh{T}^{\Gm})$.
We will denote both of these maps by $Z_{[\wh{a}]}$. \label{zwha}

\end{Emp}

\begin{Emp} \label{E:exact}
{\bf Construction of an exact sequence.}
For each maximal torus $T\subset H$, we consider the exact sequence
\begin{equation} \label{Eq:exactT}
0\to\wh{T}\overset{\mu_T}{\lra}\wh{T^2/Z(G)}\overset{\nu_T}{\lra}\wh{T/Z(G)} \to 0,
\end{equation}
dual to the exact sequence
$0\to T/Z(G) \overset{\wh{\nu}_T}{\lra}T^2/Z(G)\overset{\wh{\mu}_T}{\lra}T \to 0$,
where $\wh{\nu}_{T}$ is the diagonal morphism and $\wh{\mu}_T([t_1,t_2]):=t_1/t_2$.

Since a center of a reductive group equals to the intersection of kernels of 
all roots, we conclude that $\mu_T(Z(\wh{H}))\subset Z(\wh{H^2/Z(G)})$,  
$\nu_T(Z(\wh{H^2/Z(G)}))\subset Z(\wh{H/Z(G)})$, and the induced sequence 
\begin{equation} \label{Eq:exact}
0\to Z(\wh{H})\overset{\mu_H}{\lra}Z(\wh{H^2/Z(G)})\overset{\nu_H}{\lra} Z(\wh{H/Z(G)})
\end{equation} \label{muh} 

is $\Gm$-equivariant and exact. Furthermore, since over $\ov{E}$ all maximal 
tori of $H$ are conjugate, sequence (\ref{Eq:exact}) is independent of the 
choice of $T$. 

Observe that the composition of $\mu_H$ 
with the projection $Z(\wh{H^2/Z(G)})\to Z(\wh{H^2})=Z(\wh{H})^2$ is the map
$z\mapsto(z,z^{-1})$.
\end{Emp}

\begin{Lem} \label{L:triple}
For each $i,j\in\{1,2,3\}$, let $\mu_{i,j}: Z(\wh{H})\hra Z(\wh{H^3/Z(G)})$
be the composition of $\mu_H$ and the embedding 
$Z(\wh{H^2/Z(G)})\hra Z(\wh{H^3/Z(G)})$ corresponding to the projection 
$H^3/Z(G)\to H^2/Z(G)$ to the $i$-th and the $j$-th factors. 
Then $\mu_{1,3}=\mu_{1,2}\mu_{2,3}$. 
\end{Lem}

\begin{proof}
Consider the projection $\la_{i,j}:T^3/Z(G)\to T$ given by the rule 
$\la_{i,j}([t_1,t_2,t_3])=t_j/t_i$. Since each $\mu_{i,j}$ is the restriction of 
$\la_{i,j}:\wh{T}\hra \wh{T^3/Z(G)})$ to $Z(\wh{H})$, the equality 
$\mu_{1,3}=\mu_{1,2}\mu_{2,3}$ follows from the equality $\la_{1,3}=\la_{1,2}\la_{2,3}$.
\end{proof}

\begin{Emp} \label{E:kainv} 
{\bf Construction of two homomorphisms.} 
(a) Consider a pair $[b_i]$ of stable conjugacy classes of embeddings of maximal tori $T_i\hra H$.
Then the  $[b_i]$'s give rise to a stable conjugacy class $[b_1,b_2]$ of embeddings 
$(T_1\times T_2)/Z(G)\hra H^2/Z(G)$, hence to a $\Gm$-invariant embedding 
\begin{equation} \label{E:moriota}
\iota([b_1],[b_2]):=Z_{\wh{[b_1,b_2]}}\circ\mu_H: 
Z(\wh{H})\hra \{(T_1\times T_2)/Z(G)\}\:\wh{}. 
\end{equation} \label{ib1b2}
In its turn, $\iota([b_1],[b_2])$ induces a homomorphism
\begin{equation} \label{E:mor1}
\ka([b_1],[b_2]):\pi_0(Z(\wh{H})^{\Gm})\to 
\pi_0((\{(T_1\times T_2)/Z(G)\}\:\wh{}\,)^{\Gm}),
\end{equation}\label{kab1b2}

(b) Assume that in the notation of (a), we have $T_1=T_2=T$ and $[b_1]_G=[b_2]_G$.
Then $\iota([b_1]_G,[b_2]_G)=\iota([b_1],[b_2])|_{Z(\wh{G})}$ factors through
$\mu_T:\wh{T}\hra\wh{T^2/Z(G)}$, hence the image of $\iota([b_1]_G,[b_2]_G)$ lies in 
$\Ker\nu_T$. Thus the composition $\nu_T\circ\iota([b_1],[b_2]):Z(\wh{H})\to\wh{T/Z(G)}$
factors through $Z(\wh{H})/Z(\wh{G})$ and induces a homomorphism
\begin{equation} \label{E:mor2}
\ka\left(\frac{[b_1]}{[b_2]}\right):\pi_0(Z(\wh{H})^{\Gm}/Z(\wh{G})^{\Gm})\to\pi_0(\wh{T/Z(G)}^{\Gm}).
\end{equation}\label{ka'b1b2}

(c) For each $i=1,2$, denote by  $[\ov{b}_i]$ the stable conjugacy class of embeddings 
$T/Z(G)\hra H/Z(G)$ induced by $[b_i]$. Then the composition of the projection 
$Z(\wh{H/Z(G)})\to Z(\wh{H})$ and $\nu_T\circ\iota([b_1],[b_2]):Z(\wh{H})\to\wh{T/Z(G)}$
equals the quotient $Z_{[\wh{\ov{b}}_1]}/Z_{[\wh{\ov{b}}_2]}:Z(\wh{H/Z(G)})\to\wh{T/Z(G)}$.

This gives the following description of the map $\ka\left(\frac{[b_1]}{[b_2]}\right)$.
For each representative
$\wt{s}\in Z(\wh{H/Z(G)})$ of $\ov{s}\in \pi_0(Z(\wh{H})^{\Gm}/Z(\wh{G})^{\Gm})$, the quotient
$Z_{[\wh{\ov{b}}_1]}(\wt{s})/Z_{[\wh{\ov{b}}_2]}(\wt{s})\in\wh{T/Z(G)}$ is $\Gm$-invariant, and
$\ka\left(\frac{[b_1]}{[b_2]}\right)(\ov{s})$ equals the class of 
$Z_{[\wh{\ov{b}}_1]}(\wt{s})/Z_{[\wh{\ov{b}}_2]}(\wt{s})$. 
(Note that $\wt{s}\in Z(\wh{H/Z(G)})$ is not always $\Gm$-invariant).
\end{Emp}

\begin{Lem} \label{L:kainv}
(a) Let $[b_1]$ and $[b_2]$ be as in \re{kainv}, and let 
$\check{s}\in \pi_0(Z(\wh{H})^{\Gm})$ be a representative of
$\ov{s}\in \pi_0(Z(\wh{H})^{\Gm}/Z(\wh{G})^{\Gm})$. Then the quotient 
$Z_{[\wh{b_1}]}(\check{s})/Z_{[\wh{b_2}]}(\check{s})\in\pi_0(\wh{T}^{\Gm})$ equals 
the image of $\ka\left(\frac{[b_1]}{[b_2]}\right)(\ov{s})$.

(b) Let $[b_1],[b_2]$ and $[b_3]$ be stable conjugacy classes of embeddings of 
maximal tori $T\hra H$ such that $[b_1]_G=[b_2]_G=[b_3]_G$. Then we have  
$\ka\left(\frac{[b_1]}{[b_3]}\right)=\ka\left(\frac{[b_1]}{[b_2]}\right)\ka\left(\frac{[b_2]}{[b_3]}\right)$.

\end{Lem}
\begin{proof}
Both assertions follow from the description of  $\ka\left(\frac{[b_i]}{[b_j]}\right)$,
given in \re{kainv} (c).
\end{proof}

\subsection{Endoscopic triples: basic properties}
\begin*
\vskip 8truept
\end*
 Let  $G$ be a connected reductive group over a field $E$.
In this subsection we give basic constructions and properties of 
endoscopic triples (compare \cite[$\S$7]{Ko1}).

\begin{Def} \label{D:end}
(a) A triple $\C{E}=(H,[\eta],\ov{s})$ \label{ehetas} consisting of  

- a quasi-split reductive group  $H$ over $E$ of the same absolute rank as $G$;

- a $\Gm$-invariant $\wh{G}$-conjugacy class  $[\eta]$ of embeddings $\wh{H}\hra\wh{G}$;

- an element $\ov{s}\in \pi_0(Z(\wh{H})^{\Gm}/Z(\wh{G})^{\Gm})$, where 
the $\Gm$-equivariant embedding \\ 
\indent\; $Z(\wh{G})\hra Z(\wh{H})$ is induced by $[\eta]$ 
(see \re{pair} (b)), 

\noindent is called  an  {\em endoscopic triple} \label{endtr} for $G$, if 
for a generic representative $s\in Z(\wh{H})^{\Gm}/Z(\wh{G})^{\Gm}$ of 
$\ov{s}$ we have $\eta(\wh{H})=\wh{G}_{\eta(s)}$ for all $\eta\in [\eta]$. 
Such a representative $s$ we will call {\em $\C{E}$-compatible}. \label{ecomp} 

(b) An endoscopic triple $\C{E}$ for $G$ is called {\em elliptic}, if the group
$Z(\wh{H})^{\Gm}/Z(\wh{G})^{\Gm}$ is finite.

(c) 
An {\em isomorphism} from an endoscopic triple $\C{E}_1=(H_1,[\eta_1],\ov{s}_1)$ to 
$\C{E}_2=(H_2,[\eta_2],\ov{s}_2)$ is an isomorphism $\al:H_1\isom H_2$ such 
that the corresponding element  
$\wh{[\al]}$ of \\
$\Isom(\wh{H_2},\wh{H_1})/\Inn(\wh{H_2})$ satisfies 
$[\eta_1]\circ\wh{[\al]}=[\eta_2]$ and $\ov{s}_1=\wh{[\al]}(\ov{s}_2)$.

(d) We denote by $\Aut(\C{E})$ the group of  automorphisms of $\C{E}$
and by $\La(\C{E})$ \label{lae} the quotient $\Aut(\C{E})/H^{\ad}(E)$.

(e) An endoscopic triple $\C{E}$ is called {\em split} (resp. {\em unramified}),  
if $H$ is a split (resp. unramified) group over $E$.
\end{Def}

\begin{Rem}
(a) When $E$ is a local non-archimedean field, our notion of an endoscopic triple 
is equivalent to the standard one (\cite[7.4]{Ko1}). Indeed, let $(H,\eta,s)$ be an endoscopic
triple in the sense of \cite[7.4]{Ko1}. Since $Z(\wh{H})^{\Gm}/Z(\wh{G})^{\Gm}$ is a subgroup
of finite index in  $[Z(\wh{H})/Z(\wh{G})]^{\Gm}$, condition  \cite[7.4.3]{Ko1} asserts 
that the image of $s\in Z(\wh{H})$ in $Z(\wh{H})/Z(\wh{G})$ belongs to  
$Z(\wh{H})^{\Gm}/Z(\wh{G})^{\Gm}$, thus defines an element 
$\ov{s}\in\pi_0(Z(\wh{H})^{\Gm}/Z(\wh{G})^{\Gm})$. Moreover, the map 
$(H,\eta,s)\mapsto (H,[\eta],\ov{s})$ defines an equivalence of categories between 
endoscopic triples in the sense of \cite[7.4]{Ko1} and those in our sense.

If $(H,\eta_1,s_1)$  and  $(H,\eta_2,s_2)$ are two endoscopic 
triples in the sense of \cite[7.4]{Ko1} such that $[\eta_1]=[\eta_2]$ 
and $\ov{s_1}=\ov{s_2}$, then they are canonically isomorphic.
Therefore we have chosen to deviate from the Kottwitz' notation and 
to identify these objects.

(b) When $E$ is a number field, then our notion of an endoscopic triple 
slightly differs from the standard one (\cite[7.4]{Ko1}). 
However in this paper we will only consider the case $G=G^{\ssc}$ (see Appendix B).
In this case, $Z(\wh{G})=1$, hence both notions coincide.
\end{Rem}

\begin{Rem}
Since an inner twisting $\varphi:G\to G'$ identifies the dual groups of $G$ and $G'$, 
every endoscopic triple for $G$ defines the one for $G'$.
\end{Rem}

\begin{Rem} \label{R:isom}
 In the notation of \rd{end} (c), the map $[\al]\mapsto\wh{[\al]}$ identifies 
$\Isom(\C{E}_1,\C{E}_2)/H_1^{\ad}(E)$ with the set of all $\Gm$-invariant
$\wh{[\al]}\in\Isom(\wh{H_2},\wh{H_1})/\Inn(\wh{H_2})$ such that 
$[\eta_1]\circ\wh{[\al]}=[\eta_2]$ and $\ov{s}_1=\wh{[\al]}(\ov{s}_2)$.

In particular, for each endoscopic triple $\C{E}=(H,[\eta],\ov{s})$, the map 
$[\al]\mapsto\wh{[\al]}$ identifies $\La(\C{E})$ with the set of 
$g\in\Out(\wh{H})^{\Gm}$ such that $g(\ov{s})=\ov{s}$ and $[\eta]\circ g=[\eta]$.
\end{Rem}


\begin{Emp} \label{E:fund}
{\bf Homomorphisms, corresponding to endoscopic triples.}
 To every endoscopic triple $\C{E}=(H,[\eta],\ov{s})$ for $G$, 
we are going to associate a homomorphism 
$\pi_{\C{E}}:\La(\C{E})\to \pi_0(Z(\wh{H/Z(G)})^{\Gm})$. \label{pie}
Moreover, the image of $\pi_{\C{E}}$ lies in the image of the canonical map 
$Z_{[\ov{\eta}]}:Z(\wh{G^{\ad}})^{\Gm}\to\pi_0(Z(\wh{H/Z(G)})^{\Gm})$ from \re{pair} (c).

To define $\pi_{\C{E}}$, fix $\eta\in [\eta]$, and identify $\wh{H}$ with its image 
$\eta(\wh{H})\subset\wh{G}$. Then $\eta$ lifts to an embedding  
$\wh{H/Z(G)}\hra\wh{G^{\ad}}=\wh{G}^{\ssc}$. For each $\al\in \La(\C{E})$, choose  
$g\in\Norm_{\wh{G}^{\ad}}(\wh{H})$ inducing $\wh{\al}\in\Out(\wh{H})$. Then $g$ normalizes
$\wh{H/Z(G)}\subset\wh{G}^{\ssc}$, hence it induces an element $\wt{\wh{\al}}\in\Out(\wh{H/Z(G)})$.
Moreover, $\wt{\wh{\al}}$ is $\Gm$-invariant and independent 
of the choices of $\eta$ and $g$. 

Furthermore, $\wt{\wh{\al}}$ is trivial on the kernel 
$\Ker[\wh{H/Z(G)}\to\wh{H}/Z(\wh{G})]$. Therefore for each 
$\wt{t}\in Z(\wh{H/Z(G)})$ the quotient 
$\wt{\wh{\al}}(\wt{t})/\wt{t}\in Z(\wh{H/Z(G)})$ depends only on the image 
$t\in Z(\wh{H})/Z(\wh{G})$ of $\wt{t}$. 
Hence the map $t\mapsto \wt{\wh{\al}}(\wt{t})/\wt{t}$ defines 
a $\Gm$-equivariant homomorphism $Z(\wh{H})/Z(\wh{G})\to Z(\wh{H/Z(G)})$,
which in its turn induces a homomorphism 
$\wt{\al}:\pi_0(Z(\wh{H})^{\Gm}/Z(\wh{G})^{\Gm})\to\pi_0(Z(\wh{H/Z(G)})^{\Gm})$.

We define $\pi_{\C{E}}$ by the formula 
$\pi_{\C{E}}(\al):=\wt{\al}(\ov{s})\in\pi_0(Z(\wh{H/Z(G)})^{\Gm})$.
Since $\wh{\al}(\ov{s})=\ov{s}$, we see that $\pi_{\C{E}}$ is a homomorphism, and the image of 
$\pi_{\C{E}}$ lies in the kernel 
$\Ker[\pi_0(Z(\wh{H/Z(G)})^{\Gm})\to\pi_0(Z(\wh{H})^{\Gm}/Z(\wh{G})^{\Gm})]$. 
As the last group coincides with the image of 
$Z_{[\ov{\eta}]}:Z(\wh{G^{\ad}})^{\Gm}\to\pi_0(Z(\wh{H/Z(G)})^{\Gm})$ (see \cite[Cor. 2.3]{Ko1}), 
we get the assertion.
\end{Emp}

\begin{Lem} \label{L:hom}
For every embedding of maximal torus $b:T\hra H$ the corresponding map
$Z_{\wh{[\ov{b}]}}:\pi_0(Z(\wh{H/Z(G)})^{\Gm})]\to\pi_0(\wh{T/Z(G)}^{\Gm})$ maps
$\pi_{\C{E}}(\al)$ to $\ka\left(\frac{\al([b])}{[b]}\right)(\ov{s})$.
\end{Lem}

\begin{proof}
The assertion follows immediately from  \re{kainv} (c).
\end{proof}

\begin{Not} \label{N:ze}
Let $\C{E}=(H,[\eta],\ov{s})$ be an endoscopic triple for $G$.

(a) Denote by $\Pi_{\C{E}}$ \label{Pie} the map  
sending a stable conjugacy class $[b]$ of embeddings of maximal tori 
$T\hra H$ to a pair consisting of an $E$-rational  conjugacy class $[b]_G$ 
(see \re{pair} (c)) of embeddings $T_{\ov{E}}\hra G_{\ov{E}}$ and an element 
$\ov{\ka}_{[b]}:=Z_{\wh{[b]}}(\ov{s})\in\pi_0(\wh{T}^{\Gm}/Z(\wh{G})^{\Gm})$. \label{ovkab}

(b) Denote by $Z(\C{E})\subset Z(\wh{G^{\ad}})^{\Gm}$ \label{ze} the preimage of 
$\im\pi_{\C{E}}\subset\pi_0(Z(\wh{H/Z(G)})^{\Gm})$ under $Z_{[\ov{\eta}]}$.

(c) For each pair $([a],\ov{\ka})\in\im\Pi_{\C{E}}$, we denote by
$S_{([a],\ov{\ka})}\subset \pi_0(\wh{T/Z(G)}^{\Gm})$ \label{saka} the subgroup 
consisting of elements $\ka\left(\frac{[b_1]}{[b_2]}\right)(\ov{s})$, where $b_1$ and $b_2$ run through
the set of embeddings $b:T\hra H$ such that $\Pi_{\C{E}}([b])= ([a],\ov{\ka})$.
We denote by $Z(\C{E},[a],\ov{\ka})\subset Z(\wh{G^{\ad}})^{\Gm}$ \label{zeaka} the preimage of 
$S_{([a],\ov{\ka})}$ under the map 
$Z_{[\wh{\ov{a}}]}:Z(\wh{G^{\ad}})^{\Gm}\to\pi_0(\wh{T/Z(G)}^{\Gm})$.

(d) We call embeddings of maximal tori $a:T\hra G$ and $b:T\hra H$ {\em compatible}, 
\label{comemb} 
if $a\in [b]_G$.
\end{Not}

\begin{Emp} \label{E:conend}
{\bf Endoscopic triples, corresponding to pairs.}
(a) Following Langlands (\cite[II, 4]{La}), we associate an endoscopic triple 
$\C{E}=\C{E}_{([a],\ka)}$ \label{eaka} for $G$ to each pair $([a],\ka)$, consisting of a stable conjugacy class 
of an embedding $a:T\hra G$  of maximal torus  and an element $\ka\in\wh{T}^{\Gm}/Z(\wh{G})^{\Gm}$.
Moreover, $\C{E}_{([a],\ka)}$ is elliptic, if $a(T)$ is an elliptic torus of $G$.

For the convenience of the reader, we will recall this important construction. 
Choose an element $\wh{a}:\wh{T}\hra\wh{G}$ 
of $\wh{[a]}$, and identify $\wh{T}$ with $\wh{a}(\wh{T})$. Choose a representative 
$\wt{\ka}\in\wh{T}^{\Gm}$ of $\ka$, and put $s:=\wh{a}(\wt{\ka})\in\wh{G}/Z(\wh{G})$,
$\wh{H}:=\wh{G}^0_{s}$, and let $[\eta]$ be the conjugacy class of the inclusion 
$\wh{H}\hra\wh{G}$. 

Then the image of $\rho_T:\Gm\to\Out(\wh{T})$ lies in 
$\Norm_{\Aut(\wh{G})}(\wh{T})/\wh{T}$, and $\rho_G:\Gm\to\Out(\wh{G})$
is the composition of 
$\rho_T:\Gm\to\Norm_{\Aut(\wh{G})}(\wh{T})/\wh{T}$
with the canonical homomorphism
$\Norm_{\Aut(\wh{G})}(\wh{T})/\wh{T}\hra \Aut(\wh{G})\to \Out(\wh{G})$.

As $\wt{\ka}$ belongs to $\wh{T}^{\Gm}$, the image of $\rho_T$ lies in 
$\Norm_{\Aut(\wh{G})}(\wh{T})_s/\wh{T}$, where $\Norm_{\Aut(\wh{G})}(\wh{T})_s$ is the stabilizer
of $s$ in $\Norm_{\Aut(\wh{G})}(\wh{T})$. 
Since $\Norm_{\Aut(\wh{G})}(\wh{T})_s\subset\Norm_{\Aut(\wh{G})}(\wh{H})$, we can compose $\rho_T:\Gm\to\Norm_{\Aut(\wh{G})}(\wh{T})_s/\wh{T}$ 
with the canonical homomorphism $\Norm_{\Aut(\wh{G})}(\wh{T})_s/\wh{T}
\to\Norm_{\Aut(\wh{G})}(\wh{H})/\wh{H}\subset \Out(\wh{H})$. We denote the composition
$\Gm\to \Out(\wh{H})$ by $\rho$, and let $H$ be the unique quasi-split group over $E$ 
corresponding to the pair $(\wh{H},\rho)$ (see \re{dual} (a)).

By construction,  $[\eta]$ is $\Gm$-equivariant, and $s$ belongs to 
$Z(\wh{H})^{\Gm}$. Denote by $\ov{s}\in\pi_0(Z(\wh{H})^{\Gm}/Z(\wh{G})^{\Gm})$ the class of $s$.
Then $(H,[\eta],\ov{s})$ is an endoscopic triple for $G$, independent  
of the choices of $\wh{a}$ and $\wt{\ka}$. We denote this endoscopic triple by $\C{E}_{([a],\ka)}$. 

(b) In the notation of (a), we have $([a],\ov{\ka})\in\im\Pi_{\C{E}}$.

To find $[b]\in\Pi_{\C{E}}^{-1}([a],\ov{\ka})$, we choose $\wh{a}$ and $\wt{\ka}$
as in (a) and identify $\wh{T}$ with $\wh{a}(\wh{T})$. By construction, there exists 
$\eta\in[\eta]$ such that $\wh{T}\subset\eta(\wh{H})$ and $\eta(s)=\wt{\ka}$. Let   $\wh{[b]}$ be
the conjugacy class of the inclusion $\wh{T}\hra\eta(\wh{H})\cong\wh{H}$. Then   $\wh{[b]}$ is
is $\Gm$-equivariant, hence it gives rise to a stable conjugacy class $[b]$ of 
embeddings $T\hra H$ (see \rl{steinberg} (b)), which belongs to 
$\Pi^{-1}_{\C{E}}([a],\ov{\ka})$.
\end{Emp}


\begin{Lem} \label{L:endosc}
Let $\C{E}=(H,[\eta],\ov{s})$ be an endoscopic triple for $G$,
 $b:T\hra H$ and $a:T\hra G$ embeddings of maximal tori, $\ov{\ka}$ 
an element of $\pi_0(\wh{T}^{\Gm}/Z(\wh{G})^{\Gm})$ such that
$\Pi_{\C{E}}([b])=([a],\ov{\ka})$.

(a) There exists a representative 
$\ka\in\wh{T}^{\Gm}/Z(\wh{G})^{\Gm}$ of $\ov{\ka}$ such that 
$\C{E}\cong\C{E}_{([a],\ka)}$. 

(b) The group $Z(\C{E})$ is contained in $Z(\C{E},[a],\ov{\ka})$. Moreover, 
if $b(T)\subset H$ is elliptic, then we have  
$Z(\C{E},[a],\ov{\ka})=Z(\C{E})$. In particular, 
$Z(\C{E},[a],\ov{\ka})=Z(\C{E})$ if $a(T)\subset G$ is elliptic.
\end{Lem}

\begin{proof}
(a) Let $s\in Z(\wh{H})^{\Gm}/Z(\wh{G})^{\Gm}$ be an 
$\C{E}$-compatible representative of $\ov{s}$, 
and let $\ka\in\wh{T}^{\Gm}/Z(\wh{G})^{\Gm}$
be the image of $s$ under the embedding 
$Z(\wh{H})^{\Gm}/Z(\wh{G})^{\Gm}\hra\wh{T}^{\Gm}/Z(\wh{G})^{\Gm}$, induced by $[\wh{b}]$.
Then $\ka$ is a representative of $\ov{\ka}$, and we claim that
$\C{E}_{([a],\ka)}\cong\C{E}$. To show it, choose an embedding 
$\wh{b}:\wh{T}\hra\wh{H}$ such that $\wh{b}\in\wh{[b]}$. Then for each 
$\eta\in[\eta]$, the composition 
$\wh{a}:=\eta\circ\wh{b}:\wh{T}\hra \wh{G}$ belongs to $\wh{[a]}$ and 
satisfies $\eta(\wh{H})=\wh{G}^0_{\wh{a}(\ka)}$. Finally, since
the conjugacy class of the embedding  $\wh{b}:\wh{T}\hra\wh{H}$ is 
$\Gm$-invariant, the homomorphism
$\rho_H:\Gm\to\Out(\wh{H})$ is induced by $\rho_T:\Gm\to\Out(\wh{T})$ as in
\re{conend} (a).

(b) Since for each $\al\in\La(\C{E})$, the stable conjugacy class
$\al([b])$ belongs to $\Pi^{-1}_{\C{E}}([a],\ov{\ka})$, the inclusion 
$Z(\C{E})\subset Z(\C{E},[a],\ov{\ka})$ follows from \rl{hom}.

Assume now that $b(T)\subset H$ is elliptic. Then for the second 
assertion, it will suffice to check that for each
$[b']\in\Pi_{\C{E}}^{-1}([a],\ov{\ka})$, 
there exists $\al\in\La(\C{E})$ such that $[b']=\al([b])$.

Embed $\wh{H}$ into $\wh{G}$ by means of an element of $[\eta]$, choose
embeddings $\wh{b}:\wh{T}\hra\wh{H}\subset \wh{G}$ and 
$\wh{b'}:\wh{T}\hra\wh{H}\subset \wh{G}$ from $\wh{[b]}$ and
$\wh{[b']}$, respectively, and identify $\wh{T}$ with its image 
$\wh{b}(\wh{T})\subset\wh{H}$. Replacing $\wh{b'}$ by its $\wh{H}$-conjugate, 
we may assume that $\wh{b}'(\wh{T})=\wh{T}$. 

Choose a representative 
$s\in Z(\wh{H})^{\Gm}$ of $\ov{s}$ such that $\wh{H}=\wh{G}^0_{s}$.
Since $\Pi_{\C{E}}([b'])=\Pi_{\C{E}}([b])$, there exists an element 
$g\in \wh{G}$ such that $\wh{b}'=g\wh{b}g^{-1}$ and 
$(g^{-1}s g)s^{-1}$ belongs to $(\wh{T}^{\Gm})^0 Z(\wh{G})^{\Gm}=
(Z(\wh{H})^{\Gm})^0 Z(\wh{G})^{\Gm}\subset Z(\wh{H})$. Therefore
$g^{-1}sg\in  Z(\wh{H})$, hence $g\wh{H}g^{-1}=\wh{H}$. 
Let $\ov{g}\in\Out(\wh{H})$  be the class of $g$. 
Since the conjugacy classes of $\wh{b}:\wh{T}\hra\wh{H}$ and 
$\wh{b'}:\wh{T}\hra\wh{H}$ are $\Gm$-invariant, we get from the equality
$\wh{b}'=g\wh{b}g^{-1}$ that  $\ov{g}$ is $\Gm$-invariant.
In other words, there exists $\al\in\La(\C{E})$
such that $\wh{\al}\in \Out(\wh{H})$ equals $\ov{g}$. Then  $[b']=\al([b])$,
as claimed.  
\end{proof} 

\begin{Lem} \label{L:qisog}
Let $\C{E}=(H,[\eta],\ov{s})$ be an endoscopic triple for $G$, and 
let $\pi:\wt{G}\to G$ be a quasi-isogeny.

(a) There exists a unique pair consisting of an endoscopic triple
$\wt{\C{E}}=(\wt{H},[\wt{\eta}],\ov{\wt{s}})$ for $\wt{G}$ and a stable conjugacy class of
quasi-isogenies $\pi':\wt{H}\to H$ such that 
$[\wt{\eta}]\circ[\wh{\pi}']=[\wh{\pi}]\circ [\eta]$
and $Z_{[\wh{\pi}']}$ maps $\ov{s}$ to $\ov{\wt{s}}$.

Furthermore, the endoscopic triple $\wt{\C{E}}$ satisfies the following properties.
 
(b) For an embedding of maximal torus $a:T\hra G$ and an  
element $\ka\in\wh{T}^{\Gm}/Z(\wh{G})^{\Gm}$, denote by $\wt{a}:\wt{T}\hra \wt{G}$ the lift of 
$a$, and by $\wt{\ka}\in\wh{\wt{T}}^{\Gm}/Z(\wh{\wt{G}})^{\Gm}$ the image of $\ka$.
If  $\C{E}\cong\C{E}_{([a],\ka)}$, then $\wt{\C{E}}\cong \C{E}_{([\wt{a}],\wt{\ka})}$.

(c) We have $\La(\C{E})=\La(\wt{\C{E}})$ and $Z(\wt{\C{E}})=Z(\C{E})$.

(d) In the notation of (b), the map sending an embedding of a maximal torus $b:T\hra H$ 
to its lift $\wt{b}:\wt{T}\hra \wt{H}$ induces a bijection 
between $\Pi_{\C{E}}^{-1}([a],\ov{\ka})$ and  
$\Pi_{\wt{\C{E}}}^{-1}([\wt{a}],\ov{\wt{\ka}})$. Moreover, 
for each two embeddings $b_1$ and $b_2$, we have
$\ka\left(\frac{[b_1]}{[b_2]}\right)(\ov{s})=\ka\left(\frac{[\wt{b_1}]}{[\wt{b_2}]}\right)(\ov{\wt{s}})$. 
\end{Lem}
\begin{proof}
(a) Choose $\eta\in[\eta]$, a quasi-isogeny 
$\wh{\pi}:\wh{G}\to \wh{\wt{G}}$ corresponding to $\pi$, and an 
$\C{E}$-compatible representative $s\in\wh{H}^{\Gm}/Z(\wh{G})^{\Gm}$
of $\ov{s}$. Identify $\wh{H}$ with $\eta(\wh{H})\subset\wh{G}$ 
and put $\wt{s}:=\wh{\pi}(s)\in\wh{\wt{G}}/Z(\wh{\wt{G}})$. 

Set $\wh{\wt{H}}:=(\wh{\wt{G}})_{\wt{s}}^0$. Then $\wh{\pi}$ induces a quasi-isogeny 
$\wh{\pi}':\wh{H}\to\wh{\wt{H}}$, and we have a canonical isomorphism 
$\wh{\wt{H}}\isom\wh{H}\times_{Z(\wh{G})}Z(\wh{\wt{G}})$. 
Since both homomorphisms $\rho_{H}:\Gm\to\Out(\wh{H})$ and  
$\rho_{\wt{G}}:\Gm\to\Aut(Z(\wh{\wt{G}}))$ induce the natural Galois action on  
$Z(\wh{G})$, their product defines a homomorphism $\rho:\Gm\to\Out(\wh{\wt{H}})$. 

Denote by $\wt{H}$ the quasi-split group over $E$, 
corresponding to the pair $(\wh{\wt{H}},\rho)$ (see \re{dual} (a)), by
$\ov{\wt{s}}\in\pi_0(Z(\wh{\wt{H}})^{\Gm}/Z(\wh{\wt{G}})^{\Gm})$ the  class of
$\wt{s}$, by $[\wt{\eta}]$ the conjugacy class of the inclusion  
$\wh{\wt{H}}\hra\wh{\wt{G}}$, and by $[\pi']$ the conjugacy class of quasi-isogenies
$\wt{H}\to H$ corresponding to $\wh{\pi}'$. Then $\wt{\C{E}}:=(\wt{H},[\wt{\eta}],\ov{\wt{s}})$
is an endoscopic triple for $\wt{G}$, and the pair  
$(\wt{\C{E}}, [\pi'])$ satisfy the required properties. The proof of the uniqueness is similar.

%

(b) follow immediately from the description of $\wt{\C{E}}$ in (a). 

(c) For each $\al\in\La(\C{E})$, the corresponding element 
$\wh{\al}\in\Out(\wh{H})$ is $\Gm$-invariant and is induced by  
an element of $G^{\ad}=\wt{G}^{\ad}$. Therefore $\al$ induces a unique element  
$\wt{\al}\in\La(\wt{\C{E}})$ and vice versa. Also we have an equality 
$\pi_{\C{E}}=\pi_{\wt{\C{E}}}$ implying that $Z(\C{E})=Z(\wt{\C{E}})$. 

(d) is clear. 
\end{proof} 

\subsection{Endoscopic triples: further properties.} \label{SS:estimates} 
\begin*
\vskip 8truept
\end*

Let $\C{E}=(H,[\eta],\ov{s})$ be an elliptic endoscopic triple for $G$, and 
$b:T\hra H$ an embedding of  a maximal torus. Put $([a],\ov{\ka}):=\Pi_{\C{E}}([b])$, and
denote by $S_{[b]}$ \label{sb} the subset 
$\left\{\ka\left(\frac{[b_1]}{[b]}\right)(\ov{s})\right\}_{[b_1]\in\Pi_{\C{E}}^{-1}([a],\ov{\ka})}$ of 
$\pi_0(\wh{T/Z(G)}^{\Gm})$.

The primary goal of this subsection is to prove the following result.

\begin{Prop} \label{P:group}
The subset $S_{[b]}\subset\pi_0(\wh{T/Z(G)}^{\Gm})$ is a subgroup.
\end{Prop}

This proposition has the following corollary.

\begin{Cor} \label{C:group}
For each $z\in Z(\C{E},[a],\ov{\ka})$, there exists  
$[b_1]\in\Pi_{\C{E}}^{-1}([a],\ov{\ka})$ such that 
$\ka\left(\frac{[b_1]}{[b]}\right)(\ov{s})=Z_{[\wh{\ov{a}}]}(z)$.
\end{Cor}
\begin{proof}[Proof of the corollary]
By definition, $Z_{[\wh{\ov{a}}]}(Z(\C{E},[a],\ov{\ka}))=S_{([a],\ov{\ka})}$, 
while $S_{([a],\ov{\ka})}$ is the group generated by $S_{[b]}$. 
Since $S_{[b]}$ itself is a group, we get that $S_{([a],\ov{\ka})}=S_{[b]}$, 
implying the assertion.
\end{proof}

To prove of the proposition, we will show first several results of independent 
interest, while the proof itself will be carried out in \re{proofgroup}.

As $\C{E}$ is elliptic, we denote 
$\ov{s}\in\pi_0(Z(\wh{H})^{\Gm}/Z(\wh{G})^{\Gm})=Z(\wh{H})^{\Gm}/Z(\wh{G})^{\Gm}$
simply by $s$. Note that  $\pi_0(\wh{G}_{s})$ is a $\Gm$-invariant subgroup of 
$\Out(\wh{H})$.

\begin{Lem} \label{L:ell}
Choose $\eta\in[\eta]$ and identify $\wh{H}$ with $\eta(\wh{H})\subset\wh{G}$.

(a) There exists a natural isomorphism $\La(\C{E})\cong \pi_0(\wh{G}_{s})^{\Gm}$, 
and a $\Gm$-equivariant injection 
$\iota:\pi_0(\wh{G}_{s})\hra Z(\wh{G}^{\ssc})$. Moreover, 
$\iota$ induces an isomorphism $\La(\C{E})\isom Z(\C{E})$.

(b) If $G$ is split, then the image of $\rho_H:\Gm\to \Out(\wh{H})$ lies in 
$\pi_0(\wh{G}_{s})$, and $\La(\C{E})$ is canonically isomorphic to 
$\pi_0(\wh{G}_{s})$.
\end{Lem}
\begin{proof}
For the proof we can replace $G$ by $G^{\ssc}$ and $\C{E}$ by the 
corresponding endoscopic triple, thus we can assume that $G=G^{\ssc}$.

(a) Using the fact that $\C{E}$ is elliptic, the first assertion follows 
from \rr{isom}. Next for each $g\in \wh{G}_s$, choose a representative 
$\wt{g}\in\wh{G}^{\ssc}$ of $g$ and a representative 
$\wt{s}\in \wh{G}^{\ssc}$ of $s$. Then the element
$(\wt{g}\wt{s}\wt{g}^{-1})\wt{s}^{-1}\in Z(\wh{G}^{\ssc})$ does not depend 
on the choices, and the map $g\mapsto (\wt{g}\wt{s}\wt{g}^{-1})\wt{s}^{-1}$
defines a homomorphism $\wt{\iota}:\wh{G}_{s}\to Z(\wh{G}^{\ssc})$.
Moreover, $g\in\Ker\wt{\iota}$ if and only if $\wt{g}\in(\wh{G}^{\ssc})_{\wt{s}}$.
By \cite[$\S$8]{St}, the last group is connected, therefore $\Ker\wt{\iota}=\wh{G}^0_{s}$. 
Hence $\wt{\iota}$ induces an embedding 
$\iota:\pi_0(\wh{G}_{s})\hra Z(\wh{G}^{\ssc})$, which is clearly $\Gm$-equivariant.
 
Since $\C{E}$ is elliptic, the group $Z(\wh{H/Z(G)})^{\Gm}$ is finite.
As $\pi_{\C{E}}$ is the composition of $\iota|_{\La(\C{E})}$ with the embedding
$Z(\wh{G}^{\ssc})^{\Gm}\hra Z(\wh{H/Z(G)})^{\Gm}=\pi_0(Z(\wh{H/Z(G)})^{\Gm})$,
the last assertion follows.

(b) The first assertion follows from the definition. 
Since embedding $\iota:\pi_0(\wh{G}_{s})\hra Z(\wh{G}^{\ssc})$ is $\Gm$-equivariant and
$\Gm$ acts trivially on $Z(\wh{G}^{\ssc})$, we conclude from (a) that 
that  $\La(\C{E})\cong \pi_0(\wh{G}_{s})^{\Gm}=\pi_0(\wh{G}_{s})$, as claimed.
\end{proof}

\begin{Emp} \label{E:dynkin}
{\bf Action of $Z(\wh{G}^{\ssc})$ on the extended Dynkin 
diagram $\wt{D}_{\wh{G}}$ of $\wh{G}$}.

Let $\wh{T}^{\ad}$ and $\wh{T}^{\ssc}$ be the  abstract Cartan subgroups 
of $\wh{G}^{\ad}$ and $\wh{G}^{\ssc}$, respectively, and   
let $C\subset X_*(\wh{T}^{\ad})\otimes \B{R}$ be the fundamental alcove.
For every $\mu\in X_*(\wh{T}^{\ad})$ there exists a unique element 
$w_{\mu}\in W^{aff}$ of the affine Weyl group of $\wh{G}$ such that
$w_{\mu}(C+\mu)=C$. Then the map $c_{\mu}:x\mapsto w_{\mu}(x+\mu)$ defines 
an affine automorphism of $C$, and hence an automorphism of  
$\wt{D}_{\wh{G}}$. Moreover, the map
$c:\mu\mapsto c_{\mu}$ is a homomorphism, and 
$X_*(\wh{T}^{\ssc})\subset\Ker c$. Thus $c$ induces an action of 
 $Z(\wh{G}^{\ssc})=X_*(\wh{T}^{\ad})/X_*(\wh{T}^{\ssc})$ on $C$, hence
on $\wt{D}_{\wh{G}}$.
\end{Emp}

\begin{Lem} \label{L:split}
Let $G$ be a split simple group. Then there exists a bijection 
$[\al]\mapsto\C{E}_{[\al]}=(H_{[\al]},s_{[\al]},[\eta]_{[\al]})$ between 
the set of $Z(\wh{G}^{\ssc})$-orbits of vertexes of $\wt{D}_{\wh{G}}$ and 
the set of isomorphism classes of split elliptic endoscopic triples for $G$.
 
Moreover, for each vertex $\al\in\wt{D}_{\wh{G}}$, the stabilizer
$Z(\wh{G}^{\ssc})_{\al}$ is canonically isomorphic to $\La(\C{E}_{[\al]})$, 
and the order $\ord(s_{[\al]})$ is equal to the coefficient of 
$\al$ in the reduced linear dependence
$\sum_{\al\in\wt{D}_{\wh{G}}} n_{\al}\al=0$.
\end{Lem} 

\begin{proof}
The set of isomorphisms of split endoscopic triples for $G$ is in
bijection with the set of conjugacy classes of semisimple elements 
$s\in\wh{G}^{\ad}$ such that $\wh{G}_s^0$ is semisimple, hence
with the set of $W(\wh{G})$-orbits of elements $s\in\wh{T}^{\ad}$ such that
$[X^*(\wh{T}^{\ad}):X^*(\wh{T}^{\ad})_s]<\infty$.

Note that $X_*(\wh{T}^{\ad})\otimes \B{R}/X_*(\wh{T}^{\ad})
(=\Hom (X_*(\wh{T}^{\ad}),\B{R}/\B{Z}))$ is naturally isomorphic to  
$\wh{T}(\B{C})^1=\Hom (X_*(\wh{T}^{\ad}),S^1)$.
This isomorphism induces a bijection between $W(\wh{G})$-orbits on
$\wh{T}(\B{C})^1$ and $(X_*(\wh{T}^{\ad})\rtimes W(\wh{G}))$-orbits on 
$X_*(\wh{T}^{\ad})\otimes \B{R}$. Since $C$ is a fundamental domain 
for the action of $W^{aff}=X_*(\wh{T}^{\ssc})\rtimes W(\wh{G})$, the
latter set coincides with the set of $Z(\wh{G}^{\ssc})$-orbits on $C$. 

For each $s\in \wh{T}(\B{C})^1$ and each representative
$\wt{s}\in C\subset X_*(\wh{T}^{\ad})\otimes \B{R}$ of $s$, the set of 
roots $\al$ of $\wh{G}$ such that $\al(s)=1$ are in bijection with the set of 
affine roots $\beta$ of $\wh{G}$ such that $\beta(\wt{s})=0$.
Therefore $[X^*(\wh{T}^{\ad}):X^*(\wh{T}^{\ad})_s]<\infty$ if and only if 
$\wt{s}$ is a vertex $C$, that is, a vertex of $\wt{D}_{\wh{G}}$.  
The last assertion is clear. 
\end{proof} 

The proof of the following result is done case-by-case.

\begin{Cl}  \label{C:estimate} 
Let $G$ be an absolutely simple group over $E$ such that $(G^*)^{\ssc}$ is not 
isomorphic to $\SL_n$. For every embedding
$a:T\hra G$ of a maximal torus, and  $\ov{\ka}\in\pi_0(\wh{T}^{\Gm}/Z(\wh{G})^{\Gm})$ such that
$([a],\ov{\ka})\in\im\Pi_{\C{E}}$, we have 
\[
[Z_{[\wh{\ov{a}}]}(Z(\C{E},[a],\ov{\ka})):Z_{[\wh{\ov{a}}]}(Z(\C{E}))]\leq 2.
\] 
\end{Cl}
\begin{proof}
Replacing $G$ by $G^*$, we can assume that $G$ is quasi-split.
Replacing $G$ by $G^{\ssc}$ and $\C{E}$ by the 
corresponding endoscopic triple, we can assume that $G=G^{\ssc}$.

Assume that our assertion is false, that is,  
$[Z_{[\wh{\ov{a}}]}(Z(\C{E},[a],\ov{\ka})):Z_{[\wh{\ov{a}}]}(Z(\C{E}))]> 2$. 
Since $Z(\C{E},[a],\ov{\ka})$ is a subgroup of $Z(\wh{G}^{\ssc})^{\Gm}$, we conclude that
$|Z(\wh{G}^{\ssc})^{\Gm}|> 2$. Therefore by the classification
of simple algebraic groups, we get that $G$ (hence also 
$\wh{G}$) is of type $A$, $D$ or $E_6$. Moreover, since the group $\Out(\wh{G})$ acts 
faithfully on $Z(\wh{G}^{\ssc})$, we see case-by-case that the assumption 
$|Z(\wh{G}^{\ssc})^{\Gm}|> 2$ implies that $G$ is split.

By our assumption, $G$ is not of type $A$, therefore  $|Z(\wh{G}^{\ssc})|\leq 4$.
Since by our assumption $[Z(\wh{G}^{\ssc}):Z(\C{E})]>2$, we get that $Z(\C{E})=1$, 
hence  $\La(\C{E})=1$ (see \rl{ell} (a)). It follows from \rl{ell} (b) that 
$\im\rho_H=1$, thus $\C{E}$ is split. Therefore by \rl{split},  
$\C{E}$ corresponds to a $Z(\wh{G}^{\ssc})$-orbit $[\al]\subset\wt{D}_{\wh{G}}$. 
Moreover, since $\La(\C{E})=1$, we get that
$|[\al]|=|Z(\wh{G}^{\ssc})|>2$.

Recall that $Z_{[\wh{\ov{a}}]}(Z(\C{E},[a],\ov{\ka}))$ 
consists of images of $s=s_{[\al]}$ under certain homomorphisms 
$\ka\left(\frac{[b_1]}{[b_2]}\right):Z(\wh{H})\to\pi_0(\wh{T/Z(G)}^{\Gm})$. Therefore every
$z\in Z_{[\wh{\ov{a}}]}(Z(\C{E},[a],\ov{\ka}))$ satisfies $z^{\ord(s_{[\al]})}=1$. 
Since $Z_{[\wh{\ov{a}}]}(Z(\C{E},[a],\ov{\ka}))\neq \{1\}$, we get
$(\ord(s_{[\al]}),|Z(\wh{G}^{\ssc})|)\neq 1$.
In particular, the orbit $[\al]$ is non-special, that is, consists of non-special
vertexes.

Now it is easy to get a contradiction. Indeed, in the case of $D_n$,
there are no non-special $\Aut(\wt{D}_{\wh{G}})$-orbits of cardinality greater 
than two, while in the case of $E_6$, there is only one such orbit.
However in this case we have $\ord(s_{[\al]})=2$ and $|Z(\wh{G}^{\ssc})|=3$, 
contradicting the assumption that $(\ord(s_{[\al]}),|Z(\wh{G}^{\ssc})|)\neq 1$.
\end{proof}

\begin{Emp} \label{E:restr}
{\bf Restriction of scalars.}
(a) Let $E'$ be a finite separable extension of $E$, $\Gm'=\Gal(\ov{E}/E')$,
$G'$ a reductive group over $E'$, and  $G=R_{E'/E}G'$.
Then $(\wh{G},\rho_G)$ has the following description.
First of all, $\wh{G}=\prod_{\si\in\Hom_E(E',\ov{E})} \wh{{}^{\si}G'}$,
where ${}^{\si}G'$ is a group over $\si(E')$ induced from $G'$ by $\si$.
Every $\tau\in\Gm$ induces a canonical element of 
$\rho_{G'}(\tau,\sigma)\in
\Isom(\wh{{}^{\si}G'},\wh{{}^{\tau\si}G'})/\Inn(\wh{{}^{\si}G'})$,
which coincides with $\rho_{{}^{\si}G'}(\tau)$ if 
$\tau\in\Gal(\ov{E}/\si(E'))$.
Then for every  $\tau\in\Gm$ we have 
$\rho_G(\tau)=\prod_{\si}\rho_{G'}(\tau,\sigma)$.

(b) There is a canonical isomorphism $Z(\wh{G})^{\Gm}\isom Z(\wh{G'})^{\Gm'}$.

(c) Every endoscopic triple
$\C{E}'=(H',[\eta'],\ov{s}')$ for $G'$ gives rise to an endoscopic triple 
$\C{E}=(H,[\eta],\ov{s})$ for $G=R_{E'/E}G'$, denoted by $R_{E'/E}\C{E}'$, 
defined as follows.  
$H=R_{E'/E}H'$, $[\eta]$ is a product
$\prod_{\si} [\eta'_{\si}]$, where $[\eta'_{\si}]$ is the conjugacy class
of embeddings $\wh{{}^{\si}H'}\hra\wh{{}^{\si}G'}$ induced by $[\eta']$,
and $\ov{s}$ is the preimage of $\ov{s'}$ under the canonical isomorphism
$Z(\wh{H})^{\Gm}/Z(\wh{G})^{\Gm}\isom Z(\wh{H'})^{\Gm'}/Z(\wh{G'})^{\Gm'}$.
\end{Emp}

\begin{Lem} \label{L:restr}
(a) For every embedding of a maximal torus $a:T\hra G=R_{E'/E}G'$, 
there exists an embedding $a':T'\hra G'$ of a maximal torus such 
that $T=R_{E'/E}T'$ and $a=R_{E'/E}a'$. Moreover, the map 
$[a']\mapsto [R_{E'/E}a']$ induces a bijection between the 
sets  of stable conjugacy classes of embeddings $T'\hra G'$ and $T\hra G$.

(b) For every endoscopic triple $\C{E}=(H,[\eta],\ov{s})$ for $G$, there 
exists a unique endoscopic triple $\C{E}'=(H',[\eta'],\ov{s}')$ 
for $G'$ such that $R_{E'/E}\C{E}'\cong\C{E}$.

(c) In the above notation, for each $([a'],\ov{\ka})\in \im\Pi_{\C{E}}$
the map $b'\mapsto b:=R_{E'/E}b'$ induces a bijection  
between $\Pi_{\C{E}'}^{-1}([a'],\ov{\ka})$ and  
$\Pi^{-1}_{\C{E}}([R_{E'/E}a'],\ov{\ka})$. Moreover, 
for each two embeddings $b'_1, b'_2:T'\hra H'$, 
the isomorphism $Z(\wh{H/Z(G)})^{\Gm}\isom Z(\wh{H'/Z(G')})^{\Gm'}$ maps
$\ka\left(\frac{[b_1]}{[b_2]}\right)(\ov{s})$ to $\ka\left(\frac{[b'_1]}{[b'_2]}\right)(\ov{s}')$.
\end{Lem}

\begin{proof}
(a) Note first that $G'$ is a direct factor of $G_{E'}$, and denote by $T'$
the image of the composition map $T_{E'}\overset{a}{\hra}G_{E'}\to G'$.
Then    $R_{E'/E}T'\cong T$, and the embedding $a':T'\hra G'$
satisfies $R_{E'/E}a'=a$. The second assertion follows from the first one.

(b) Choose $a:T\hra G$ and $\ka\in\wh{T}^{\Gm}/Z(\wh{G})^{\Gm}$ such that
$\C{E}\cong \C{E}_{([a],\ka)}$. Let $a':T'\hra G'$ be the embedding as in (a), 
and let $\ka'\in \wh{T'}^{\Gm'}/Z(\wh{G'})^{\Gm'}$ be the image of $\ka$.
Then $\C{E}':=\C{E}_{([a'],\ka')}$ satisfies
$R_{E'/E}\C{E}'\cong\C{E}$. The uniqueness is clear.

(c) follows immediately from (a).
\end{proof}

\begin{Emp} \label{E:proofgroup} 
\begin{proof}[Proof of \rp{group}]
The proof will be carried out in two steps. First we will treat the case 
$G=\SL_n$, and then reduce the general case to it.

{\bf Step 1: The case  $G=\SL_n$.}
In this case, we will describe all the objects involved explicitly.

First of all, there exists a divisor $m|n$, 
a cyclic Galois extension $K\subset\ov{E}$ of $E$ of degree $m$ such that $H$ is 
isomorphic to 
$(R_{K/E}\GL_{\frac{n}{m}})^1=\{g\in R_{K/E}\GL_{\frac{n}{m}}\,|\,N_{K/E}(\det g)=1\}$.
Next $s\in Z(\wh{H})\cong(\B{C}\m)^{\Gal(K/E)}/\B{C}\m$ is a class 
$[(\mu(\si))_{\si\in \Gal(K/E)}]$ for a certain isomorphism 
$\mu:\Gal(K/E)\isom \{z\in\B{C}\m\,|\,z^m=1\}$.  

Moreover, if we embed  $H\hra\SL_n$ by means of any $E$-linear isomorphism $K\isom E^m$, 
then for every embedding of maximal torus $b:T\hra H$, we get that
$[b]_G$ is the stable conjugacy class of the composition $T\overset{b}{\hra} H\hra\SL_n$.

Every maximal torus $T\subset H\subset SL_n$ is of the form 
$(\prod_{i=1}^l R_{K_i/E}\B{G}_m)^1$, where  
$(\prod_{i=1}^l R_{K_i/E}\B{G}_m)^1= 
\{t_i\in \prod_{i=1}^l R_{K_i/E}\B{G}_m\,|\,\prod_i N_{K_i/E}(t_i)=1\}$, 
for certain finite extensions $K_i/K$ with $\sum_i[K_i:K]=\frac{n}{m}$. 
Denote by $b$ the inclusion $T\hra H$. Then embeddings 
$b':T\hra H$ such that $[b']=[b]$ (resp. $[b']_G=[b]_G$) are in canonical bijection with
$K$-linear (resp. $E$-linear) algebra embeddings 
$\oplus_{i=1}^l K_i\hra\Mat_{\frac{n}{m}}(K)$.

In particular, we get a bijection $[\iota]\mapsto b_{[\iota]}$ from the set of   
$l$-tuples $\ov{\iota}=(\iota_1,\ldots,\iota_l)$ of $E$-algebra embeddings $\iota_i:K\hra K_i$
to that of stable conjugacy classes $[b']$ of embeddings $T\hra H$  such that $[b']_G=[b]_G$.
Therefore both sets are principal homogeneous spaces for the action of the group $\Gal(K/E)^l$.

The dual torus $\wh{T}$ is  
$[\prod_{i=1}^l (\B{C}\m)^{\Hom_E(K_i,\ov{E})}]/\B{C}\m$, and 
$(\wh{T}^{\Gm})^0$ is the image of the diagonal map $[(\B{C}\m)^l]/\B{C}\m\hra\wh{T}$.
Also $\wh{T/Z(G)}$ consists of 
$\{c_{i,\si_i}\}_{i,\si_i}\in\prod_i (\B{C}\m)^{\Hom_E(K_i,\ov{E})}$ such that 
$\prod_i\prod_{\si_i} c_{i,\si_i}=1$.

For each stable conjugacy class $[b_{\ov{\iota}}]$ of embeddings $T\hra H$, the 
corresponding embedding $Z_{[\wh{b_{\ov{\iota}}}]}:Z(\wh{H})\hra \wh{T}$ sends 
$s\in Z(\wh{H})$ to an element
$[(c_{i,\si_i})_{i,\si_i}]$, given by the rule $c_{i,\si_i}:=\mu(\si_i\circ\iota_i)$
(here $\si_i\circ\iota_i\in\Hom_E(K,\ov{E})=\Gal(K/E)$). 
When $\ov{\iota}$ is replaced by $\ov{\iota}\circ\ov{\tau}$ for certain 
$\ov{\tau}=(\tau_1,\ldots,\tau_l)$, then each $c_{i,\si_i}$ is multiplied by 
$\mu(\tau_i)$. It follows that 
$\ka\left(\frac{[b_{\ov{\iota}\circ\ov{\tau}}]}{[b_{\ov{\iota}}]}\right)(s)$ is the class  
$[(\mu(\tau_i))_{i,\si_i}]\in\wh{T/Z(G)}^{\Gm}$. In particular, the image of each 
$\ka\left(\frac{[b_{\ov{\iota}\circ\ov{\tau}}]}{[b_{\ov{\iota}}]}\right)(s)$ in  
$\pi_0(\wh{T}^{\Gm})$ is trivial. Thus (see \rl{kainv})
we get that $\Pi_{\C{E}}([b_{\ov{\iota}}])$ is independent of $\ov{\iota}$ (hence 
$\Pi_{\C{E}}([b_{\ov{\iota}}])=([a],\ov{\ka})$ for every $\ov{\iota}$).
As a result, the subset $S_{[b]}$ consists of classes of elements 
$[(\mu(\tau_i))_{i,\si_i}]$ where $\ov{\tau}$ runs through $\Gal(K/E)^l$.  
Hence $S_{[b]}$ is a group, as claimed.

{\bf Step 2: The general case.}
It follows from \rl{qisog} (d) that the subset $S_{[b]}$ will not change if we replace
$G$ by $G^{\ssc}$, $\C{E}$ by the corresponding endoscopic triple for 
$G^{\ssc}$ (see \rl{qisog} (a)) and $b$ by its lifting $b^{\ssc}:T^{\ssc}\hra G^{\ssc}$.
Thus we are reduced to the case when $G$ is semisimple and simply connected. 

Then $G$ is the product of its simple factors $G=\prod_{i} G_i$, and there exist 
embeddings of maximal tori $b_i:T_i\hra G_i$ such that $T=\prod_{i} T_i$ and  
$b=\prod_{i} b_i$. Then $S_{[b]}$ decomposes as a product
$\prod_{i} S_{[b_i]}\subset\prod_i \pi_0(\wh{T_i/Z(G_i)}^{\Gm})=
\pi_0(\wh{T/Z(G)}^{\Gm})$. Thus it will suffice to show that each
$S_{[b_i]} $ is a subgroup, thus reducing us to the case when 
$G$ is simple and simply connected. 

There exists  a finite separable extension $E'$ of $E$ and 
an absolutely simple simply connected algebraic group $G'$ over $E'$ 
such that $G\cong R_{E'/E}G'$. Using \rl{restr}, the subset  
$S_{[b]}$ will not change if we replace $E$ by $E'$, $G$ by $G'$, $\C{E}$ by
$\C{E}'$ and $[b]$ by $[b']$. Thus we can assume that $G$ is 
absolutely simple. Replacing $G$ by $G^*$, we can assume that $G$ is 
quasi-split, not isomorphic to $SL_n$.

For each embedding 
$b_1:T\hra H$ such that $\Pi_{\C{E}}([b_1])=\Pi_{\C{E}}([b])$ and each 
$\al\in\La(\C{E})$ we have  $\Pi_{\C{E}}(\al([b_1]))=\Pi_{\C{E}}([b])$.
By \rl{kainv} and \rl{hom} we have 
$\ka\left(\frac{\al([b_1])}{[b]}\right)(\ov{s})=
\ka\left(\frac{[b_1]}{[b]}\right)(\ov{s})Z_{[\wh{\ov{a}}]}(\pi_{\C{E}}(\al))$. 
In other words, $S_{[b]}$ is invariant under the multiplication 
by elements from $Z_{[\wh{\ov{a}}]}(Z(\C{E}))$.  

On the other hand, by definition, 
$S_{[b]}\subset Z_{[\wh{\ov{a}}]}(Z(\C{E},[a],\ov{\ka}))$. 
Since by \rcl{estimate}, we have   
$[Z_{[\wh{\ov{a}}]}(Z(\C{E},[a],\ov{\ka})):Z_{[\wh{\ov{a}}]}(Z(\C{E}))]\leq 2$, we conclude 
that $S_{[b]}$ is equal either to $Z_{[\wh{\ov{a}}]}(Z(\C{E},[a],\ov{\ka}))$ 
or to $Z_{[\wh{\ov{a}}]}(Z(\C{E}))$. Hence it is a group, as claimed.  
\end{proof}
\end{Emp}

\subsection{The case of local fields} \label{SS:local}
\begin*
\vskip 8truept
\end*

In this subsection we will apply the results from 1.1--1.3 to the case of endoscopic triples for  
a reductive group $G$ over a  local non-archimedean field $E$.

\begin{Emp} \label{E:TN}
{\bf Tate--Nakayama duality.} 
For every torus $T$ over $E$, Tate--Nakayama duality provides us 
with a functorial isomorphism $\C{D}_T:H^1(E,T)\isom\pi_0(\wh{T}^{\Gm})^D$
of finite abelian groups. 

Kottwitz showed that for every connected reductive group $G$ over $E$ its cohomology 
group $H^1(E,G)$ has a unique structure of a finite abelian group 
such that for every maximal torus $T\subset G$ the natural map 
$H^1(E,T)\to H^1(E,G)$ is a group homomorphism (see \cite[Thm 1.2]{Ko2}).
Moreover, there exists a group isomorphism  
$\C{D}_G:H^1(E,G)\isom\pi_0(Z(\wh{G})^{\Gm})^D$ \label{dg} such that for every 
maximal torus  $T$ of $G$, the embedding $T\hra G$
induces a commutative diagram:
\[
\CD
H^1(E,T)@>>> H^1(E,G)\\
@V{\C{D}_T}VV          @V{\C{D}_G}VV\\
\pi_0(\wh{T}^{\Gm})^D @>>> \pi_0(Z(\wh{G})^{\Gm})^D.
\endCD
\]
In particular, 
we have a canonical surjection
\[
\wh{T}^{\Gm}/Z(\wh{G})^{\Gm}\to \Coker[\pi_0(Z(\wh{G})^{\Gm})\to\pi_0(\wh{T}^{\Gm})]\isom
(\Ker\,[H^1(E,T)\to H^1(E,G)])^D.
\]
\end{Emp}

\begin{Rem} \label{R:Borovoi}
Borovoi \cite{Bo} showed that for every reductive group $G$ over $E$ 
there is a functorial group isomorphism $H^1(E,G)\isom (\pi_1(G)_{\Gm})_{\tor}$, where
$(\cdot)_{\tor}$ means for torsion. 
In particular, for every homomorphism of reductive groups $f:G_1\to G_2$, the induced map 
$H^1(E,G_1)\to H^1(E,G_2)$ is a group homomorphism as well. 
Now the existence of Kottwitz' isomorphism $\C{D}_G$ follows from the 
$\Gm$-equivariant isomorphism $\pi_1(G)\isom X^*(Z(\wh{G}))$.
\end{Rem}

\begin{Lem} \label{L:elltor}
Let  $a:T\hra G$ be an embedding of a maximal elliptic torus.
Then for every  inner twisting $\varphi:G\to G'$, 
there exists an embedding $a':T\hra G'$ stably conjugate to $a$.
\end{Lem}

\begin{proof}
By assumption, $T/Z(G)$ is anisotropic, therefore $\wh{T/Z(G)}^{\Gm}$ is 
finite. Hence the canonical map 
$\pi_0(Z(\wh{G^{\ad}})^{\Gm})= Z(\wh{G^{\ad}})^{\Gm}\hra{\wh{T/Z(G)}}^{\Gm}=
\pi_0({\wh{T/Z(G)}}^{\Gm})$ is injective. By duality, the canonical map 
$H^1(E,T/Z(G))\to H^1(E,G^{\ad})$ is surjective (see \re{TN}). This implies the assertion 
(see \re{inv} (b)). 
\end{proof}

\begin{Lem} \label{L:pairing}
Assume that we are in the situation of \re{coh} with $k=2$. Let $\C{E}=(H,[\eta],\ov{s})$ be 
an endoscopic triple for $G$, $b_1:T_1\hra H$ and $b_2:T_2\hra H$ are embeddings of maximal tori 
compatible with $a_1$ and $a_2$, respectively.
Then for every $z\in\pi_0(Z(\wh{G})^{\Gm})$ its image
$\wt{z}:=\ka([b_1],[b_2])(z)\in\pi_0(([(T_1\times T_2)/Z(G)]\;\wh{ }\,)^{\Gm})$ 
(see \re{kainv}) satisfies 
$\lan\ov{\inv}((a_1,a'_1);(a_2,a'_2),\wt{z}\ran=1$.
\end{Lem}

\begin{proof}
Recall that  $\wt{z}$ is the image of $\mu_G(z)\in\pi_0(Z(\wh{G^2/Z(G)})^{\Gm})$
under the map $Z_{\wh{[a_1,a_2]}}:\pi_0(Z(\wh{G^2/Z(G)})^{\Gm})\to
\pi_0(([(T_1\times T_2)/Z(G)]\;\wh{ }\,)^{\Gm})$.
Therefore by the commutative diagram of \re{TN}, we have 
\[
\lan\ov{\inv}((a_1,a'_1);(a_2,a'_2)),\wt{z}\ran=\lan\Dt(\inv(G,G')),\mu_G(z)\ran.
\] 
Moreover, the latter expression equals  
$\lan\inv(G,G'),\nu_G(\mu_G(z))\ran=\lan\inv(G,G'),0\ran=1$, as claimed.
\end{proof}

\begin{Not} \label{N:inv}
Let $\C{E}=(H,[\eta],\ov{s})$ be an endoscopic triple for $G$,
$\varphi:G\to G'$ an inner twisting, and $(a_i,a'_i,[b_i])$ 
be two triples consisting of  
stably conjugate embeddings of maximal tori $a_i:T_i\hra G$ and $a'_i:T_i\hra G'$, and 
stably conjugate classes $[b_i]$ of embeddings of maximal tori $T_i\hra H$, compatible with $a_i$.
To these data one associates elements 
$\ov{\inv}((a_1,a'_1);(a_2,a'_2))\in H^1(E,(T_1\times T_2)/Z(G))$ 
(see \re{coh}) and 
$\ka([b_1],[b_2])(\check{s})\in\pi_0((\{(T_1\times T_2)/Z(G)\}\:\wh{}\,)^{\Gm})$ 
for every representative $\check{s}\in\pi_0(Z(\wh{H})^{\Gm})$ of $\ov{s}$
(see \re{kainv}). By \rl{pairing}, the pairing
\[
\lan\frac{a_1,a'_1;[b_1]}{a_2,a'_2;[b_2]}\ran=\lan\frac{a_1,a'_1;[b_1]}{a_2,a'_2;[b_2]}\ran_{\C{E}}:=
\lan\ov{\inv}((a_1,a'_1);(a_2,a'_2)),\ka([b_1],[b_2])(\check{s})\ran\in\B{C}\m
\] \label{brac}
is independent of the choice of $\check{s}$.
\end{Not}

\begin{Rem} \label{R:delta1}
This invariant is essentially the term $\Dt_1$
of Langlands--Shelstad (\cite[(3.4)]{LS}).
\end{Rem}

\begin{Lem} \label{L:inv}
(a) For any three triples $(a_i,a'_i,[b_i]), (i=1,2,3)$, we have 
\[
\lan\frac{a_1,a'_1;[b_1]}{a_3,a'_3;[b_3]}\ran=\lan\frac{a_1,a'_1;[b_1]}{a_2,a'_2;[b_2]}\ran
\lan\frac{a_2,a'_2;[b_2]}{a_3,a'_3;[b_3]}\ran.
\]

(b) Assume that $T_1=T_2=T$, $a_1=a_2=a$ and $a'_1=a'_2=a'$. Then
$\lan\frac{a,a';[b_1]}{a,a';[b_2]}\ran=
\lan\ov{\inv}(a,a'),\ka\left(\frac{[b_1]}{[b_2]}\right)(\ov{s})\ran$.
If, moreover, $\ka\left(\frac{[b_1]}{[b_2]}\right)(\ov{s})=Z_{[\wh{\ov{a}}]}(z)$ 
for some $z\in Z(\wh{G^{\ad}})^{\Gm}$, then 
$\lan\frac{a,a';[b_1]}{a,a';[b_2]}\ran=\lan\inv(G,G'),z\ran$.

(c)  Assume that $T_1=T_2=T$ and $b_1=b_2=b$. Then 
\[
\lan\frac{a_1,a'_1;[b]}{a_2,a'_2;[b]}\ran=\lan\inv(a_1,a_2),\ov{\ka}_{[b]}\ran
\lan\inv(a'_1,a'_2),\ov{\ka}_{[b]}\ran^{-1}.
\]

(d) Assume that $\varphi=\Id_G$. Then 
$\lan\frac{a_1,a'_1;[b_1]}{a_2,a'_2;[b_2]}\ran=
\lan\inv(a_1,a'_1),\ov{\ka}_{[b_1]}\ran\lan\inv(a_2,a'_2),\ov{\ka}_{[b_2]}\ran^{-1}$.

(e) Let $\pi:\wt{G}\to G$ be a quasi-isogeny, $\wt{\C{E}}$ the 
endoscopic triple for $\wt{G}$ induced from $\C{E}$, 
$\wt{\varphi}:\wt{G}\to\wt{G}'$ the inner twisting
induced from  $\varphi$, and $\wt{a_i}$, 
$\wt{a_i}'$ and $[\wt{b_i}]$ the lifts of $a_i,a'_i$ and $[b_i]$, respectively.
Then $\lan\frac{\wt{a_1},\wt{a_1}';[\wt{b_1}]}{\wt{a_2},\wt{a_2}';[\wt{b_2}]}
\ran_{\wt{\C{E}}}=\lan\frac{a_1,a'_1;[b_1]}{a_2,a'_2;[b_2]}\ran_{\C{E}}$.

(f) Let $\varphi':G'\to G''$ be another inner twisting, and for each $i=1,2$, let
$a''_i:T\hra G''$ be a stable conjugate of $a_i$ and $a'_i$. Then
\[
\lan\frac{a_1,a''_1;[b_1]}{a_2,a''_2;[b_2]}\ran=\lan\frac{a_1,a'_1;[b_1]}{a_2,a'_2;[b_2]}\ran
\lan\frac{a'_1,a''_1;[b_1]}{a'_2,a''_2;[b_2]}\ran.
\]
\end{Lem}

\begin{proof}
All assertions follow from the functoriality of the Tate--Nakayama duality \re{TN}.

(a) By \rl{propinv} (c), $\lan\frac{a_i,a'_i;[b_i]}{a_j,a'_j;[b_j]}\ran$ equals the pairing
of $\ov{\inv}((a_1,a'_1);(a_2,a'_2);(a_3,a'_3))$ with the image of
$\ka([b_i],[b_j])(\check{s})$ in $\pi_0((\{(\prod_{i=1}^3 T_i)/Z(G)\}\:\wh{}\,)^{\Gm})$. 
Thus  the assertion follows from \rl{triple}.

(b) The first assertion follows from \rl{propinv} (d), while 
the second one follows from the fact that the image of $\ov{\inv}(a,a')$ in 
$H^1(E,G^{\ad})$ equals $\inv(G,G')$. 

(c) Since $\ka([b],[b])$ equals $\mu_T\circ Z_{[\wh{b}]}$, we get that
\[
\lan\frac{a_1,a'_1;[b]}{a_2,a'_2;[b]}\ran=
\lan \wh{\mu_T}(\ov{\inv}((a_1,a'_1);(a_2,a'_2)), Z_{[\wh{b}]}(\check{s})\ran.
\]
Using \rl{propinv} (d) and (e), we conclude that 
\[
\wh{\mu_T}(\ov{\inv}((a_1,a'_1);(a_2,a'_2)))=\inv(a_1,a_2)
\inv(a'_1,a'_2)^{-1},
\]
implying the assertion.

(d) By \rl{propinv} (b), $\ov{\inv}((a_1,a'_1);(a_2,a'_2))$
is the image of  $(\inv(a_1,a'_1),\inv(a_2,a'_2))$. 
Since the image of $\ka([b_1],[b_2])(\check{s})$ in 
$\pi_0(\wh{T_1}^{\Gm}/Z(\wh{G})^{\Gm})\times\pi_0(\wh{T_2}^{\Gm}/Z(\wh{G})^{\Gm})$
equals $(\ov{\ka}_{[b_1]},\ov{\ka}_{[b_2]}^{-1})$, the assertion follows.

(e) By \re{qisog} (d), 
$\ov{\inv}((a_1,a'_1);(a_2,a'_2))\in H^1(E,(T_1\times T_2)/Z(G))$ is the image of 
$\ov{\inv}((\wt{a_1},\wt{a_1}');(\wt{a_2},\wt{a_2}'))\in 
H^1(E,(\wt{T_1}\times \wt{T_2})/Z(\wt{G}))$.  Choose a representative 
$\check{s}\in\pi_0(Z(\wh{H})^{\Gm})$ of $\ov{s}$, and let 
$\wt{s}\in\pi_0(Z(\wh{\wt{H}})^{\Gm})$ be the image of $\check{s}$. Then 
$\ka([\wt{b_1}],[\wt{b_2}])(\wt{s})$ is the image of 
$\ka([b_1],[b_2])(\check{s})$, and the assertion follows.

(f) Follows from \rl{propinv} (f).
\end{proof}

\begin{Def} \label{D:adm}
Let $\C{E}$ be an endoscopic triple for $G$, and $([a],\ov{\ka})$ a 
pair belonging to  $\im \Pi_{\C{E}}$.
An inner twisting $\varphi:G\to G'$ is called {\em $\C{E}$-admissible} \label{eadm}
(resp. {\em $(\C{E},[a],\ov{\ka})$-admissible}), \label{eakaadm}
if the corresponding class
$\inv(G,G')\in H^1(E,G^{\ad})\cong (Z(\wh{G^{\ad}})^{\Gm})^D$ is orthogonal to 
$Z(\C{E})\subset Z(\wh{G^{\ad}})^{\Gm}$ 
(resp. orthogonal to $Z(\C{E},[a],\ov{\ka})\subset Z(\wh{G^{\ad}})^{\Gm}$).
\end{Def}

\begin{Lem} \label{L:adm}
(a) If $\varphi$ is $(\C{E},[a],\ov{\ka})$-admissible, then $\varphi$ is 
$\C{E}$-admissible. The converse is true, if $a(T)\subset G$ is elliptic.

(b) For each $i=1,2$, the function 
$[b_i]\mapsto \lan\frac{a_1,a'_1;[b_1]}{a_2,a'_2;[b_2]}\ran$ 
is constant on the fiber $\Pi_{\C{E}}^{-1}([a_i],\ov{\ka})$ 
if and only if $\varphi$ is $(\C{E},[a_i],\ov{\ka})$-admissible. 
\end{Lem}
\begin{proof}
(a) The assertion is a translation of \rl{endosc} (b).

(b) We will show the assertion for $i=1$, while the case $i=2$ is similar.
For each $[b_1], [b'_1]\in\Pi^{-1}_{\C{E}}([a_1],\ov{\ka})$, the quotient
 $\lan\frac{a_1,a'_1;[b'_1]}{a_2,a'_2;[b_2]}\ran\big/\lan\frac{a_1,a'_1;[b_1]}{a_2,a'_2;[b_2]}\ran$
equals the pairing $\lan\ov{\inv}(a_1,a'_1),\ka\left(\frac{[b_1]}{[b'_1]}\right)(\ov{s})\ran$ (use \rl{inv} (a) and (b)).
Since by definition, elements $\ka\left(\frac{[b_1]}{[b'_1]}\right)(\ov{s})$
run through $Z_{[\wh{\ov{a}}]}(Z(\C{E},[a_1],\ov{\ka}))$, our assertion follows from
the last assertion of \rl{inv} (b).
\end{proof}

\begin{Not} 
(a) Assume that in \rn{inv}, $\varphi$ is $(\C{E},[a_1],\ov{\ka}_{[b_1]})$-admissible
(resp. $(\C{E},[a_2],\ov{\ka}_{[b_2]})$-admissible), then we will denote 
$\lan\frac{a_1,a'_1;[b_1]}{a_2,a'_2;[b_2]}\ran$ by 
$\lan\frac{a_1,a'_1;\ov{\ka}_{[b_1]}}{a_2,a'_2;[b_2]}\ran$ 
(resp. $\lan\frac{a_1,a'_1;[b_1]}{a_2,a'_2;\ov{\ka}_{[b_2]}}\ran$). \label{brac2}
This notion is well defined by \rl{adm} (b).

(b) If in addition, $\varphi_X:X\to X'$ is an inner twisting induced by $\varphi$,
$a_1=a_x$ and $a'_1=a_{x'}$ (see \re{inner} (c) and \re{toremb} (c)), then we will denote 
$\lan\frac{a_1,a'_1;\ov{\ka}_{[b_1]}}{a_2,a'_2;[b_2]}\ran$ simply by 
$\lan\frac{x,x';\ov{\ka}_{[b_1]}}{a_2,a'_2;[b_2]}\ran$. \label{brac3}
\end{Not}

\subsection{Definitions of stability and equivalence}
\begin*
\vskip 8truept
\end*
In this subsection we will define the notions of stability and equivalence of invariant 
generalized functions.
 
\begin{Emp} \label{E:setup}
{\bf Set up.} (a) Let $G$ be a connected reductive group over a local non-archimedean
field $E$, $\om_G$ \label{omg} a non-zero translation invariant differential form on $G$ of 
the top degree, and $dg:=|\om_G|$ the corresponding
Haar measure on $G(E)$. Let $\om_{\C{G}}$ a non-zero translation 
invariant differential form on $\C{G}$ of the top degree such that 
the identification $T_1(G)=\C{G}=T_0(\C{G})$ identifies $\om_G|_{g=1}$ with 
$\om_{\C{G}}|_{x=0}$.

We also assume that $\C{G}$ contain regular semisimple elements, 
which hold automatically if the characteristic of $E$ is different from two.

(b) Let $(X,\om_X)$ be a pair consisting of a smooth algebraic variety $X$ 
over $E$, equipped with an $G^{\ad}$-action, and 
a non-vanishing $G$-invariant top degree differential form $\om_X$ \label{omx} on $X$, and let 
$dx:=|\om_X|$ be the corresponding measure on $X(E)$. 
We assume that $X^{\sr}\subset X$ is Zariski dense. This condition automatically implies that  
$X^{\sr}\subset X$ is open (see \rl{quotient} (a), below). In order to avoid dealing with algebraic 
spaces in \rss{l1}, we assume that $X$ is quasi-projective.

The results from this and the next subsections will be later used in 
two particular cases: $(X,\om_X)=(G,\om_G)$ and $(X,\om_X)=(\C{G},\om_{\C{G}})$ 
with the actions of $\Int G=G^{\ad}$ and $\Ad G=G^{\ad}$, respectively.
 
(c) Let $\varphi:G\to G'$ be an inner twisting. The inner twist $X'$ of $X$
 is smooth, and  the differential form $\om_{X'}:=(\varphi_X^{-1})^*(\om_X)$ \label{omx'} 
on $X'$ is $E$-rational and $G'^{\ad}$-invariant.

We also denote by $G^*$  \label{g*} the quasi-split inner twist of $G$, and by  $X^*$ the corresponding twist of 
$X$.

(d) Let $\C{E}=(H,[\eta],\ov{s})$ be an endoscopic triple for $G$.
\end{Emp}

%


\begin{Not} \label{N:meas}
(a) let $C_c^{\infty}(X(E))$ \label{ccinfty} be the space of locally constant functions 
on $X(E)$ with compact support, and let 
$\C{S}(X(E))$ \label{sxe} be the space of locally constant measures
on $X(E)$ with compact support, that is, measures of the form $\phi=fdx$, where
$f\in C_c^{\infty}(X(E))$. 

Denote by $\C{D}(X(E))=(\C{S}(X(E))^*)^{G(E)}$ \label{dxe} the space
of $G(E)$-invariant linear functionals on $\C{S}(X(E))$, which we call 
{\em (invariant)  generalized functions}.

(b) Let $U\subset X(E)$ be an open and closed subset. 
For each $\phi=fdx\in\C{S}(X(E))$, put $\phi|_U:=(f|_U)dx\in\C{S}(X(E))$. \label{phiu}
Moreover, if $U$ is $G(E)$-invariant, 
then for each $F\in\C{D}(X(E))$, the generalized function $F|_U$ given by the formula
$F|_U(\phi):=F(\phi|_U)$ \label{fu} belongs to $\C{D}(X(E))$.

(c) For every smooth  morphism $\pi:X_1\to X_2$, the integration along fibers $\pi_!$ maps
$\C{S}(X_1(E))$ to $\C{S}(X_2(E))$. 
Moreover, if $\pi$ is $G$-equivariant, then the dual of $\pi_!$ induces
a map $\pi^*:\C{D}(X_2(E))\to\C{D}(X_1(E))$.
\end{Not}

\begin{Rem} \label{R:notation}
The map $\phi\mapsto \frac{\phi}{dx}$ identifies $\C{S}(X(E))$ with 
$C_c^{\infty}(X(E))$, hence the space $\C{D}(X(E))$ with the space of 
invariant distributions on $X(E)$. Below we list several reasons why 
$\C{S}(X(E))$ and $\C{D}(X(E))$ are more convenient to work with.

i) The space $C_c^{\infty}(X(E))$ is not functorial with respect to non-proper
maps.

ii) Characters of admissible representations of $G(E)$ belong to $\C{D}(X(E))$.

iii) Orbital integrals behave better (see \rr{orbital} below).
\end{Rem}

\begin{Not} \label{N:orbital}
For each $x\in X^{\sr}(E)$ and $\ov{\xi}\in\pi_0(\wh{G_{x}}^{\Gm}/Z(\wh{G})^{\Gm})$,

(i) fix an invariant measure $dg_{x}$ on $G_{x}(E)$ such that
the total measure of the maximal compact subgroup of $G_x(E)$ is $1$,
and define orbital integral $O_{x}\in \C{D}(X(E))$ \label{ox} by the formula
\[
O_{x}(\phi):=\int_{G(E)/G_{x}(E)}\left(\frac{\phi}{dx}\right)(g(x))\frac{dg}{dg_x}
\]
for each $\phi\in\C{S}(X(E))$. 

(ii) denote by $O^{\ov{\xi}}_{x}\in \C{D}(X(E))$ \label{oxix} the sum 
$\sum_{x'}\lan\inv(x,x'),\ov{\xi}\ran O_{x'}$,
taken over a set of representatives $x'\in X(E) $ of conjugacy classes stably 
conjugate to $x$.

(iii) When $\ov{\xi}=1$, we will write $SO_{x}$ \label{sox} instead of  
$O^{\ov{\xi}}_{x}$. More generally, for each  $x^*\in (X^*)^{\sr}(E)$
we define $SO_{x^*}\in \C{D}(X(E))$ be zero unless there exists a stable 
conjugate $x\in X(E)$ of $x^*$, in which case, $SO_{x^*}:=SO_{x}$ 
(compare \rco{steinberg}).
\end{Not}

\begin{Rem} \label{R:orbital}
If $(X,\om_X)$ is either $(G,\om_G)$ or $(\C{G},\om_{\C{G}})$, then
measure $dx$ is induced by $dg$, and orbital integrals $O_{x}$ are independent of a choice of $dg$.
\end{Rem}

\begin{Def} \label{D:stab}

(i) A measure $\phi\in\C{S}(X(E))$ is called {\em $\C{E}$-unstable} \label{eunstable} if 
$O^{\ov{\xi}}_{x}(\phi)=0$ for all  $x\in X^{\sr}(E)$ and 
$\ov{\xi}\in\pi_0(\wh{G_{x}}^{\Gm}/Z(\wh{G})^{\Gm})$
such that $([a_{x}],\ov{\xi})\in \im\Pi_{\C{E}}$.

(ii) A generalized function $F\in\C{D}(X(E))$ is called {\em $\C{E}$-stable} \label{estable} if 
$F(\phi)=0$ for all $\C{E}$-unstable $\phi\in\C{S}(X(E))$.
\end{Def}

\begin{Rem} \label{R:sreg}
 Denote by $\C{D}^0(X(E))\subset\C{D}(X(E))$ \label{d0xe} the closure of the linear span of
$\{O_x\}_{x\in X^{\sr}(E)}$. Our notion of $\C{E}$-stability 
(and of $(a,a';[b])$-equivalence below) seems to be "correct" only for generalized 
functions belonging to $\C{D}^0(X(E))$. 

However, all generalized functions considered in this paper belong to $\C{D}^0(X(E))$. 
Indeed, if $(X,\om_X)$ is either $(G,\om_G)$ or $(\C{G},\om_{\C{G}})$, 
it follows from the results of Harish-Chandra \cite[Thm. 3.1]{HC2}
(at least when the characteristic of $E$ is zero) that 
$\C{D}^0(X(E))=\C{D}(X(E))$ (see also Remarks \ref{R:l1} and \ref{R:SO} below).
\end{Rem}

\begin{Not} \label{N:kaint}
Fix a triple $(a,a';[b])$, consisting of stably conjugate 
embeddings of maximal tori $a:T\hra G, a':T\hra G'$, and a stable conjugacy class $[b]$ of  
embeddings $T\hra H$, compatible with $a$ and $a'$.
For every  $\phi'\in\C{S}(X'(E))$, $x\in X^{\sr}(E)$ and embedding $c:G_x\hra H$ compatible with 
$a_x:G_x\hra G$, we define 
\[
(O^{[c]}_{x})_{(a,a';[b])}:=\sum_{x'}\lan \frac{a_{x},a_{x'};[c]}{a,a';[b]}
\ran O_{x'}\in \C{D}(X'(E)),
\]
where the sum is taken over a set of representatives $x'\in X'(E)$ 
of conjugacy classes stably conjugate to $x$.
\end{Not}

\begin{Rem} \label{R:kaint}
If $\varphi=\Id_G$ (and $a'=a$), then 
$(O^{[c]}_{x})_{(a,a;[b])}=O^{\ov{\ka}_{[c]}}_{x}$ (by \rl{inv} (d)).
In general,  $(O^{[c]}_{x})_{(a,a';[b])}$ vanishes unless there exists 
a stable conjugate $x'\in X'(E)$ of $x$, in which case, 
$(O^{[c]}_{x})_{(a,a';[b])}=\lan \frac{a_{x},a_{x'};[c]}{a,a';[b]}\ran 
O^{\ov{\ka}_{[c]}}_{x'}$ (by \rl{inv} (a) and (c)). 
\end{Rem}

\begin{Def} \label{D:main}
Let $(a,a';[b])$ be as in \rn{kaint}. 
By \rl{steinberg} (b), we can choose a stably conjugate embedding $a^*:T\hra G^*$ of $a$ and $a'$.

(a) Measures $\phi\in\C{S}(X(E))$ and $\phi'\in\C{S}(X'(E))$ are called 
{\em $(a,a';[b])$-indistinguishable}, \label{aa'bin} if for each $x^*\in (X^*)^{\sr}(E)$ and each embedding 
$c:G^*_{x^*}\hra H$ compatible with $a_{x^*}:G^*_{x^*}\hra G^*$, we have 
\begin{equation} \label{E:equiv}
(O^{[c]}_{x^*})_{(a^*,a;[b])}(\phi)=(O^{[c]}_{x^*})_{(a^*,a';[b])}(\phi').
\end{equation}

(b) Generalized functions $F\in\C{D}(X(E))$ and $F'\in\C{D}(X'(E))$ are 
called  {\em $(a,a';[b])$-equivalent} \label{aa'beq} if $F(\phi)=F'(\phi')$ for every  
pair of $(a,a';[b])$-indistinguishable 
measures $\phi\in\C{S}(X(E))$ and $\phi'\in\C{S}(X'(E))$.
\end{Def}

\begin{Lem} \label{L:indist}
 Measures $\phi\in\C{S}(X(E))$ and $\phi'\in\C{S}(X'(E))$ are 
$(a,a';[b])$-indistin-guishable if and only if the following conditions are
satisfied:

$(i)$ For each $x\in X^{\sr}(E)$ and
$\ov{\xi}\in\pi_0(\wh{G_{x}}^{\Gm}/Z(\wh{G})^{\Gm})$
such that $([a_{x}],\ov{\xi})\in \im\Pi_{\C{E}}$ and $x$ does not have
a stable conjugate element in $X'(E)$, we have
$O^{\ov{\xi}}_{x}(\phi)=0$.

$(ii)$ Condition $(i)$ holds if $x,X,G,\phi$ are replaced by 
$x',X',G',\phi'$.

$(iii)$  For each stable conjugate $x\in X^{\sr}(E)$ and $x'\in X'^{\sr}(E)$
and each $\ov{\xi}\in\pi_0(\wh{G_{x}}^{\Gm}/Z(\wh{G})^{\Gm})$
such that $([a_{x}],\ov{\xi})\in \im\Pi_{\C{E}}$ we have

$(iii)'$ $O^{\ov{\xi}}_{x}(\phi)=O^{\ov{\xi}}_{x'}(\phi')=0$, if $\varphi$ is 
not $(\C{E},[a_x],\ov{\xi})$-admissible, and 

$(iii)''$ $O^{\ov{\xi}}_{x}(\phi)=\lan \frac{x,x';\ov{\xi}}{a,a';[b]}\ran
O^{\ov{\xi}}_{x'}(\phi')$, if $\varphi$ is 
$(\C{E},[a_x],\ov{\xi})$-admissible.
\end{Lem}

\begin{proof}
Fix $x^*\in (X^*)^{\sr}(E)$ and $\ov{\xi}\in\pi_0(\wh{G^*_{x^*}}^{\Gm}/Z(\wh{G})^{\Gm})$
such that $([a_{x^*}],\ov{\xi})\in \im\Pi_{\C{E}}$. Using \rr{kaint}, we see that equalities 
(\ref{E:equiv}) for all $[c]\in\Pi^{-1}_{\C{E}}([a_{x^*}],\ov{\xi})$ are equivalent to equalities

-  $0=0$, if there are no stable conjugates of $x^*$ neither in $X(E)$ nor in
$X'(E)$;

- $O^{\ov{\xi}}_{x}(\phi)=0$, if there exists a stable conjugate $x\in X(E)$ 
of $x^*$ but there is not such conjugate in $X'(E)$;

- $O^{\ov{\xi}}_{x'}(\phi')=0$, if there exists a stable conjugate 
$x'\in X'(E)$ of $x^*$ but there is not such conjugate in $X(E)$;

-  $O^{\ov{\xi}}_{x}(\phi)=\lan \frac{a_{x},a_{x'};[c]}{a,a';[b]}\ran
O^{\ov{\xi}}_{x'}(\phi')$ for all $[c]\in\Pi^{-1}_{\C{E}}([a_x],\ov{\xi})$, 
if there exist stable conjugates $x\in X(E)$ and $x'\in X'(E)$ of $x^*$ (use \rl{inv} (f)).

Moreover, by \rl{adm} (b), the last equalities are equivalent to the equalities
$(iii)'$ and $(iii)''$. Now the assertion follows from \rco{steinberg}.
\end{proof}


\begin{Cor} \label{C:end}
(a) The notion of $(a,a';[b])$-equivalence is independent of the choice of $a^*$.

(b) Every two $(a,a';[b])$-equivalent generalized functions $F$ and $F'$ are 
$\C{E}$-stable. 

(c) Assume that  $\varphi$ is not $\C{E}$-admissible. Then every $\C{E}$-stable 
$F\in\C{D}(X(E))$ and $F'\in\C{D}(X'(E))$ are $(a,a';[b])$-equivalent.

(d) Assume that $a(T)$ is elliptic. If $F\in\C{D}(X(E))$ and $F'\in\C{D}(X'(E))$ 
are  $(a,a';[b])$-equivalent, then they are $(a,a';[b'])$-equivalent
for all $b':T\hra H$ such that $\Pi_{\C{E}}([b'])=\Pi_{\C{E}}([b])$.
\end{Cor}

\begin{proof} 
 All assertions follow almost immediately from  \rl{indist}. 

(a) is clear.

(b) By duality, we have to check that every $\C{E}$-unstable measures
$\phi\in \C{S}(X(E))$ and  $\phi'\in \C{S}(X'(E))$ are 
$(a,a';[b])$-indistinguishable, which is clear.

(c) By duality, we have to check that every  $(a,a';[b])$-indistinguishable 
$\phi\in\C{S}(X(E))$ and $\phi'\in\C{S}(X'(E))$
are $\C{E}$-unstable. Hence the assertion follows from the first assertion of
\rl{adm} (a).

(d) When $\varphi$ is not $\C{E}$-admissible, the assertion was proved in (c).
When $\varphi$ is $\C{E}$-admissible, the assertion follows from \rl{adm}.
\end{proof}

\begin{Cor} \label{C:comp}
Let $\pi:\wt{G}\to G$ be a quasi-isogeny, and let $\wt{\varphi}:\wt{G}\to\wt{G}', 
\wt{\C{E}}, (\wt{a},\wt{a}';[\wt{b}])$ be the corresponding objects for $\wt{G}$.
Generalized functions $F\in\C{D}(X(E))$ and $F'\in\C{D}(X'(E))$ are 
$(\wt{a},\wt{a}';[\wt{b}])$-equivalent if and only if they are $(a,a';[b])$-equivalent. 
\end{Cor}
\begin{proof}
By duality, we have to show that measures $\phi\in\C{S}(X(E))$ and $\phi'\in\C{S}(X'(E))$
are $(a,a';[b])$-indistinguishable if and only if they are
$(\wt{a},\wt{a}';[\wt{b}])$-indistinguishable. It follows from \rl{qisog}, 
that for each $x\in X^{\sr}(E)$ and $\ov{\wt{\xi}}\in
\pi_0(\wh{\wt{G}_{x}}^{\Gm}/Z(\wh{\wt{G}})^{\Gm})$
such that $([\wt{a}_{x}],\ov{\wt{\xi}})\in \im\Pi_{\wt{\C{E}}}$,
there exists $\ov{\xi}\in\pi_0(\wh{G_{x}}^{\Gm}/Z(\wh{G})^{\Gm})$
such that $(a_{x},\ov{\xi})\in \im\Pi_{\C{E}}$ and $\ov{\wt{\xi}}$ 
is the image of $\ov{\xi}$. 
Moreover, we have $O_x^{\ov{\wt{\xi}}}=c O_x^{\ov{{\xi}}}$,
where $c\in\B{C}\m$ is such that measure $\frac{d\wt{g}}{d\wt{g}_x}$ on
$(\wt{G}/\wt{G}_x)(E)=(G/G_x)(E)$ equals  $c\frac{dg}{d{g}_x}$.
Therefore the assertion follows from Lemmas \ref{L:indist}, \ref{L:inv} (e), and \ref{L:qisog} (c), (d).  
\end{proof}

\begin{Def} \label{D:def3}
Let $a:T\hra G$ and $a':T\hra G'$ be stable conjugate embeddings of 
maximal elliptic 
tori, and let $\ka$ be an element of $\wh{T}^{\Gm}/Z(\wh{G})^{\Gm}$ such that
$([a],\ka)\in\im\Pi_{\C{E}}$. We say that $F\in\C{D}(X(E))$ and 
$F'\in\C{D}(X'(E))$ are {\em $(a,a';\ka)$-equivalent} \label{aa'kaeq} if they are $(a,a';[b])$-equivalent
for some or, equivalently, for all  
$[b]\in\Pi_{\C{E}}^{-1}([a],\ka)$ (see \rco{end} (d)).
\end{Def}

\begin{Lem} \label{L:pullback}
(a) Let $\pi:X_1\to X_2$ be a smooth $G$-equivariant morphism. For every $\C{E}$-unstable
$\phi\in\C{S}(X_1(E))$, its push-forward $\pi_!(\phi)\in\C{S}(X_2(E))$ is $\C{E}$-unstable.

(b) Let $\varphi:G\to G'$ be an inner twisting, and  $\pi':X'_1\to X'_2$ the corresponding
inner twisting of $\pi$.  For every  $(a,a';[b])$-indistinguishable
$\phi\in\C{S}(X_1(E))$ and $\phi'\in\C{S}(X'_1(E))$, their  push-forwards
$\pi_!(\phi)\in\C{S}(X_1(E))$ and $\pi'_!(\phi')\in\C{S}(X'_1(E))$ are
$(a,a';[b])$-indistinguishable.
\end{Lem}

\begin{proof}
Recall that $\om_{X_1}$ and $\om_{X_2}$ are global nowhere vanishing sections 
of sheaves of top degree differential forms $\Om^{\dim\,X_1}_{X_1}$ and $\Om^{\dim\,X_2}_{X_2}$, 
respectively. 
Therefore $\om_{X_1}\otimes \pi^*(\om_{X_2})^{-1}$  is a global nowhere vanishing 
section of $\Om^{\dim\,X_1}_{X_1}\otimes\pi^*(\Om^{\dim\,X_2}_{X_2})^{-1}$, which induces a
measure, denoted by $dy:=\frac{dx_1}{dx_2}$, on all fibers of 
$\pi(E):X_1(E)\to X_2(E)$.

For each $x\in X_2^{\sr}(E)$ and 
$y\in\pi(E)^{-1}(x)$, we have $y\in X_1^{\sr}(E)$ and 
$G_x=G_y$. Moreover, for all $\ov{\xi}\in\pi_0(\wh{G_x}^{\Gm}/Z(\wh{G})^{\Gm})$, we have 
$O_x^{\ov{\xi}}(\pi_!(\phi))=\int_{\pi(E)^{-1}(x)} O_y^{\ov{\xi}} (\phi)dy$.
From this the assertion follows.
\end{proof}

\rl{pullback} has the following corollary.

\begin{Cor} \label{C:pullback}
(a) Let $\pi:X_1\to X_2$ be a smooth $G$-equivariant morphism. For every  $\C{E}$-stable
$F\in\C{D}(X_2(E))$, its pullback  $\pi^*(F)\in\C{D}(X_1(E))$ is $\C{E}$-stable.

(b) Let $\varphi:G\to G'$ be an inner twisting, and  $\pi':X'_1\to X'_2$ the inner 
twist of $\pi$. For every $(a,a';[b])$-equivalent $F\in\C{D}(X_2(E))$ and 
$F'\in\C{D}(X'_2(E))$, their pullbacks $\pi^*(F)\in\C{D}(X_1(E))$ and 
$\pi'^*(F')\in\C{D}(X'_1(E))$ are $(a,a';[b])$-equivalent.
\end{Cor}

\subsection{Locally $L^1$ functions} \label{SS:l1}
\begin*
\vskip 8truept
\end*
The goal of this subsection is to write down explicitly the condition for 
$\C{E}$-stability and $(a,a';[b])$-equivalence of generalized functions coming from 
invariant locally $L^1$ functions.


\begin{Not}
(a) Denote by $L^1_{loc}(X(E))$ \label{l1loc} the space of $G(E)$-invariant locally $L^1$ functions on $X(E)$, 
whose restriction to its open subset $X^{\sr}(E)$ is locally constant. 

(b) We have a canonical embedding $L^1_{loc}(X(E))\hra \C{D}(X(E))$, 
which sends each $F\in  L^1_{loc}(X(E))$ to a generalized function 
$\phi\mapsto\int_{X(E)}F\phi$.
For simplicity of notation, we identify 
functions from $L^1_{loc}(X(E))$ with the corresponding generalized functions from $\C{D}(X(E))$.
\end{Not}

\begin{Rem} \label{R:l1}
For each $F\in L^1_{loc}(X(E)))$, the  corresponding generalized function is contained in 
$\C{D}^0(X(E))$.
\end{Rem} 

\begin{Not} \label{N:coinv}
 For a $G(E)$-invariant function $F:X(E)\to\B{C}$ and a pair 
$x\in X^{\sr}(E)$ and $\ov{\xi}\in \pi_0(\wh{G_x}^{\Gm}/Z(\wh{G})^{\Gm})$, 
we put
$F(x,\ov{\xi}):=\sum_{x'}\lan\inv(x,x'),\ov{\xi}\ran^{-1} F(x')$, \label{fxxi}
where $x'\in X(E)$ runs over a set of representatives of $G(E)\bs [x]\subset G(E)\bs X(E)$.
\end{Not}

\begin{Prop} \label{P:reform}
In the notation of \rd{main}, 

(a) $F\in L^1_{loc}(X(E))\subset \C{D}(X(E))$ is $\C{E}$-stable if and only if for each pair 
$x\in X^{\sr}(E)$ and $\ov{\xi}\in\pi_0(\wh{G_x}^{\Gm}/Z(\wh{G})^{\Gm})$
such that 
$([a_x],\ov{\xi})\notin \im\Pi_{\C{E}}$, we have $F(x,\ov{\xi})=0$.

(b) $F\in L^1_{loc}(X(E))$ and $F'\in L^1_{loc}(X'(E))$ are $(a,a';[b])$-equivalent 
if and only if the following two conditions are satisfied:

$(i)$ $F$ and $F'$ are $\C{E}$-stable;

$(ii)$ for all  stably conjugate $x\in X^{\sr}(E)$ and $x'\in X'^{\sr}(E)$
and all $\ov{\xi}\in\pi_0(\wh{G_x}^{\Gm}/Z(\wh{G})^{\Gm})$ such that 
$([a_x],\ov{\xi})\in\im\Pi_{\C{E}}$ and 
$\varphi$ is $(\C{E},[a_x],\ov{\xi})$-admissible, we have 
\[
F'(x',\ov{\xi})=\lan\frac{x,x';\ov{\xi}}{a,a';[b]}\ran F(x,\ov{\xi}).
\]
\end{Prop}

After certain preparations, the proof of the proposition will be carried out 
in \re{proofreform}.

\begin{Not} \label{N:Tor}
Denote by $\un{Tor}$ \label{tor} the variety of all maximal tori in $G$.
\end{Not}

\begin{Lem} \label{L:quotient}
(a) The subset $X^{\sr}\subset X$ is open, and there exists a smooth morphism
$\pi:X^{\sr}\to\un{Tor}$ such that $\pi(x)=G_x$ for each $x\in X^{\sr}$.

(b) There exists a geometric quotient $Y=G\bs X^{\sr}$. 
Moreover, the canonical projection $f:X^{\sr}\to Y$ is smooth, 
the restriction of $f$ to each fiber of $\pi$ is \'etale, 
and the induced map
$f(E):X^{\sr}(E)\to Y(E)$ is a (locally) trivial fibration.
\end{Lem}

\begin{proof}
(a) Denote by $X^{\reg}$ the set of $x\in X$ such that 
$\dim\C{G}_x=\rk_{\ov{E}}(G)$. Then $X^{\reg}$ contains $X^{\sr}$, therefore
$X^{\reg}$ is Zariski dense in $X$. Our first step will be to show that 
$X^{\reg}$ is open in $X$, and the map $x\mapsto\C{G}_x$ gives an algebraic morphism 
$\ov{\pi}$ from $X^{\reg}$ to the Grassmannian 
$\Gr_{\C{G},\rk_{\ov{E}}(G)}$, classifying linear subspaces of $\C{G}$ of dimension
$\rk_{\ov{E}}(G)$.

 Observe that the action $\mu:G\times X\to X$ induces a map 
$T(\mu):T(G)\times T(X)\to T(X)$ of tangent bundles. The restriction of $T(\mu)$ to 
$\C{G}\times X$, where $\C{G}=T_1(G)\subset T(G)$, and $X\subset T(X)$ is the zero 
section, is a map of vector bundles $f:\C{G}\times X\to T(X)$ such that 
for each $x\in X$, the kernel of $f_x:\C{G}\to T_x(X)$ is $\C{G}_x$.
In other words, $X^{\reg}$ can be described as the set of $x\in X$ such that 
$\rk f_x=\dim\C{G}-\rk_{\ov{E}}(G)$. Since $X^{\reg}$ is dense in $X$, 
we get that  $\rk f_x\leq\dim\C{G}-\rk_{\ov{E}}(G)$ for each $x\in X$, 
and $X^{\reg}\subset X$ is open. Moreover, the restriction
$\Ker f|_{X^{\reg}}$ is a vector subbundle of $\C{G}\times X^{\reg}$, 
therefore it gives rise to a morphism 
$\ov{\pi}:X^{\reg}\to\Gr_{\C{G},\rk_{\ov{E}}(G)}$ such that $\ov{\pi}(x)=\C{G}_x$.

Next consider a subset $X^{\rgss}$ of $X^{\reg}$ consisting of points $x$ such that 
$\C{G}_x\subset\C{G}$ is a Cartan subalgebra of $\C{G}$ (hence $G_x^0\subset G$ is 
a maximal torus). Since we assumed that $\C{G}^{\rgss}\neq\emptyset$, 
every Cartan subalgebra of 
$\C{G}$ has a non-zero intersection with $\C{G}^{\rgss}$. 
Hence $X^{\rgss}$ equals the set of $x\in X^{\reg}$ such that 
$\C{G}_x\cap \C{G}^{\rgss}\neq\emptyset$. Since $\C{G}^{\rgss}$ is open in $\C{G}$, 
we conclude that $X^{\rgss}$ is open in $X^{\reg}$, hence in $X$.

Note that the map $T\mapsto\C{T}=\Lie T$ identifies $\un{Tor}$ with the variety 
of Cartan subalgebras of $\C{G}$, which is a locally closed subvariety of 
$\Gr_{\C{G},\rk_{\ov{E}}(G)}$. Therefore the restriction of $\ov{\pi}|_{X^{\rgss}}$ 
can be viewed as a morphism $\pi:X^{\rgss}\to\un{Tor}$ such that $\pi(x)=G_x^0$ 
for each $x\in X^{\rgss}$. 

We claim that $\pi$ is smooth. Since both $X^{\rgss}$ and $\un{Tor}$ are smooth, 
we only have to check that the differential $d\pi_x:T_x(X)\to T_{\pi(x)}(\un{Tor})$ 
is surjective for each $x\in X^{\rgss}$. Put $T:=\pi(x)$. Then $G_x^0=T$, 
$G_x\subset \Norm_G(T)$, $G(x)\cong G/G_x$ and $\un{Tor}\cong  G/\Norm_G(T)$. 
Hence $\pi|_{G(x)}$ is \'etale, thus $d\pi_x|_{T_x(G(x))}$ is surjective, which 
implies the surjectivity of  $d\pi_x$. 

It remains to show that $X^{\sr}$ is open in $X^{\rgss}$. Fix $T\in\un{Tor}$, and 
put $W=\Norm_G(T)/T$. Then $W$ acts on $\pi^{-1}(T)$, and $Z_T:=\pi^{-1}(T)\cap X^{\sr}$ 
consists of points $x\in\pi^{-1}(T)$ such that $w(x)\neq x$ for all $w\neq 1$.
Hence $Z_T\subset\pi^{-1}(T)$ is open, and $W$ acts freely on $Z_T$.
In particular, $Z_T$ is smooth. 

Consider the natural map $\iota:(G/T)\times_W Z_T\to X^{\rgss}:[g,x]\mapsto g(x)$. 
This is a map between  smooth spaces, which induces an isomorphism between tangent 
spaces, hence $\iota$ is \'etale. Since $\iota$ induces a bijection between 
$(G/T)\times_W Z_T$ and  $X^{\sr}\subset X^{\rgss}$, we get that 
$X^{\sr}\cong (G/T)\times_W Z_T$ is open in $X^{\rgss}$.

(b) Since $X$ is quasi-projective,  $Z_T$ is quasi-projective as well. 
Hence a quasi-projective scheme $Y:=W\bs Z_T$ is a geometric quotient
$G\bs X^{\sr}=G\bs [(G/T)\times_W Z_T]$.  
Moreover, the projection $f:X^{\sr}\to Y$ is smooth, and  $f|_{Z_T}$ is \'etale.

To show the last assertion, choose $x\in X^{\sr}(E)$ and put
$T:=G_x$. Since the projection $f|_{Z_T}$ is \'etale, there exist 
open neighborhoods $U\subset Z_T(E)$ of $x$ and $V\subset Y(E)$ of $f(x)$ 
such that  $f$ induces a 
homeomorphism $U\isom V$. Then the map $(G/T)(E)\times U\to X^{\sr}(E)$
sending $([g],u)$ to $g(u)$ induces a $G(E)$-equivariant isomorphism between
$(G/T)(E)\times U\cong (G/T)(E)\times V$ and $f(E)^{-1}(V)$.
\end{proof}

\begin{Con} \label{C:measure}

(a) For each $T\in\un{Tor}(E)$ and an open and compact subset 
$U\subset\pi^{-1}(T)(E)$, there exists a measure 
$\phi_{U}\in\C{S}(X(E))$ such that $O_{x}(\phi_{U})=1$ for each $x\in G(E)(U)$, 
and $O_{x}(\phi_{U})=0$ otherwise.

Explicitly, for each open and compact subgroup $K\subset G(E)$, the measure 
$\phi_{U}:=\frac{|dt|(K\cap T(E))}{|dg|(K)}\chi_{K(U)} dx$, where
$\chi_{K(U)}$ is the characteristic function of
$K(U)\subset X^{\sr}(E)$ and $dt$ is an invariant measure on $T(E)$ such that
the total measure of the maximal compact subgroup of $T(E)$ is $1$,
satisfies the required properties.

(b) For each stable conjugates $x\in X^{\sr}(E)$ and  $x'\in X'^{\sr}(E)$, there 
exists a natural isomorphism $\varphi_{x,x'}$ \label{varphixx'} between 
$\pi^{-1}(G_x)\subset X^{\sr}$ and $\pi'^{-1}(G'_{x'})\subset X'^{\sr}$.

Explicitly, choose $g\in G(E^{\sep})$ such that $x'=\varphi_X(g(x))$.
Then the map $X_{E^{\sep}}\to X'_{E^{\sep}}:y\mapsto\varphi_X(g(y))$ maps 
$\pi^{-1}(G_x)$ to $\pi'^{-1}(G'_{x'})$, and the corresponding morphism 
$\varphi_{x,x'}:\pi^{-1}(G_x)\to\pi'^{-1}(G'_{x'})$ is $E$-rational, independent
of $g$ and $\varphi_{x,x'}(x)=x'$.
\end{Con}

\begin{Cor} \label{C:constr}
(a) Given $x\in X^{\sr}(E)$ and $\ov{\xi}\in\pi_0(\wh{G_x}^{\Gm}/Z(\wh{G})^{\Gm})$ such that
$([a_x],\ov{\xi})\notin\im\Pi_{\C{E}}$, there exists an 
$\C{E}$-unstable measure $\phi\in\C{S}(X(E))$
such that $O^{\ov{\xi}}_x(\phi)\neq 0$ and  $O^{\ov{\xi}'}_x(\phi)=0$ for each 
$\ov{\xi}'\neq\ov{\xi}$. 

(b) Let $x\in X^{\sr}(E)$ and $\ov{\xi}\in\pi_0(\wh{G_x}^{\Gm}/Z(\wh{G})^{\Gm})$ be 
such that $([a_x],\ov{\xi})\in\im\Pi_{\C{E}}$, $\varphi$ is 
$(\C{E},[a_x],\ov{\xi})$-admissible and there exists a stable conjugate $x'\in X'(E)$ 
of $x$. Then there exist 
$(a,a';[b])$-indistinguishable measures $\phi\in\C{S}(X(E))$ and $\phi'\in\C{S}(X'(E))$
such that $O^{\ov{\xi}}_x(\phi)\neq 0$, $O^{\ov{\xi}}_{x'}(\phi')\neq 0$ and  
$O^{\ov{\xi}'}_x(\phi)=O^{\ov{\xi}'}_{x'}(\phi')=0$ for each $\ov{\xi}'\neq\ov{\xi}$. 
\end{Cor}
\begin{proof}

(a) Let $x_1=x,\ldots,x_n\in X(E)$ be a set of representatives of conjugacy classes 
stably conjugate to $x$. Choose an open neighborhood $U\subset\pi^{-1}(G_x)$
of $x$, and for each $i=1,\ldots,n$ put 
$U_i:=\varphi_{x,x_i}(U)\subset\pi^{-1}(G_{x_i})(E)$, and let $\phi_{U_i}\in \C{S}(X(E))$
be as in \rc{measure} (a).
Then measure $\phi:=\sum_{i}\lan\inv(x,x_i),\ov{\xi}\ran^{-1}\phi_{U_i}$
satisfies the required property. 

(b) Now choose a set of representatives $x'_1,\ldots,x'_n\in X'(E)$ 
of conjugacy classes stably conjugate to $x$, and put 
$U'_i:=\varphi_{x,x'_i}(U)\subset\pi'^{-1}(G'_{x_i})(E)$. 
Then  measures $\phi:=\sum_{i}\lan\inv(x,x_i),\ov{\xi}\ran^{-1}\phi_{U_i}$ and
 $\phi':=\sum_{i}\lan\frac{x,x'_i;\ov{\xi}}{a,a';[b]}\ran^{-1}\phi_{U_i}$ 
satisfy the required property (use \rl{indist}).
\end{proof}

\begin{Lem} \label{L:orbit}
(a) Let $F\in\C{D}(X(E))$ be of the form $F=\sum_{\ov{\xi}} c_{x,\ov{\xi}}O^{\ov{\xi}}_{x}$, 
where $x\in X^{\sr}(E),\, c_{x,\ov{\xi}}\in\B{C}$ and $\ov{\xi}$ runs over
$\pi_0(\wh{G_{x}}^{\Gm}/Z(\wh{G})^{\Gm})$. Then 
$F$ is $\C{E}$-stable if and only if $c_{x,\ov{\xi}}=0$ for each $\ov{\xi}$ with  
$([a_{x}],\ov{\xi})\notin \im\Pi_{\C{E}}$.

(b) Let $F\in\C{D}(X(E))$ and $F'\in\C{D}(X'(E))$ be of the form
 $F=\sum_{\ov{\xi}} c_{x,\ov{\xi}}O^{\ov{\xi}}_{x}$ and 
 $F'=\sum_{\ov{\xi}} c_{x',\ov{\xi}}O^{\ov{\xi}}_{x'}$ for some stable conjugate
 $x\in X^{\sr}(E)$ and  $x'\in X'^{\sr}(E)$. Then 
$F$ and $F'$ are $(a,a';[b])$-equivalent if and only if they satisfy the following 
two conditions:

$(i)$ $F$ and $F'$ are $\C{E}$-stable;

$(ii)$ for each $\ov{\xi}\in\pi_0(\wh{G_x}^{\Gm}/Z(\wh{G})^{\Gm})=
\pi_0(\wh{G'_{x'}}^{\Gm}/Z(\wh{G'})^{\Gm})$ such that  $([a_{x}],\ov{\xi})\in \im\Pi_{\C{E}}$
and $\varphi$ is $(\C{E},[a_x],\ov{\xi})$-admissible, we have
$c_{x',\ov{\xi}}=\lan \frac{x,x';\ov{\xi}}{a,a';[b]}\ran c_{x,\ov{\xi}}$.
\end{Lem}
\begin{proof}
(a) The ``if'' assertion is clear. The ``only if'' assertion follows from the equality 
$F(\phi)=0$ applied to measure $\phi$ from \rco{constr} (a).

(b) Assume that $F$ and $F'$ satisfy assertions $(i)$ and $(ii)$. 
Then it follows from \rl{indist} and assertion (a), that 
for each  $(a,a';[b])$-indistinguishable $\phi\in\C{S}(X(E))$ and 
$\phi'\in\C{S}(X'(E))$ we have
$F(\phi)=F'(\phi')$. Conversely, assume that $F$ and $F'$ are 
$(a,a';[b])$-equivalent. Then condition $(i)$
was proved in \rco{end} (b) and condition $(ii)$ follows from the equality
$F(\phi)=F'(\phi')$ applied to measures $\phi$ and $\phi'$ from \rco{constr} (b).
\end{proof}

The following result is clear. 

\begin{Lem} \label{L:deltadist}
Let $f:Z\to Y$ be a morphism of smooth algebraic varieties over $E$ such that  the induced map 
$f(E):Z(E)\to Y(E)$ is a locally trivial fibration.
Fix a measure $\mu$ on $Y(E)$, and let $U_1\supset U_2\supset\ldots$ be a basis of open and compact 
neighborhoods of $y\in Y(E)$. Then for every locally constant function $F$ on $Z(E)$ and every 
$\phi\in\C{S}(Z(E))$, the sequence $\frac{1}{\mu(U_i)}F|_{f^{-1}(U_i)}(\phi)$ stabilizes.
\end{Lem}

\begin{Not}
For each $x\in X^{\sr}(E)$ and $F\in L^1_{loc}(X(E))$, denote by $F_x\in \C{D}(X(E))$ \label{fx} 
the generalized function $\phi\mapsto SO_{x}(F\phi)$. For each 
$x^*\in (X^*)^{\sr}(E)$, $F\in L^1_{loc}(X(E))$ and $F'\in L^1_{loc}(X'(E))$, 
we denote by $F_{x^*}\in\C{D}(X(E))$ and $F'_{x^*}\in\C{D}(X'(E))$ the generalized 
functions $\phi\mapsto SO_{x^*}(F\phi)$ and $\phi'\mapsto SO_{x^*}(F'\phi')$, 
respectively (see \rn{orbital} (iii)).
\end{Not}

\begin{Rem} \label{R:SO}
Clearly, $F_x$ and $F_{x^*}$ belong to $\C{D}^0(X(E))$.
\end{Rem}

\begin{Cl} \label{C:points}
(a) $F\in L^1_{loc}(X(E))$ is $\C{E}$-stable if and only if 
$F_x\in \C{D}(X(E))$ is $\C{E}$-stable for all  $x\in X^{\sr}(E)$.

(b) $F\in L^1_{loc}(X(E))$ and $F'\in L^1_{loc}(X'(E))$ are $(a,a';[b])$-equivalent 
if and only if  $F_{x^*}\in \C{D}(X(E))$ and  $F'_{x^*}\in \C{D}(X'(E))$
are $(a,a';[b])$-equivalent for all $x^*\in (G^*)^{\sr}(E)$.
\end{Cl}

\begin{proof}

(a) Since $X^{\sr}$ is Zariski dense in $X$, the complement 
$X(E)\sm X^{\sr}(E)$ is nowhere dense. As $F\in L^1_{loc}(X(E))$, we get that
$F$ is $\C{E}$-stable if and only if the restriction $F|_{X^{\sr}(E)}$ is $\C{E}$-stable.

Consider the map $f:X^{\sr}\to Y$ from \rl{quotient}. Then for each 
$x\in X^{\sr}(E)$ and $\phi\in\C{S}(X(E))$, the value $F_{x}(\phi)$ is the limit of the 
stabilizing sequence $\frac{1}{\mu(U_i)}F|_{f^{-1}(U_i)}(\phi)$, where $U_i$ is any basis 
of open and compact neighborhoods of $f(x)\in Y(E)$ (use \rl{deltadist}).
From this the assertion follows.

Indeed, if $F|_{X^{\sr}(E)}$ is $\C{E}$-stable, then each  $F|_{f^{-1}(U_i)}$ is 
 $\C{E}$-stable. In particular, for every $\C{E}$-unstable
$\phi$, we have  $F|_{f^{-1}(U_i)}(\phi)=0$ for each $i$, hence $F_{x}(\phi)=0$.
Conversely, assume that each $F_{x} $ is $\C{E}$-stable, and pick an 
$\C{E}$-unstable $\phi$. Then there exists an open disjoint covering $\{U_{\al}\}_{\al}$ of 
$f(Supp\, \phi)$ such that each $F|_{f^{-1}(U_{\al})}(\phi)=0$ for each $\al$, hence 
$F(\phi)=\sum_{\al}F|_{f^{-1}(U_{\al})}(\phi)=0$. This shows that $F$ is $\C{E}$-stable.

(b) As in (a),  $F$ and $F'$ are $(a,a';[b])$-equivalent 
if and only if  $F|_{X^{\sr}(E)}$ and  $F'|_{X'^{\sr}(E)}$ are 
$(a,a';[b])$-equivalent.
Next since $G$ acts trivially on $Y$, we get identifications $Y'=Y^*=Y$, 
and  the projection $f:X^{\sr}\to Y$ induce maps $f':X'^{\sr}\to Y$ 
and  $f^*:(X^*)^{\sr}\to Y$.
Moreover, by \rco{steinberg}, both $\im f(E)$ and $\im f'(E)$ are contained in $\im f^*(E)$.
Now the assertion follows from \rl{deltadist} by exactly the same arguments as (a).
\end{proof}

\begin{Emp} \label{E:proofreform}
\begin{proof}[Proof of \rp{reform}]
(a) For each $x\in X^{\sr}(E)$, we denote by  $N_x$ the cardinality of 
$\pi_0(\wh{G_{x}}^{\Gm}/Z(\wh{G})^{\Gm})$.
By \rcl{points} (a), $F$ is $\C{E}$-stable if and only if 
each $F_x$ is stable. Since $F_x=\frac{1}{N_x}\sum_{\ov{\xi}}F(x,\xi)O_x^{\ov{\xi}}$,
the assertion then is just a reformulation of \rl{orbit} (a). 

(b) The proof of (b) is similar. By  \rcl{points} (b), 
 $F\in L^1_{loc}(X(E))$ and $F'\in L^1_{loc}(X'(E))$ are $(a,a';[b])$-equivalent 
if and only if each $F_{x^*}\in \C{D}(X(E))$ and  $F'_{x^*}\in \C{D}(X'(E))$
are $(a,a';[b])$-equivalent. By definition this means that each $F_x$ and $F'_{x'}$ are
$\C{E}$-stable, and for every stable conjugate $x\in X^{\sr}(E)$ and
$x'\in X'^{\sr}(E)$, $F_x$ and  $F'_{x'}$ are  $(a,a';[b])$-equivalent.
But by  \rcl{points} (a) and \rl{orbit} (b), these conditions are equivalent to 
conditions $(i)$ and $(ii)$, respectively.
\end{proof}
\end{Emp}

\begin{Cor} \label{C:func}
(a) Let $\pi:X_1\to X_2$ be a smooth $G$-equivariant morphism of 
varieties as in \re{setup} (b) and $U\subset X_1(E)$  
an open $G(E)$-invariant subset. For each $F\in L^1_{loc}(X_2(E))$ such 
that $\pi^*(F)|_{U}$ is $\C{E}$-stable, the restriction $F|_{\pi(U)}$ is 
$\C{E}$-stable. 

(b) Let $\pi':X'_1\to X'_2$ be the inner twisting of $\pi$ and 
$U'\subset X'_1(E)$ an open 
$G'(E)$-invariant subset. For each $F'\in L^1_{loc}(X'_2(E))$ such that 
$\pi^*(F)|_{U}$ and  $\pi'^*(F')|_{U'}$ are
$(a,a';[b])$-equivalent, the restrictions $F|_{\pi(U)}$ and  $F'|_{\pi'(U')}$ are
$(a,a';[b])$-equivalent.
\end{Cor}
\begin{proof}
(a) Since $\pi$ is smooth and $G$-equivariant,
the subset $\pi(U)\subset X_2(E)$  
is open and $G(E)$-invariant. Thus  $F|_{\pi(U)}$ 
is defined and belongs to  $L^1_{loc}(X_2(E))$. Let $(x,\ov{\xi})$
be as in \rp{reform} (a). If 
$[x]\cap \pi(U)=\emptyset$, then $F|_{\pi(U)}$ vanishes on $[x]$, hence 
$F|_{\pi(U)}(x,\ov{\xi})=0$. Assume now that $[x]\cap \pi(U)\neq\emptyset$. Replacing $x$ by
a stable conjugate, we can assume that 
$x=\pi(x_1)$ for some $x_1\in X_1(E)$. Then $x_1\in X_1^{\sr}(E)$, $G_{x_1}= G_x$,
$a_{x_1}=a_{x}$, and $\pi$ induces a $G(E)$- and $\Gm$-equivariant isomorphism 
$[x_1]\isom [x]$. Hence 
$F|_{\pi(U)}(x,\ov{\xi})=\pi^*(F)|_U(x_1,\ov{\xi})=0$. Thus the 
assertion follows from \rp{reform} (a).

(b) Follows from \rp{reform} (b) by exactly the same arguments as (a).
\end{proof}

\subsection{Quasi-logarithm maps} \label{SS:qlog}
\begin*
\vskip 8truept
\end*
Starting from \re{BT}, $E$ will be a local non-archimedean field, $\C{O}$ the ring of integers
of $E$,  $\frak{m}$ the maximal ideal of $\C{O}$, $\fq$ the residue field of 
$E$, $p$ the characteristic of $\fq$, and $G$ a reductive group over $E$  
split over $E^{\nr}$.

\begin{Def} \label{D:qlog}
Let $G$ be an algebraic group over a field $k$. By a  {\em quasi-logarithm} \label{qlog} we call 
a $G^{\ad}$-equivariant algebraic morphism $\Phi:G\to \C{G}$ \label{phi} such that $\Phi(1)=0$ and 
$d\Phi_1:\C{G}=T_1(G)\to\C{G}$ is the identity map.
\end{Def}

\begin{Ex}
Let $\rho:G\to\Aut V$ be a representation such that the corresponding 
$G$-invariant pairing $\lan a,b\ran_{\rho}:=\Tr(\rho(a)\rho(b))$ \label{lanranrho} on $\C{G}$ 
is non-degenerate. Denote by $\pr_{\rho}:\End V\to\C{G}$ be the projection given by
the rule $\Tr(\pr_{\rho}(A)\rho(b))=\Tr(A\rho(b))$ for each $b\in\C{G}$.
Then the map $\Phi_{\rho}:g\mapsto\pr_{\rho}(\rho(g)-\Id_V)$ \label{phirho} is a 
quasi-logarithm $G\to \C{G}$.
\end{Ex}

\begin{Lem} \label{L:qlogf}
Let $\Phi:G\to\C{G}$ be a quasi-logarithm map.

(a) For every Borel subgroup $B$ of $G$, we have $\Phi(B)\subset \Lie B$;

(b) If a Cartan subgroup of $G$ is a maximal torus, then 
$\Phi$ induces a quasi-logarithm map $G_{\red}:=G/R_u(G)\to \C{G}_{\red}$
 \end{Lem}

\begin{proof}


 Let $T\subset B$ be a maximal torus, and $C:=\Cent_G(T)$ the corresponding Cartan
subgroup. 

(a) Since $\Phi$ is $G^{\ad}$-equivariant, 
$\Phi(C)$ is contained in the set of fixed points of $\Ad T$ in $\C{G}$, 
that is, $\Phi(C)\subset\Lie C$. Therefore 
\[
\Phi(\Inn B(C))=\Ad B(\Phi(T))\subset \Ad B(\Lie C)\subset\Lie B .\]
Since $C$ is a Cartan subgroup of $B$, $\Inn B(C)\subset B$ is Zariski dense, 
hence $\Phi(B)$ is contained in $\Lie B$.

(b) We have to show that for each $g\in G$ and $u\in R_u(G)$, we have
\begin{equation} \label{E:ru}
\Phi(gu)-\Phi(g)\in \Lie R_u(G).
\end{equation}
 Since we assumed that $T$ is a Cartan subgroup, $\Inn G(T)$ is Zariski dense in $G$. Therefore 
it is enough to check the equality
(\ref{E:ru}) only for $g\in\Inn G (T)$, hence (since $\Phi$ is $G^{\ad}$-equivariant and
$R_u(G)$ is normal in $G$), only for $g=t\in T$. Consider subgroup $H=TR_u(G)\subset G$. 
Then $T$ is
a Cartan subgroup of $H$, hence $\Inn H(T)$ is Zariski dense in $H$. Since 
$H=T\ltimes R_u(G)$, it will therefore suffice to check (\ref{E:ru}) 
under the additional assumption that $tu=vtv^{-1}$ for certain $v\in R_u(G)$. 
In this case, 
\[
\Phi(tu)-\Phi(t)=(\Ad v-\Id)(\Phi(t))\in (\Ad v-\Id)(\C{T})\subset\Lie R_u(G),
\] 
as claimed.
\end{proof}

\begin{Rem}
If we do not assume that Cartan subgroup of $G$ is a maximal torus, then the 
assertion (b) of the lemma is obviously false. For example, it is false for abelian 
groups.
\end{Rem} 

From now on,  $G$ is a reductive group over a local non-archimedean field $E$, which
is split over $E^{\nr}$.

\begin{Emp} \label{E:BT}
{\bf Bruhat-Tits building.} 
(a) Denote by $\C{B}(G)$ \label{Bg} the (non-reduced) Bruhat-Tits building 
of $G$. For every point $x\in\C{B}(G)$, 
we denote by $G_{x}\subset G(E)$ \label{gx} (resp. $\C{G}_x\subset\C{G}$) 
be corresponding  parahoric subgroup (resp. subalgebra), and let
$G_{x^+}\subset G_{x}$ \label{gx+}(resp. $\C{G}_{x^+}\subset \C{G}_{x}$)  be the pro-unipotent 
(resp. pro-nilpotent) radical of $G_x$ (resp. of $\C{G}_x$) (compare \cite{MP1}). 

(b) For each $x\in \C{B}(G)$, denote by $\un{G}_x$ \label{ungx} the canonical
smooth connected group scheme over $\C{O}$ whose generic fiber is $G$ and $\un{G}_x(\C{O})=G_x$, 
and let $\ov{G}_x$ \label{ovgx} be the special fiber of $\un{G}_x$. Then  $\ov{G}_x$ is a connected group over $\fq$,
whose  Cartan subgroup is a maximal torus. (Here we use the assumption that  $G$ splits over $E^{\nr}$).

(c) For each $x\in\C{B}(G)$, denote by $L_x$  \label{lx} the quotient $(\ov{G}_x)_{\red}=\ov{G}_x/R_u(\ov{G}_x)$.
We have canonical identifications $L_x(\fq)=G_x/G_{x^+}$ and
$\C{L}_x:=\Lie L_x=\C{G}_x/\C{G}_{x^+}$. For every $g\in G_x$ and $a\in\C{G}_x$, we put
$\ov{g}:=gG_{x^+}\in L_x(\fq)$ \label{ovg} and $\ov{a}:=a+\C{G}_{x^+}\in\C{L}_x(\fq)$.

(d) If $G=T$ is a torus, then the group scheme $\un{T}_x$ is independent of $x\in\C{B}(T)$, and 
coincides with the canonical $\C{O}$-structure $T_{\C{O}}$ \label{to} of $T$. 
We denote by $\ov{T}$ the special fiber of $T_{\C{O}}$, and will write $T(\C{O})$ instead of
$T_{\C{O}}(\C{O})$.
\end{Emp}

\begin{Not} \label{N:qlog}

(a) We call an invariant pairing $\lan\cdot,\cdot\ran$ on $\C{G}$ 
{\em non-degenerate at $x\in\C{B}(G)$}, \label{nondegatx}
if it is non-degenerate over $E$ and the dual lattice 
\[
(\C{G}_x)^{\perp}:=\{x\in \C{G}(E)\,|\,\lan x,y\ran\in \frak{m}\text { for each } y\in\C{G}_x\}
\] 
equals
$\C{G}_{x^+}$.

(b) We call a quasi-logarithm $\Phi:G\to\C{G}$ {\em defined  at $x\in\C{B}(G)$},  \label{defatx} if 
$\Phi$ extends to the morphism $\Phi_x:\un{G}_x\to\C{G}_x$ \label{phix} of schemes over $\C{O}$.
\end{Not}

\begin{Lem} \label{L:qlogx}
(a) If an invariant pairing $\lan\cdot,\cdot\ran$ on $\C{G}$ is  
non-degenerate at $x$ for some $x\in\C{B}(G)$, then it is non-degenerate at $x$ for all 
$x\in\C{B}(G)$.
In this case, $\lan\cdot,\cdot\ran$ defines an invariant non-degenerate 
pairing $\lan\cdot,\cdot\ran_x$ \label{lanranx} on $\C{L}_x$ for all $x\in\C{B}(G)$.

(b) If  a quasi-logarithm  $\Phi:G\to\C{G}$ over $E$ is defined  at $x$ for some $x\in\C{B}(G)$, 
then it is defined at $x$ for all $x\in\C{B}(G)$. In this case,
$\Phi$ gives rise to a quasi-logarithm map $\ov{\Phi}_x:L_x\to\C{L}_x$ \label{ovphix} for all $x\in\C{B}(G)$.
\end{Lem}

\begin{proof}

(a) Assume  that  $(\C{G}_x)^{\perp}=\C{G}_{x^+}$. Since 
$\frak{m}\C{G}_x\subset\C{G}_{x^+}$, we get 
that $\lan a,b\ran\in\C{O}$ for every $a,b\in\C{G}_x$, and the pairing on $\C{L}_x$ given by the 
formula $\lan \ov{a},\ov{b}\ran_x:=\ov{\lan a,b\ran}\in\fq$ for every $a,b\in\C{G}_x$ is 
well-defined, $L_x$-invariant and non-degenerate.  

It remains to show that for every $x,y\in\C{B}(G)$, the equalities 
$(\C{G}_x)^{\perp}=\C{G}_{x^+}$ and $(\C{G}_y)^{\perp}=\C{G}_{y^+}$ are equivalent.
For this we can extend scalars to $E^{\nr}$.
Also we can assume that $\C{G}_y$ is an Iwahori subalgebra of $\C{G}_x$.

  Assume first that  $(\C{G}_x)^{\perp}=\C{G}_{x^+}$. Since 
$\C{G}_{x^+}\subset\C{G}_y\subset\C{G}_x$, we get that
$\C{G}_{x^+}\subset(\C{G}_y)^{\perp}\subset\C{G}_x$. Moreover, 
$(\C{G}_y)^{\perp}/\C{G}_{x^+}\subset \C{L}_x$ is the orthogonal complement of 
the Borel subalgebra $\C{G}_y/\C{G}_{x^+}\subset\C{L}_x$ with respect to the
non-degenerate pairing $\lan \cdot,\cdot\ran_x$. Hence $(\C{G}_y)^{\perp}/\C{G}_{x^+}$
is the nilpotent radical of $\C{G}_y/\C{G}_{x^+}$, thus $(\C{G}_y)^{\perp}=\C{G}_{y^+}$. 
Conversely, assume that  $(\C{G}_y)^{\perp}=\C{G}_{y^+}$.
Since $\C{G}_x=\cup_{g\in G_x}\Ad g(\C{G}_y)$ and 
$\C{G}_{x^+}=\cap_{g\in G_x}\Ad g(\C{G}_{y^+})$,
we get that  $\C{G}_{x^+}=(\C{G}_x)^{\perp}$, as claimed.

(b) The strategy will be similar to that of (a). Assume that $\Phi$ extends to a morphism 
$\Phi_x:\un{G}_x\to\C{G}_x$. Then the special fiber of  $\Phi_x$ is a quasi-logarithm 
$\ov{G}_x\to\ov{\C{G}}_x$. Since $L_x=(\ov{G}_x)_{\red}$, the existence $\ov{\Phi}_x$ follows from 
\rl{qlogf} (b) and the observation of \re{BT} (b).

It remains to show that for every  $x,y\in\C{B}(G)$, the existence of $\Phi_x$ is equivalent to that 
of $\Phi_y$. Notice first that the existence of $\Phi_x$ is equivalent to the fact that 
$\Phi(\un{G}_x(\C{O}_{E^{\nr}}))\subset \C{G}_x\otimes_{\C{O}} \C{O}_{E^{\nr}}$ 
(see \cite[Prop. 1.7.6]{BT}).
Thus we can extend scalars to $E^{\nr}$, and we are required to check that the inclusions  
$\Phi(G_x)\subset \C{G}_x$ and $\Phi(G_y)\subset \C{G}_y$ are equivalent.
Also we can assume that $G_y$ is an Iwahori subgroup of $G_x$. 

 Assume that $\Phi(G_x)\subset \C{G}_x$. As we have shown, $\Phi$ induces a quasi-logarithm
$\ov{\Phi}_x:L_x\to\C{L}_x$. By \rl{qlogf} (a), we get 
$\ov{\Phi}_x(G_y/G_{x^+})\subset \C{G}_y/\C{G}_{x^+}$, thus $\Phi(G_y)\subset \C{G}_y$.
Conversely, assume that $\Phi(G_y)\subset\C{G}_y$. As $\Phi$ is 
$G^{\ad}$-equivariant, the inclusion $\Phi(G_x)\subset\C{G}_x$ follows from equalities 
$G_x=\cup_{g\in G_x}gG_y g^{-1}$ and $\C{G}_x=\cup_{g\in G_x}\Ad g(\C{G}_y)$.
\end{proof}

\rl{qlogx} allows us to give the following definition.

\begin{Def} \label{D:defO}
(a) We call an invariant pairing $\lan\cdot,\cdot\ran$ on $\C{G}$ {\em non-degenerate over $\C{O}$}, \label{nondego}
if it is non-degenerate at $x$ for some (or, equivalently, for all) $x\in\C{B}(G)$.

(b) We call a quasi-logarithm $\Phi:G\to\C{G}$ {\em defined  over $\C{O}$} \label{defovero} if it is defined at $x$
for some (or, equivalently, for all) $x\in\C{B}(G)$.
\end{Def}

\begin{Lem} \label{L:qlog}
(a) Let $\varphi:G\to G'$ be an inner twisting defined over $E^{\nr}$. Every non-degenerate over 
$\C{O}$ invariant pairing $\lan\cdot,\cdot\ran$ on $\C{G}$ gives rise to the corresponding 
pairing $\lan\cdot,\cdot\ran'$ on $\C{G}'$. Every  quasi-logarithm  $\Phi:G\to\C{G}$
defined over $\C{O}$ gives rise to the corresponding  quasi-logarithm  $\Phi':G'\to\C{G}'$.

(b) If the pairing $\lan \cdot,\cdot\ran_{\rho}$ on $\C{G}$ corresponding to a representation
$\rho:G\to\Aut(V)$ is non-degenerate over $\C{O}$, then the corresponding quasi-logarithm
$\Phi_{\rho}:G\to\C{G}$ is defined over $\C{O}$.
\end{Lem}

\begin{proof}
(a) Recall that $\varphi$ induces isomorphisms $G_{E^{\nr}}\isom G'_{E^{\nr}}$ and  
$\C{G}_{E^{\nr}}\isom \C{G}'_{E^{\nr}}$.
Hence $\Phi$ and $\lan\cdot,\cdot\ran$ give rise to a quasi-logarithm 
$\Phi':G'_{E^{\nr}}\to\C{G}'_{E^{\nr}}$ and a pairing  
$\lan\cdot,\cdot\ran':\C{G}_{E^{\nr}}\times \C{G}_{E^{\nr}}\to E^{\nr}$, respectively.
Furthermore, since the twisting $\varphi$ is inner, while $\Phi$ and  $\lan\cdot,\cdot\ran$
are $G^{\ad}$-equivariant, the quasi-logarithm $\Phi'$ and the inner twisting 
$\lan\cdot,\cdot\ran'$ are defined over $E$. Finally, to show that  $\Phi'$ and 
$\lan\cdot,\cdot\ran'$ are defined over $\C{O}$, we can extend scalars to $E^{\nr}$.
Then the assertion follows from the corresponding assertion for  $\Phi$ and 
$\lan\cdot,\cdot\ran$.

(b) For the proof we can replace $E$ by its finite unramified extension, so that $G$ is 
split over $E$. Fix a hyperspecial vertex $x\in\C{B}(G)$, and we have to show that
$\Phi_{\rho}(G_x)\subset\C{G}_x$. As $(\C{G}_x)^{\perp}=\C{G}_{x^+}$, it is enough to show that
$\Tr(\rho(g)\rho(a))\in \frak{m}$ for every $g\in G_x$ and $a\in\C{G}_{x^+}=\frak{m}\C{G}_{x}$. 
Choose any $\rho(G_x)$-invariant $\C{O}$-lattice $V_x\subset V$. Then $V_x$ is 
$\rho(\C{G}_x)$-invariant, hence $\rho(g)\rho(a)(V_x)\subset \frak{m}V_x$. Thus
$\Tr(\rho(g)\rho(a))\in \frak{m}$, as claimed. 
\end{proof}

%

\begin{Def} \label{D:vg}
We say that the group $G$ over $E$ {\em satisfies property $(vg)$}, \label{vg}
if $G^{\ssc}$ admits a quasi-logarithm map $G^{\ssc}\to\C{G}^{\ssc}$ 
defined over $\C{O}$, $\C{G}^{\ssc}$ admits an invariant pairing non-degenerate 
over $\C{O}$, and $p$ does not divide the order of $Z(G^{\ssc})$.
\end{Def}

\begin{Rem} \label{R:vginner}
By \rl{qlog} (a), for every inner twisting $\varphi:G\to G'$, the group
$G$ satisfies property $(vg)$ if and only if $G'$ satisfies property $(vg)$.
\end{Rem}

\begin{Lem} \label{L:vg}
Write $(G^*)^{\ssc}$ in the form $\prod_i R_{E_i/E}H_i$, where each $H_i$ is a 
quasi-split absolutely simple over a finite unramified extension $E_i$ of $E$. 
Then $G$ satisfies property $(vg)$, if the 
following conditions are satisfied: 

$(i)$ $p$ is good for each $H_i$ in the sense of \cite[I, $\S$4]{SS};

$(ii)$ $p$ does not divide the order of each $Z(H_i)$;

$(iii)$ $p$ does not divide $[E_i[H_i]:E_i]$, where $E_i[H_i]$ is the splitting field of $H_i$.
\end{Lem}

\begin{proof}
Assume that $p$ satisfies assumptions $(i)-(iii)$ of the lemma. 
By \rl{qlog} (a), we may replace $G$ by $(G^*)^{\ssc}$ hence by each $R_{E_i/E}H_i$, 
thus assuming that $G$ is quasi-split simple and simply connected.
By \rl{qlog} (b), it will suffice to construct a representation 
$\rho$ of $G$ such that the corresponding pairing $\lan\cdot,\cdot\ran_{\rho}$ on 
$\C{G}$ is non-degenerate over $\C{O}$. We will construct such a $\rho$ in three steps.

Assume first that $G=H_i$ is split. In this case, take $\rho$ be the standard representation, 
if $G$ is classical, and the adjoint representation, if $G$ is exceptional. Then  
assumptions $(i)$ and $(ii)$ imply (as in \cite[I, Lem. 5.3]{SS}) that the pairing 
$\lan\cdot,\cdot\ran_{\rho}$ is non-degenerate at every hyperspecial vertex of $\C{B}(G)$, 
hence non-degenerate over $\C{O}$.

Next we assume that $G=H_i$ is absolutely simple, set $E':=E[G]$, and put $G':=G_{E'}$. 
Then by the proven above, there exists a representation $\rho':G'\to\Aut_{E'}(V)$ such that the pairing
$\lan\cdot,\cdot\ran_{\rho'}$ on $\C{G}'$ is non-degenerate over $\C{O}_{E'}$. Take
$\rho$ be the restriction of  $R_{E'/E}\rho':R_{E'/E}G'\to\Aut_{E}(V)$ to $G$.
Extending scalars to $E'$, we get $\rho_{E'}\cong(\rho')^{[E':E]}$, hence
$\lan\cdot,\cdot\ran_{\rho_{E'}}=[E':E]\lan\cdot,\cdot\ran_{\rho'}$. Therefore by assumption $(iii)$,
$\lan\cdot,\cdot\ran_{\rho_{E'}}$ is non-degenerate over $\C{O}_{E'}$, thus
$\lan\cdot,\cdot\ran_{\rho}$ is non-degenerate over $\C{O}$.

In the general case, choose a representation $\rho':H_i\to \Aut_{E_i}(V)$ such that
$\lan\cdot,\cdot\ran_{\rho'}$ is non-degenerate over $\C{O}_{E_i}$. Then the representation
$\rho:=R_{E_i/E}\rho':G\to\Aut_{E}(V)$ satisfies the required property.
\end{proof}

\begin{Rem} \label{R:vg}
The name $(vg)$ was chosen to indicate the fact that it is closely related to the notion of a 
very good prime. (Recall that $p$ is called very good for $H_i$, if it satisfies properties 
$(i)$ and $(ii)$ of \rl{vg}.)
\end{Rem} 

\begin{Not} \label{N:tnil}
(a) For an algebraic group $H$, we denote by $\C{U}(H)\subset H$ \label{uh} and  
$\C{N}(\C{H})\subset \C{H}$ \label{nh} the subvarieties of unipotent elements of $H$ and of
nilpotent elements of $\C{H}$, respectively.
 
(b)  For every $x\in \C{B}(G)$, we denote by $G_{x,\tu}\subset G_x$ and $\C{G}_{x,\tn}\subset\C{G}_x$ 
\label{gxtu}
the preimages of $\C{U}(L_x)(\fq)\subset L_x(\fq)$ and $\C{N}(\C{L}_x)(\fq)\subset \C{L}_x$, 
respectively.  

(c) Put $G(E)_{\tu}:=\cup_{x\in\C{B}(G)} G_{x,\tu}$ \label{getu}
and  $\C{G}(E)_{\tn}:=\cup_{x\in\C{B}(G)} \C{G}_{x,\tn}$.\label{getn}

\end{Not}

\begin{Lem} \label{L:tnil}
For every $x\in\C{B}(G)$, 

(a) $G_{x,\tu}=\cup_y G_{y^+}$ and $\C{G}_{x,\tn}=\cup_y \C{G}_{y^+}$, where $y$ runs over 
the union of alcoves in $\C{B}(G)$, whose closures contain $x$.

(b) $G_{x,\tu}=G_x\cap G(E)_{\tu}\subset G(E)$ and 
$\C{G}_{x,\tn}=\C{G}_x\cap\C{G}(E)_{\tn}\subset\C{G}(E)$.
\end{Lem}

\begin{proof}
(a) is clear.

(b) The first assertion follows from the equality 
$G_{x,\tu}=\{g\in G_x\,|\,g^{p^n}\underset{n\to\be}{\lra} 1\}$. 
For the second equality, we have to show 
that for each $z\in\C{B}(G)$, there exists $y$ as in (a), such that 
$\C{G}_x\cap\C{G}_{z^+}\subset\C{G}_{y^+}$.
But every $y$, lying in the segment $[x,z]\subset\C{B}(G)$, satisfies this property.
\end{proof}

\begin{Prop} \label{P:ident}
Let  $\Phi:G\to\C{G}$ be a quasi-logarithm defined over $\C{O}$. 

(a) For each $x\in\C{B}(G)$, $\Phi$ induces measure preserving analytic isomorphisms 
$G_{x,\tu}\isom\C{G}_{x,\tn}$ and $G_{x^+}\isom\C{G}_{x^+}$ (with respect to measures
$|\om_G|$ and $|\om_{\C{G}}|$ chosen in \re{setup}).

(b) $\Phi$ induces a  measure preserving analytic isomorphism 
$G(E)_{\tu}\isom\C{G}(E)_{\tn}$. 

(c) Let ${}^0G\subset G$ \label{0g} be the biggest open subset $U\subset G$ such that $\Phi|_{U}$
is \'etale. Then  ${}^0G(E)$ contains $G(E)_{\tu}$.
\end{Prop}

\begin{proof}
(a) For the proof, one can replace $E$ by a finite unramified extension, so we can assume that
$G$ splits over $E$. By \rl{qlogx} (b), we have  $\Phi(G_{x^+})\subset\C{G}_{x^+}$ for each 
$x\in\C{B}(G)$,
therefore by \rl{tnil} (a),   $\Phi(G_{x,\tu})\subset\C{G}_{x,\tn}$.
Since $|\om_G|(G_{x^+})=|\om_{\C{G}}|(\C{G}_{x^+})$, the second assertion follows from 
the first one.

Let us first show the assertion for a hyperspecial vertex $x\in\C{B}(G)$.
By \rl{qlogx} (b), $\Phi$ extends to the morphism $\Phi_x:\un{G}_x\to\C{G}_x$ of schemes over $\C{O}$, 
whose special fiber is a quasi-logarithm map $\ov{\Phi}_x:L_x\to\C{L}_x$. 
By \cite[9.1, 9.2, 9.3.3 and 6.3]{BR} (compare \cite[9.4]{BR}), 
$\ov{\Phi}_x$ induces an isomorphism $\C{U}(L_x)\isom\C{N}(\C{L}_x)$.
Moreover, there exists an open affine neighborhood $V\subset L_x$ of $\C{U}(L_x)$ such that
$\ov{\Phi}_x|_{V}$ is \'etale (see \cite[Thm 6.2 and 9.1]{BR}). Therefore by Hensel's lemma, 
$\Phi_x$ induces an analytic isomorphism $G_{x,\tu}\isom \C{G}_{x,\tn}$.
 
Since $\Phi_x$ is an algebraic morphism over $\C{O}$, we get that
$\Phi(g G_{x,r})=\Phi(g)+\C{G}_{x,r}$ for each $g\in G_{x,\tu}$ and $r\in\B{N}$. 
But $|\om_G|(g G_{x,r})=|\om_\C{G}|(\Phi(g)+\C{G}_{x,r})$, and $\{g G_{x,r}\}_{g,r}$ 
form a basis of open neighborhoods of $G_{x,\tu}$. Hence 
the analytic isomorphism $G_{x,\tu}\isom\C{G}_{x,\tn}$ is measure preserving.

It remains to show that for every $x,y\in\C{B}(G)$, the assertions for $x$ 
and $y$ are equivalent. Moreover, we can assume that $G_y$ is an Iwahori subgroup of 
$G_x$ (compare the proof of  \rl{qlogx}). Then $G_{y,\tu}=G_{x,\tu}\cap G_y$ and 
$\C{G}_{y,\tn}=\C{G}_{x,\tn}\cap \C{G}_y$ (see \rl{tnil}), 
so the assertion for $x$ implies that for $y$. The opposite direction follows from 
equalities $G_{x,\tu}=\cup_{g\in G_x}gG_{y,\tu} g^{-1}$ and 
$\C{G}_{x,\tn}=\cup_{g\in G_x}\Ad g(\C{G}_{y,\tn})$.

(b) By (a), we get that $\Phi(G(E)_{\tu})=\C{G}(E)_{\tn}$, and that the induced map  
$G(E)_{\tu}\to\C{G}(E)_{\tn}$ is open. Thus we have to check that the restriction of 
$\Phi$ to $G(E)_{\tu}$ is one-to-one. 

Assume that  $g_1,g_2\in G(E)_{\tu}$ 
satisfy $\Phi(g_1)=\Phi(g_2)$. Choose $x,y\in\C{B}(G)$ such that 
$g_1\in G_{x^+}$ and $g_2\in G_{y^+}$ (use \rl{tnil}). 
By (a), $\Phi$ induces a measure preserving embedding
$G_{x^+}\cap G_{y^+}\hra\C{G}_{x^+}\cap\C{G}_{y^+}$. As measures of 
both sides are equal, the last embedding is surjective. But
$\Phi(g_1)=\Phi(g_2)$ belongs to $\C{G}_{x^+}\cap\C{G}_{y^+}$, hence there exists 
$g_3\in G_{x^+}\cap G_{y^+}$ such that  
$\Phi(g_1)=\Phi(g_3)=\Phi(g_2)$. Since $\Phi|_{G_{x^+}}$ and  $\Phi|_{G_{y^+}}$
are injective, we get that $g_1=g_3=g_2$, as claimed.

(c) We have to check that for each $g\in G(E)_{\tu}$, the differential
$d\Phi_g:T_g(G)\to T_{\Phi(g)}(\C{G})$ is an isomorphism. But this follows from (b).
\end{proof} 

We finish this subsection with a result, which we will need later. 

\begin{Lem} \label{L:isog} 
(a) Let $\pi:\wt{G}\to G$ be an isogeny (that is, a finite surjective quasi-isogeny) 
of order prime to $p$. Then $\pi(\wt{G}(E)_{\tu})=G(E)_{\tu}$.

(b) Let $\pi:\wt{G}\to G$ be a surjective quasi-isogeny such that $S=\ker\pi$ is a torus split over 
$E^{\nr}$. Then for every $x\in \C{B}(\wt{G})$, we have $\pi(\wt{G}_x)=G_{\pi(x)}$. 
\end{Lem} 

\begin{proof}
(a) Since $\pi$ is an isogeny, it identifies $\C{B}(\wt{G})$ with $\C{B}(G)$.  
Thus we have to check that for each  
$x\in\C{B}(\wt{G})=\C{B}(G)$, we have   
$\pi(\wt{G}_{x,\tu})=G_{x,\tu}$. Since the order of $\pi$ is prime to  
$p$, the corresponding map $\pi_x:\un{\wt{G}}_x\to\un{G}_x$ of group schemes  
over $\C{O}$ is \'etale. Since the special fiber  
$\ov{\pi}_x:\ov{\wt{G}}_x\to\ov{G}_x$ induces an isomorphism   
$\C{U}(\ov{\wt{G}}_x)\isom\C{U}(\ov{G}_x)$, the assertion 
follows from Hensel's lemma.  

(b) Homomorphism $\pi$ gives rise to an exact sequence
$1\to S_{\C{O}}\to \un{\wt{G}}_x\overset{\pi_x}{\lra} \un{G}_{\pi(x)}\to 1$
of group schemes over $\C{O}$. In particular, $\pi_x$ is smooth.
Passing to special fibers, we get an exact sequence
$1\to \ov{S}\to \ov{\wt{G}}_x\overset{\ov{\pi}_x}{\lra} \ov{G}_{\pi(x)}\to 1$
of groups over $\fq$. Since $\ov{S}$ is connected, we get that $H^1(\fq,\ov{S})=1$, hence
the map $\ov{\pi}_x(\fq)=\pi_x(\fq)$
is surjective. Therefore the surjectivity of $\pi_x(\C{O})$ follows from Hensel's lemma.
\end{proof}  

\begin{Cor} \label{C:unrtori} 
Let $\iota:G^{\ssc}\to G$ be a canonical map. Then for every $x\in\C{B}(G)$ and an 
unramified maximal torus $T\subset G$, we have $G_x\subset \iota(G^{\ssc})(E)\cdot T(\C{O})$. 
\end{Cor}

\begin{proof}
Assume first that $G^{\der}=G^{\ssc}$. 
Denote by $q$ the projection $G\to G^{\ab}$. Then we have to check that $q(G_x)\subset q(T(\C{O}))$. 
Since $q(G_x)\subset G^{\ab}(\C{O})$, the assertion follows from part (b) of
the lemma applied to the morphism $q|_T:T\to G^{\ab}$. 

For a general $G$, there exists a surjective quasi-isogeny $\pi:\wt{G}\to G$ 
such that $\wt{G}^{\der}=\wt{G}^{\ssc}(=G^{\ssc})$, and $\Ker\pi$ is an induced torus splitting 
over $E^{\nr}$ (see \cite[Prop. 3.1]{MS}). Then for every $\wt{x}\in \C{B}(\wt{G})$ such that 
$\pi(\wt{x})=x$ we have $\pi(\wt{G}_{\wt{x}})=G_{x}$, so the assertion for $G_x$ and $T$
follows from that for $\wt{G}_{\wt{x}}$ and $\pi^{-1}(T)\subset\wt{G}$.
\end{proof}
\section{Endoscopic decomposition}

\subsection{Main Theorem}
\begin*
\vskip 8truept
\end*
In this subsection we will give two equivalent formulations of the main result of the paper.

\begin{Emp} \label{E:finite}
{\bf Deligne--Lusztig representations.}  \label{DLr}
 Let $L$ be a connected reductive group over $\fq$, $\ov{a}:\ov{T}\hra L$ an embedding of a maximal 
torus of $L$, and $\ov{\theta}:\ov{T}(\fq)\to\B{C}\m$ a character. To this data Deligne and Lusztig
\cite{DL} associate a virtual representation $R^{\theta}_{\ov{a}(\ov{T})}$ of $L(\fq)$. Moreover,
if the torus $\ov{a}(\ov{T})\subset L$ is elliptic and  the character $\ov{\theta}$ is non-singular,
 then $\rho_{\ov{a},\ov{\theta}}:=(-1)^{\rk_{\fq}(L)-\rk_{\fq}(\ov{T})}R^{\theta}_{\ov{a}(\ov{T})}
(=e(L)R^{\theta}_{\ov{a}(\ov{T})})$ \label{rhoatheta} is a cuspidal representation 
(see \cite[Prop. 7.4 and Thm. 8.3]{DL}). In particular, $\rho_{\ov{a},\ov{\theta}}$ is a genuine 
representation and not a virtual one. Moreover, $\rho_{\ov{a},\ov{\theta}}$ is irreducible, if 
$\ov{\theta}$ is in general position. 
\end{Emp}

\begin{Emp}
Recall that there is an equivalence of categories $T\mapsto \ov{T}$
between tori over $E$ splitting over $E^{\nr}$ and tori over $\fq$. Moreover, 
every such $T$ has a canonical $\C{O}$-structure.
We denote by $T(\C{O})^+$ and $\C{T}(\C{O})^+$ the kernels of the reduction maps
$\Ker[T(\C{O})\to\ov{T}(\fq)]$ and $\Ker[\C{T}(\C{O})\to\ov{\C{T}}(\fq)]$, respectively.
\end{Emp}

\begin{Not} \label{N:repr}
(a) Let $G$ be a reductive group over $E$,  $T$ a torus over $E$ splitting over $E^{\nr}$,
and $a:T\hra G$ \label{atg} an embedding of a maximal elliptic torus of $G$.

For every vertex $x$ of $\C{B}(T)$, $a(x)$ is a vertex of $\C{B}(G)$.
Moreover, since $T$ is elliptic, we have $\C{B}(T)=\C{B}(Z(G)^0)$. 
Thus the $Z(G)^0(E)$-orbit of $a(x)$, hence also the parahoric subgroup $G_{a(x)}$ does not 
depend on $x$. Therefore we can denote  $G_{a(x)}$ by $G_{a}$, \label{ga} and similarly for 
$\C{G}_{a(x)}$, $G_{a(x)^+}$, $\C{G}_{a(x)^+}$, $L_{a(x)}$ and $\C{L}_{a(x)}$.
We also set $\wt{G_a}:=Z(G)(E)G_a$. \label{wtga}

An embedding $a:T\hra G$ induces an embedding $\ov{a}:\ov{T}\hra L_a$  \label{ovatg} of a maximal elliptic torus of $L_a$.

(b) Let $\theta:T(E)\to\B{C}\m$ \label{th} be a non-singular character (that is, $\theta$ is not orthogonal 
to any coroot of $(G,T)$), trivial on $T(\C{O})^+$. Denote by $\ov{\theta}:\ov{T}(\fq)\to\B{C}\m$ \label{ovth}
the character of $\ov{T}(\fq)$ defined by $\theta$. 
Then there exists a unique representation $\rho_{a,\theta}$ \label{rhoath} of $\wt{G_a}$,
whose central character is the restriction of $\theta$, extending the inflation to $G_a$
of the Deligne--Lusztig representation $\rho_{\ov{a},\ov{\theta}}$ of $L_a(\fq)$.
We denote by $\pi_{a,\theta}$ the induced representation $\Ind_{\wt{G_a}}^{G(E)}\rho_{a,\theta}$
of $G(E)$. Since for each irreducible factor $\rho'\subset\rho_{a,\theta}$, the
 induced representation $\Ind_{\wt{G_a}}^{G(E)}\rho'$  is cuspidal and irreducible
(see \cite[Prop. 6.6]{MP2}), we get that $\pi_{a,\theta}$ \label{piath}
is a semisimple cuspidal representation of finite length, which is irreducible, if 
$\theta$ is in general position.
\end{Not} 

\begin{Def} \label{D:vga}

(a) Let $a:T\hra G$ be an embedding of a maximal torus split over $E^{\nr}$.
We say that an element $\ov{t}\in\ov{\C{T}}(\fq)$ is {\em $a$-strongly regular}, \label{astrreg} 
if $\ov{t}$ is not fixed by a non-trivial element of the Weyl group $W(G,a(T))\subset
\Aut(T_{E^{\nr}})=\Aut(\ov{T}_{\ov{\fq}})$.

(b) Let $a^{\ssc}:T^{\ssc}\hra G^{\ssc}$ be the lift of $a$. We say that $G$ 
{\em satisfies property $(vg)_a$}, \label{vga} if $G$ satisfies property $(vg)$ (see \rd{vg})
and there exists an  {\em $a^{\ssc}$-strongly regular} element $\ov{t}\in\ov{\C{T}}^{\ssc}(\fq)$.
\end{Def} 

\begin{Not} \label{N:char}
To each  $\ka\in \wh{T}^{\Gm}/Z(\wh{G})^{\Gm}$, an embedding  $a_0:T\hra G$, and 
a character $\theta$ of $T(E)$ as in \rn{repr}, we associate an invariant generalized function 
\[
\chi_{a_0,\ka,\theta}:=e(G)\sum_{a}\lan \inv(a_0,a),\ka\ran \chi(\pi_{a,\theta})\in\C{D}(G(E)).
\] \label{chia0kath}
Here $a$ runs over a set of representatives of conjugacy classes of embeddings 
$T\hra G$ which are 
stably conjugate to $a_0$, and $\chi(\pi_{a,\theta})$ \label{chipiath} denotes the character of $\pi_{a,\theta}$.
\end{Not}

Now we are ready to formulate our main result of the paper.

\begin{Thm} \label{T:main}
Let $(a_0,\ka,\theta)$ be as in \rn{char}. Assume that the characteristic of $E$ is zero and
$G$ satisfies property $(vg)_{a_0}$. Then 

(a) The generalized function $\chi_{a_0,\ka,\theta}$ is $\C{E}_{([a_0],\ka)}$-stable.

(b) For each  inner twisting $\varphi:G\to G'$ and each embedding $a'_0:T\hra G'$, 
stably conjugate to  $a_0$, the generalized functions $\chi_{a_0,\ka,\theta}$ on $G(E)$ and  
$\chi_{a'_0,\ka,\theta}$ on $G'(E)$ are $(a_0,a'_0;\ka)$-equivalent.
\end{Thm}

\begin{Rem}
By \rco{end} (b), assertion (a) is a particular case of (b).
Moreover, if $\varphi$ is not $\C{E}_{([a_0],\ka)}$-admissible, then assertions (a) and (b)  
are equivalent (by \rco{end} (c)).

\end{Rem} 

\begin{Not} \label{N:deff}
(a) To each $a:T\hra G$ and $\theta:T(E)\to\B{C}\m$ as in \rn{repr} we associate a function 
$t_{a,\theta}$ \label{tath} on $G(E)$, supported on $\wt{G_a}$ and equal to $\Tr\,\rho_{a,\theta}$ there.

(b) Since $t_{a,\theta}$ is cuspidal, it follows from \cite[Lem. 23]{HC} that  
for each  $\gm\in G^{\sr}(E)$ and each compact open subgroup $K\subset G(E)$, the sum
$\sum_{g\in D_b}t_{a,\theta}(g\gm g^{-1})$, where $D_b:=\wt{G_a}\bs \wt{G_a}bK$, 
does not vanish only for finitely many $b\in \wt{G_a}\bs G(E)/K$.
Therefore the sum
\[
F_{a,\theta}(\gm):=\sum_{b\in \wt{G_a}\bs G(E)/K}
\left[\sum_{g\in D_b}t_{a,\theta}(g\gm g^{-1})\right]
\]\label{fatheta}
stabilizes, and the resulting value is independent of $K$. 

Explicitly,
$F_{a,\theta}(\gm)=\sum_{g\in \wt{G_a}\bs\Om}t_{a,\theta}(g\gm g^{-1})$ for each sufficiently 
large compact modulo center subset $\Om=\wt{G_a}\Om K\subset G(E)$.
In particular, $F_{a,\theta}$ is a locally constant invariant function on $G^{\sr}(E)$.

(c) For each  $\ka\in \wh{T}^{\Gm}/Z(\wh{G})^{\Gm}$, 
put 
\[
F_{a_0,\ka,\theta}:=e(G)\sum_{a}\lan \inv(a_0,a),\ka\ran F_{a,\theta},
\]\label{fa0katheta}
where $a$ runs over a set of representatives of conjugacy classes of embeddings $T\hra G$
which are stably conjugate to $a_0$.
\end{Not}

\begin{Lem} \label{L:HCH}
Assume that the characteristic of $E$ is zero. Then for each $a$ and $\theta$ as in \rn{repr}, 
$F_{a,\theta}$ belongs to $L^1_{loc}(G(E))$, and 
the corresponding generalized function is equal to $\chi({\pi_{a,\theta}})$.
\end{Lem} 
\begin{proof}
Since $\pi_{a,\theta}$ is cuspidal, the assertion is a combination of the theorem of 
Harish-Chandra (\cite[Thm. 16]{HC}) and 
a formula for characters of induced representations.
\end{proof}

For the next result, we will use \rn{coinv}.

\begin{Thm} \label{T:main'}
Under the assumptions of \rt{main}, let $\gm\in G^{\sr}(E)$ and 
$\ov{\xi}\in\pi_0(\wh{G_{\gm}}^{\Gm}/Z(\wh{G})^{\Gm})$ be such
that $F_{a_0,\ka,\theta}(\gm,\ov{\xi})\neq 0$ . Then 

 $(i)$ $([a_{\gm}],\ov{\xi})\in \im\Pi_{\C{E}};$

 $(ii)$ if $\varphi:G\to G'$ is $(\C{E},[a_{\gm}],\ov{\xi})$-admissible, 
then for every stable conjugate $\gm'\in G'(E)$ of $\gm$ we have 
$F'_{a'_0,\ka,\theta}(\gm',\ov{\xi})=\lan\frac{\gm,\gm';\ov{\xi}}{a,a';\ka}\ran 
F_{a_0,\ka,\theta}(\gm,\ov{\xi})$.
\end{Thm}

\begin{Lem} \label{L:equiv}
\rt{main'} is equivalent to \rt{main}. 
\end{Lem}

\begin{proof}
The equivalence follows from \rl{HCH} and \rp{reform}. 
More precisely, \rp{reform} (a) implies the equivalence between 
\rt{main} (a) and \rt{main'} $(i)$, while \rp{reform} (b) implies the equivalence
between \rt{main} (b) and a combination of \rt{main} (a) and \rt{main'} $(ii)$.
\end{proof}

\begin{Rem} \label{R:pos}
If the characteristic of $E$ is positive, then it is not known that $\chi(\pi_{a,\theta})$ belongs 
to $L^1_{loc}(G(E))$. However the restriction $\chi(\pi_{a,\theta})|_{G^{\sr}(E)}$ belongs to
$L^1_{loc}(G^{\sr}(E))$, therefore \rp{reform} implies that \rt{main'} for $E$ is equivalent
to an analog of \rt{main} for restrictions $\chi_{a_0,\ka,\theta}|_{G^{\sr}(E)}$ 
and $\chi_{a'_0,\ka,\theta}|_{G'^{\sr}(E)}$.

Moreover, \rt{main'} for local fields of positive characteristic 
follows from that for local fields of characteristic zero by approximation arguments of 
\cite{Ka3} and \cite{De} (see \cite{KV2}).
\end{Rem}

\subsection{Stability of the restriction to $G(E)_{\tu}$} \label{SS:unip} 
\begin*
\vskip 8truept
\end*
Starting from this subsection we will assume that the characteristic of $E$ is zero.
In this subsection we will strongly use definitions and results from \rss{qlog}.

\begin{Emp} \label{E:ass}
{\bf Assumptions.} Assume that $G$ admits a quasi-logarithm map $\Phi:G\to\C{G}$
defined over $\C{O}$, $\C{G}$ admits an invariant pairing  $\lan\cdot,\cdot\ran$
non-degenerate over $\C{O}$, and there exists $t\in\C{T}(\C{O})$, whose 
reduction $\ov{t}\in\ov{\C{T}}(\fq)$ is $a_0$-strongly regular (see \rd{vga} (a)).
\end{Emp}

\begin{Not}
(a) For every generalized function $F\in\C{D}(G(E))$, denote by $F_{\tu}$ \label{ftu} the
restriction of $F|_{G(E)_{\tu}}$ (see \rn{meas} and \rn{tnil}).  
Since $G(E)_{\tu}\subset {}^0G(E)$ (see \rp{ident} (c)), we can consider $F_{\tu}$ 
as an element either of $\C{D}(G(E))$ or of $\C{D}({}^0 G(E))$. 

(b) Denote by ${}^0\Phi:{}^0G\to\C{G}$ \label{0Phi} the restriction of $\Phi$ to ${}^0G$.
\end{Not}

The goal of this subsection is to prove the following particular case of 
\rt{main}.

\begin{Thm} \label{T:unip}
Let $(a_0,\ka,\theta)$ be as in \rn{char}. Under the assumptions of \re{ass}, 
the generalized functions $(\chi_{a_0,\ka,\theta})_{\tu}$ and  
$(\chi_{a'_0,\ka,\theta})_{\tu}$ are $(a_0,a'_0;\ka)$-equivalent. 
In particular, each $(\chi_{a_0,\ka,\theta})_{\tu}$ is $\C{E}_{([a_0],\ka)}$-stable.
\end{Thm}

\rt{unip} will be deduced in \re{proofunip} 
from the corresponding statement about generalized functions on Lie algebras.

\begin{Not} 
For every $a:T\hra G$ as in \rn{repr}, we denote by $\ov{\Om}_{a,t}\subset \C{L}_a(\fq)$ \label{omat}
the $\Ad L_a(\fq)$-orbit of $\ov{a}(\ov{t})$, by $\Om_{a,t}\subset \C{G}_{a}\subset \C{G}(E)$
the preimage of $\ov{\Om}_{a,t}$, and let $\dt_{a,t}$ and $\ov{\dt}_{a,t}$ \label{dtat} be the characteristic 
functions of $\Om_{a,t}$ and $\ov{\Om}_{a,t}$, respectively.
\end{Not}

\begin{Lem} \label{L:lie}
(a) For each $y\in {\Om}_{a,t}$, its stabilizer $G_y\subset G$ is 
$G_a$-conjugate to $a(T)$; 

(b) for each $y\in{\Om}_{a,t}$ and $g\in G(E)$ such that $\Ad g(y)\in {\Om}_{a,t}$, we have
$g\in \wt{G_a}$.

\end{Lem}

\begin{proof}
(a) Since $\ov{t}\in\ov{\C{T}}(\fq)$ is $a_0$-strongly regular, 
we see that $a(t)\in a(\C{T}(E))\subset \C{G}(E)$ is strongly regular, hence $G_{a(t)}=a(T)$.

First we will show that for every $y\in a(t)+\C{G}_{a^+}$, we have 
$y\in\C{G}^{\sr}(E)$, and  $G_{y}$ is $G_{a^+}$-conjugate to $a(T)$. By \cite[Lem. 2.2.2]{DB}, 
it will suffice to prove that $y\in\C{G}(E)$ is $G$-regular, and 
$G_{y}$ splits over $E^{\nr}$. 
For this we can replace $E$ by an unramified extension, so we may assume that $T$ splits
over $E$. Under this assumption we will show that $y$ is $G(E)$-conjugate to an element of 
$a(t+\C{T}(\C{O})^+)$. 

Choose an Iwahori subgroup $I\subset {G}_a$, containing $a(T(\C{O}))$, and let 
$\C{I}, I^+$ and $\C{I}^+$ be the corresponding  Iwahori subalgebra, the pro-unipotent radical of 
$I$ and the pro-nilpotent radical of $\C{I}$, respectively.
Since $\al(a(t))\in\C{O}\m$ for each root $\al$ of $(G,a(T))$, it follows from  
direct calculations that every element of $a(t)+\C{I}^+$ is $I^+$-conjugate to an element of
$a(t+\C{T}(\C{O})^+)$. But $y\in a(t)+\C{G}_{a^+}\subset a(t)+\C{I}^+$, therefore the get the 
assertion in this case.

For an arbitrary $y\in\Om_{a,t}$, there exists $h\in G_a$ such that 
$\Ad h(y)\in a(t)+\C{G}_{a^+}$. So the general case follows from the previous one. 

(b) Replacing  $y$  and $\Ad g(y)$ by their $G_a$-conjugates, we can assume that 
$y\in a(t)+\C{G}_{a^+}$ and $\Ad g(y)\in a(t)+\C{G}_{a^+}$. Then by the shown in (a), 
one can further replace  $y$  and $\Ad g(y)$ by their $G_{a^+}$-conjugates, so that 
both $G_y$ and $G_{\Ad g(y)}=gG_yg^{-1}$ equal $a(T)$. Thus $g\in\Norm_G(a(T))$.

Since $\ov{\Ad g(y)}=\ov{a}(\ov{t})=\ov{y}$ is not fixed by a non-trivial element of the Weyl 
group $W(G,a(T))$,
we get that $g\in a(T)(E)$. By \rco{tilde} (a) below, 
$g$ therefore belongs to $\wt{G_a}$, as claimed.
\end{proof}
\begin{Lem} \label{L:torus}
Let $T$ be an unramified torus over $E$, and $S\subset T$ a maximal split subtorus. 
Then $T(E)=T(\C{O})S(E)$.
\end{Lem}

\begin{proof}
By a very particular case of \rl{isog} (b), the projection $T(\C{O})\to (T/S)(\C{O})$ is surjective, therefore
we have to check that $(T/S)(E)=(T/S)(\C{O})$. 
Since $T/S$ is anisotropic over $E$, the group $(T/S)(E)$ is compact. Hence
$(T/S)(E)$ is contained in $(T/S)(E)\cap (T/S)(\C{O}_{E^{\sep}})=(T/S)(\C{O})$, 
as claimed.
\end{proof}

\begin{Cor} \label{C:tilde}
(a) $\wt{G_a}=a(T)(E)G_a$;

(b) $G_a$ is the unique maximal compact subgroup of $\wt{G_a}$. 
\end{Cor}
\begin{proof}
(a) Since $Z(G)\subset a(T)$, we get the inclusion $\wt{G_a}\subset a(T)(E)G_a$.
It remains to show that $a(T)(E)$ is contained in $\wt{G_a}$. 
Let $S\subset T$ be the maximal split subtorus. Since $a(T)\subset G$
is elliptic, we get $a(S)\subset Z(G)$. Now the assertion follows from 
the inclusion $a(T(\C{O}))\subset G_a$ and \rl{torus}.

(b) Assume that $g\in\wt{G_a}$ belongs to a compact subgroup. Choose $g_a\in G_a$ and 
$z\in Z(G)(E)\subset a(T)(E)$ such that $g=g_a z$.  Since $G_a$ is compact and the sequence 
$\{g^n\}_{n}=\{g_a^n z^n\}_{n}\subset\wt{G_a}$ has a convergent subsequence, 
the sequence $\{z^n\}_{n}\subset a(T)(E)$ has a convergent subsequence.
Hence $z$ is contained in $a(T)(\C{O})\subset G_a$, thus $g\in G_a$.
\end{proof}

\begin{Not} 
It follows from \rl{lie} (b), that for each $x\in \C{G}(E)$ there exists at most one coset
$g\in\wt{G_a}\bs G(E)$ such that $\dt_{a,t}(\Ad g(x))\neq 0$. Therefore
\[
\Dt_{a,t}(x):=\sum_{g\in \wt{G_a}\bs G(E)}\dt_{a,t}(\Ad g(x)) 
\]\label{Dtat}
is the characteristic function of an open and closed subset 
$\Ad G(E)(\Om_{t,a})\subset\C{G}(E)$. In particular, $\Dt_{a,t}$ lies in 
$L^1_{loc}(\C{G}(E))\subset \C{D}(\C{G}(E))$. 
Similarly to \rn{char}, we define elements
${\Dt}_{a_0,\ka,t}:=e(G)\sum_{a}\lan \inv(a_0,a),\ka\ran \Dt_{a,t}\in 
L^1_{loc}(\C{G}(E))$ \label{Dta0kat} and ${\Dt}'_{a'_0,\ka,t}:=e(G')\sum_{a'}
\lan \inv(a'_0,a'),\ka\ran \Dt_{a',t}\in L^1_{loc}(\C{G'}(E))$.
\end{Not}

\begin{Lem} \label{L:stab}
$e(G){\Dt}_{a_0,\ka,t}$ is $(a_0,a'_0;\ka)$-equivalent to 
$e(G'){\Dt}'_{a'_0,\ka,t}$. 
\end{Lem}
\begin{proof}
By \rp{reform}, we have to show that for each $x_0\in\C{G}^{\sr}(E)$ and
$\ov{\xi}\in\pi_0(\wh{G_{x_0}}^{\Gm}/Z(\wh{G})^{\Gm})$ such that 
${\Dt}_{a_0,\ka,t}(x_0,\ov{\xi})\neq 0$, we have

$(i)$ $([a_{x_0}],\ov{\xi})\in\im\Pi_{\C{E}}$;

$(ii)$  if $\varphi:G\to G'$ is $(\C{E},[a_{x_0}],\ov{\xi})$-admissible, 
then for every stable conjugate $x'_0\in \C{G}'(E)$ of $x_0$, we have 
$
e(G')\Dt'_{a'_0,\ka,\theta}(x'_0,\ov{\xi})=\lan\frac{x_0,x'_0;\ov{\xi}}{a_0,a'_0;\ka}\ran 
e(G)\Dt_{a_0,\ka,\theta}(x_0,\ov{\xi})$.

Recall that
\[ 
e(G){\Dt}_{a_0,\ka,t}(x_0,\ov{\xi})=\sum_{x}\sum_a\lan 
\inv(x_0,x),\ov{\xi}\ran ^{-1}\lan \inv(a_0,a),\ka\ran{\Dt}_{a,t}(x),
\]
where $x$ runs over a set of representatives of $G(E)\bs [x_0]\subset G(E)\bs\C{G}(E)$. 
We identify $T$ with $a_0(T)\subset G$ and $\C{T}$ with $a_0(\C{T})\subset \C{G}$.
By \rl{lie} (a), the support of each $\Dt_{a,t}$ consists of elements, stably conjugate 
to $\C{T}(E)$. Hence replacing $x_0$ by a stable conjugate, 
we can assume that $x_0\in\C{T}(E)$. Then $G_{x_0}=T$, $a_{x_0}=a_0$, and thus 
$\ov{\xi}$ is an element of $\wh{T}^{\Gm}/Z(\wh{G})^{\Gm}$.

For a stable conjugate $x$ of $x_0$, we have 
${\Dt}_{a,t}(x)=1$ if and only if $\inv(x_0,x)=\inv(a_0,a)$. Therefore
$e(G){\Dt}_{a_0,\ka,t}(x_0,\ov{\xi})=\sum_a
\lan \inv(a_0,a),\ka\ov{\xi}^{-1}\ran$. Since the latter sum is non-zero, we get that
$\ov{\xi}=\ka$,  and $e(G){\Dt}_{a_0,\ka,t}(x_0,\ov{\xi})=|\wh{T}^{\Gm}/Z(\wh{G})^{\Gm}|$.
Since $\C{E}=\C{E}_{([a_0],\ka)}$, we get that $([a_{x_0}],\ov{\xi})=([a_0],\ka)\in\im
\Pi_{\C{E}}$, showing the assertion $(i)$.

To show $(ii)$, we can replace  $x'_0$ by a stably conjugate 
$a'_0(a_0^{-1}(x_0))\in a'_0(\C{T}(E))\subset\C{G}'(E)$.  
Then the same arguments show that 
$e(G'){\Dt}'_{a'_0,\ka,t}(x'_0,\ov{\xi})=|\wh{T}^{\Gm}/Z(\wh{G})^{\Gm}|$. Now the required assertion
$\lan\frac{x_0,x'_0;\ov{\xi}}{a_0,a'_0,\ka}\ran=1$ follows from equalities
 $a_{x_0}=a_0$, $a_{x'_0}=a'_0$ and $\ov{\xi}=\ka$.
\end{proof}

\begin{Emp} \label{E:four}
{\bf Fourier transform.} \label{ft} 
Fix an additive character $\psi:E\to\B{C}\m$ such that $\psi|_{\C{O}}$ is non-trivial, but 
$\psi|_{\frak{m}}$ is trivial. 

(a) The pairing $\lan\cdot,\cdot\ran$ and the measure $dx=|\om_{\C{G}}|$ on 
$\C{G}(E)$ give rise to the Fourier transform $\C{F}=\C{F}(\psi,\lan\cdot,\cdot\ran, dx)$
on $C_c^{\infty}(\C{G}(E))$.
Then $\C{F}$ \label{f} induces Fourier transforms on 
 $\C{S}(\C{G}(E))$ and $\C{D}(\C{G}(E))$ given by the formulas 
$\C{F}(fdx):=\C{F}(f)dx$ for each $f\in C_c^{\infty}(\C{G}(E))$ and
$\C{F}(F)(\phi):=F(\C{F}(\phi))$ for each  $\phi\in\C{S}(\C{G}(E))$ and 
$F\in\C{D}(\C{G}(E))$.

(b) For each $f_1,f_2\in C_c^{\infty}(\C{G}(E))$ we have 
$\int_{\C{G}(E)}f_1\C{F}(f_2)dx=\int_{\C{G}(E)}\C{F}(f_1)f_2dx$. Therefore
the embedding $C_c^{\infty}(\C{G}(E))\hra \C{D}(\C{G}(E))$ commutes with 
the Fourier transform.

(c) For each parahoric subalgebra $\C{G}_a\subset\C{G}(E)$, we denote by 
$\ov{\C{F}}=\ov{\C{F}}(\ov{\psi},\lan\cdot,\cdot\ran_a, \mu)$ the Fourier transform 
on $\C{L}_a(\fq)$, where the character $\ov{\psi}:\fq\to\B{C}\m$ is induced by $\psi$,
pairing $\lan\cdot,\cdot\ran_a$ is induced by $\lan\cdot,\cdot\ran$ (see \rl{qlogx} (a))
and $\mu(l)=1$ for each $l\in L_a(\fq)$.
\end{Emp}

\begin{Lem} \label{L:four}
Denote by $\C{I}^+$ \label{I+} the pro-nilpotent radical of an Iwahori subalgebra of $\C{G}$.
Then for each $u\in G(E)_{\tu}$, we have
\[
t_{a,\theta}(u)=\C{F}(\dt_{a,t})(\Phi(u))|\om_{\C{G}}|(\C{I}^+)^{-1}.
\]
\end{Lem}
\begin{proof} 
First we claim that $\C{F}(\dt_{a,t})(x)=0$ for each $x\in\C{G}(E)\sm\C{G}_{a}$, and 
$\C{F}(\dt_{a,t})(x)=\ov{\C{F}}(\ov{\dt}_{a,t})(\ov{x})|\om_{\C{G}}|(\C{G}_{a^+})$ for
each $x\in\C{G}(E)\sm\C{G}_{a}$. Indeed, since $\dt_{a,t}$ vanishes outside of $\C{G}_a$, 
we have an equality  
\[
\C{F}(\dt_{a,t})(x)=\int_{\C{G}(E)}\psi(\lan x, y\ran)\dt_{a,t}(y) dy=
\int_{\C{G}_a}\psi(\lan x, y\ran)\dt_{a,t}(y) dy
\]
for each $x\in\C{G}(E)$. Since  $\dt_{a,t}(y+y')=\dt_{a,t}(y)$ 
for each $y'\in \C{G}_{a^+}$, we conclude that $\C{F}(\dt_{a,t})(x)$ equals
\[
\sum_z\psi(\lan x, z\ran)\dt_{a,t}(z)\int_{\C{G}_{a^+}}\psi(\lan x, y\ran)dy,
\]
where $z$ runs over a set of representatives of $\C{G}_a/\C{G}_{a^+}$ in $\C{G}_a$.

The assumptions on $\psi$ and $\lan\cdot,\cdot\ran$  imply that $\C{G}_{a}$ is the orthogonal 
complement of $\C{G}_{a^+}$ with respect to the pairing $(x,y)\mapsto\psi(\lan x,y\ran)$.
Therefore $\C{F}(\dt_{a,t})(x)=0$ for each $x\notin\C{G}_a$, and 
\[
\C{F}(\dt_{a,t})(x)=\sum_{\ov{z}\in\C{L}_a(\fq)}\psi(\lan \ov{x}, \ov{z}\ran)
\ov{\dt}_{a,t}(\ov{z})|\om_{\C{G}}|(\C{G}_{a^+})=\ov{\C{F}}(\ov{\dt}_{{a},{t}})(\ov{x})
|\om_{\C{G}}|(\C{G}_{a^+})
\]
for each $x\in\C{G}_a$.

Now we are ready to prove the lemma. Assume first that $u\in G_{a,\tu}$. 
It follows from \rl{qlog} (b) that $\Phi(u)\in \C{G}_{a}$ and that $\Phi$ induces a quasi-logarithm 
$\ov{\Phi}_a:L_a\to\C{L}_a$ satisfying $\ov{\Phi(u)}=\ov{\Phi}_a(\ov{u})$.
Therefore $\C{F}(\dt_{a,t})(\Phi(u))$ equals
$\ov{\C{F}}(\ov{\dt}_{{a},{t}})(\ov{\Phi}_a(\ov{u}))|\om_{\C{G}}|(\C{G}_{a^+})$.
It now follows from a combination of \rt{Spr} (see Appendix A) and the equality   
$|\om_{\C{G}}|(\C{G}_{a^+})=q^{-\frac{1}{2}\dim (L_a/\ov{T})}|\om_{\C{G}}|(\C{I}^+)$
that $\C{F}(\dt_{a,t})(\Phi(u))$
equals $|\om_{\C{G}}|(\C{I}^+)\Tr\,\rho_{\ov{a},\ov{\theta}}(\ov{u})=
|\om_{\C{G}}|(\C{I}^+)t_{a,\theta}(u)$. 

Finally, assume that $u\in G(E)_{\tu}\sm G_{a,\tu}$. In this case, $t_{a,\theta}(u)=0$.
On the other hand, by \rp{ident}, $\Phi$ induces bijections 
$G(E)_{\tu}\isom \C{G}(E)_{\tn}$ and $G_{a,\tu}\isom\C{G}_{a,\tn}$, therefore 
$\Phi(u)\in\C{G}(E)_{\tn}\sm\C{G}_{a,\tn}$. Using the equality 
$\C{G}_{a,\tn}=\C{G}(E)_{\tn}\cap \C{G}_a$ from \rl{tnil} (b), we conclude that
$\Phi(u)\notin \C{G}_a$, hence $\C{F}(\Dt_{a,t})(\Phi(u))=0$. 
This completes the proof of the lemma.
\end{proof}

\begin{Cor} \label{C:four}
For each $a$, we have $\chi(\pi_{a,\theta})_{\tu}=
{}^0\Phi^*(\C{F}(\Dt_{a,t}))_{\tu}|\om_{\C{G}}|(\C{I}^+)^{-1}$.
\end{Cor}

\begin{proof}
Note first that ${}^0\Phi$ is \'etale, hence smooth, therefore 
the pullback ${}^0\Phi^*(\C{F}(\Dt_{a,t}))$ is defined.
Consider generalized functions $t_{a,\theta}\in C_c^{\infty}(G(E))\subset\C{D}(G(E))$
and $\dt_{a,t}\in C_c^{\infty}(\C{G}(E))\subset\C{D}(\C{G}(E))$.
In light of \re{four} (b), \rl{four} implies the equality of generalized functions  
\begin{equation} \label{Eq:four1}
(t_{a,\theta})_{\tu}={}^0\Phi^*(\C{F}(\dt_{a,t}))_{\tu}|\om_{\C{G}}|(\C{I}^+)^{-1}.
\end{equation}
Since $\pi_{a,\theta}=\Ind_{\wt{G_a}}^{G(E)}\rho_{a,\theta}$ is admissible, 
it follows from the formula for characters of induced representations 
that 
\begin{equation} \label{Eq:four2}
\chi(\pi_{a,\theta})=\sum_{g\in\wt{G_a}\bs G(E)}(\Inn g)^*(t_{a,\theta}).
\end{equation}
Explicitly, $\chi(\pi_{a,\theta})(\phi)=\sum_{g\in\wt{G_a}\bs G(E)}(\Inn g)^*(t_{a,\theta})(\phi)$
for each $\phi\in\C{S}(G(E))$, where only finitely many terms in the sum are non-zero.
On the other hand, by the very definition of $\Dt_{a,t}$, we have
\begin{equation} \label{Eq:four3}
\Dt_{a,t}=\sum_{g\in\wt{G_a}\bs G(E)}(\Ad g)^*(\dt_{a,t}).
\end{equation}
Since $\langle\cdot,\cdot\rangle$ and $\Phi$ are $G$-equivariant,
we get the equality
\begin{equation} \label{Eq:four4}
{}^0\Phi^*(\C{F}((\Ad g)^*\dt_{a,t}))={}^0\Phi^*(\Ad g)^*(\C{F}(\dt_{a,t}))=
(\Inn g)^*{}^0\Phi^*(\C{F}(\dt_{a,t})).
\end{equation}
Now our corollary is an immediate consequence of equalities
(\ref{Eq:four1})--(\ref{Eq:four4}).
\end{proof}

\begin{Emp} \label{E:proofunip}
\begin{proof}[Proof of \rt{unip}]
Let $\Phi':G'\to\C{G}'$ and  $\lan\cdot,\cdot\ran'$ be the  quasi-logarithm
map and the pairing on $\C{G}'$ induced by $\Phi$ and $\lan\cdot,\cdot\ran$, respectively
(see \rl{qlog} (a)). We denote by $\C{F}=\C{F}(\psi,\lan\cdot,\cdot\ran',|\om_{\C{G}'}|)$
the corresponding  Fourier transform $\C{G}'(E)$ and by $\C{I}'^+$
the pro-nilpotent radical of an Iwahori subalgebra of $\C{G}'(E)$.

By \rl{stab},   $e(G){\Dt}_{a_0,\ka,t}$ is $(a_0,a'_0;\ka)$-equivalent to 
$e(G'){\Dt}'_{a'_0,\ka,t}$. 
Using the equality $e(G)e(G')=e'(G)e'(G')$, it follows from \rt{Wa} (see Appendix B) that  
$\C{F}(\Dt_{a_0,\ka,t})$ is  $(a_0,a'_0;\ka)$-equivalent to 
$\C{F}(\Dt'_{a'_0,\ka,t})$. Hence by \rco{pullback}, the pullback
${}^0\Phi^*(\C{F}(\Dt_{a_0,\ka,t}))$ is $(a_0,a'_0;\ka)$-equivalent to 
${}^0\Phi'^*(\C{F}(\Dt'_{a'_0,\ka,t}))$.
By \rco{four}, we thus get that 
$|\om_{\C{G}}|(\C{I}^+)(\chi_{a_0,\ka,\theta})_{\tu}$ and 
$|\om_{\C{G}'}|(\C{I'}^+)(\chi_{a'_0,\ka,\theta})_{\tu}$ are $(a_0,a'_0;\ka)$-equivalent.
Since $|\om_{\C{G}}|(\C{I}^+)=|\om_{\C{G}'}|(\C{I'}^+)$ (see, for example, \cite[p. 632]{Ko4}),
the assertion follows.
\end{proof}
\end{Emp}

\subsection{Reduction formula} 
\begin*
\vskip 8truept
\end*
In this subsection we will assume that $G^{\der}=G^{\ssc}$. Our goal is to rewrite 
character $\chi(\pi_{a,\theta})$ in terms of restrictions
to topologically unipotent elements of the corresponding characters of the centralizers
$G_{\dt}(E)$. 

\begin{Lem} \label{L:assumptions}
 Assume that $G^{\der}=G^{\ssc}$. Then

(a) For each semisimple element $\dt\in G$, the centralizer $G_{\dt}$ is connected.

(b) The stabilizer in $G(E)$ of each $x\in\C{B}(G)$ is $G_x$.
\end{Lem}

\begin{proof} 
(a) was shown in \cite[Cor. 8.5]{St} when $G$ is semisimple, and in \cite[pp. 788--789]{Ko3} 
in the general case. 

(b) When $G$ is semisimple, the result was proved in \cite[Prop. 4.6.32]{BT}. 
For a general $G$, we can replace $E$ by an unramified extension so that $G$ is split over $E$.
Choose a split maximal torus $T\subset G$ such that $x\in\C{B}(T)$.
Since $Stab_{G(E)}(x)$ is compact, we see as in \rco{unrtori} that 
$Stab_{G(E)}(x)$ is contained in $G^{\der}(E)T(\C{O})$. Since
$T(\C{O})\subset Stab_{G(E)}(x)$, we get that $Stab_{G(E)}(x)$ is contained
in $Stab_{G^{\der}(E)}(x)T(\C{O})$, hence (use \cite[Prop. 4.6.32]{BT}) in
$(G^{\der})_xT(\C{O})=G_x$, as claimed.
\end{proof}


\begin{Not} 
(a) We will call an element $\gm\in G(E)$ {\em compact}, \label{compact} if it generates a relatively 
compact subgroup of $G(E)$. 

(b) We will call an element $\gm\in G(E)$ {\em topologically unipotent}, \label{topun} if the sequence   
$\{\gm^{p^n}\}_n$ converges to $1$. 

\end{Not} 

\begin{Cor} \label{C:compact}
(a) The set of compact elements of $G(E)$ is $\cup_{x\in\C{B}(G)}G_x$.

(b) The set of topologically unipotent elements of $G(E)$ is $G(E)_{\tu}$ (see \rn{tnil}).
\end{Cor}
\begin{proof}
As each $G_x$ is compact and every compact element of $G(E)$ stabilizes a point 
of $\C{B}(G)$, (a) follows from \rl{assumptions} (b). 
Since every topologically unipotent element is compact, (b) follows from (a) and the fact that 
$G_{x,\tu}$ is set of all topologically unipotent elements of $G_x$.
\end{proof}

The following result is a straightforward generalization of \cite[Lem 2, p. 226]{Ka2}. 

\begin{Lem} \label{L:jor} 
For every compact element $\gm\in G(E)$, there exists a unique decomposition $\gm={\dt}u=u\dt$ 
such that ${\dt}$ is of finite order prime to $p$, and $u$ is topologically unipotent. 
In particular, this decomposition is compatible with conjugation and field extensions.  
\end{Lem} 

\begin{Not} \label{N:jor}
The decomposition  $\gm={\dt}u$ from \rl{jor} is called the {\em topological
Jordan decomposition} \label{tjd} of $\gm$.
\end{Not}

\begin{Rem} \label{R:strss} 
If $\dt\in G(E)$ is an element of finite order prime to $p$, then $\dt$ is automatically semisimple.
\end{Rem}

The goal of this subsection is to prove the following result.

\begin{Prop} \label{P:reduction}
For every embedding $a:T\hra G$ and a compact element $\gm\in G(E)$ with topological
Jordan decomposition $\gm=\dt u$, we have the following formula
\begin{equation} \label{Eq:reduction}
e(G)F_{a,\theta}(\gm)=e(G_{\dt})\sum_{b}\theta(b^{-1}({\dt}))F_{b,\theta}(u).
\end{equation}
Here $b$ runs over a set of representatives of conjugacy classes of embeddings $T\hra G_{\dt}$,
whose composition with the inclusion $G_{\dt}\hra G$ is conjugate to $a$.
\end{Prop}

First we need to prove two preliminary results.

\begin{Lem} \label{L:cent}
(a) For each vertex $x$ of $\C{B}(G)$, we have $e(L_x)=e(G)$.

(b) Let $\dt\in G(E)$ be an element of finite order, prime to $p$. Then the centralizer
$G_{\dt}$ splits over $E^{\nr}$, and the building $\C{B}(G_{\dt})$ is canonically identified with
the set of invariants $\C{B}(G)^{\dt}\subset \C{B}(G)$.

(c) For each $x\in\C{B}(G_{\dt})\subset\C{B}(G)$, the parahoric subgroup
$(G_{\dt})_x$ is a subgroup of finite index in  $(G_x)_{\dt}$, and the canonical  
map $(G_x)_{\dt}\hra G_x\to L_x(\fq)$ induces isomorphisms
$(G_x)_{\dt}/(G_{\dt})_{x^+}\isom(L_x)_{\ov{\dt}}(\fq)$ and 
$(G_{\dt})_x/(G_{\dt})_{x^+}\isom(L_x)^{0}_{\ov{\dt}}(\fq)$.

(d) Let $\dt,\dt'\in G_x$ be two elements of finite orders prime to $p$. 
Then $\dt$ and $\dt'$ are $G_x$-conjugate if and only if their reductions 
$\ov{\dt},\ov{\dt}'\in L_x(\fq)$ are $L_x(\fq)$-conjugate.
\end{Lem}

\begin{proof}

(a) For each $x\in\C{B}(G)$, the maximal split torus of $L_x$ is the reduction of that of $G$, 
therefore $\rk_{\fq}(L_x)=\rk_E(G)$. If, moreover, $x$ is a vertex, then  
$\rk_{\fq}(Z(L_x)^0)=\rk_E(Z(G)^0)$, hence $e(L_x)=e(G)$.

(b) The second assertion is shown in \cite{PY}. For the first, recall that $G_{E^{\nr}}$ 
splits, hence there exists a split maximal torus $T\subset G_{E^{\nr}}$. 
Choose $g\in G(\ov{E})$ such that $g\dt g^{-1}\in T(\ov{E})$. Then $g\dt g^{-1}$ is  
of finite order, prime to $p$, therefore it follows from Hensel's lemma that 
$g\dt g^{-1}\in T(\C{O}_{E^{\nr}})\subset T(E^{\nr})$. Hence $g$ gives rise to a cocycle 
$\sigma\mapsto g^{-1}{}^{\sigma}g\in G_{\dt}(\ov{E})$ over $E^{\nr}$.
Since $H^1(E^{\nr},G_{\dt})=0$, there exists $h\in G_{\dt}(\ov{E})$ such that
$h^{-1}{}^{\sigma}h=g^{-1}{}^{\sigma}g$ for each $\sigma\in\Gal(\ov{E}/E^{\nr})$.
It follows that $gh^{-1}\in G(E^{\nr})$, and $(gh^{-1})\dt (gh^{-1})^{-1}=g\dt g^{-1}\in T(E^{\nr})$.
Therefore  $(gh^{-1})^{-1}T(gh^{-1})$ is a split maximal torus of 
$(G_{\dt})_{E^{\nr}}$.

(c) As $(G_x)_{\dt}\subset G_{\dt}(E)$ is the stabilizer of $x\in\C{B}(G_{\dt})$,
it is compact. Therefore $(G_{\dt})_x$ is a subgroup of finite index in  $(G_x)_{\dt}$, thus 
the corresponding group scheme $\un{G_{\dt}}_x$ over $\C{O}$ is the connected component 
of $(\un{G}_x)_{\dt}$. In particular,  $(\un{G}_x)_{\dt}$ is smooth over $\C{O}$. 
Therefore by Hensel's lemma, the reduction map  $G_x\to\ov{G}_x(\fq)$
 surjects $(G_x)_{\dt}=(\un{G}_x)_{\dt}(\C{O})$ onto 
$(\ov{G}_x)_{\dt}(\fq)=(\un{G}_x)_{\dt}(\fq)$. 

As $L_x$ is the quotient $\ov{G}_x/R_{u}(\ov{G}_x)$, 
we see that $(L_x)_{\ov{\dt}}$ is the quotient $(\ov{G}_x)_{\dt}/R_{u}(\ov{G}_x)_{\dt}$. 
Since $R_{u}(\ov{G}_x)_{\dt}$ is a connected group over $\fq$, Lang's theorem implies that 
the projection $\ov{G}_x(\fq)\to L_x(\fq)$ surjects 
$(\ov{G}_x)_{\ov{\dt}}(\fq)$ onto $(L_x)_{\ov{\dt}}(\fq)$. Therefore the projection 
$G_x\to L_x(\fq)$ induces a surjection $(G_x)_{\dt}\to (L_x)_{\ov{\dt}}(\fq)$, whose kernel is 
$(G_x)_{\dt}\cap G_{x^+}=(G_{\dt})_{x^+}$. This shows the first isomorphism, while the proof 
of the second one is similar but easier.

(d) The ``only if'' assertion is clear. Assume now that 
$\ov{\dt}$ and $\ov{\dt}'$ are $L_x(\fq)$-conjugate. 
Let us first show that $\dt$ and $\dt'$ are $\un{G}_x(\C{O}_{E^{\nr}})$-conjugate.
For this we can replace $E$ by an unramified extension, so that $G$, $G_{\dt}$ and $G_{\dt'}$  
are split over $E$ (use (b)). Since $\dt$ lies in $G_x$, we get that $x$ belongs to $\C{B}(G)^{\dt}=\C{B}(G_{\dt})$ 
(use (b)). Therefore there exists a split maximal torus $T\subset G_{\dt}\subset G$ such that 
 $x\in\C{B}(T)$. Similarly,  there exists a split maximal torus $T'\subset G_{\dt'}\subset G$ such that 
 $x\in\C{B}(T')$. By a property of buildings, there exists $g\in G_x$ such that 
$gTg^{-1}=T'$. Replacing $\dt'$ by $g^{-1}\dt' g$, we may assume that $\dt,\dt'\in T(E)$.

Next we observe that the projection 
$\Norm_{G_x}(T)\to\Norm_{L_x(\fq)}(\ov{T})$ is surjective. Indeed, for each
$\ov{g}\in \Norm_{L_x(\fq)}(\ov{T})\subset L_x(\fq)$, choose a representative
$g\in G_x$. Then $\ov{g T g^{-1}}=\ov{T}$, hence by \cite[Lem 2.2.2]{DB}, there exists 
$h\in G_{x^+}$ such that $h(g T g^{-1})h^{-1}=T$. In other words, 
$hg\in  \Norm_{G_x}(T)$ is a preimage of $\ov{g}$.

By the assumption, $\ov{\dt},\ov{\dt'}\in\ov{T}(\fq)$ are conjugate in $L_x(\fq)$, therefore
they are conjugate in $\Norm_{L_x(\fq)}(\ov{T})$. Hence 
there exists $g\in\Norm_{G_x}(T)$ such that $\ov{g^{-1}\dt' g}=\ov{\dt}$.
But the projection $T(\C{O})\to\ov{T}({\fq})$ defines a bijection between 
elements of $T(\C{O})$ of finite order prime to $p$ and elements of $\ov{T}({\fq})$.
Hence $g^{-1}\dt' g=\dt$, implying that 
$\dt$ and $\dt'$ are conjugate by an element of $\un{G}_x(\C{O}_{E^{\nr}})$.

To show that $\dt$ and $\dt'$ are conjugate by $G_x$, consider the closed subscheme $Z$ (resp. $Z'$)
of $\un{G}_x$ (resp. $L_x$) consisting of elements $g$ (resp. $\ov{g}$) such that 
$g\dt g^{-1}=\dt'$ (resp. $\ov{g}\ov{\dt}\ov{g}^{-1}=\ov{\dt}'$). By the assumption, 
$Z'(\fq)\neq\emptyset$, and we have to show that
$Z(\C{O})\neq\emptyset$. By the shown above, $Z(\C{O}_{E^{\nr}})\neq\emptyset$. 
Thus  $Z$ and $(\un{G}_x)_{\dt}$ are isomorphic over $\C{O}_{E^{\nr}}$. In particular, 
$Z$ is smooth over $\C{O}$, thus by Hensel's lemma it suffice to show that 
the projection $Z(\fq)\to Z'(\fq)$ is surjective.

Denote by $\ov{Z}\subset\ov{G}_x$ the special fiber of $Z$. Since all fibers of the projection
$\ov{Z}\to Z'$ are principal homogeneous spaces for the connected group $R_u(\ov{G}_x)_{\dt}$, 
the surjectivity of the projection $Z(\fq)=\ov{Z}(\fq)\to Z'(\fq)$ follows from Lang's theorem.
\end{proof}

\begin{Lem} \label{L:DL}
For every $\gm\in G_a$ with topological Jordan decomposition $\gm={\dt}u$, we have 
an equality
\begin{equation} \label{Eq:DL}
e(G) t_{a,\theta}(\gm)=e(G_{\dt})\sum_{b}
\sum_{h\in (G_{\dt})_b\bs (G_a)_{\dt}}\theta(b^{-1}({\dt}))t_{b,\theta}(huh^{-1}),
\end{equation}
where the $b$ runs over a set of representatives of conjugacy classes of embeddings
$T\hra G_{\dt}$, which are $G_a$-conjugate to $a:T\hra G$.
\end{Lem}

\begin{proof}
We start from the following claim
\begin{Cl} \label{C:bij}
The correspondence $b\mapsto \ov{b}$ induces a bijection 
between the set of conjugacy classes of embeddings $T\hra G_{\dt}$ 
which are $G_a$-conjugate to $a:T\hra G$ and the set of conjugacy classes of embeddings 
$\ov{b}:\ov{T}\hra (L_a)_{\ov{{\dt}}}$, which are $L_a(\fq)$-conjugate to 
$\ov{a}:\ov{T}\hra L_a$.
\end{Cl}
\begin{proof}
For simplicity of notation, we identify $T$ with $a(T)$ and $\ov{T}$ with $\ov{a}(\ov{T})$. Then 
the maps $b\mapsto t:=b^{-1}({\dt})$ and
$\ov{b}\mapsto \ov{t}:=\ov{b}^{-1}(\ov{{\dt}})$ identify our sets with the sets of 
elements $t\in T(E)$, which are $G_a$-conjugate to $\dt$, and elements $\ov{t}\in \ov{T}(\fq)$, 
which are $L_a(\fq)$-conjugate to $\ov{{\dt}}$, respectively.
Since the reduction map $t\mapsto \ov{t}$ induces a bijection between elements of $T(E)$ of 
finite order prime to $p$ and elements of $\ov{T}(\fq)$, we get the 
injectivity. The surjectivity follows from \rl{cent} (d).
\end{proof}

Now we are ready to prove the lemma.
By \rl{cent} and \rcl{bij}, the right hand side of (\ref{Eq:DL})
equals
\begin{equation} \label{Eq:DL1}
e((L_a)^0_{\ov{\dt}})\sum_{\ov{b}}
\sum_{\ov{h}\in (L_a)^0_{\ov{\dt}}(\fq)\bs (L_a)_{\ov{\dt}}(\fq)}
\ov{\theta}(\ov{b}^{-1}(\ov{\dt}))\Tr\,\rho_{\ov{b},\ov{\theta}}(\ov{h}\ov{u}\ov{h}^{-1}),
\end{equation}
where $\ov{b}$ runs over a set of representatives of $(L_a)_{\ov{\dt}}(\fq)$-conjugacy classes of 
embeddings $\ov{T}\hra (L_a)_{\ov{\dt}}$, which are $L_a(\fq)$-conjugate to $\ov{a}$. 

Next note that (\ref{Eq:DL1}) can be rewritten as 
\begin{equation} \label{Eq:DL2}
e((L_a)^0_{\ov{\dt}})\sum_{\ov{b}:\ov{T}\hra (L_a)^0_{\ov{\dt}}}
\theta(\ov{b}^{-1}(\ov{\dt}))\Tr\,\rho_{\ov{b},\ov{\theta}}(\ov{u}),
\end{equation}
where $\ov{b}$ runs over a set of representatives of conjugacy classes of 
embeddings, $L_a(\fq)$-conjugate to $\ov{a}$. By  the 
formula of Deligne--Lusztig \cite[Thm 4.2]{DL},  (\ref{Eq:DL2}) equals  
$e(L_a)\Tr\rho_{\ov{a},\ov{\theta}}(\ov{\gm})$. Hence it is equal to $e(G) t_{a,\theta}(\gm)$,
as claimed.
\end{proof} 
\begin{Emp}
\begin{proof}[Proof of \rp{reduction}]


Notice first that the map $b\mapsto b^{-1}({\dt})$ embeds the set of conjugacy classes of 
embeddings $b:T\hra G_{\dt}$,  conjugate to $a:T\hra G$, into the finite set 
$\{t\in T(E)\,|\,t^{\ord{\dt}}=1\}$. Therefore the sum in (\ref{Eq:reduction}) is finite.

Fix a set of representatives $J\subset G(E)$ of double classes 
$\wt{G_a}\bs G(E)/G_{\dt}(E)$. 
For every $h\in J$, put $\gm_{h}={h}\gm h^{-1}$, and let $\gm_h={\dt}_h u_h$ be the
topological Jordan decomposition of $\gm_h$.

By \cite[Lem. 23]{HC} (compare \rn{deff} (b)), for each sufficiently large 
compact modulo center $\wt{G_a}$-bi-invariant subset $\Om\subset G(E)$, 
we have
\[
F_{a,\theta}(\gm)= e(G)\sum _{g\in \wt{G_a}\bs \Om}t_{a,\theta}(g\gm g^{-1}).
\]

Since $G(E)$ decomposes as a disjoint union $\sqcup_{h\in J}
\wt{G_a} h G_{\dt}(E)= \sqcup_{h\in J} \wt{G_a} G_{\dt_h}(E)h$,
we have a finite decomposition
$\wt{G_a}\bs \Om=\sqcup_{h\in J}\wt{G_a}\bs 
[\wt{G_a} G_{\dt_h}(E) h\cap\Om]$. Therefore
\[
F_{a,\theta}(\gm)= e(G)\sum_{h\in J} \sum _{g\in \wt{G_a}\bs 
[\wt{G_a} G_{\dt_h}(E) h\cap\Om]}t_{a,\theta}(g\gm g^{-1}).
\]
Using the identifications
\[
\wt{G_a}\bs [\wt{G_a} G_{\dt_h}(E) h\cap\Om]=
(\wt{G_a})_{\dt_h}\bs [G_{\dt_h}(E) h\cap\Om]=
(\wt{G_a})_{\dt_h}\bs [G_{\dt_h}(E)\cap\Om h^{-1}]h,
\]
we get that $F_{a,\theta}(\gm)$ equals 
\begin{equation} \label{Eq:red2}
\sum_{h\in J} e(G)
\sum_{g\in (\wt{G_a})_{\dt_h}\bs [G_{\dt_h}(E)\cap\Om h^{-1}]} 
t_{a,\theta}(g\gm_h g^{-1}).
\end{equation}

Using \rco{tilde}, we see that for each embedding $b:T\hra G_{{\dt}_h}$ which is 
$G_a$-conjugate to $a$, the group $\wt{(G_{\dt_h})_b}$ is contained in $(\wt{G_a})_{\dt_h}$
and we have a natural isomorphism $(G_{\dt_h})_b\bs (G_a)_{\dt_h}\cong
\wt{(G_{\dt_h})_b}\bs (\wt{G_a})_{\dt_h}$. Then by \rl{DL} applied to 
$\gm_h=\dt_h u_h$, the contribution of each $h\in J$ to
(\ref{Eq:red2}) equals
\begin{equation} \label{Eq:red3}
e(G_{{\dt}_h})\sum_{b:T\hra G_{{\dt}_h}}\theta(b^{-1}({\dt}_h))
\sum_{g\in \wt{(G_{\dt_h})_b}\bs [G_{\dt_h}(E)\cap\Om h^{-1}]} t_{b,\theta}(g u_h g^{-1}),
\end{equation} 
where $b$ runs over the a set of representatives of conjugacy classes of embeddings, which are 
$G_a$-conjugate to $a:T\hra G$.  

Conjugating by $h^{-1}$, we can rewrite (\ref{Eq:red3}) in the form
\begin{equation} \label{Eq:red4}
e(G_{{\dt}})\sum_{b:T\hra G_{{\dt}}}\theta(b^{-1}({\dt}))
\sum_{g\in \wt{(G_{\dt})_b}\bs [G_{\dt}(E)\cap h^{-1}\Om]} t_{b,\theta}(gu g^{-1}),
\end{equation}
where $b$ runs over a set of representatives of conjugacy classes of embeddings 
such that $a=gbg^{-1}$ for some $g\in \wt{G_a} h G_{\dt}(E)$. 
In particular, $b$ has a non-trivial contribution to (\ref{Eq:red4}) only for a unique $h\in J$,
which we denote by $h_b$. It follows that $F_{a,\theta}(\gm)$ equals  
\begin{equation} \label{Eq:red5}
e(G_{{\dt}})\sum_{b:T\hra G_{\dt}}\theta(b^{-1}({\dt})) 
\sum_{g\in \wt{(G_{\dt})_b}\bs [G_{\dt}(E)\cap h^{-1}_b\Om]} 
t_{b,\theta}(gu g^{-1}),
\end{equation}
where $b$ runs over a (finite) set of representatives of conjugacy classes to embeddings,
which are $G(E)$-conjugate to $a$.

Replacing $b$'s by their $G_{\dt}(E)$-conjugates, we can assume that $a=gbg^{-1}$ for some  
$g\in\wt{G_a} h_{b}$. Then $h_b {(G_\dt)_b}h_b^{-1}\subset{G_a}$, hence  
$h_b \wt{(G_\dt)_b}h_b^{-1}\subset\wt{G_a}$ (use \rco{tilde}). It follows that the subset 
$G_{\dt}(E)\cap h_b^{-1}\Om\subset G_{\dt}(E)$
is  compact modulo center, $\wt{(G_\dt)_b}$-invariant from the left, and
$h_b \wt{(G_\dt)_b}h_b^{-1}$-invariant from the right.
Since the number of $b$'s is finite, it follows from \cite[Lem. 23]{HC}, 
that for each sufficiently large $\Om\subset G(E)$, the contribution of each $b$ to 
(\ref{Eq:red5}) equals $\theta(b^{-1}({\dt}))F_{b,\theta}(u)$. This completes the proof.
\end{proof}
\end{Emp}

\subsection{Endoscopy for $G$ and $G_{\dt}$} \label{SS:technical}
\begin*
\vskip 8truept
\end*

Let $G$ be a connected reductive group over $E$ such that 
$G^{\der}=G^{\ssc}$, $\dt\in G(E)$ a semisimple element, and 
$\iota:G_{\dt}\hra G$ the canonical embedding. 
In this subsection we will compare endoscopic triples for $G$ and for $G_{\dt}$.

\begin{Emp} \label{E:technical}
(a) Similarly to \re{pair} (b), there exists a natural embedding 
$Z(\wh{G})\hra Z(\wh{G_{\dt}})$. 
Indeed, every maximal torus $T$ of $G_{\dt}$ is a maximal torus of $G$, and
the set of roots of $(G_{\dt},T)$ equals the set of those roots of $(G,T)$ 
which vanish on $\dt$. Hence the set of coroots (hence also of roots) of 
$(\wh{G_{\dt}},\wh{T})$ is naturally a subset of those of $(\wh{G},\wh{T})$. 
Thus $Z(\wh{G})$ is naturally a subgroup of $Z(\wh{G_{\dt}})$.

(b) Fix an embedding $a_{\dt}:T\hra G_{\dt}$ of a maximal torus and 
$\ka\in\wh{T}^{\Gm}/Z(\wh{G})^{\Gm}$. Set
$a:=\iota\circ a_{\dt}:T\hra G$, let $\ov{\ka}\in\pi_0(\wh{T}^{\Gm}/Z(\wh{G})^{\Gm})$
and $\ov{\ka}'\in\pi_0(\wh{T}^{\Gm}/Z(\wh{G_{\dt}})^{\Gm})$ be the classes of $\ka$, 
and put $\C{E}:=\C{E}_{([a],\ka)}=(H,[\eta],\ov{s})$ and 
$\C{E}':=\C{E}_{([a_{\dt}],\ka)}=(H',[\eta'],\ov{s}')$. 
  
Also we fix embeddings $c':T\hra H'$ and $c:T\hra H$ 
such that $\Pi_{\C{E}}([c])=([a],\ov{\ka})$ and  
$\Pi_{\C{E}'}([c'])=([a_{\dt}],\ov{\ka}')$ (see \re{conend} (b)).
\end{Emp}

\begin{Lem} \label{L:embtori}

(a) There is a natural $\Gm$-equivariant embedding $Z(\wh{H})\hra Z(\wh{H}')$ 
mapping $Z(\wh{G})$ into $Z(\wh{G_{\dt}})$ and an embedding $W(H')\hra W(H)$,
both of which depend on $c$ and $c'$. The induced map 
$\pi_0( Z(\wh{H})^{\Gm}/Z(\wh{G})^{\Gm})\to 
\pi_0( Z(\wh{H}')^{\Gm}/Z(\wh{G_{\dt}})^{\Gm})$ sends $\ov{s}$ to $\ov{s}'$.

(b) There is a natural map $[b']\mapsto [b]$, depending on $c$ and $c'$, from the set of 
stable conjugacy classes of embeddings of maximal tori $S\hra H'$ to those of 
embeddings $S\hra H$. 

(c) In the notation of (b), we have $[b]_G=\iota\circ[b']_{G_{\dt}}$, and 
$Z_{\wh{[b]}}:Z(\wh{H})\hra\wh{S}$ is the restriction of 
$Z_{\wh{[b']}}:Z(\wh{H}')\hra\wh{S}$ (see (a)).
In particular, $\ov{\ka}_{[b']}\in\pi_0(\wh{S}^{\Gm}/Z(\wh{G_{\dt}})^{\Gm})$ is the image of
$\ov{\ka}_{[b]}\in\pi_0(\wh{S}^{\Gm}/Z(\wh{G})^{\Gm})$.

(d) Let $[b'_i]$ be two stable conjugacy classes of embeddings of 
maximal tori $T_i\hra H'$, and let $[b_i]$ be the corresponding 
stable conjugacy classes of embeddings $T_i\hra H$. Then the following 
diagram is commutative
\[
\CD
Z(\wh{H}) @>{\iota([b_1],[b_2])}>>\{(T_1\times T_2)/Z(G)\}\:\wh{}\\
@VVV          @AAA\\
Z(\wh{H'}) @>{\iota([b'_1],[b'_2])}>> \{(T_1\times T_2)/Z(G_{\dt})\}\:\wh{}.
\endCD
\]
(Here the left vertical map was defined in (a), the right one is induced by the inclusion 
$Z(G)\hra Z(G_{\dt})$, and the horizontal maps are the homomorphisms (\ref{E:moriota}) from
\re{kainv}.)
\end{Lem}
\begin{proof}
(a) Embed $T$ into $G, G_{\dt}, H$ and $H'$ by $a,a_{\dt},c$ and $c'$, 
respectively. Then the set of roots of $(\wh{G_{\dt}},\wh{T})$ 
(resp. of $(\wh{H},\wh{T})$, resp. of $(\wh{H}',\wh{T})$) is the set of those roots 
$\hat{\al}$ of $(\wh{G},\wh{T})$ such that ${\al}(\dt)=1$ 
(resp. $\hat{\al}(\ka)=1$, resp. ${\al}(\dt)=1$ and $\hat{\al}(\ka)=1$). 
In particular, the set of roots of $(\wh{H}',\wh{T})$ is canonically a subset
of those of $(\wh{H},\wh{T})$. This gives us the required embeddings $W(H')\hra W(H)$ and 
$Z(\wh{H})\hra Z(\wh{H}')$.
The last assertion follows from the fact 
that both $\ov{s}$ and $\ov{s}'$ are the classes of $\ka$.

(b) Choose an embedding $b':S\hra H'$ from $[b']$,
and identify $S$ with $b'(S)\subset H'$ and $T$ with $c'(T)\subset H'$.
Choose $g\in H'(\ov{E})$ such that $gSg^{-1}=T$. Then $\Inn g$ defines an 
isomorphism $S_{\ov{E}}\isom T_{\ov{E}}$.
Let $b:S_{\ov{E}}\hra H$ be the composition $c\circ\Inn g$. We claim that
the $H(\ov{E})$-conjugacy class $[b]$ of $b$ is $\Gm$-invariant and independent of 
the choices of $g$ and $b'$. 

If $g'\in H'(\ov{E})$ is another element such that  $g'Sg'^{-1}=T$, then
$g^{-1}g'\in\Norm_{H'}(S)$, and $b':=c\circ\Inn g'$ equals
$b\circ \Inn(g^{-1}g')$. But $\Inn(g^{-1}g'):S_{\ov{E}}\to S_{\ov{E}}$
is induced by an element of $W(H')\subset W(H)$, therefore $b'$ is conjugate 
to $b$. Thus $[b]$ is independent of the choice of $g$.  For each $\si\in\Gm$, we have 
${}^{\si}b=c\circ\Inn({}^{\si}g)$ and ${}^{\si}gS({}^{\si}g)^{-1}=T$.
Hence from the shown above, ${}^{\si}b$ is conjugate to $b$. 
Finally, if $b'$ is replaced by $\Inn h\circ b'$ and $g$ by $gh^{-1}$ for some $h\in H'(\ov{E})$, 
then the resulting isomorphism $S_{\ov{E}}\isom T_{\ov{E}}$ 
(and, therefore, $b$) do not change.

(c) Since the assertion is over $\ov{E}$, we can identify  
$S_{\ov{E}}\isom T_{\ov{E}}$,
as in (b), and thus replace $[b']$  by $[c']$ and $[b]$ by $[c]$.
Now the first assertion follows from the fact that 
$\iota\circ [c']_{G_{\dt}}=\iota\circ[a_{\dt}]=[a]=[c]_G$, while the  second
one was the definition of the embedding $Z(\wh{H})\hra Z(\wh{H}')$.

(d)  Mimicking the proof of (a) and (b), we see that   
there exists a $\Gm$-equivariant embedding 
$ Z(\wh{H^2/Z(G)})\hra Z(\wh{(H')^2/Z(G)})$ 
characterized by the following property:
For each pair $[c'_i]$ of stable conjugacy classes of embeddings of maximal 
tori $S_i\hra H'$ with the corresponding stable conjugacy classes $[c_i]$
of embeddings of maximal tori $S_i\hra H$, the embedding
$Z_{\wh{[c_1,c_2]}}:Z(\wh{H^2/Z(G)})\hra\{(S_1\times S_2)/Z(G)\}\:\wh{}\:\:$ is 
the restriction of 
$Z_{\wh{[c'_1,c'_2]}}:Z(\wh{(H')^2/Z(G)})\hra\{(S_1\times S_2)/Z(G)\}\:\wh{}$.
Then our diagram extends to the following diagram
\[
\CD
Z(\wh{H}) @>{\mu_H}>> Z(\wh{H^2/Z(G)})@>Z_{\wh{[b_1,b_2]}}>> 
\{(T_1\times T_2)/Z(G)\}\:\wh{}\\
@VVV          @VVV    @|\\
Z(\wh{H'}) @>{\mu_{H'}}>> Z(\wh{(H')^2/Z(G)})@>Z_{\wh{[b'_1,b'_2]}}>> 
\{(T_1\times T_2)/Z(G)\}\:\wh{}\\
@|          @AAA    @AAA\\
Z(\wh{H'}) @>{\mu_{H'}}>> Z(\wh{(H')^2/Z(G_{\dt})})@>Z_{\wh{[b'_1,b'_2]}}>> 
\{(T_1\times T_2)/Z(G_{\dt})\}\:\wh{}.
\endCD
\]
It remains to show that each inner square of the diagram is commutative. 
The commutativity of the top right square follows from the characterization of the embedding 
 $ Z(\wh{H^2/Z(G)})\hra Z(\wh{(H')^2/Z(G)})$. The commutativity of the 
top left square follows
 from the characterization of the vertical maps and the fact that both $\mu_H$ and $\mu_{H'}$ 
are restrictions of $\mu_T:\wh{T}\hra\wh{T^2/Z(G)}$. The commutativity of the two bottom 
squares is clear.
\end{proof}

\begin{Cor} \label{C:embtori}
 Let $[b'_1]$ and $[b'_2]$ be stable conjugacy classes of embeddings of maximal tori 
$S\hra H'$ such that $[b'_1]_{G_{\dt}}=[b'_2]_{G_{\dt}}$, and let $[b_1]$ and $[b_2]$ be the 
corresponding stable conjugacy classes of embeddings $S\hra H$.

Then $[b_1]_{G}=[b_2]_{G}$, and the image of $\ka\left(\frac{[b'_1]}{[b'_2]}\right)(\ov{s}')\in 
\pi_0(\wh{S/Z(G_{\dt})}^{\Gm})$ in $\pi_0(\wh{S/Z(G)}^{\Gm})$
equals $\ka\left(\frac{[b_1]}{[b_2]}\right)(\ov{s})$.
\end{Cor}
\begin{proof}
The first assertion follows from \rl{embtori} (c).
For the second one, choose a representative $\check{s}\in\pi_0(Z(\wh{H})^{\Gm})$ of 
$\ov{s}$, and let $\check{s}'\in\pi_0(Z(\wh{H}')^{\Gm})$ be the image of $\check{s}$. 
Then $\ka(\frac{[b_1]}{[b_2]})(\ov{s})=\nu_S(\iota([b_1],[b_2])(\check{s}))$
and  $\ka(\frac{[b'_1]}{[b'_2]})(\ov{s}')=\nu_S(\iota([b'_1],[b'_2])(\check{s}'))$.
So the assertion follows from \rl{embtori} (d). 
\end{proof}

\begin{Emp}
Let $\varphi:G\to G'$ be an inner twisting such that $\dt':=\varphi(\dt)$ belongs to $G'(E)$.
 For each $i=1,2$,  
let $(a_{\dt})_i:T_i\hra G_{\dt}$ and $(a'_{\dt})_i:T_i\hra G'_{\dt'}$ be stably conjugate 
embeddings of maximal tori, and let 
$b'_i:T_i\hra H'$ be an embedding of a maximal torus compatible with $(a_{\dt})_i$.

Let $\iota':G'_{\dt'}\hra G'$ be the natural embedding, and for each  $i=1,2$,
set  $a_i:=\iota\circ(a_{\dt})_i:T_i\hra G$ and 
$a'_i:=\iota'\circ(a'_{\dt})_i:T_i\hra G'$, and denote by 
$[b_i]$ the stable conjugacy class of 
embeddings of maximal tori $T_i\hra H'$ corresponding to $[b'_i]$ 
(and compatible with $a_i$ by \rl{embtori} (b)).
\end{Emp}

\begin{Lem} \label{L:inv2}
(a)  The image of 
 $\ov{\inv}((a_1,a'_1);(a_2,a'_2))\in H^1(E,(T_1\times T_2)/Z(G))$ in 
$H^1(E,(T_1\times T_2)/Z(G_{\dt}))$ equals 
$\ov{\inv}((a_{\dt})_1,(a'_{\dt})_1);((a_{\dt})_2,(a'_{\dt})_2))$.

(b) We have an equality $\lan \frac{a_1,a'_1;[b_1]}{a_2,a'_2;[b_2]}\ran_{\C{E}}=
\lan\frac{(a_{\dt})_1,(a'_{\dt})_1;[b'_1]}{(a_{\dt})_2,(a'_{\dt})_2;[b'_2]}
\ran_{\C{E}'}$.
\end{Lem}

\begin{proof}
(a) Follows immediately from the definition of the invariant.

(b) Let $\check{s}$ and $\check{s}'$ be as in the proof of \rco{embtori}. Then 
by \rl{embtori} (a) and (d), we obtain that  
$\ka([b_1],[b_2])(\check{s})\in \pi_0((\{(T_1\times T_2)/Z(G)\}\:\wh{}\:)^{\Gm})$ 
is the image of  $\ka([b'_1],[b'_2])(\check{s}')\in 
\pi_0((\{(T_1\times T_2)/Z(G_{\dt})\}\:\wh{}\:)^{\Gm})$. Now the assertion follows 
from (a) and the functoriality of the Tate--Nakayama duality. 
\end{proof}

\begin{Emp} \label{E:compatib}
Let $d_{\dt}:S\hra G_{\dt}$ be an embedding of a maximal torus, $\ov{\xi}$ an element of 
$\pi_0(\wh{S}^{\Gm}/Z(\wh{G})^{\Gm})$, $d:=\iota\circ d_{\dt}:S\hra G$,
$\ov{\xi}_{\dt}\in\pi_0(\wh{S}^{\Gm}/Z(\wh{G_{\dt}})^{\Gm})$ the class of $\ov{\xi}$, and
$\varphi:G\to G'$ an $(\C{E},[d],\ov{\xi})$-admissible inner twisting.
\end{Emp}
\begin{Lem} \label{L:compatib}
Assume that there exists 
$[b']\in\Pi_{\C{E}'}^{-1}([d_{\dt}],\ov{\xi}_{\dt})$ such that the corresponding stable 
conjugacy class $[b]$ of embeddings $S\hra H$ satisfies $\Pi_{\C{E}}([b])=([d],\ov{\xi})$.

If there exists a stably conjugate embedding $d':S\hra G'$ of $d$, then there 
exists a stable conjugate $d'$ of $d$ for which  $G'_{d'(d^{-1}(\dt))}$ is an 
$(\C{E}',[d_{\dt}],\ov{\xi}_{\dt})$-admissible inner form of $G_{\dt}$.
\end{Lem}

\begin{proof}
For shortness, we will denote $Z(\C{E},[d],\ov{\xi})\subset Z(\wh{G^{\ad}})^{\Gm}$ by  $Z$
and $Z(\C{E}',[d_{\dt}],\ov{\xi}_{\dt})\subset Z(\wh{G_{\dt}^{\ad}})^{\Gm}$
by $Z_{\dt}$. Embedding $d_{\dt}:S\hra G_{\dt}\subset G$ 
induce homomorphisms
\[
H^1(E,(G_{\dt})^{\ad})\overset{g}{\lla}H^1(E,S/Z(G))\overset{f}{\lra}
H^1(E,G^{\ad}),
\] 
hence dual homomorphisms $Z(\wh{(G_{\dt})^{\ad}}^{\Gm})
\overset{g^D}{\lra}\pi_0(\wh{S/Z(G)}^{\Gm})
\overset{f^D}{\lla}Z(\wh{G^{\ad}})^{\Gm}$. 

First we claim that the assertion of the lemma is equivalent to the inclusion
\begin{equation} \label{E:inclusion}
g^D(Z_{\dt})\cap \im f^D\subset f^D(Z).
\end{equation}
Indeed, put $x:=\inv(G,G')$. By our assumptions, $x\in  Z^{\perp}\cap\im f$, 
and the lemma is equivalent to the assertion that there exists $y\in f^{-1}(x)$ such that 
$g(y)\in (Z_{\dt})^{\perp}$.  Equivalently, we have to show that  
$f^{-1}(x)\cap g^{-1}((Z_{\dt})^{\perp})\neq\emptyset$. Since 
$g^{-1}((Z_{\dt})^{\perp})=[g^D(Z_{\dt})]^{\perp}$, we have to check that  
$x$ belongs to $f([g^D(Z_{\dt})]^{\perp})=[(f^D)^{-1}(g^D(Z_{\dt}))]^{\perp}$. 
In other words, the lemma asserts that 
$\im f\cap Z^{\perp}\subset [(f^D)^{-1}(g^D(Z_{\dt}))]^{\perp}$, or by duality that
$(f^D)^{-1}(g^D(Z_{\dt}))\subset Z+\Ker f^D$.
But the last inclusion is equivalent to (\ref{E:inclusion}).


To show (\ref{E:inclusion}), take  any element 
$y\in g^D(Z_{\dt})\cap \im f^D\subset \pi_0(\wh{S/Z(G)}^{\Gm})$. Note that
$g^D$ factors through $Z_{[\wh{\ov{d_{\dt}}}]}:Z(\wh{(G_{\dt})^{\ad}})^{\Gm}
\to\pi_0(\wh{S/Z(G_{\dt})}^{\Gm})$. Therefore it follows from 
\rco{group} that there exists $[b'_1]\in\Pi^{-1}_{\C{E}'}([d_{\dt}],\ov{\xi}_{\dt})$ such that 
$y$ is the  image of $\ka\left(\frac{[b'_1]}{[b']}\right)(\ov{s}')\in\pi_0(\wh{S/Z(G_{\dt})}^{\Gm})$.
Let $[b_1]$ be the stable conjugacy class of embeddings $S\hra H$
corresponding to $[b'_1]$. Then  $[b_1]_G=[b]_G$ and
$y=\ka\left(\frac{[b_1]}{[b]}\right)(\ov{s})$ (by \rco{embtori}).
Since $y\in\im f^D$, the image of $\ka\left(\frac{[b_1]}{[b]}\right)(\ov{s})$ in 
$\pi_0(\wh{S}^{\Gm}/Z(\wh{G})^{\Gm})$ is trivial. 
Since this image equals $\ov{\ka}_{[b_1]}/\ov{\ka}_{[b]}$ (use \rl{kainv} (a)), we conclude that
$\Pi_{\C{E}}([b_1])=\Pi_{\C{E}}([b])=([d],\ov{\xi})$. Therefore
$y=\ka\left(\frac{[b_1]}{[b]}\right)(\ov{s})$ belongs to  $f^D(Z)$, as claimed.
\end{proof}

\subsection{Preparation for the proof of the Main Theorem}
\begin*
\vskip 8truept
\end*

\begin{Lem}\label{L:characters}
Let $\pi:\wt{G}\to G$ be a quasi-isogeny such that $\wt{G}$ splits
over $E^{\nr}$, $\wt{a}_0:\wt{T}\hra\wt{G}$ the lift of $a_0:T\hra G$, and 
$\wt{\theta}$ the composition 
$\wt{T}(E)\to T(E)\overset{\theta}{\lra}\B{C}\m$.

For each stable conjugate $a$ of $a_0$, 
the representation $\pi_{a,\theta}\circ\pi$
of $\wt{G}(E)$ is isomorphic to the direct sum 
$\sum_{\wt{a}}\pi_{\wt{a},\wt{\theta}}$, taken over the set of all 
conjugacy classes of embeddings $\wt{a}:\wt{T}\hra \wt{G}$ such that
$\pi\circ\wt{a}:T\hra G$ is conjugate to $a$.
\end{Lem}
\begin{proof}
Observe first that for each quasi-isogeny $\ov{\pi}:\wt{L}\to L_a$, the representation
$\rho_{\ov{a},\ov{\theta}}\circ \ov{\pi}$ of $\wt{L}(\fq)$
is isomorphic to the Deligne--Lusztig representation 
$\rho_{\wt{\ov{a}},\wt{\ov{\theta}}}$, where $\wt{\ov{a}}:\wt{\ov{T}}\hra \wt{L}$
is the lift of $\ov{a}$, and   $\wt{\ov{\theta}}$ is the composition
$\wt{\ov{T}}(\fq)\to \ov{T}(\fq)\overset{\ov{\theta}}{\lra}\B{C}\m$.

For each $\wt{a}$ as in the lemma, denote by 
$\pi_{\wt{a}}:\wt{\wt{G}_{\wt{a}}}\to\wt{G_a}$ the restriction of $\pi$.
Then by the above observation, the representation 
$\rho_{a,\theta}\circ \pi_{\wt{a}}$ is isomorphic to 
$\rho_{\wt{a},\wt{\theta}}$. It follows that each $\pi_{\wt{a},\wt{\theta}}$ is a 
subrepresentation of $\pi_{a,\theta}\circ\pi$.

Since conjugacy classes of the 
$\wt{a}$'s  are naturally identified with the double coset $\pi(\wt{G}(E))\bs G(E)/a(T(E))$, while 
the set of irreducible factors of $\pi_{a,\theta}\circ\pi$ 
is naturally identified with  $\pi(\wt{G}(E))\bs G(E)/\wt{G_a}$, 
it remains to check that these two double cosets coincide. 

By \rco{tilde}, we have $\wt{G_a}=a(T(E))G_a$, therefore it will suffice to show 
that $G_a\subset \pi(\wt{G}(E)) a(T(\C{O}))$. But this inclusion
follows from \rco{unrtori}.
\end{proof}

\begin{Cor}\label{C:characters}
In the notation of \rl{characters}, let 
$\wt{\ka}\in\wh{\wt{T}}^{\Gm}/Z(\wh{\wt{G}})^{\Gm}$ be the image of
$\ka\in\wh{T}^{\Gm}/Z(\wh{{G}})^{\Gm}$. Then $\pi^*(\chi_{a_0,\ka,\theta})=\chi_{\wt{a}_0,\wt{\ka},\wt{\theta}}$.
\end{Cor}
\begin{proof}
Since each $\wt{a}$ as in the lemma satisfies 
$\lan\inv(\wt{a}_0,\wt{a}),\wt{\ka}\ran=\lan\inv(a_0,a),\ka\ran$,
the assertion follows from the lemma.
\end{proof}

\begin{Lem} \label{L:sc}
It will suffice to prove \rt{main} under the assumption that 
the derived group of $G$ is simply connected.
\end{Lem} 
\begin{proof}
Let $G$ be an arbitrary group satisfying the assumptions of \rt{main}.
Since $G$ splits over $E^{\nr}$,
there exists a surjective quasi-isogeny $\pi:\wt{G}\to G$ 
 such that $\wt{G}^{\der}=\wt{G}^{\ssc}$, and $\Ker\pi$ is an induced torus splitting 
over $E^{\nr}$ (use \cite[Prop. 3.1]{MS}). (Such a quasi-isogeny Kottwitz calls $z$-extension.) 
Let $\wt{a}_0:\wt{T}\hra \wt{G}$ be the 
lift of $a_0$. Then $\wt{T}$ splits over $E^{\nr}$, hence 
$\wt{G}$ satisfies all the assumptions of \rt{main}. 

Let $\wt{\ka}\in\wh{\wt{T}}^{\Gm}/Z(\wh{\wt{G}})^{\Gm}$ be the image of
$\ka$, $\wt{\theta}$ the composition $\wt{T}(E)\to T(E)\overset{\theta}{\lra}\B{C}\m$,
$\pi':\wt{G}'\to G'$ the inner twist of $\pi$, induced by $\varphi$, and  
$\wt{a}'_0:\wt{T}\hra\wt{G}'$ the lift of $a'_0:T\hra G'$.
By the assumption, generalized functions 
$\chi_{\wt{a}_0,\wt{\ka},\wt{\theta}}$ and 
$\chi_{\wt{a}'_0,\wt{\ka},\wt{\theta}}$ are 
$(\wt{a}_0,\wt{a}'_0;\wt{\ka})$-equivalent. Since $H^1(E,\Ker\pi)=0$, 
we get that $\pi(\wt{G}(E))=G(E)$ and $\pi'(\wt{G}'(E))=G'(E)$.
Therefore it follows from Corollaries \ref{C:characters} and \ref{C:func}
that generalized functions $\chi_{a_0,\ka,\theta}$ and $\chi_{a'_0,\ka,\theta}$
are $(\wt{a}_0,\wt{a}'_0;\wt{\ka})$-equivalent.
Thus by \rco{comp}, they are  $({a}_0,{a}'_0;{\ka})$-equivalent, as claimed.
\end{proof}

From now on we will assume that $G^{\der}=G^{\ssc}$. 

\begin{Lem} \label{L:red} 
Let $G$ and $a_0$ be as in \rt{main}, $\dt\in G(E)$ an element of finite order 
prime to $p$, and $b_0:T\hra G_{\dt}$ an embedding stably conjugate to $a_0$. Then 

(a) The conclusion of \rt{unip} holds for $G_{\dt}$ and $b_0$. 

(b) The set of topologically unipotent elements of $G_{\dt}(E)$ equals  
$G_{\dt}(E)_{\tu}$.
\end{Lem} 
\begin{proof}
(a) First of all, it follows from \rl{cent} (b) that $G_{\dt}$ splits over $E^{\nr}$. 
Consider the canonical isogeny $\pi:G^{\ssc}\times Z(G)^0\to G$. 
It induces the isogeny   
$\pi_{\dt}:(G^{\ssc})_{\dt}\times Z(G)^0\to G_{\dt}$, 
where $(G^{\ssc})_{\dt}:=\{g\in G^{\ssc}\,|\,\Inn \dt(g)=g\}$ is connected by 
\cite[Thm. 8.1]{St}.

Since $G$ satisfies property $(vg)$, the order of  
$Z(G^{\ssc})$ is prime to $p$. Therefore $\pi$ and hence $\pi_{\dt}$ are of order 
prime to $p$. As in the proof of \rl{sc}, we see 
(using \rl{isog} (a) and Corollaries \ref{C:characters}, \ref{C:func} and \ref{C:comp}) 
that we can replace 
$G_{\dt}$ by $(G^{\ssc})_{\dt}\times Z(G)^0$ and $b_0$ by its lift.
Thus it will suffice to show that $(G^{\ssc})_{\dt}$ and the lift 
$b_0^{\ssc}:T^{\ssc}\hra(G^{\ssc})_{\dt}$ of $b_0$ satisfy the assumptions of 
\re{ass}.

Since $G$ satisfies property $(vg)$, $G^{\ssc}$ admits a quasi-logarithm 
$\Phi:G^{\ssc}\to \C{G}^{\ssc}$ defined over $\C{O}$, and  $\C{G}^{\ssc}$ admits
a  non-degenerate over $\C{O}$ invariant pairing $\lan\cdot,\cdot\ran$. As
$\Lie (G^{\ssc})_{\dt}=(\C{G}^{\ssc})_{\dt}$, the quasi-logarithm  
$\Phi$ induces a quasi-logarithm $\Phi_{\dt}$ for $(G^{\ssc})_{\dt}$, and 
$\lan\cdot,\cdot\ran$ induces an invariant pairing $\lan\cdot,\cdot\ran_{\dt}$ on 
$\Lie(G^{\ssc})_{\dt}$. Furthermore, as $\dt$ is of finite order prime to   
$p$, we get that  $\Phi_{\dt}$ is defined over $\C{O}$ and 
$\lan\cdot,\cdot\ran_{\dt}$ is non-degenerate over $\C{O}$. 
Finally, since $b_0^{\ssc}:T^{\ssc}\hra(G^{\ssc})_{\dt}$ is stably conjugate to 
$a_0^{\ssc}:T^{\ssc}\hra G^{\ssc}$ and since $G$ satisfies property $(vg)_{a_0}$,
there exists a $b_0^{\ssc}$-strongly regular element of $\ov{T}^{\ssc}(\fq)$, 
implying the last assumption of \re{ass}.

(b) The proof is a generalization of that of \rco{compact}. Each topologically unipotent element 
$u\in G_{\dt}(E)$ stabilizes some point $x\in\C{B}(G_{\dt})\subset\C{B}(G)$. Hence 
$u$ belongs to $G_{\dt}(E)\cap G_{x,\tu}$ (by \rl{assumptions} (b)). Therefore 
$\ov{u}\in L_x(\fq)$ belongs to $(L_x)_{\ov{\dt}}\cap \C{U}(L_x)$ (by \rl{cent} (b)).
By \rl{connected} below, $\ov{u}$ belongs to  
$(L_x)^0_{\ov{\dt}}\cap \C{U}(L_x)=\C{U}((L_x)^0_{\ov{\dt}})$, hence $u$ belongs
to $(G_{\dt})_{x,\tu}\subset G_{\dt}(E)_{\tu}$, as claimed.
\end{proof}

\begin{Lem} \label{L:connected}
Let $L$ be a connected reductive group, and $s\in L$ a semisimple element. Then 
$\C{U}(L)_s=\C{U}(L)\cap L_s$ is contained in $L^0_s$.
\end{Lem}
\begin{proof}
Recall that the canonical homomorphism $\iota:L^{\ssc}\to L$ induces an $L$-equivariant 
isomorphism $\C{U}(L^{\ssc})\isom \C{U}(L)$, hence an isomorphism 
$\C{U}(L^{\ssc})_s\isom \C{U}(L)_s$. Therefore $\C{U}(L)_s=\iota(\C{U}(L^{\ssc})_s)$ 
is contained in $\iota ((L^{\ssc})_s)$. 
Since $(L^{\ssc})_s$ is connected (by \cite[Thm. 8.1]{St}), the latter group is contained in
$L^0_s$, as claimed.
\end{proof}

\begin{Not}
(a) For a compact element $\gm_0\in G^{\sr}(E)$ with a topological Jordan decomposition
$\gm_0=\dt_0 u_0$, we say that $t\in T(E)$ is 
{\em $(G,a_0,\gm_0)$-relevant}, \label{ga0gm0} if there exists an embedding 
$b_0:T\hra  G_{\dt_0}\subset G$ stably conjugate to $a_0$ such that $b_0(t)=\dt_0$. 

(b) Assume that  $t\in T(E)$ is $(G,a_0,\gm_0)$-relevant. Since 
$b_0(T)\subset G_{\dt_0}$ is elliptic, Kottwitz' theorem (see \re{TN}) implies that 
$H^1(E,T)\to H^1(E,G_{\dt_0})$ is surjective (compare the proof of \rl{elltor}).
Hence for each $\dt\in G(E)$ stably conjugate to $\dt_0$ there exists an embedding 
$b_{t,\dt}:T\hra G_{\dt}\subset G$ \label{btdt} stably conjugate to $a_0$ such that $b_{t,\dt}(t)=\dt$.
Furthermore, $b_{t,\dt}$ is unique up to a stable conjugacy, and the endoscopic triple
$\C{E}_{t}:=\C{E}_{([b_{t,\dt}],\ka)}=(H_t,[\eta_t],\ov{s}_t)$ \label{et} of $G_{\dt_0}$ is 
independent of $\dt$.
\end{Not}

\begin{Lem} \label{L:formula}

For each compact element $\gm_0\in G^{\sr}(E)$ with topological Jordan decomposition
$\gm_0=\dt_0 u_0$ and each $\ov{\xi}\in\pi_0(\wh{G_{\gm_0}}^{\Gm}/Z(\wh{G})^{\Gm})$, 
we have

\begin{equation} \label{E:formula} 
F_{a_0,\ka,\theta}(\gm_0,\ov{\xi})=\sum_{t}\theta(t)\sum_{{\dt}}I_{t,\dt},
\end{equation}
where

 (i) $t$ runs over the set of $(G,a_0,\gm_0)$-relevant elements of $T(E)$;

 (ii) ${\dt}$ runs over a set of representatives of the conjugacy classes in $G(E)$
stably conjugate to ${\dt}_0$, for which there exists a stably conjugate 
$\gm$ of $\gm_0$ with topological Jordan decomposition $\gm={\dt}u$;

(iii) $I_{t,\dt}$ equals 
\begin{equation} \label{E:for1}
\lan \inv(a_0,b_{t,\dt}),\ka\ran \lan \inv(\gm_0,\gm),\ov{\xi}\ran ^{-1}
F_{b_{t,\dt},\ka,\theta}(u,\ov{\xi})
\end{equation} 
for each $\gm=\dt u$ as in (ii).
\end{Lem}

\begin{proof} 
Recall that
\[
F_{a_0,\ka,\theta}(\gm_0,\ov{\xi})=e(G)\sum_{a}\sum_{\gm} 
\lan \inv(a_0,a),\ka\ran \lan \inv(\gm_0,\gm),\ov{\xi}\ran ^{-1}F_{a,\theta}(\gm),
\]
where $a:T\hra G$ and $\gm\in G(E)$ run over sets of representatives of the conjugacy classes within 
the stable conjugacy classes of $a_0$ and $\gm_0$, respectively. 

Using \rp{reduction}, we see that $F_{a_0,\ka,\theta}(\gm_0,\ov{\xi})$ equals the 
triple sum 
\begin{equation} \label{E:for21}
\sum_{\gm=\dt u} \sum_{a}\lan\inv(a_0,a),\ka\ran 
\lan \inv(\gm_0,\gm),\ov{\xi}\ran ^{-1}e(G_{\dt})
\sum_{b}\theta(b^{-1}({\dt}))F_{b,\theta}(u),
\end{equation}
where $\gm$ and $a$ are as above, and $b$ runs over conjugacy classes of 
embeddings $T\hra G_{\dt}$, which are conjugate to $a:T\hra G$.

Then (\ref{E:for21}) can be rewritten in the form 
\begin{equation} \label{E:for2}
\sum_{\gm=\dt u} e(G_{\dt}) \sum_{b}\lan\inv(a_0,b),\ka\ran 
\lan \inv(\gm_0,\gm),\ov{\xi}\ran ^{-1}\theta(b^{-1}({\dt}))F_{b,\theta}(u),
\end{equation}
where $b$ runs over the set of conjugacy classes of embeddings $T\hra G_{\dt}$, 
whose composition with the inclusion $G_{\dt}\hra G$ is stably conjugate to $a_0$. 

Furthermore, each $t:=b^{-1}({\dt})\in T(E)$, appearing in the sum, is $(G,a_0,\gm_0)$-relevant,
and the contribution of each such $t$ to (\ref{E:for2}) is 

\begin{equation} \label{E:for3}
\theta(t)\sum_{\gm=\dt u}\lan \inv(a_0,b_{t,\dt}),\ka\ran \lan \inv(\gm_0,\gm),\ov{\xi}\ran ^{-1}
e(G_{\dt})\sum_{b}\lan\inv(b_{t,\dt},b),\ka\ran F_{b,\theta}(u),
\end{equation}
where $b$ runs over a set of representatives of conjugacy classes of embeddings 
$T\hra G_{\dt}$ stably conjugate to $b_{t,\dt}$. 
But (\ref{E:for3}) coincides with the sum $\theta(t)\sum_{{\dt}}I_{t,\dt}$ 
as in the lemma.
\end{proof}

\begin{Emp}
We fix  $\gm_0\in G^{\sr}(E)$ and 
$\ov{\xi}\in\pi_0(\wh{G_{\gm_0}}^{\Gm}/Z(\wh{G})^{\Gm})$ such that 
$F_{a_0,\ka,\theta}(\gm_0,\ov{\xi})\neq 0$, and we are 
going to show that $(\gm_0,\ov{\xi})$ satisfies the conditions $(i),(ii)$ of 
\rt{main'}. 
\end{Emp}

\subsection{Proof of \rt{main'} $(i)$} \label{SS:i}
\begin*
\vskip 8truept
\end*

\begin{Emp}
Since the support of each $t_{a,\theta}$ is contained in $\wt{G_a}$,
the assumption on $\gm_0$ implies that there exists $z\in Z(G)(E)$ 
such that $z\gm_0$ is compact.
But $F_{a_0,\ka,\theta}(z\gm_0,\ov{\xi})=\theta(z)
F_{a_0,\ka,\theta}(\gm_0,\ov{\xi})$
for each  $z\in Z(G)(E)$, therefore the assertions of \rt{main'} for $\gm_0$ are 
equivalent to those for $z\gm_0$. 
Hence we can and will assume that $\gm_0$ is compact. In particular,
\rl{formula} holds for $\gm_0$.
 
Every stably conjugate $\gm$ of $\gm_0$ 
is compact as well, and we denote by $\gm_0=\dt_0 u_0$ and  $\gm=\dt u$ their 
topological Jordan decompositions. We also let $\ov{\xi}_{\dt_0}$ be the image of $\ov{\xi}$  in 
$\pi_0(\wh{G_{\gm_0}}^{\Gm}/Z(\wh{G_{\dt_0}})^{\Gm})$. 
\end{Emp}

\begin{Not} \label{N:eq}
We will say that stable conjugates $\dt_1,\dt_2\in G(E)$ of $\dt_0$ are 
{\em $(\gm_0,\ov{\xi})$-equivalent}, \label{gmoxieq} if there exist stable conjugate $\gm_1$ and $\gm_2$
of $\gm_0$ with topological Jordan decompositions $\gm_i=\dt_i u_i$, and 
$G_{\dt_2}$ is an $(\C{E}_t,[a_{\gm_0}],\ov{\xi}_{\dt_0})$-admissible inner form of 
$G_{\dt_1}$. In this case we will write 
$\dt_1\sim_{(\gm_0,\ov{\xi})}\dt_2$. \label{dt1simdt2}
\end{Not}

\begin{Emp}
Fix $t$ which has a non-zero contribution to (\ref{E:formula}). Since the set of 
conjugacy classes of $\dt$ in \rl{formula} (ii) decomposes as a union of 
$(\gm_0,\ov{\xi})$-equivalent classes, we can replace  
$\gm_0$ by a stably conjugate element, so that  
$\sum_{\dt\sim_{(\gm_0,\ov{\xi})}\dt_0}I_{t,\dt}\neq 0$.
Further replacing $\gm_0$, we can moreover assume that 
$I_{t,\dt_0}\neq 0$, thus $F_{b_{t,\dt_0},\ka,\theta}(u_0,\ov{\xi}_{\dt_0})\neq 0$. 

We also fix embeddings of maximal tori $c:T\hra H$ and $c':T\hra H_t$ 
such that $\Pi_{\C{E}}([c])=([a_0],\ov{\ka})$ and 
$\Pi_{\C{E}_t}([c'])=([b_{t,\dt_0}],\ov{\ka}')$ (see \re{technical} (b)). 
This enables us to apply the results of \rss{technical}.
\end{Emp}

\begin{Emp}
Since $G^{\der}=G^{\ssc}$, we get that $u_0\in G_{\dt_0}(E)_{\tu}$ and
the conclusion of \rt{unip} holds for $G_{\dt_0}$ and $b_{t,\dt_0}$ (see \rl{red}). 
Therefore as in \rl{equiv}, there exists an
embedding $b':G_{\gm_0}=(G_{\dt_0})_{u_0}\hra H_t$ such that
$\Pi_{\C{E}_{t}}([b'])=([a_{\gm_0}],\ov{\xi}_{\dt_0})$, where 
$a_{\gm_0}:G_{\gm_0}\hra G_{\dt_0}$ is the natural inclusion.

Let $[b]$ be the stable conjugacy class of embeddings $G_{\gm_0}\hra H$ corresponding 
to
$[b']$ (see \rl{embtori} (b)), and put 
$\ov{\xi}_{[b']}:=\ov{\ka}_{[b]}\in\pi_0(\wh{G_{\gm_0}}^{\Gm}/Z(\wh{G})^{\Gm})$. 

To prove the assertion $(i)$ of \rt{main'}, it will suffice to show the
existence of $[b']\in\Pi_{\C{E}_{t}}^{-1}([a_{\gm_0}],\ov{\xi}_{\dt_0})$ such that 
$\ov{\xi}_{[b']}=\ov{\xi}$. 
\end{Emp}

\begin{Emp}
For each $[b']\in\Pi_{\C{E}_{t}}^{-1}([a_{\gm_0}],\ov{\xi}_{\dt_0})$,  
we denote by $z_{[b']}$ the image of the quotient 
$\ov{\xi}_{[b']}/\ov{\xi}\in\pi_0(\wh{G_{\gm_0}}^{\Gm}/Z(\wh{G})^{\Gm})$ in 
$\pi_0(\wh{G_{\gm_0}}^{\Gm}/Z(\wh{G})^{\Gm}Z(\C{E}_{t},[a_{u_0}],\ov{\xi}_{\dt_0}))$, where  
$Z(\C{E}_{t},[a_{u_0}],\ov{\xi}_{\dt_0})\subset Z(\wh{(G_{\dt_0})^{\ad}})^{\Gm}$ is mapped into 
$\wh{G_{\gm_0}}^{\Gm}$ via the homomorphism $Z(\wh{(G_{\dt_0})^{\ad}})\to 
\wh{G_{\gm_0}/Z(G_{\dt_0})}\to \wh{G_{\gm_0}}$.

We claim that $z_{[b']}$ does not depend on $[b']$. Indeed, for each 
$[b'_1],[b'_2]\in\Pi_{\C{E}_{t}}^{-1}([a_{\gm_0}],\ov{\xi}_{\dt_0})$,   the quotient 
$\ov{\xi}_{[b'_1]}/\ov{\xi}_{[b'_2]}$ is the image of 
$\ka\left(\frac{[b'_1]}{[b'_2]}\right)(\ov{s}_t)\in\pi_0(\wh{G_{\gm_0}/Z(G_{\dt_0})}^{\Gm})$
(by \rco{embtori} and \rl{kainv} (a)).
Thus $\ov{\xi}_{[b'_1]}/\ov{\xi}_{[b'_2]}$ belongs to the image of 
$Z(\C{E}_{t},[a_{u_0}],\ov{\xi}_{\dt_0})$. It follows that
$z:=z_{[b']}$ is independent of $[b']$, as claimed.
\end{Emp}

\begin{Emp}
Our next goal is to show that $z=1$. By definition, $z$ belongs to the image  
\begin{equation} \label{Eq:group}
\im\left[\pi_0(Z(\wh{G_{\dt_0}})^{\Gm}/Z(\wh{G})^{\Gm})\to
\pi_0(\wh{G_{\gm_0}}^{\Gm}/Z(\wh{G})^{\Gm}Z(\C{E}_{t},[a_{u_0}],\ov{\xi}_{\dt_0}))\right].
\end{equation}
Denote by $V\subset\Ker[H^1(E,G_{\dt_0})\to H^1(E,G)]$ the intersection of 
the image of $\Ker [H^1(E,G_{\gm_0})\to H^1(E,G)]$ and the preimage of 
$Z(\C{E}_{t},[a_{u_0}],\ov{\xi}_{\dt_0})^\perp\subset H^1(E,(G_{\dt_0})^{\ad})$. 
Then for a stable conjugate $\dt\in G(E)$ of $\dt_0$, we have 
$\dt\sim_{(\gm_0,\ov{\xi})}\dt_0$ if and only if the invariant
$\inv(\dt_0,\dt)\in \Ker[H^1(E,G_{\dt_0})\to H^1(E,G)]$ lies in $V$.

By Kottwitz' theorem (\re{TN}), the dual group $V^D$ of $V$ is naturally identified 
with the group (\ref{Eq:group}). 
In particular, $z$ belongs to $V^D$. Therefore 
for each $\dt\sim_{(\gm_0,\ov{\xi})}\dt_0$ one can form a pairing
$\lan \inv(\dt_0,\dt),z\ran$. Explicitly, 
\begin{equation} \label{E:invz}
\lan \inv(\dt_0,\dt),z\ran= \lan \inv(\gm_0,\gm),\wt{z}\ran
\end{equation}
for every stably conjugate $\gm$ of $\gm_0$   with topological Jordan decomposition
$\gm=\dt u$, and every representative $\wt{z}\in\pi_0(\wh{G_{\gm_0}}^{\Gm}/Z(\wh{G})^{\Gm})$ of $z$.
\end{Emp}

\begin{Cl} \label{C:eq}
For each $\dt\sim_{(\gm_0,\ov{\xi})}\dt_0$, we have 
$I_{t,\dt}=\lan \inv(\dt_0,\dt),z\ran I_{t,\dt_0}$.
\end{Cl}

\begin{proof}
Let $\gm\in G(E)$ be a stable conjugate of $\gm_0$ with topological Jordan decomposition
$\gm=\dt u$. Then  $u\in G_{\dt}(E)$ is a stable conjugate of $u_0\in G_{\dt_0}(E)$. 
Since by assumption $F_{b_{t,\dt_0},\ka,\theta}(u_0,\ov{\xi}_{\dt_0})\neq 0$ and 
$u_0\in G_{\dt_0}(E)$ is topologically unipotent, we conclude from \rl{red} (as in \rl{equiv}) that  
\[
F_{b_{t,\dt},\ka,\theta}(u,\ov{\xi}_{\dt_0})=
\lan \frac{u_0,u;\ov{\xi}_{\dt_0}}{b_{t,\dt_0} ,b_{t,\dt};\ka}\ran_{\C{E}_t}
F_{b_{t,\dt_0},\ka,\theta}(u_0,\ov{\xi}_{\dt_0}).
\]
Hence the quotient $I_{t,\dt}/I_{t,\dt_0}$ equals
\[
\lan \frac{u_0,u;\ov{\xi}_{\dt_0}}{b_{t,\dt_0} ,b_{t,\dt};\ka}\ran_{\C{E}_t}
\lan \inv(b_{t,\dt_0},b_{t,\dt}),\ka\ran\lan \inv(\gm_0,\gm),\ov{\xi}\ran^{-1}.
\]
Thus our claim is equivalent to the equality
\begin{equation} \label{E:invariant}
\lan \frac{u_0,u;\ov{\xi}_{\dt_0}}{b_{t,\dt_0} ,b_{t,\dt};\ka}\ran_{\C{E}_t}=
\lan \inv(\dt_0,\dt),z\ran
\lan \inv(\gm_0,\gm),\ov{\xi}\ran\lan \inv(b_{t,\dt_0},b_{t,\dt}),\ka\ran^{-1}.
\end{equation}
The left hand side of (\ref{E:invariant}) 
equals $\lan \frac{a_{u_0},a_u;[b']}{b_{t,\dt_0} ,b_{t,\dt};[c']}\ran_{\C{E}_t}$ for 
each $[b']\in\Pi_{\C{E}_{t}}^{-1}([a_{\gm_0}],\ov{\xi}_{\dt_0})$.
By  \rl{inv2} (b) and \rl{inv} (d), it therefore equals   
\[
\lan \frac{a_{\gm_0},a_{\gm};[b]}{b_{t,\dt_0},b_{t,\dt};[c]}\ran_{\C{E}}=
\lan \inv(\gm_0,\gm),\ov{\ka}_{[b]}\ran\lan 
\inv(b_{t,\dt_0},b_{t,\dt}),\ov{\ka}_{[c]}\ran^{-1}.
\] 
But $\ov{\ka}_{[c]}=\ka$, 
$\ov{\ka}_{[b]}=(\ov{\xi}_{[b']}/\ov{\xi})\ov{\xi}$ and
$\ov{\xi}_{[b']}/\ov{\xi}$ is a representative of $z$. Therefore equality 
(\ref{E:invariant}) follows from (\ref{E:invz}).
\end{proof}

\begin{Emp}
By \rcl{eq}, the sum 
$\sum_{\dt\sim_{(\gm_0,\ov{\xi})}\dt_0}I_{t,\dt_0}$ equals
$I_{t,\dt_0}(\sum_{v\in V}\lan v,z\ran)$. It follows that 
$\sum_{v\in V}\lan v,z\ran\neq 0$, hence $z=1$. 

Choose now an arbitrary  $[b']\in\Pi_{\C{E}_{t}}^{-1}([a_{\gm_0}],\ov{\xi}_{\dt_0})$.
Since $z=1$, the quotient 
$\ov{\xi}_{[b']}/\ov{\xi}\in\pi_0(\wh{G_{\gm_0}}^{\Gm}/Z(\wh{G})^{\Gm})$
lies in the image of $Z(\C{E}_{t},[a_{u_0}],\ov{\xi}_{\dt_0})$. Hence by \rco{group},
there exists $[b'_1]\in\Pi_{\C{E}_{t}}^{-1}([a_{\gm_0}],\ov{\xi}_{\dt_0})$ 
such that $\ov{\xi}_{[b']}/\ov{\xi}$ equals the image of 
$\ka\left(\frac{[b']}{[b'_1]}\right)(\ov{s}_t)$. Since by \rco{embtori} and \rl{kainv} (a),
the image of $\ka\left(\frac{[b']}{[b'_1]}\right)(\ov{s}_t)$ in $\pi_0(\wh{G_{\gm_0}}^{\Gm}/Z(\wh{G})^{\Gm})$
equals  $\ov{\xi}_{[b']}/\ov{\xi}_{[b'_1]}$, 
 we get that $\ov{\xi}_{[b'_1]}=\ov{\xi}$, completing the proof of $(i)$.
\end{Emp}

\subsection{Proof of \rt{main'} $(ii)$}
\begin*
\vskip 8truept
\end*

\begin{Emp}
Let $\gm'_0\in G'(E)$ be a stable conjugate of $\gm_0$. Since  $\gm_0$ is compact, so is
$\gm'_0$, and we denote by $\gm'_0=\dt'_0 u'_0$ its topological Jordan decomposition.
By \rl{formula}, we can write $F_{a'_0,\ka,\theta}(\gm'_0,\ov{\xi})$ in the form
\[
F_{a'_0,\ka,\theta}(\gm'_0,\ov{\xi})=\sum_{t'}\theta(t')\sum_{\dt'}I_{t',\dt'},
\]
where $t'$, $\dt'$ and $I_{t',\dt'}$ have the same meaning as in \rl{formula}.

First we claim that  an element $t\in T(E)$ is $(G,a_0,\gm_0)$-relevant 
if and only if it is  $(G',a'_0,\gm'_0)$-relevant. Indeed, assume that
 $t$ is $(G,a_0,\gm_0)$-relevant, and let $b_0:T\hra G_{\dt_0}\subset G$ be the 
corresponding embedding. Since $T/Z(G)$ is anisotropic, Kottwitz' theorem implies that
the map $H^1(E,T/Z(G))\to H^1(E,G_{\dt_0}/Z(G))$ is surjective (compare the proof of \rl{elltor}).  
Hence there exists a stable conjugate $b'_0:T\hra G'_{\dt'_0}\subset G'$ of $b_0$ such that 
$b'_0(t)=\dt'_0$, thus  $t$ is $(G',a'_0,\gm'_0)$-relevant.

Therefore it will suffice to show that 
for each  $(G,a_0,\gm_0)$-relevant $t$, we have
\begin{equation} \label{E:final} 
\sum_{\dt'}I_{t,\dt'}=\lan\frac{\gm_0,\gm'_0;\ov{\xi}}{a_0,a'_0,\ka}\ran_{\C{E}} 
\sum_{\dt}I_{t,\dt}.
\end{equation}
\end{Emp}

\begin{Emp}
Fix  $(G,a_0,\gm_0)$-relevant $t$. Generalizing \rn{eq}, we will say that stable 
conjugates $\dt\in G(E)$ and  $\dt'\in G'(E)$ of $\dt_0$ are 
{\em $(\gm_0,\ov{\xi})$-equivalent} \label{gmoxieq2} (and will write  
$\dt\sim_{(\gm_0,\ov{\xi})}\dt'$), \label{dt'simdt} if there exist 
stable conjugates $\gm\in G(E)$ and $\gm'\in G'(E)$ of $\gm_0$ with 
topological Jordan decompositions $\gm=\dt u$ and $\gm'=\dt' u'$ such that
$G'_{\dt'}$ is an $(\C{E}_t,[a_{u_0}],\ov{\xi}_{\dt_0})$-admissible inner form of 
$G_{\dt}$. 

Assume that $\varphi:G\to G'$ is $(\C{E},[a_{\gm_0}],\ov{\xi})$-admissible.
We claim that for every stable conjugate 
$\gm\in G(E)$ of $\gm_0$ with topological Jordan decomposition $\gm=\dt u$,
there exists a stable conjugate $\gm'\in G'(E)$ of $\gm'_0$
with topological Jordan decomposition $\gm'=\dt' u'$ such that 
$\dt\sim_{(\gm_0,\ov{\xi})}\dt'$. 

Indeed, we have shown in \rss{i} that there exists 
$[b']\in \Pi_{\C{E}_t}^{-1}([a_{\gm_0}],\ov{\xi}_{\dt_0})$ such that 
the corresponding stable conjugacy class $[b]$ of embeddings $G_{\gm_0}\hra H$
satisfies  $\Pi_{\C{E}}([b])=([a_{\gm_0}],\ov{\xi})$. Since the inclusion 
$a_{\gm}:G_{\gm}\hra G$
has a stable conjugate $a_{\gm'_0}:G_{\gm}\cong G'_{\gm'_0}\hra G'$, 
the assertion follows from \rl{compatib}.

Therefore equality (\ref{E:final}) follows from the following generalization
of \rcl{eq}. 
\end{Emp}

\begin{Cl} \label{C:eq1}
For each $\dt'\sim_{(\gm_0,\ov{\xi})}\dt$, we have 
$I_{t,\dt'}=\lan\frac{\gm_0,\gm'_0;\ov{\xi}}{a_0,a'_0;\ka}\ran_{\C{E}} I_{t,\dt}$.
\end{Cl}
\begin{proof}
The proof is very similar to that of \rcl{eq}. 
By \rt{unip} for the inner twisting $G_{\dt}\to G'_{\dt'}$ (use \rl{red}), we have
\[
F_{b_{t,\dt'},\ka,\theta}(u',\ov{\xi}_{\dt_0})=
\lan \frac{u,u';\ov{\xi}_{\dt_0}}{b_{t,\dt} ,b_{t,\dt'};\ka}\ran_{\C{E}_t}
F_{b_{t,\dt},\ka,\theta}(u,\ov{\xi}_{\dt_0}).
\]
Hence the quotient $I_{t,\dt'}/I_{t,\dt}$ equals
\[
\lan \frac{u,u';\ov{\xi}_{\dt_0}}{b_{t,\dt} ,b_{t,\dt'};\ka}\ran_{\C{E}_t}
\lan \inv(a'_0,b_{t,\dt'}),\ka\ran\lan \inv(a_0, b_{t,\dt}),\ka\ran^{-1}
\lan \inv(\gm'_0,\gm'),\ov{\xi}\ran^{-1}\lan \inv(\gm_0,\gm),\ov{\xi}\ran.
\]
Thus we have to check that 
$\lan \frac{u,u';\ov{\xi}_{\dt_0}}{b_{t,\dt} ,b_{t,\dt'};\ka}\ran_{\C{E}_t}=
\lan \frac{a_{\gm},a_{\gm'};[b']}{b_{t,\dt} ,b_{t,\dt'};[c']}\ran_{\C{E}_t}$ equals
\[
\lan\frac{\gm_0,\gm'_0;\ov{\xi}}{a_0,a'_0;\ka}\ran_{\C{E}}
\lan \inv(a'_0,b_{t,\dt'}),\ka\ran^{-1}\lan \inv(a_0,b_{t,\dt}),\ka\ran
\lan \inv(\gm'_0,\gm'),\ov{\xi}\ran\lan \inv(\gm_0,\gm),\ov{\xi}\ran ^{-1}.
\]
Since the latter expression equals 
$\lan\frac{a_{\gm},a_{\gm'};[b]}{b_{t,\dt} ,b_{t,\dt'};[c]}\ran_{\C{E}}$
(use \rl{inv} (a), (c)), the assertion follows from \rl{inv2} (b).
\end{proof}

This completes the proof of Theorems \ref{T:main'} and  \ref{T:main}.

%

\appendix
\section{Springer Hypothesis} \label{S:Spr}
The goal of this appendix is to prove the following result, 
conjectured by Springer and playing a crucial role in \rss{unip}.
In the case of large characteristic this result was first proved in  
\cite{Ka1}.

\begin{Thm'} \label{T:Spr}
Let  $L$ be a reductive group over a finite field $\fq$, $\Phi:L\to \C{L}$
a quasi-logarithm (see \rd{qlog}),  $\lan\cdot,\cdot\ran$ a non-degenerate invariant 
pairing on $\C{L}$, $T\subset L$ a maximal torus, $\theta$ a character of $T(\fq)$, 
$\psi$ a character of $\fq$, and $t$ an element of $\C{T}(\fq)\cap \C{L}^{\sr}(\fq)$.

Denote by $\dt_t$ the characteristic function of the $\Ad L(\fq)$-orbit of $t$, 
by $\C{F}(\dt_t)$ its Fourier transform, and by 
$(-1)^{\rk_{\fq}(L)-\rk_{\fq}(T)}\rho_{T,\theta}$
the Deligne--Lusztig \cite{DL} virtual representation of $L(\fq)$ corresponding to $T$ and $\theta$.

For every unipotent $u\in L(\fq)$, we have
$$
\Tr\rho_{T,\theta}(u)=q^{-\frac{1}{2}\dim(L/T)}\C{F}(\dt_t)(\Phi(u)).
$$
\end{Thm'}

\begin{Not'}
For each Weil sheaf $\C{A}$ over a variety $X$ over a finite field $\fq$, we 
denote by $Func(\C{A})$ the corresponding function on $X(\fq)$.
\end{Not'}

\begin{proof}
The theorem is a consequence of the results of Lusztig (\cite{Lu}) and 
Springer (\cite{Sp}), and seems to be well-known to experts. Set 
$\C{U}:=\C{U}(L)\subset  L$ and $\C{N}:=\C{N}(\C{L})\subset  \C{L}$. 
To carry out the proof, we will construct a Weil sheaf $\C{A}$ on $\C{L}\times \C{T}$ 
such that the restrictions $\C{A}_t$ to $\C{L}\times\{t\}\cong\C{L}$ are perverse
for all $t\in\C{T}(\fq)$ and satisfy the following properties:

$(i)$ The restriction $\C{A}_t|_{\C{N}}$ is independent of $t$.

$(ii)$ If $t\in\C{L}^{\sr}(\fq)$, then $Func(\C{A}_t)=q^{-\frac{1}{2}\dim(L/T)}\C{F}(\dt_t)$.

$(iii)$ $Func(\Phi^*(\C{A}_0))|_{\C{U}}=\Tr\rho_{T,\theta}|_{\C{U}}$.

The existence of such an $\C{A}$ implies the Theorem. Indeed, fix $u\in \C{U}(\fq)$.
By $(iii)$, $\Tr\rho_{T,\theta}(u)$ equals  $Func(\C{A}_0)(\Phi(u))$. 
Since $\Phi(\C{U})\subset\C{N}$ (see \cite[9.1, 9.2]{BR}), we have $\Phi(u)\in \C{N}$.
Therefore the assertion follows from $(i)$ and $(ii)$.

\begin{Emp'} \label{E:spr}
{\bf Construction of $\C{A}$.}
Let $T'$ be the abstract Cartan subgroup of $L$, and $W\subset \Aut(T')$ the 
Weyl group of $L$ (see \cite[1.1]{DL}).
Denote by $\wt{\C{L}}$ the Springer resolution of $\C{L}$ classifying 
pairs $(\C{B},x)$, where $\C{B}\subset\C{L}$ is a Borel subalgebra
and $x$ is an element of $\C{B}$. Consider the diagram 
\[
\C{L}\times\C{T}'\overset{\pi\times\Id }{\lla}\wt{\C{L}}\times\C{T}'
\overset{\al\times\Id}{\lra}\C{T}'\times\C{T}'\overset{\lan\cdot,\cdot\ran'}{\lra}\B{A}^1,
\]
where $\pi$ and $\al$ send $(\C{B},x)$ to $x$ and 
$\ov{x}\in\C{B}/[\C{B},\C{B}]=\C{T}'$, respectively, and 
$\lan\cdot,\cdot\ran'$ is the form on $\C{T}'$ induced by $\lan\cdot,\cdot\ran$.
Put 
\[
\C{A}':=(\pi\times\Id)_!(\al\times\Id)^*\lan\cdot,\cdot\ran'^*(\C{L}_{\psi})[\dim L],
\] 
where $\C{L}_{\psi}$ is the Artin-Schreier local system on $\B{A}^1$ corresponding to $\psi$.

To construct $\C{A}$, we will show first that for every $w\in W$, there exists a canonical 
isomorphism $(\Id\times w)^*\C{A}'\isom\C{A}'$.
Denote by upper index $(\cdot)^{\sr}$ the restriction to (the preimage of) 
$\C{L}^{\sr}$. First we will show that $(\Id\times w)^*\C{A}^{\sr}$ 
is canonically isomorphic to $\C{A}^{\sr}$.
As $\pi^{\sr}:\wt{\C{L}}^{\sr}\to\C{L}^{\sr}$ is an unramified Galois covering with 
the Galois group $W$, the functor $\pi^{\sr}_!$ is isomorphic to 
$\pi^{\sr}_!\circ w^*$.
Since $\al^{\sr}$ is $W$-equivariant, and $\lan\cdot,\cdot\ran'=\lan\cdot,\cdot\ran'\circ (w\times w)$, 
we have a canonical isomorphism of functors
\[
(\Id\times w)^*(\pi^{\sr}\times\Id)_!(\al^{\sr}\times\Id)^* \lan\cdot,\cdot\ran'^*\cong 
(\pi^{\sr}\times\Id)_!(\Id\times w)^*(\al^{\sr}\times\Id)^* \lan\cdot,\cdot\ran'^*\cong
\]\[
(\pi^{\sr}\times\Id)_!(w\times w)^* (\al^{\sr}\times\Id)^*\lan\cdot,\cdot\ran'^*\cong 
(\pi^{\sr}\times\Id)_!(\al^{\sr}\times\Id)^*\lan\cdot,\cdot\ran'^*,
\]
implying the isomorphism  $(\Id\times w)^*\C{A}'^{\sr}\isom\C{A}'^{\sr}$.
Since $\al$ and $\lan\cdot,\cdot\ran'$ are smooth morphisms, while $\pi$ is small, 
we see that $\C{A}'[\dim T]$ is a semisimple perverse sheaf, which is
the intermediate extension of $\C{A}'^{\sr}[\dim T]$.
Thus the constructed above isomorphism  
$(\Id\times w)^*\C{A}'^{\sr}\isom\C{A}'^{\sr}$ uniquely extends to an
isomorphism $(\Id\times w)^*\C{A}'\isom\C{A}'$.

Denote by $\Fr:T'\to T'$ the geometric Frobenius morphism corresponding to the 
$\fq$-structure of $T'$, and
choose an isomorphism between $T$ and $T'$ over $\ov{\fq}$. Then there exists
$w\in W$ such that the $\fq$-structure of $T$ corresponds to the morphism
$\Fr_w:=w\circ\Fr:T'\to T'$. Denote by $\C{A}$ the Weil sheaf
$\C{A}$ on $\C{L}\times\C{T}$, which is isomorphic to $\C{A}'$ over 
$\ov{\fq}$, and the Weil structure corresponds to the composition
$\Fr^*(\Id\times w)^*\C{A}'\isom\Fr^*\C{A}'\isom\C{A}'$, where the first isomorphism
was constructed above, and the second one comes from the Weil structure of $\C{A}'$.
\end{Emp'}

It remains to show that $\C{A}$ satisfies  properties $(i)$--$(iii)$.

\begin{Emp'}
{\bf Proof of properties $(i)$--$(iii)$.}

$(i)$ Put $\wt{\C{N}}:=\pi^{-1}(\C{N})\subset \wt{\C{L}}$. Then  $\al(\wt{\C{N}})=0$, hence  
$(\al\times\Id)^* \lan\cdot,\cdot\ran'^*(\C{L}_{\psi})|_{\wt{\C{N}}\times\C{T}'}\cong\qlbar$. 
This implies the assertion for $\C{A}'$. To show the assertion for $\C{A}$, notice that
$\wt{\C{N}}$ is smooth of dimension $\dim(L/T)$, 
and the projection $\wt{\C{N}}\to \C{N}$ is semi-small. It follows that
$\C{A}'|_{\C{N}\times \C{T}}$  is a semisimple perverse sheaf.
Therefore for each $t\in\C{T}(\fq)$ the restriction map 
$\Hom(\Fr_w^*\C{A}'|_{\C{N}\times \C{T}},\C{A}'|_{\C{N}\times \C{T}})\to
\Hom(\Fr_w^*\C{A}'|_{\C{N}\times\{t\}},\C{A}'|_{\C{N}\times\{t\}})$ 
is an isomorphism. Thus the Weil structure 
of $\C{A}_t|_{\C{N}}$ is independent of $t$ as well. 

$(ii)$ Denote by $\IC_t$ the constant perverse sheaf 
$\qlbar(\frac{\dim (L/T)}{2})[\dim(L/T)]$ 
on the orbit $\Ad L(t)\subset\C{L}$, and let $\C{F}(\IC_t)$ be the 
Fourier--Deligne transform of $\IC_t$.  As  $\dim (L/T)$ is even,
we get that $Func(\C{F}(\IC_t))=q^{-\frac{\dim (L/T)}{2}}\C{F}(\dt_t)$. Therefore
the assertion follows from well-known equality $\C{F}(\IC_t)=\C{A}_t$
(see for example \cite{Sp}).

$(iii)$ By \cite[Thm. 4.2]{DL}, $\Tr\rho_{T,\theta}|_{\C{U}(\fq)}$ does not depend on 
$\theta$, hence we can assume that $\theta=1$. It was proved by Lusztig 
(see \cite[Thm. 1.14 (a) and Prop. 8.15]{Lu} and compare \cite{BP})
that there exists a perverse sheaf $\C{K}_{T}$ on $L$ such that 
$\Tr\rho_{T,1}=Func(\C{K}_{T})$. More precisely, in \cite[Thm. 1.14 (a)]{Lu}
Lusztig showed the corresponding result for general character sheaves if $q$ is
sufficiently large, while by \cite[Prop. 8.15]{Lu} the restriction on $q$ 
is unnecessary in our situation.

Thus it will suffice  to check that $\C{K}_{T}|_{\C{U}}=\Phi^*(\C{A}_0)|_{\C{U}}$.
The description of $\C{K}_T$ is very similar to that of $\C{A}$. 
Let $\pi_L:\wt{L}\to L$ be the Springer resolution, put 
$\wt{\C{U}}:=\pi_L^{-1}(\C{U})$, and we denote by 
$\pi^{\sr}_L:\wt{L}^{\sr}\to L^{\sr}$ the restriction of $\pi$ to the preimage of 
$L^{\sr}$. Then the semisimple perverse sheaf $\C{K}_T$ is equal to 
$(\pi_L)_!(\qlbar)[\dim L]$ over $\fqbar$, while the Weil structure of $\C{K}_{T}$  
is induced (as in  \re{spr}) by isomorphism of functors  
$\Fr_w^*(\pi^{\sr}_L)_!\isom (\pi^{\sr}_L)_!$.

By \cite[Thm. 6.2 and 9.1]{BR}, there exists a $L^{\ad}$-invariant open affine
neighborhood $V\supset\C{U}$ in $L$ such that $\Phi|_V:V\to\C{L}$ is \'etale.
We claim that  $\C{K}_{T}|_V$ is isomorphic to $\Phi^*(\C{A}_0)|_V=(\Phi|_V)^*(\C{A}_0)$.
As both $\C{K}_{T}|_V$ and $(\Phi|_V)^*(\C{A}_0)$ are semisimple perverse sheaves, 
which are immediate extensions of their restrictions to $V\cap \Phi^{-1}(\C{L}^{\sr})$, 
it will suffice to show that
$\Phi^*(\C{A}_0)|_{\Phi^{-1}(\C{L}^{\sr})}\cong\C{K}_{T}|_{\Phi^{-1}(\C{L}^{\sr})}$.

Note that $\Phi^{-1}(\C{L}^{\sr})\subset L^{\sr}$ and that
$\Phi$ gives rise to the commutative diagram 
\[
\CD
\wt{L}@>\wt{\Phi}>> \wt{\C{L}}\\
@V\pi_LVV         @V\pi VV\\
L@>\Phi>> \C{L}
\endCD
\]
(use \rl{qlogf} (a)), whose the restriction to $\C{L}^{\sr}$ is Cartesian and $W$-equivariant.
Therefore the required isomorphism $\Phi^*(\C{A}_0)|_{\Phi^{-1}(\C{L}^{\sr})}\isom 
\C{K}_{T}|_{\Phi^{-1}(\C{L}^{\sr})}$ follows from the proper 
base change theorem.
\end{Emp'}
\end{proof}

\section{$(a,a';[b])$-equivalence and Fourier transform} \label{S:Wa}

\subsection{Formulation of the result}

\begin*
\vskip 8truept
\end*

\begin{Emp}
Let $G$ be a reductive group over a local field $E$ of characteristic zero, 
$\C{E}=(H,[\eta],\ov{s})$ an endoscopic triple for $G$, and 
$\varphi:G\to G'$ an inner twisting. 
Fix a triple $(a,a';[b])$, where 
$a:T\hra G$ and $a':T\hra G'$ are stably conjugate embeddings of 
maximal tori, and $[b]$ is a stable conjugacy class of embeddings of maximal
tori $T\hra H$, compatible with $a$ and $a'$.

Fix a non-trivial character $\psi: E\to \B{C}\m$, a non-degenerate 
$G$-invariant pairing $\lan\cdot,\cdot\ran$ on $\C{G}$, and a non-zero 
translation invariant top degree differential form $\om_G$ on $G$.
Denote by  $\lan\cdot,\cdot\ran'$ the $G'$-invariant pairing on $\C{G}'$, 
induced $\lan\cdot,\cdot\ran$, and let $dx=|\om_{\C{G}}|$ and 
$dx'=|\om_{\C{G}'}|$ be the 
invariant measures on $\C{G}(E)$ and $\C{G}'(E)$ induced by $\om_G$.
These data define Fourier transforms $F\mapsto \C{F}(F)$ on $\C{G}(E)$ and 
 $F'\mapsto \C{F}(F')$ on $\C{G}'(E)$ (see \re{four}).
\end{Emp}

The following result generalizes both the theorem of 
Waldspurger \cite{Wa2} (who treated the case $\phi'=0$)   
and of Kazhdan--Polishchuk \cite[Thm. 2.7.1]{KP} (where the stable case is treated).

\begin{Thm} \label{T:Wa}
Generalized functions $F\in\C{D}(\C{G}(E))$ and $F'\in\C{D}(\C{G}'(E))$ 
are  $(a,a';[b])$-equivalent if and only if $e'(G)\C{F}(F)$ and 
$e'(G')\C{F}(F')$ are  $(a,a';[b]))$-equivalent.
\end{Thm}

\begin{Rem} \label{R:ch}
When this work was already written, we have 
learned that our \rt{Wa} seems to follow from the recent work of Chaudouard \cite{Ch}.
\end{Rem}

By duality, \rt{Wa} follows from the following result.

\begin{Thm} \label{T:Wa1}
Measures $\phi\in\C{S}(\C{G}(E))$ and $\phi'\in\C{S}(\C{G}'(E))$ are  
$(a,a';[b])$-indistin-guishable
if and only if $e'(G)\C{F}(\phi)$ and $e'(G')\C{F}(\phi')$ are 
$(a,a';[b])$-indistinguishable.
\end{Thm}

For the proof, we will combine arguments from \cite{Wa1,Wa2} 
with those from \cite{KP}. 

\begin{Lem} \label{L:reduction}
(a) The validity of \rt{Wa1} is independent of the choice of 
$\om_G$, $\psi$ and $\lan .,.\ran$. 

(b) It will suffice to show \rt{Wa1} under the assumption
that $G=G^{\ssc}$.

(c) It will suffice to show \rt{Wa1} under the assumption
that $\C{E}$ is elliptic. (Note that this is the only case used in this paper). 
\end{Lem}
\begin{proof}
The proof follows by essentially the same arguments as \cite[II]{Wa2}.

(a) Another choice of the data results in replacing the Fourier transform $\C{F}$ by 
$B^*\circ\C{F}$ for a certain linear automorphism $B$ of $\C{G}$ commuting with 
$\Ad G$. Thus the assertion follows from \rl{pullback}.

(b) Let $\wt{\C{E}}$ be the endoscopic triple for $G^{\ssc}$ induced by $\C{E}$, and
let $(\wt{a},\wt{a}';[\wt{b}])$ be the lift of $(a,a';[b])$
(see \rl{qisog}). Fix a pair of $(a,a';[b])$-indistinguishable measures 
$\phi\in\C{S}(\C{G}(E))$ and $\phi'\in\C{S}(\C{G}'(E))$.

Denote by  $\C{Z}$ the Lie algebra of $Z(G)=Z(G')$. Then $\C{G}=\C{G}^{\ssc}\oplus\C{Z}$
and $\C{G}'=\C{G}'^{\ssc}\oplus\C{Z}$, hence there exist measures
$h_i\in\C{S}(\C{Z}(E))$, $f_i\in\C{S}(\C{G}^{\ssc}(E))$ and 
$f'_i\in\C{S}(\C{G}'^{\ssc}(E))$ such that  the $h_i$'s are linearly independent, 
$\phi=\sum_i f_i\times h_i$ and $\phi'=\sum_i f'_i\times h_i$.

For each $x\in(\C{G}^{\ssc})^{\sr}(E)$, $z\in\C{Z}(E)$ and 
$\ov{\ka}\in\pi_0(\wh{G_{x+z}}^{\Gm}/Z(\wh{G})^{\Gm})$, we have  
 $O^{\ov{\ka}}_{x+z}(\phi)=\sum_i O^{\ov{\wt{\ka}}}_{x}(f_i) O_z(h_i)$, where
$\ov{\wt{\ka}}\in\pi_0(\wh{G^{\ssc}_{x}}^{\Gm})$ is the image of $\ov{\ka}$, 
and similarly for $\phi'$. Since the $h_i$'s are linearly independent, 
it follows from Lemmas \ref{L:indist} and \ref{L:qisog} that $f_i$ and 
$f'_i$ are $(\wt{a},\wt{a}';[\wt{b}])$-indistinguishable for each $i$
(compare \rco{comp} and its proof).

Let $\lan .,.\ran$ be a direct sum of pairings on $\C{G}^{\ssc}$ and $\C{Z}$.
Then $\C{F}(\phi)=\sum_i \C{F}(f_i)\times\C{F}(h_i)$ and
$\C{F}(\phi')=\sum_i \C{F}(f'_i)\times \C{F}(h_i)$. 
By our assumptions,  $e'(G)\C{F}(f_i)$ and $e'(G')\C{F}(f'_i)$ are 
$(\wt{a},\wt{a}';[\wt{b}])$-indistinguishable for each $i$, therefore 
$e'(G)\C{F}(\phi)$ and $e'(G')\C{F}(\phi')$ are  $(a,a';[b])$-indistinguishable, as claimed.

(c) The assertion follows from the arguments of \cite[II. 3]{Wa2}.
Since we do not use this result in the main body of the paper, we omit the details.
\end{proof} 

From now on, we assume that $\C{E}$ is elliptic, $G$ is semisimple and simply 
connected, and $\lan.,.\ran $ is the Killing form.

\begin{Not}
(a) Consider the natural map $[y]\mapsto[y]_{G}$ \label{yg} from the set of 
stable conjugacy classes of elements of $\C{H}^{\sr}(E)$ to the set of $E$-rational
conjugacy classes in $\C{G}^{\sr}(\ov{E})$, defined as follows.
For each $y\in\C{H}^{\sr}(E)$, denote by $b_y$ the inclusion $H_y\hra H$. Each
embedding $a:(H_y)_{\ov{E}}\hra G_{\ov{E}}$ from $[b_y]_G$ defines an embedding  
$da:(\C{H}_y)_{\ov{E}}\hra \C{G}_{\ov{E}}$, and 
we denote by  $[y]_{{G}}$ be the conjugacy class of $da(y)\in\C{G}(\ov{E})$.

(b) We say that $y\in \C{H}^{\sr}(E)$ and $x\in\C{G}^{\sr}(E)$ are {\em compatible} \label{comp2}
if $x\in [y]_{G}$. 

(c) For each  $x\in\C{G}^{\sr}(E)$, the map $b\mapsto y:=db(x)$ defines a 
bijection between embeddings of maximal tori $b:G_x\hra H$ compatible with 
$a_x:G_x\hra G$
and elements  $y\in \C{H}^{\sr}(E)$ compatible with $x$. Let 
$y\mapsto b_y$ be the inverse map. 

We will write $\ov{\ka}_{[y]}\in\pi_0(\wh{G_x}^{\Gm})$ \label{kay} instead of $\ov{\ka}_{[b_y]}$, $O^{[y]}_x$ \label{oyx} instead of 
$O^{\ov{\ka}_{[y]}}_x$, and $\lan\frac{x,x';[y]}{a,a';[b]}\ran$, \label{brac4} where $x'\in\C{G}'(E)$ is a stable conjugate of $x$, instead 
of $\lan\frac{a_{x},a_{x'};[b_y]}{a,a';[b]}\ran$.
 
\end{Not}

\subsection{Local calculations}
\begin*
\vskip 8truept
\end*
The primary goal of this subsection is to construct $(a,a';[b])$-indistinguishable 
measures whose Fourier transforms are $(a,a';[b])$-indistinguishable in some region.
We mostly follow \cite{KP}.

\begin{Emp} \label{E:difforms}
(a) For each $t\in\C{G}^{\sr}(E)$, fix a top degree
form $\om_{\C{G}_t}\neq 0$ on the vector space $\C{G}_{t}$ and identify it with the
corresponding top degree translation invariant differential form. Then  $\om_{\C{G}_t}$ defines a 
$G_t$-invariant top degree form $\nu=\nu_t:=\om_{\C{G}}\otimes(\om_{\C{G}_t})^{-1}$ on
$\C{G}/\C{G}_t$, which uniquely extends to a non-zero top degree $G$-invariant form on 
$G/G_t$, which we will also denote by $\nu$.

(b) We denote by $d\ov{g}$ and  $du$ the measures $|\nu|$ on $(G/G_t)(E)$  and  $|\om_{\C{G}_t}|$
on $\C{G}_t(E)$, respectively.

(c) Consider the map $\Pi:(G/G_t)\times\C{G}_t^{\sr}\to \C{G}^{\sr}$  given by the rule
$\Pi(\ov{g},x)=\Ad\ov{g}(x)$. Then $\Pi$ is \'etale, and we have an equality
$\Pi^*(\om_{\C{G}})=\nu\wedge\om_{\C{G}_t}$.
\end{Emp}

\begin{Emp} \label{E:qform}
{\bf Preliminaries on quadratic forms over local fields}.
(a) To every non-degenerate quadratic form $q$ on an $E$-vector space $V$ 
one associates a rank $\rk q=\dim V$, a determinant $\det q\in (\det V)^{\otimes(-2)}$
and a Hasse--Witt invariant $e(q)\in\{-1,1\}$. Any trivialization $\det V\isom E$
associates to $\det q$ an element of $E\m$. Moreover, 
the class of  $\det q$ in $E\m/(E\m)^2$ is independent on the trivialization. 

To each isomorphism class $(q,\psi)$, where $\psi$ is a 
non-trivial additive character of $E$, Weil \cite{We} associated an $8$th root of unity 
$\gm(q,\psi)$. For each non-zero top degree form $\nu$ on $V$, we set
$c(q,\nu,\psi):=\gm(q,\psi)|\det(q)/\nu^2|^{-1/2}$. 

(b) Weil proved that for every non-degenerate quadratic forms $q$ and $q'$ 
satisfying 
$\rk q=\rk q'$ and  $\det q \equiv\det q'\mod (E\m)^2$, we have 
$\gm(q,\psi)/\gm(q,\psi)=e(q)e(q')$. 

(c) To each $t\in \C{G}^{\sr}(E)$ and $y,z\in\C{G}^{\sr}_t(E)$ we associate 
a non-degenerate quadratic form $q=q_{y,z}:\ov{x}\mapsto \lan y,(\ad \ov{x})^2(z)\ran$ on 
$V:=\C{G}/\C{G}_t$. Then the form $\nu$ on $V$, chosen in \re{difforms} (a) gives rise to
an invariant $c(q_{y,z},\nu,\psi)$.

(d) Let $t'\in\C{G}'^{\sr}(E)$ be a stable conjugate of $t$,  $\nu'$ the  
form on $G'/G'_{t'}$ induced by $\nu$, and  
$\varphi_{t,t'}$ the canonical  isomorphism $\C{G}_t\isom\C{G}'_{t'}$ 
from \rc{measure} (b).
For each $y,z\in\C{G}_t^{\sr}$, we set $y'=\varphi_{t,t'}(y)$ and $z'=\varphi_{t,t'}(z)$.

It follows from results of \cite{KP} that 
$c(q_{y,z},\nu,\psi)/c(q_{y',z'},\nu',\psi)=e'(G)e'(G')$.
In particular, $c(q_{y,z},\nu,\psi)=c(q_{y',z'},\nu',\psi)$ if $G'=G$.
Indeed, since $\det q_{y,z} /\nu^2=\det q_{y',z'} /\nu'^2$ and $\rk q=\rk q'$, 
the assertion is a combination of the result of Weil (see (b)) and  
\cite[Lem. 2.7.5 and the remark following it]{KP}. 
\end{Emp}

\begin{Not} \label{N:constr}
Fix $t\in\C{G}^{\sr}(E)$ and sufficiently small open compact subgroups $K\subset G(E)$
and $U\subset \C{G}_t(E)$ such that $t+U\subset \C{G}^{\sr}_t(E)$ and 
$K\cap\Norm_{G(E)}(\C{G}_t)=K\cap G_t(E)$.

For each $a\in E\m$, put 
$S_a:=\Ad K(a t+U)\subset\C{G}(E)$, and we denote by $\chi_a=\chi_{a,t,K,U}$ 
be the characteristic function of $S_a$. 
\end{Not}

\begin{Emp} \label{E:stationary}
{\bf Stationary phase principle.}
For each $u\in\C{G}^{\sr}(E)$ and $x\in\C{G}(E)$, we define function 
$f_{x,u}:(G/G_u)(E)\to E$ by the rule $f_{x,u}(\ov{g}):=\lan x,\Ad \ov{g}(u)\ran$. 
Then $\ov{g}\in (G/G_u)(E)$ is a critical point for 
$f_{x,u}$ if and only if $x\in\Ad\ov{g}(\C{G}_{u})=\C{G}_{\Ad\ov{g}(u)}$. In this case, 
the corresponding quadratic form on $T_{\ov{g}}(G/G_u)=\C{G}/\C{G}_{u}$ is 
$q_{\Ad\ov{g}^{-1}(x),u}$.

By the stationary phase principle (see \cite[Lem 2.5.1]{KP}), for each compact
subset $C\subset G^{\sr}(E)$ there exists $N_0=N_0(t,K,U,c)\in\B{N}$ such that for each 
$x\in C$, $u\in t+U$ and $a\in E\m$ with $|a|>N_0$, the integral 
$\int_{K/K\cap G_t(E)}\psi(a\lan x,\Ad\ov{g}(u)\ran) d\ov{g}$ equals
\begin{equation} \label{E:int}
c(q_{\Ad\ov{k}(x),u},\nu,\psi)\psi( a\lan x, u\ran )|a|^{-\frac{1}{2}\dim \C{G}/\C{G}_t},
\end{equation}
if there exists (a unique) element $\ov{k}\in K/K\cap G_t(E)$ such that
$\Ad\ov{k}(x)\in\C{G}^{\sr}_t(E)$, and vanishes otherwise.

Indeed, the map $f_{x,u}|_{K/K\cap G_t(E)}$ has a unique non-degenerate critical point
$\ov{g}=\ov{k}^{-1}$ in the former case and has no critical points in the latter one.
Thus the assertion for $a\in (E\m)^2$ follows from \cite[Lem 2.5.1]{KP}.
Since the quotient $E\m/(E\m)^2$ is finite, the general case 
now follows from the previous one applied to the compact set $\sqcup_{b} bC$, 
where $b\in E\m$ runs over a set of representatives of $E\m/(E\m)^2$.
\end{Emp}

\begin{Lem} \label{L:constr}
Let $t,K$ and $U$ be as in \rn{constr}, and let $[x]\subset\C{G}^{\sr}(E)$ be 
a stable conjugacy class. Then there exists 
$N=N(t,K,U,[x])\in\B{N}$ such that for each $x\in[x]$ and $a\in E\m$ with $|a|>N$, 
the Fourier transform $\C{F}(\chi_a)(x)$ equals 
\begin{equation} \label{E:fourier}
c(q_{\Ad\ov{k}(x),t},\nu,\psi)\psi(a\lan x,t\ran )|a|^{\frac{1}{2}\dim \C{G}/\C{G}_t}
\int_U \psi(\lan x,u\ran )du,
\end{equation}
if there exists (a unique) element $\ov{k}\in K/K\cap G_t(E)$ such that
$\Ad\ov{k}(x)\in\C{G}^{\sr}_t(E)$, and vanishes otherwise.
\end{Lem}
\begin{proof}  

First we claim that there exists a compact set $C_0\subset \C{G}(E)$ containing the support
of $\C{F}(\chi_a)$ for all $a\in E\sm\C{O}$. 
To show this we will find an $\C{O}$-lattice $L\subset \C{G}$ such that
$S_a+L=S_a$ for all $a\in E\sm\C{O}$. Then the dual lattice
$C_0:=\{x\in\C{G}\,|\,\psi(\lan x,L\ran)=1\}$ satisfies the required property.

Recall that the map $(\ov{k},u)\mapsto\Ad k(t+u)$ gives an analytic isomorphism
$F:(K/K\cap G_t(E))\times U\isom S_1$, and we denote by $\pi:S_1\to U$
the composition $\pr_2\circ F^{-1}$. Since $S_1$ is compact, $\pi$ has a bounded derivative.
Therefore there exists a lattice $L\subset\C{G}(E)$ such that for each $b\in E$ 
and $x\in S_1$, we have $\pi(x+bL)\subset \pi(x)+ bU$. Shrinking $L$ if necessary, 
we can moreover assume that $S_1+L=S_1$.

Fix $a\in E\sm\C{O}$. Then $a^{-1}S_a=\Ad K(t+a^{-1}U)\subset S_1$ and, moreover, $a^{-1}S_a$ 
is the preimage of $a^{-1}U\subset U$ under $\pi$. Therefore $a^{-1}S_a+a^{-1}L\subset S_1+L=S_1$ 
and $\pi(a^{-1}S_a+a^{-1}L)\subset \pi(a^{-1}S_a)+a^{-1}U=a^{-1}U$.
Hence $a^{-1}S_a+a^{-1}L\subset\pi^{-1}(a^{-1}U)=a^{-1}S_a$, thus $S_a+L=S_a$, as claimed.

As the intersection $\C{G}^{\sr}_t\cap[x]$ is finite, the set   
$\Ad K(\C{G}^{\sr}_t\cap[x])$ and therefore also $C:=C_0\cup \Ad K(\C{G}^{\sr}_t\cap[x])$ is compact.
Since $[x]\subset\C{G}(E)$ is closed, the intersection $C\cap[x]$ is compact as well.
Take any $N\geq N_0(t,K,U,C\cap[x])$ (see \re{stationary}) such that
quadratic forms  $q_{x,u}$ and $q_{x,t}$ are isomorphic for each element $x$ 
of the finite set $\C{G}_t(E)\cap[x]$ and each $u\in t+a^{-1}U$ with $|a|>N$. We claim that
this $N$ satisfies the required properties.

Indeed, the Fourier transform 
$\C{F}(\chi_a)(x)=\int_{S_a} \psi(\lan x,y\ran) dy$ equals 
\[
|a|^{\dim\C{G}}\int_{a^{-1} S_a} \psi(a \lan x,y\ran) dy=
|a|^{\dim\C{G}}\int_{t+a^{-1}U}du\int_{K/K\cap G_t(E)}
\psi(x\lan a,\Ad\ov{g}(u)\ran)d\ov{g} 
\]
(by \re{difforms} (c) and our assumptions in \rn{constr}). 
Therefore the assertion  for $x\in [x]\cap C$ follows from 
\re{stationary}. Finally, if $x\in[x]\sm C$, then $x\notin \Ad K(\C{G}^{\sr}_t)$ and 
$\C{F}(\chi_a)(x)=0$, implying the assertion in the remaining case.
\end{proof}

\begin{Cor} \label{C:constr2}
For each triple $(y,x,x')$, where $x\in\C{G}^{\sr}(E)$ and 
$x'\in\C{G}'^{\sr}(E)$ are stably conjugate, and  $y\in\C{H}^{\sr}(E)$ is compatible with $x$, 
there exist measures $\phi\in \C{S}(\C{G}(E))$ and $\phi'\in \C{S}(\C{G}'(E))$ 
satisfying the following properties:

(i) $\phi$ and $\phi'$ are supported on elements stably conjugate to $\C{G}^{\sr}_x(E)$.

(ii) $O^{[y]}_{u}(\phi)=\lan\frac{x,x';[y]}{a,a';[b]}\ran O^{[y]}_{u'}(\phi')$
for each $u\in \C{G}^{\sr}_x(E)$ and $u'=\phi_{x,x'}(u)\in \C{G}'^{\sr}_{x'}(E)$.

(iii) $O^{\ov{\xi}}_{u}(\phi)=O^{\ov{\xi}}_{u'}(\phi')=0$ 
for each $u\in \C{G}^{\sr}_x(E)$, $u'=\phi_{x,x'}(u)$ and
$\ov{\xi}\neq \ov{\ka}_{[y]}$.

(iv) $O^{\ov{\xi}}_{x}(\C{F}(\phi))=O^{\ov{\xi}}_{x'}(\C{F}(\phi'))=0$ 
for each $\ov{\xi}\neq \ov{\ka}_{[y]}$.

(v) $e'(G)O^{[y]}_{x}(\C{F}(\phi))=e'(G') \lan\frac{x,x';[y]}{a,a';[b]}\ran 
O^{[y]}_{x'}(\C{F}(\phi'))\neq 0$.
\end{Cor}
\begin{proof}
(compare \rco{constr}). Let $x_1=x, x_2,\ldots,x_l$ be all elements of $\C{G}_x(E)\cap [x]$. 
Pick $t\in\C{G}^{\sr}_x(E)$ such that all $\lan x_j,t\ran$ are distinct, and put 
$t':=\phi_{x,x'}(t)\in \C{G}'^{\sr}_{x'}(E)$. Choose sufficiently small subgroups
$K\subset G(E)$, $K'\subset G'(E)$ and $U\subset\C{G}_t(E)$ satisfying the assumptions
of \rn{constr} and such that $\psi(\lan x_j,u\ran)=1$ for each $u\in U$ and each $j=1,\ldots,l$.
 
Let $t_1=t,t_2,\ldots,t_k\in\C{G}(E)$ and $t'_1=t',t'_2,\ldots,t'_k\in\C{G}'(E)$ be sets
of representatives of conjugacy classes stably conjugate to $t$ and $t'$. 
For each $i=1,\ldots,k$, put $U_i:=\varphi_{t,t_i}(U)\subset\C{G}_{t_i}(E)$ and 
$U'_i:=\varphi_{t,t'_i}(U)\subset\C{G}'_{t'_i}(E)$.

For each $i=1,\ldots,k$, let $\nu_i$ and $\nu'_i$ be the forms on $G/G_{t_i}$ and 
$G'/G'_{t'_i}$ respectively, induced by $\nu$. For each $a\in E\sm\C{O}$, put 
$\phi_i=\frac{1}{|\nu_i|(K/K\cap G_{t_i}(E))}\chi_{a,t_i,U_i,K}dx$ and
$\phi'_i=\frac{1}{|\nu'_i|(K'/K'\cap G'_{t'_i}(E))}\chi_{a,t'_i,U'_i,K'}dx'$. Define 
$\phi$ and $\phi'$ by the formulas 
$\phi:=\sum_i\lan\inv(t,t_i),\ov{\ka}_{[y]}\ran^{-1}\phi_i$ and 
$\phi':=\sum_i\lan\frac{a_t,a_{t'_i};[b_y]}{a,a';[b]}\ran^{-1}\phi'_i$.
Then $\phi$ and $\phi'$ clearly satisfy properties (i)--(iii) for each $a\in E\sm\C{O}$, 
so it remains to show the existence of $a$ for which properties (iv) and (v) are satisfied.

Let $N$ be the maximum of the $N(t_i,U_i,K,[x])$'s and the $N(t'_i,U'_i,K',[x'])$'s.
Then \rl{constr} implies that for each $\ov{\xi}\in\pi_0(\wh{G_x}^{\Gm})$ and 
$a\in E\m$ with $|a|>N$, we have 
\begin{equation} \label{E:vanish}
O^{\ov{\xi}}_x(\C{F}(\phi_1))=|a|^{\frac{1}{2}\dim\C{G}/\C{G}_t}|\om_{\C{G}_t}|(U)
\sum_{j=1}^l\lan\inv(x,x_j),\ov{\xi}\ran
c(q_{x_j,t},\nu,\psi)\psi(a\lan x_j,t\ran).
\end{equation}
Moreover, for each $i=1,\ldots,n$, observation \re{qform} (d) and 
the equality $\varphi_{t,t_i}^*(\om_{\C{G}_{t_i}})=\om_{\C{G}_t}$ imply that
\begin{equation} \label{E:eq}
O^{\ov{\xi}}_x(\C{F}(\phi_i))=\lan\inv(t,t_i),\ov{\xi}\ran O^{\ov{\xi}}_x(\C{F}(\phi_1)).
\end{equation}
Then (\ref{E:vanish}) and (\ref{E:eq}) together with similar formulas for the $\phi'_i$'s 
imply the equalities from (iv) and (v) and that 
$O^{[y]}_{x}(\C{F}(\phi))=kO^{[y]}_{x}(\C{F}(\phi_1))$.

It remains to show the existence of $a\in E\m$ with $|a|>N$ such that 
$O^{[y]}_{x}(\C{F}(\phi_1))\neq 0$.
Since all $\lan x_j,t\ran$ are distinct, the functions 
$a\mapsto \psi(a\lan x_j,t\ran)$ are linearly independent. The assertion now follows from 
(\ref{E:vanish}).
\end{proof}


Later we will need the following result

\begin{Lem} \label{L:unit}
Let $G$ be an unramified semisimple simply connected group over $E$, 
$K\subset G(E)$ a hyperspecial subgroup,  
$\C{K}\subset\C{G}(E)$ 
the corresponding subalgebra, $1_{\C{K}}$ the characteristic function of 
$\C{K}$, and $\C{E}=(H,[\eta],\ov{s})$ an endoscopic triple for $G$.

(a) If $H$ is ramified, then $O^{[y]}_x(1_{\C{K}}dx)=0$ for all 
compatible elements $y\in\C{H}^{\sr}(E)$ and $x\in\C{G}^{\sr}(E)$.

(b) If $H$ is unramified, then there exists an open neighborhood of zero $\Om\subset\C{H}(E)$
such that $O^{[y]}_x(1_{\C{K}}dx)\neq 0$ for all compatible elements $y\in\Om\cap\C{H}^{\sr}(E)$
and $x\in\C{G}^{\sr}(E)$.

(c) If $x\in\C{K}$ has a regular reduction modulo $\frak{m}$, 
then every stably conjugate
$x'\in\C{K}$ of $x$ is $K$-conjugate. In particular, $O^{\ov{\xi}}_x(1_{\C{K}}dx)\neq 0$ for each 
$\ov{\xi}\in\pi_0(\wh{G_x}^{\Gm})$.
\end{Lem}
\begin{proof}
Since in the notation of \cite{Wa1,Wa2}, $O^{[y]}_x(1_{\C{K}}dx)$ is a non-zero multiple of 
$J^{G,H}(y,1_{\C{K}})$, the assertions follow from 
\cite[7.2 and 7.4]{Wa1} and \cite[III, Prop.]{Wa2} (compare \cite[Prop. 7.1 and 7.5]{Ko2}).
\end{proof} 

\subsection{Global results}
\begin*
\vskip 8truept
\end*
For the proof of \rt{Wa1} we will use global methods. In this subsection we will recall 
necessary notation and results.

\begin{Emp} \label{E:glob}
(a) Let $\un{E}$ be a number field, which we will always assume to be totally 
imaginary, $\un{\Gm}$ the absolute Galois group of $\un{E}$, and 
$\B{A}$ the ring of ad\'eles of $\un{E}$. We denote by $V,V_{\infty}$ and $V_f$ the set of 
all places, all infinite places and all finite places of $\un{E}$, respectively. 
 For each $v\in V_f$,
we have a natural conjugacy class of embeddings $\un{\Gm}_v\hra\un{\Gm}$.
For every object $\un{S}$ over $\un{E}$ and 
each $v\in V$, we will denote by $\un{S}_v$ the corresponding object over $\un{E}_v$.

 
(b) For every reductive group $\un{G}$ over $\un{E}$, consider a sequence
\begin{equation} \label{E:kot}
H^1(\un{E},\un{G})\to\bigoplus_{v\in V_f}H^1(\un{E}_v,\un{G})\overset{*}{\lra}
\pi_0(Z(\wh{\un{G}})^{\un{\Gm}})^D,
\end{equation} 
where the restriction of $*$ to $H^1(\un{E}_v,\un{G})$ 
is the composition of the isomorphism $\C{D}_G:H^1(\un{E}_v,\un{G})\isom
\pi_0(Z(\wh{\un{G}})^{\un{\Gm}_v})^D$ from \re{TN} and the projection 
$\pi_0(Z(\wh{\un{G}})^{\un{\Gm}_v})^D\to\pi_0(Z(\wh{\un{G}})^{\un{\Gm}})^D$.
Kottwitz proved (see \cite[Prop 2.6]{Ko2}) that the sequence (\ref{E:kot}) is exact.

(c) For each $c\in H^1(\un{E},\un{G})$, $\ka\in\pi_0(Z(\wh{\un{G}})^{\un{\Gm}})$ and 
$v\in V_f$, we denote by $c_v\in H^1(\un{E}_v,\un{G})$ and 
$\ka_v\in\pi_0(Z(\wh{\un{G}})^{\un{\Gm}_v})$ the images of $c$ and $\ka$, respectively.
Then Kottwitz' theorem from (b) asserts that $\prod_{v\in V_f}\lan c_v,\ka_v\ran=1$.
\end{Emp}

\begin{Lem} \label{L:globalcoh}
Let $\un{G}$ be a reductive group over $\un{E}$, and $u\in V_f$.

(a) If $u$ is inert in the splitting field $\un{E}[\un{G}^*]$ of the quasi-split inner form 
$\un{G}^*$ of $\un{G}$, then the diagonal map 
 $H^1(\un{E},\un{G})\to\bigoplus_{v\neq u}H^1(\un{E}_v,\un{G})$
is surjective.

Assume in addition that $\un{G}$ is either semisimple and simply connected or adjoint. Then 

(b) The map from (a) is an isomorphism.

(c) Let $\un{T}$ be a maximal torus of $\un{G}$, and let $c\in H^1(\un{E},\un{G})$ be such that
$c_v\in H^1(\un{E}_v,\un{G})$ belongs to $\im[H^1(\un{E}_v,\un{T})\to  H^1(\un{E}_v,\un{G})]$ 
for each $v\neq u$. 
Then $c\in \im[H^1(\un{E},\un{T})\to  H^1(\un{E},\un{G})]$ in each of the following cases:

 $(i)$  $u$ is inert in the splitting field $\un{E}[\un{T}]$ of $\un{T}$;

$(ii)$  $u$ inert in $\un{E}[\un{G}^*]$, and  $\un{T}_u\subset\un{G}_u$ is elliptic.
\end{Lem}
\begin{proof}
(a), (b) By assumption, we have $Z(\wh{\un{G}})^{\un{\Gm}}=Z(\wh{\un{G}})^{\un{\Gm}_u}$, so assertion 
(a) follows from the exactness of (\ref{E:kot}) while assertion (b) follows 
from the Hasse principle.

(c) Consider commutative diagram
\[
\CD
H^1(\un{E},\un{T})@>A>>\bigoplus_{v\neq u} H^1(\un{E}_v,\un{T})\\
@VVV          @VDVV\\
H^1(\un{E},\un{G})@>B>> \bigoplus_{v\neq u}H^1(\un{E}_v,\un{G}).
\endCD
\]
In both cases, $u$ is inert in $\un{E}[\un{G}^*]$, hence $B$ is injective (by (b)).
Since by our assumption, $B(c)=(c_v)_{v\neq u}$ belongs to $\im D$, it will suffice to show that
$A$ is surjective. In the case $(i)$, the surjectivity of $A$ follows from (a). In the case $(ii)$, 
the canonical map $\pi_0(\wh{\un{T}}^{\un{\Gm}})=
\wh{\un{T}}^{\un{\Gm}}/Z(\wh{\un{G}})^{\un{\Gm}}\hra\wh{\un{T}}^{\un{\Gm}_v}/
Z(\wh{\un{G}})^{\un{\Gm}_v}=\pi_0(\wh{\un{T}}^{\un{\Gm}_u})$
is injective. Therefore surjectivity of $A$ follows from the exactness of (\ref{E:kot}).
\end{proof}

From now on, $\un{G}$ is a semisimple and simply connected group over $\un{E}$.

\begin{Emp} \label{E:globend}
(a) Let $\C{E}'=(H',[\eta'],\ov{s}')$ be an endoscopic triple for $\un{G}$. For each $v\in V_f$,
$\C{E}'$ gives rise to an (isomorphism class of an) endoscopic triple 
$\C{E}'_v=(H'_v,[\eta'_v],\ov{s}'_v)$ for $\un{G}_v$. 

In particular, if $\C{E}'\cong\C{E}_{([\un{a}],\un{\ka})}$ 
for a certain pair $(\un{a},\un{\ka})$ consisting of an embedding of a
maximal torus $\un{a}:\un{T}\hra\un{G}$ and $\un{\ka}\in\wh{T}^{\un{\Gm}}$, 
then $\C{E}'_v\cong\C{E}_{([\un{a}_v],\un{\ka}_v)}$.

(b) Let $\C{E}'_i=(H'_i,[\eta'_i],\ov{s}'_i)$, $i=1,2$ be a pair of endoscopic triples 
for $\un{G}$. Assume that there exists $v\in V_f$, inert in both splitting fields
$\un{E}[H'_1]$ and $\un{E}[H'_2]$, such that $(\C{E}'_1)_v\cong(\C{E}'_2)_v$.
Then $\C{E}'_1\cong \C{E}'_2$. 

 Indeed,  the image of 
$\rho_{H'_i}:\un{\Gm}\to\Out(\wh{H'_i})$ coincides with that of 
$\rho_{(H'_i)_v}:\un{\Gm}_v\to\Out(\wh{H'_i})$ for each $i=1,2$, therefore by \rr{isom},
the map 
\[
\Isom(\C{E}'_1,\C{E}'_2)/(H'_1)^{\ad}(\un{E})\to 
\Isom((\C{E}'_1)_v,(\C{E}'_2)_v)/(H'_1)^{\ad}(\un{E}_v)\] 
is bijective.


(c) Let $\C{E}'=(H',[\eta'],\ov{s}')$ be an endoscopic triple for $G$,
$\un{\varphi}:\un{G}\to \un{G}'$ an inner twisting, and $(\un{a}'_i,\un{a}_i,[\un{b}_i])$ 
be two triples consisting of  
stably conjugate embeddings of maximal tori $\un{a}_i:\un{T}_i\hra \un{G}$ and 
$\un{a}'_i:\un{T}_i\hra \un{G}'$, and 
stably conjugate classes $[\un{b}_i]$ of embeddings of maximal tori $\un{T}_i\hra H'$, compatible with
$\un{a}_i$. Then we have the following product formula
\[
\prod_{v\in V_f}\lan\frac{(\un{a}_1)_v,(\un{a}'_1)_v;[\un{b}_1]_v}
{(\un{a}_2)_v,(\un{a}'_2)_v;[\un{b}_2]_v}\ran_{\C{E}'_v}=1.
\]
Indeed, consider elements 
$c:=\ov{\inv}((\un{a}_1,\un{a}'_1);(\un{a}_2,\un{a}'_2))\in 
H^1(\un{E},(\un{T}_1\times \un{T}_2)/Z(\un{G}))$ 
(see \re{coh}) and 
$\ka:=\ka([b_1],[b_2])(\ov{s}')\in\pi_0((\{(\un{T}_1\times \un{T}_2)/Z(\un{G})\}\:\wh{}\,)^{\un{\Gm}}$ 
(see \re{kainv}). Since $\lan\frac{(\un{a}_1)_v,(\un{a}'_1)_v;[\un{b}_1]_v}
{(\un{a}_2)_v,(\un{a}'_2)_v;[\un{b}_2]_v}\ran_{\C{E}'_v}=\lan c_v,\ka_v\ran$ 
for each $v\in V_f$, 
the product formula follows from Kottwitz' theorem (see \re{glob} (c)).
\end{Emp}

\begin{Emp}
Denote by $dg$ and $dx=\prod_v dx_v$ the Tamagawa measures on $\un{G}(\B{A})$ and 
$\un{\C{G}}(\B{A})$, respectively (defined by a non-zero translation invariant 
top degree differential form $\om_{\un{G}}$ on $\un{G}$).

For each $v\in V_{\be}$, denote by $\C{S}(\un{G}(\un{E}_v))$ the space of measures
on $\un{G}(\un{E}_v)$ of the form $f_v d x_v$, where $f_v$ is a Schwartz function.
Put $\C{S}(\un{\C{G}}(\B{A})):=\otimes'_{v\in V}\C{S}(\un{\C{G}}(\un{E}_v))$, and 
fix a non-trivial character $\un{\psi}=\prod_v \psi_v:\B{A}/\un{E}\to \B{C}\m$.
Then $\un{\psi}$ gives rise to a Fourier transform 
$\C{F}:\C{S}(\un{\C{G}}(\B{A}))\to \C{S}(\un{\C{G}}(\B{A}))$ such that 
$\C{F}(\otimes_v \phi_v)=\otimes_v\C{F}(\phi_v)$, where $\C{F}(\phi_v)$ is the Fourier
transform of $\phi_v$ corresponding to a measure $dx_v$ and
a character $\psi_v$.

For each $\un{x}\in\un{\C{G}}^{\sr}(\un{E})$, $\un{\ka}\in \wh{\un{G}_{\un{x}}}^{\un{\Gm}}$, 
and $\phi=\otimes_v \phi_v\in \C{S}(\un{\C{G}}(\B{A}))$ put
$O^{\un{\ka}}_{\un{x}}(\phi)=:\prod_{v\in V}O^{\un{\ka}_v}_{\un{x}_v}(\phi_v)$, where
$O^{\un{\ka}_v}_{\un{x}_v}:=O_{\un{x}_v}$ for each $x\in V_{\be}$.
It follows from Kottwitz' theorem (see \re{glob} (c)) that generalized function $O^{\un{\ka}}_{\un{x}}$ 
depends only on the stable conjugacy class of $\un{x}$.
\end{Emp}

The main technical tool for the proof of \rt{Wa1} is the following simple version of the 
trace formula.
Let $\theta$ be a generalized function on $\un{\C{G}}(\B{A})$ defined by the rule 
$\theta(fdx)=\sum_{x\in\C{G}(F)}f(x)$. For each $g\in \un{G}(\B{A})$, put
$\theta^g:=(\Ad g)^*(\theta)$.

\begin{Prop} \label{P:trfor} 
(a) Let $\phi=\otimes_v\phi_v$ be an element of $\C{S}(\un{\C{G}}(\B{A}))$ 
such that $Supp(\phi)\cap \Ad\un{G}(\B{A})(\un{\C{G}}(\un{E}))$ 
consists of regular elliptic elements. 
Then the integral
\[
\Theta(\phi):=\int_{\un{G}(\B{A})/\un{G}(\un{E})}\theta^g(\phi)dg
\]
converges absolutely. 
Furthermore, 
\[
\Theta(\phi)=\sum_{\un{x}}\sum_{\un{\ka}\in \wh{\un{G}_{\un{x}}}^{\un{\Gm}}} 
O^{\un{\kappa}}_{\un{x}}(\phi),
\]
where $\un{x}$ runs over the set of regular elliptic stably conjugacy classes of 
$\un{\C{G}}(\un{E})$.

(b) If $\C{F}(\phi)$ also satisfies the support assumption of (a), then 
$\Theta(\phi)=\Theta(\C{F}(\phi))$.
\end{Prop}
\begin{proof}

(a) The first assertion (see \cite[10.8]{Wa1}) is a direct analog of the corresponding
result of Arthur, while the second one  (see, for example, \cite[Thm 3.2.1]{KP}) is an analog
of a result of Kottwitz.

(b) By the Poisson summation formula, we have  $\C{F}(\theta)=\theta$. Therefore the assertion 
follows from the absolute convergence of  $\Theta(\phi)$ and $\Theta(\C{F}(\phi))$.
\end{proof}

To apply the trace formula, we will embed our local data into global ones. 

\begin{Cl} \label{C:global}
There exist a totally imaginary number field $\un{E}$, two finite places $w$ and $u$ of 
$\un{E}$, a semisimple simply connected
group $\un{G}$ over $\un{E}$,  an inner twisting 
$\un{\varphi}:\un{G}\to \un{G}'$, an endoscopic triple
$\un{\C{E}}=(\un{H},[\un{\eta}],\ov{\un{s}})$ for $\un{G}$, a tori $\un{T}$ over $\un{E}$,
a pair of stably conjugate embeddings of maximal tori $\un{a}:\un{T}\hra\un{H}$ and 
$\un{a}':\un{T}\hra\un{G}'$, and an embedding
$\un{b}:\un{T}\hra\un{G}$ compatible with $\un{a}$ satisfying the following conditions:

(a) $\un{E}_{w}\cong E$, $\un{G}_w\cong G$, $\un{\varphi}_w\cong \varphi$, 
$\un{\C{E}}_w\cong\C{E}$, $\un{T}_w\cong T$, $[\un{b}]_w=[b]$, while
$\un{a}_w$ are $\un{a}'_w$ are conjugate
to $a$ and $a'$, respectively.

(b) For each $v\neq u,w$, the groups $\un{G}_v$ and $\un{G}'_v$ are quasi-split, 
and $\un{\varphi}_v$ is trivial. Moreover, after we identify 
$\un{G}_v$ with $\un{G}'_v$ by means of some $\un{G}(\ov{\un{E}_v})$-conjugate
of $\un{\varphi}_v$, embeddings $\un{a}_v$ and $\un{a}'_v$ are conjugate.

(c) $u$ is inert in $\un{E}[\un{T}]$, and $\un{\C{E}}_u$ is elliptic.
\end{Cl}
\begin{proof}

(I) Put $E':=E[T]$, and set $\Gm':=\Gal(E'/E)$. Choose a dense subfield $F'$ of $E'$, which is 
a finite extension of $\B{Q}$. Increasing $F'$, we may assume that $F'$ is $\Gm'$-invariant. Set 
$F:=(F')^{\Gm'}$, and let $w_0$ be the prime of $F$, corresponding to the 
embedding $F\hra (E')^{\Gm'}=E$. In particular, $F_{w_0}\cong E$.
Choose a totally imaginary quadratic extension $\un{E}$ of $F_0$ such that 
$w_0$ splits in $\un{E}$, and let $w$ and $u$ be the primes of $\un{E}$ lying over $w_0$. 
Finally, let $\un{E}'$ be the composite field $\un{E}\cdot F'$. 
We have natural identifications $\un{E}_w\cong\un{E}_u\cong E$ and 
$\Gal(\un{E}'/\un{E})\cong\Gal(E'/E)$, and both $w$ and $u$ are inert in $\un{E}'$.

(II) Let $\varphi^*:G\to G^*$ be the inner twisting such that $G^*$ is quasi-split. 
Since $G,H$ and $T$ split
over $E'$, the homomorphisms $\rho_G, \rho_H$ and $\rho_T$ factor through 
$\Gal(E'/E)$. We denote by $\un{G}^*$ (resp. $\un{H}$, resp. $\un{T}$) the 
quasi-split group over $\un{E}$ such that
$\wh{\un{G}^*}=\wh{G}$ (resp. $\wh{\un{H}}=\wh{H}$, resp. $\wh{\un{T}}=\wh{T}$) 
and $\rho_{\un{G}^*}$ (resp. $\rho_{\un{H}}$, resp. $\rho_{\un{T}}$)
is the composition of the projection $\un{\Gm}\to\Gal(\un{E}'/\un{E})\cong\Gal(E'/E)$
and the homomorphism $\rho_G:\Gal(E'/E)\to\Out(\wh{G})$ 
(resp. $\rho_H:\Gal(E'/E)\to\Out(\wh{H})$, resp.  $\rho_T:\Gal(E'/E)\to\Out(\wh{T})$).

By construction, we have $\un{G}^*_w\cong G^*$, $\un{H}_w\cong H$ and $\un{T}_w\cong T$.
Moreover, the conjugacy classes of embeddings $\wh{T}\hra\wh{G}$ and $\wh{T}\hra\wh{H}$
corresponding to $a:T\hra G$ and $[b]$ are $\un{\Gm}$-invariant.
Therefore they come from stable conjugacy classes of embeddings
$\un{a}^*:\un{T}\hra \un{G}^*$ and $\un{b}:\un{T}\hra\un{H}$.
Furthermore, $(\un{a}^*)_w$ is stably conjugate to $a$, and  $[\un{b}]_w=[b]$.

(III) Since $u$ is inert in $\un{E}'=\un{E}[\un{T}]\supset\un{E}[\un{G}^*]$, 
the canonical map   $H^1(\un{E},(\un{G}^*)^{\ad})\to\bigoplus_{v\neq u}H^1(\un{E}_v,(\un{G}^*)^{\ad})$ is an isomorphism (by \rl{globalcoh} (b)).
Hence there exist unique inner twistings $\un{\varphi}^*:\un{G}^*\to\un{G}$ and 
 $\un{\varphi}:\un{G}\to\un{G}'$ such that $\un{\varphi}^*_w\cong(\varphi^*)^{-1}$, 
$\un{\varphi}_w\cong\varphi$, while $\un{\varphi}^*_v$ and $\un{\varphi}_v$ are trivial 
for all $v\neq w,u$.  

Applying \rl{globalcoh} (c) for the embedding
$\un{a}^*:\un{T}/Z(\un{G})\hra (\un{G}^*)^{\ad}$ we conclude from 
\re{inv} (b) that there exist embeddings $\un{a}:\un{T}\hra\un{G}$ and  
$\un{a}':\un{T}\hra\un{G}'$ stably conjugate to $\un{a}^*$. Applying now
\rl{globalcoh} (c) for $\un{T}$, we can further replace $\un{a}$ and $\un{a}'$ so that 
$\un{a}_w$ is conjugate to $a$, $\un{a}'_w$ is conjugate to $a'$, while 
$\un{a}_v$ and $\un{a}'_v$ are conjugate for all $v\neq w,u$.

(IV) Since $w$ is inert in $\un{E}[\un{T}]\supset \un{E}[\un{H}]$, we get that 
$Z(\wh{\un{H}})^{\un{\Gm}}=Z(\wh{H})^{{\Gm}}$, and the conjugacy class 
$[\eta]$ of embeddings $\wh{\un{H}}=\wh{H}\hra\wh{G}=\wh{\un{G}}$ is $\un{\Gm}$-invariant.
Hence the triple $\un{\C{E}}:=(\un{H},[\eta],\ov{s})$ is an endoscopic triple for $\un{G}$.
Moreover, 
$\un{\C{E}}_w\cong\C{E}$ and  $\un{\C{E}}_u$ is elliptic. Indeed,
$u$ is inert in $\un{E}[\un{H}]$, therefore 
$Z(\wh{\un{H}})^{\un{\Gm}_u}=Z(\wh{\un{H}})^{\un{\Gm}}=Z(\wh{H})^{\Gm}$ is finite.
\end{proof}

\subsection{Proof of \rt{Wa1}}
\begin*
\vskip 8truept
\end*

\begin{Emp} \label{E:1}
Fix $(a,a';[b])$-indistinguishable $\phi\in\C{S}(\C{G}(E))$ and $\phi'\in\C{S}(\C{G}'(E))$.
We want to show that $e'(G)\C{F}(\phi)$ and $e'(G')\C{F}(\phi')$ are 
$(a,a';[b])$-indistinguishable. 
Using \rl{indist} and the symmetry between $G$ and $G'$,
it will suffice to check that for each compatible  $x\in\C{G}^{\sr}(E)$ and $y\in\C{H}^{\sr}(E)$
we have:

$(i)$  $O^{[y]}_{x}(\C{F}(\phi))=0$, if $x$ does not have a stable conjugate in $\C{G}'(E)$;

$(ii)$ $e'(G)O^{[y]}_{x}(\C{F}(\phi))=
e'(G')\lan\frac{x,x';[y]}{a,a';[b]}\ran O^{[y]}_{x'}(\C{F}(\phi'))$, for each 
stable conjugate  $x'\in\C{G}'^{\sr}(E)$ of $x$.
\end{Emp}


\begin{Emp} \label{E:Omega} 
Fix $x$ and $y$ as in \re{1}. Let $\Om_x\subset\C{G}_x^{\sr}(E)$ 
be an open neighborhood of $x$ such that 
$O^{\ov{\ka}_[y]}_{\wt{x}}(\C{F}(\phi))=O^{\ov{\ka}_[y]}_{x}(\C{F}(\phi))$
for each $\wt{x}\in \Om_x$, and  
$O^{\ov{\ka}_[y]}_{\wt{x}'}(\C{F}(\phi'))=O^{\ov{\ka}_[y]}_{x'}(\C{F}(\phi'))$
for each stable conjugate $x'\in\C{G}'(E)$ of $x$ and each 
$\wt{x}'\in\varphi_{x,x'}(\Om_x)\subset\C{G}'^{\sr}_{x'}(E)$.
Denote by $\Om_y\subset \C{H}_y^{\sr}(E)$ the image of $\Om_x$ under the natural isomorphism
$\C{G}_x\isom\C{H}_y$, sending $x$ to $y$, and choose
an open neighborhood $\Om\subset\C{H}^{\sr}(E)$ of $y$ contained in $\Ad H(E)(\Om_y)$. 

By construction, for each  $\wt{y}\in\Om$ there exists $\wt{x}\in\Om_x\subset\C{G}^{\sr}(E)$ compatible 
with $y$. Moreover, since conditions $(i), (ii)$ from \re{1} do not change when 
$y$ and $x$ are replaced by stable conjugates, it will suffice to show  $(i), (ii)$ for 
some pair of compatible elements $\wt{y}\in\Om\subset\C{H}^{\sr}(E)$ and $\wt{x}\in\C{G}^{\sr}(E)$.
\end{Emp}

\begin{Emp} \label{E:2}
{\bf Strategy of the proof.} 

Choose $\un{E}$, $u$, $w$, $\un{G}$, $\un{\varphi}$, $\un{\C{E}}$, $\un{T}$, 
$\un{a}$, $\un{a}'$, $[\un{b}]$ as in \rcl{global}, 
and identify $\un{G}_v$ with $\un{G}'_v$ for each $v\neq w,u$ as in
\rcl{global} (b). Finally, we fix a
non-trivial character $\un{\psi}:\B{A}_F/F\to \B{C}\m$.


Our strategy will be to construct measures 
$\un{\phi}=\otimes_{v\in V}\phi_v\in \C{S}(\C{G}(\B{A}))$ and 
$\un{\phi'}=\otimes_{v\in V} \phi'_v\in\C{S}(\C{G}'(\B{A}) )$ and compatible elements
$\un{y}\in\un{\C{H}}^{\sr}(\un{E})$ and $\un{x}\in\un{\C{G}}^{\sr}(\un{E})$ 
satisfying the following properties:

$(A)$ $\un{y}_w\in\Om$,  $\phi_w=\phi$ and $\phi'_w=\phi'$.

$(B)$ $\un{x}$ has a stable conjugate in $\un{\C{G}'}(\un{E})$  
if $\un{x}_w$ has a stable conjugate in $\C{G}'({E})$.

$(C)$ both $\un{\phi}$ and $\un{\phi'}$ 
satisfy the support assumption of \rp{trfor} (a), and we have $\Theta(\un{\phi})=\Theta(\un{\phi'})$. 

$(D)$ $\C{F}(\un{\phi})$ and  $\C{F}(\un{\phi}')$ satisfy the support
assumption of \rp{trfor} (a), and we have 
$\Theta(\C{F}(\un{\phi}))=O^{[\un{y}]}_{\un{x}}(\C{F}(\un{\phi}))$ and 
$\Theta(\C{F}(\un{\phi}'))=O^{[\un{y}]}_{\un{x}}(\C{F}(\un{\phi}'))$.
(We define $O^{[\un{y}]}_{\un{x}}(\C{F}(\un{\phi}'))$ to be zero unless 
there exists 
a stable conjugate $\un{x}'\in\un{\C{G}}'(\un{E})$ of $\un{x}$, in which case 
we define $O^{[\un{y}]}_{\un{x}}(\C{F}(\un{\phi}'))$ to be 
$O^{[\un{y}]}_{\un{x}'}(\C{F}(\un{\phi}'))$.)

$(E)$ For each $v\neq w$, we have $O^{[\un{y}]_v}_{\un{x}_v}(\C{F}(\phi_v))\neq 0$.

$(F)$ For each $v\in V_f\sm w$, there exists a stable conjugate $x'_v\in\un{\C{G}}'(\un{E}_v)$ of $\un{x}_v$, and we have 
\[
e'(\un{G}_v)O^{[\un{y}]_v}_{\un{x}_v}(\C{F}(\phi_v))=e'(\un{G}'_v)
\lan\frac{\un{x}_v,x'_v;[\un{y}]_v}{\un{a}_v,\un{a}'_v;[\un{b}]_v}\ran_{\un{\C{E}}_v}
O^{[\un{y}]_v}_{x'_v}(\C{F}(\phi'_v)).
\]

Once these data are constructed, the result follows. 
Indeed, by $(A)$ and the observation at the end of \re{Omega}, it will suffice 
to check that $\un{y}_w$ and $\un{x}_w$ satisfy conditions $(i), (ii)$ of \re{1}.
Next $(C)$, $(D)$ and \rp{trfor} (b) imply that 
$O^{[\un{y}]}_{\un{x}}(\C{F}(\un{\phi}))=O^{[\un{y}]}_{\un{x}}(\C{F}(\un{\phi}'))$.
Assume first that $\un{x}_w$ does not have a stable conjugate in $\C{G}'(E)$. Then 
$\un{x}$ does not have a stable conjugate in $\un{\C{G}}'(\un{E})$, thus 
$O^{[\un{y}]}_{\un{x}}(\C{F}(\un{\phi}))=\prod_v O^{[\un{y}]_v}_{\un{x}_v}(\C{F}(\phi_v))=0$.
Hence the vanishing of $O^{[\un{y}]_w}_{\un{x}_w}(\C{F}(\phi_w))$ follows from 
$(E)$. 

Assume now that that $\un{x}_w$ has a stable conjugate in $\C{G}'(E)$,
then by $(B)$, there exists a stable conjugate $\un{x}'\in\un{\C{G}}'(\un{E})$ of $\un{x}$.
Then $O^{[\un{y}]}_{\un{x}}(\C{F}(\un{\phi}))=O^{[\un{y}]}_{\un{x}'}(\C{F}(\un{\phi}'))$ and
each $\un{x}'_v\in\un{\C{G}}(\un{E})$ is a stable conjugate of $\un{x}_v$.
Using product formulas $\prod_{v\in V} e'(\un{G}_v)=\prod_{v\in V} e'(\un{G}'_v)=1$
(see \cite{Ko5}) and $\prod_{v\in V_f}\lan\frac{\un{x}_v,\un{x}'_v;[\un{y}]_v}
{\un{a}_v,\un{a}'_v;[\un{b}]_v}\ran_{\un{\C{E}}_v}=1$
(\re{globend} (c)), the required equality 
$e'(G)O^{[\un{y}]_w}_{\un{x}_w}(\C{F}(\phi))=e'(G')\lan\frac{\un{x}_w,\un{x}'_w;[\un{y}]_w}
{a,a';[b]}\ran O^{[\un{y}]_w}_{\un{x}'_w}(\C{F}(\phi'))$ follows  from $(E)$ and $(F)$.
\end{Emp}

\begin{Emp}  \label{E:constr} 
{\bf Construction of $\un{\phi},\un{\phi}',\un{y}$ and $\un{x}$.}

(a) Choose an $\un{\C{O}}$-subalgebra $\C{K}\subset\un{\C{G}}(\un{E})$, and  
let $S_1\subset V$ be a finite subset containing $V_{\be}\cup \{w,u\}$ such 
that for each $v\notin S_1$ we have

-  $\un{H}_v$ and $\un{\C{G}}_v$ are unramified;

-  the $\un{\C{O}}_v$-subalgebra $\C{K}_v\subset\un{\C{G}}(\un{E}_v)$, spanned by $\C{K}$, 
is hyperspecial and satisfies $\C{F}(1_{\C{K}_v})=1_{\C{K}_v}$.

(b) Let $\C{A}$ be the set of isomorphisms classes of those endoscopic triples for $\un{G}$, 
which are unramified outside of $S_1$. Then $\C{A}$ is finite (see \cite[Lem. 8.12]{La}), 
and $\un{\C{E}}\in\C{A}$. Let $\C{A}'$ be the subset of $\C{A}$ consisting of triples 
$(H_{a},[\eta_a],\ov{s}_a)$ such that $\un{E}[H_a]$ is not contained in $\un{E}[\un{H}]$. 
For each $a\in\C{A}'$, we fix a prime $v_a\in V_f\sm S_1$ which splits in 
$\un{E}[\un{H}]$ but does not
split in $\un{E}[H_a]$. Put $S_2:=\{v_a\, |\,a\in\C{A}'\}$, and set $S:=S_1\cup S_2$.

(c) Choose $\un{y}\in\un{\C{H}}^{\sr}(\un{E})$ such that $\un{y}_w\in\Om$, 
$\un{y}_u\in\C{H}^{\sr}(\un{E}_u)$ is elliptic 
(that is, $\un{H}_{\un{y}_{u}}\subset \un{H}_{u}$ is elliptic) and 
$\un{y}_v\in\C{H}^{\sr}(\un{E}_v)$ is split  (that is, 
$\un{H}_{\un{y}_{v}}\subset \un{H}_v$ is split) for each $v\in S_2$. 

Choose an element $\un{x}^*\in\un{\C{G}}^*(\un{E})$ compatible with $\un{y}$
(exists by \rl{steinberg} (b)).
Then $(\un{G}^*_{\un{x}^*})_u\subset \un{G}^*_u$ is an elliptic torus.  
Since $\un{y}_w\in\Om$, element $\un{x}^*_w$ has a stable conjugate in 
$\Om_x\subset\C{G}^{\sr}_x(E)$. Since $\un{\varphi}^*_v$ is trivial for each $v\neq w,u$,
it follows from \rl{globalcoh} (c) (ii) (as in the proof \rcl{global} (III) 
that there exists a stably conjugate 
$\un{x}\in \un{\C{G}}(\un{E})$ of $\un{x}^*$.

(d) Choose a stably conjugate $x'_u\in\un{\C{G}}'(\un{E}_u)$ of $\un{x}_u$ 
(which exists by \rl{elltor}), and set 
$x'_v:=\un{x}_v$ for each $v\neq w,u$.  
For each $v\in S\sm (V_{\be}\cup w)$, choose measures  
$\phi_v\in\C{S}(\un{G}(\un{E}_v))$ and $\phi'_v\in\C{S}(\un{G}'(\un{E}_v))$ 
constructed in \rco{constr2} for the triple $(\un{y}_v,\un{x}_v,{x}'_v)$. 
In particular, $\phi_v=\phi'_v$ for each $v\in S\sm (V_{\be}\cup\{u,w\})$.

For each $v\in S\sm V_{\be}$, let $\om_v\subset\un{E}_v\m$ be an open 
neighborhood of the identity such that $\Om$ is  invariant under the multiplication by $\om_w$, 
while $\C{F}(\phi_v)$ and $\C{F}(\phi'_v)$ are invariant under the multiplication by $\om_v$ 
if $v\neq w$.

For each $v\in V\sm S$, put $\phi_v=\phi'_v=1_{\C{K}_v}dx_v$. Finally, put $\phi_w:=\phi$ and 
$\phi'_w:=\phi'$.

(e) Choose a finite set $S_3\supset S$ such that for each  
$v\notin S_3$, we have $\un{x}_v\in\C{K}_v$, and the reduction  of $\un{x}_v$ 
modulo $v$ is regular. Choose $\la\in \un{E}\m$ such that  

(i)  $\val_v(\la)=0$ for each $v\notin S_3$;  

(ii) $\la\in\om_v$ for each $v\in S\sm V_{\be}$;  

(iii) For each $v\in S_3\sm S$, $\val_v(\la)$ is so large that $(\la\ov{y})_v$ 
belongs to the open neighborhood of zero prescribed in \rl{unit} (b).

Finally, we replace $\un{y}$ and $\un{x}$ constructed in (c) by 
$\la\un{y}$ and $\la\un{x}$, respectively.

(f) Recall that $\C{F}(\phi_v)$ 
and $\C{F}(\phi'_v)$ are compactly supported for each $v\in V_f$ and
$\C{F}(\phi_v)=\C{F}(\phi'_v)=1_{\C{K}_v}dx_v$ for each $v\notin S$.
Since $\un{E}\subset \B{A}$ is discrete, one can choose a
compact neighborhood $C_v\subset \un{\C{G}}(\un{E}_v)=\un{\C{G}}'(\un{E}_v)$ of 
$\un{x}_v$ for each $v\in V_{\be}$ such that all elements of 
$\un{\C{G}}(\un{E})\cap
(\prod_{v\in V_{\be}}C_v\times\prod_{v\in V_f}Supp(\C{F}(\phi_v)))$ 
and $\un{\C{G}}'(\un{E})\cap
(\prod_{v\in V_{\be}}C_v\times\prod_{v\in V_f}Supp(\C{F}(\phi'_v)))$ 
are stable conjugate to $\un{x}$.

For each $v\in V_{\be}$, choose a measure $\phi_v=\phi'_v$ of 
$\C{S}(\un{\C{G}}(\un{E}_v))$ such that 
$\C{F}(\phi_v)=f_v d{x}_v$ for a smooth non-negative function $f_v$ on 
$\un{\C{G}}(\un{E}_v)$ supported on $C_v$ such that $f_v(\un{x}_v)\neq 0$.
Put $\un{\phi}:=\otimes_{v\in V}\phi_v$ and $\un{\phi}':=\otimes_{v\in V} \phi'_v$.
\end{Emp}

To complete the proof of \rt{Wa1}, it remains to show that the constructed above 
$\un{\phi},\un{\phi}',\un{y}$ and $\un{x}$ satisfy conditions  $(A)-(F)$ of \re{2}.

\begin{Emp} 
{\bf Proof of conditions $(A)-(F)$.}

$(A)$ is clear (see \re{constr} (c) and (e)). 

$(B)$ Since $\un{x}_u$ is elliptic, the assertion follows from \rl{globalcoh} (c) (ii)
(as in the proof of \rcl{global} (III)).

$(C)$ Since the support of $\phi_u$ and $\phi'_u$ is regular elliptic, 
both $\un{\phi}$ and $\un{\phi}'$ satisfy the support assumption of \rp{trfor} (a).
Because of symmetry between $\un{\phi}$ and $\un{\phi}'$, it will therefore
suffice to check that for every $\un{z}\in\un{\C{G}}(\un{E})$ and  
$\un{\ka} \in \wh{\un{G}_{\un{z}}}^{\un{\Gm}}$ such that 
$O^{\un{\kappa}}_{\un{z}}(\un{\phi})\neq 0$, there exists a stable conjugate 
$\un{z}'\in\un{\C{G}'}(\un{E})$ of $\un{z}$, and we have 
$O^{\un{\kappa}}_{\un{z}}(\un{\phi})=O^{\un{\kappa}}_{\un{z}'}(\un{\phi}')$. 

Fix $\un{z}\in\un{\C{G}}(\un{E})$ and  $\un{\ka} \in \wh{\un{G}_{\un{z}}}^{\un{\Gm}}$ 
such that $O^{\un{\kappa}}_{\un{z}}(\un{\phi})\neq 0$. 
Consider the endoscopic triple $\C{E}':=\C{E}_{([a_{\un{z}}],\un{\ka})}=(H',[\eta'],\ov{s}')$ 
for $\un{G}$. Following \cite[10.9]{Wa1}, we will show that $\C{E}'\cong\un{\C{E}}$. 

By the assumption, $O^{\un{\kappa}_v}_{\un{z}_v}(\phi_v)\neq 0$ for each $v\in V$. 
Since $\C{K}_v\subset \un{\C{G}}(\un{E}_v)$ is a hyperspecial subalgebra and 
$\phi_v=1_{\C{K}_v}dx_v$ for each $v\in V\sm S$, we conclude from 
\rl{unit} (a) that $\C{E}'$ is unramified outside of $S$. 

For each $v\in S_2$, measure $\phi_v$ is supported on split elements (by \rco{constr2} (i)) and 
satisfies $O^{\un{\ka}_v}_{\un{z}_v}(\phi_v)\neq 0$. Therefore $\un{G}_{\un{z}_v}$ is split, 
hence $\C{E}'_v\cong \C{E}_{([a_{\un{z}_v}],\un{\ka}_v)}$ is split. In particular,  
 $\C{E}'$ is unramified outside of $S_1$. Moreover, $H'$ splits at each $v\in S_2$, hence
$\un{E}[H']\subset \un{E}[\un{H}]$. 

Since $O^{\un{\ka}_u}_{\un{z}_u}(\phi_u)\neq 0$, 
we get by \rco{constr2} (i), (iii) that $\un{z}_u$ is stable conjugate to an element of
$\un{G}^{\sr}_{\un{x}_u}(\un{E}_u)$, and the class of 
$\un{\ka}_u\in\wh{\un{G}_{\un{z}}}^{\un{\Gm}_u}\cong 
\wh{\un{G}_{\un{x}}}^{\un{\Gm}_u}$ equals $\un{\ka}_{[\un{y}]_u}$.
Since $\un{x}_u$ is elliptic, we get that $\C{E}'_u\cong
\C{E}_{([a_{\un{x}_u}],\ov{\ka}_{[\un{y}]_u})}$ is isomorphic to
$\un{\C{E}}_u$ (use \rl{endosc} (a)).
  As $u$ is inert in $\un{E}[\un{H}]\supset \un{E}[H']$, we conclude from this that  
$\C{E}'\cong\un{\C{E}}$ (use \re{globend} (b)), as claimed.

By the proven above, there exists $\un{y}'\in\un{\C{H}}^{\sr}(\un{E})$ compatible with $\un{z}$ 
such that $\un{\kappa}=\ov{\ka}_{[\un{y}']}$, thus 
$O^{\un{\kappa}}_{\un{z}}(\un{\phi})=O^{[\un{y}']}_{\un{z}}(\un{\phi})\neq 0$. 
In particular, $O^{[\un{y}']_w}_{\un{z}_w}(\phi_w)\neq 0$. Since $\phi_w$ and $\phi'_w$ are 
$(a,a';[b])$-indistinguishable, it follows from \rl{indist} that there exists
a stable conjugate $z'_w\in\C{G}'(E)$ of $\un{z}_w$. 
Since $\un{G}_{\un{z}_u}\subset \un{G}_u$ is elliptic, there exists 
$\un{z}'\in \un{G}'(\un{E})$ stably conjugate to $\un{z}$ such that
$\un{z}'_v$ is conjugate to $\un{z}_v$ for each $v\neq u,w$ 
(use \rl{globalcoh} (c) and \re{inv}). 

It now remains to show that 
$O^{[\un{y}']}_{\un{z}}(\un{\phi})=O^{[\un{y}']}_{\un{z}'}(\un{\phi}')$.
By the product formula (\re{globend} (c)), it will suffice to check that for each $v\in V_f$, 
we have
\[
O^{[\un{y}']_v}_{\un{z}_v}(\phi_v)=
\lan\frac{\un{z}_v,\un{z}'_v;{y}'_v}{\un{a}_v,\un{a}'_v;[\un{b}]_v}\ran_{\un{\C{E}}_v} 
O^{[\un{y}']_v}_{\un{z}'_v}(\phi'_v).
\]
The assertion for $v=w$ follows from \rl{indist}, while the assertion
for $v=u$ follows from \rco{constr2} (ii).
Finally, the assertion for $v\neq u,w$ follows from the fact that under the identification
$\un{G}'_v=\un{G}_v$, we have $\phi'_v=\phi_v$, $\un{z}'_v$ is conjugate to $\un{z}_v$, and
$\un{a}'_v$ is conjugate to $\un{a}_v$. 

$(D)$ By \re{constr} (f),  
$Supp(\C{F}(\un{\phi}))\cap \Ad\un{G}(\B{A})(\C{G}(\un{E}))$ 
consists of elements stably conjugate to $\un{x}$. 
Since $\un{x}_u$ is elliptic, $\un{x}$ is elliptic, thus $\C{F}(\un{\phi})$ satisfies the 
support assumption of \rp{trfor} (a). Therefore 
$\Theta(\C{F}(\un{\phi}))=\sum_{\un{\ka}\in \wh{\un{G}_{\un{x}}}^{\un{\Gm}}} 
O^{\un{\ka}}_{\un{x}} (\C{F}(\un{\phi}))$. Let $\un{\ka}\in \wh{\un{G}_{\un{x}}}^{\un{\Gm}}$ be such 
that $O^{\un{\ka}}_{\un{x}} (\C{F}(\un{\phi}))\neq 0$.
Then $O^{\un{\kappa}_u}_{\un{x}_u}(\C{F}({\phi}_u))\neq 0$. 
Since $\un{x}_u$ is elliptic, it follows from \rco{constr2} (iv) that   
$\un{\ka}_u=\ov{\ka}_{[\un{y}]_u}$. But the map
$\wh{\un{G}_{\un{x}}}^{\un{\Gm}}\hra\wh{\un{G}_{\un{x}}}^{\un{\Gm}_u}=
\pi_0(\wh{\un{G}_{\un{x}}}^{\un{\Gm}_u})$ is injective, therefore 
$\un{\ka}=\ov{\ka}_{[\un{y}]}$. This shows that 
$\Theta(\C{F}(\un{\phi}))=O^{[\un{y}]}_{\un{x}}(\C{F}(\un{\phi}))$.
The proof for $\un{\phi}'$ is similar. 

$(E),(F)$ For $v\in S\sm (V_{\be}\cup w)$, the assertions follows from \rco{constr2} (v).
For $v\neq u,w$, $(F)$ follows from the fact that under the identification 
$\un{G}_v=\un{G}'_v$, we have
$\phi'_v=\phi_v$, and $\un{a}'_v$ is conjugate to $\un{a}_v$. 
It remains to show $(E)$ for $v\in V_{\be}\cup(V\sm S)$. 
If $v\in V_{\be}$, the assertion follows from the fact that $\C{F}(\phi_v)=f_v d{x}_v$,
while $f_v$ is non-negative and satisfies $f_v(\un{x}_v)\neq 0$.
Assume now that $v\notin S$. Then $\C{F}(\phi_v)=1_{\C{K}_v}dx_v$.
The assertion now follows from \rl{unit} (c) if $v\notin S_3$ and from the choice
of $\la$ in \re{constr} (e) (and  \rl{unit} (b)) if $v\in S_3$.
\end{Emp}


























\section*{List of main terms and symbols}


\begin{multicols}{3}

\e $\C{E}$-admissible, \pg{eadm}

\e $(\C{E},[a],\ov{\ka})$-admissible, \pg{eakaadm}

\e compact element, \pg{compact} 

\e compatible, \pg{comemb}, \pg{comp2}

\e $\C{E}$-compatible, \pg{ecomp}
 
\e conjugate, \pg{conj}

\e $E^{\sep}$-conjugate, \pg{esepconj}

\e defined at $x$, \pg{defatx}

\e defined over $\C{O}$, \pg{defovero}

\e  Deligne--Lusztig representation, \pg{DLr}

\e endoscopic triple, \pg{endtr}

\e $(a,a';[b])$-equivalent, \pg{aa'beq}

\e $(a,a';\ka)$-equivalent, \pg{aa'kaeq}

\e $(\gm_0,\ov{\xi})$-equivalent, \pg{gmoxieq}, \pg{gmoxieq2}

\e Fourier transform, \pg{ft} 

\e $(a,a';[b])$-indistinguishable, \pg{aa'bin}

\e inner twisting, \pg{it} 

\e non-degenerate at $x$, \pg{nondegatx}

\e non-degenerate over $\C{O}$, \pg{nondego}

\e property $(vg)$, \pg{vg}

\e property $(vg)_a$, \pg{vga}

\e $(G,a_0,\gm_0)$-relevant, \pg{ga0gm0}

\e $\C{E}$-stable, \pg{estable}

\e stably conjugate, \pg{stcon}

\e strongly regular,  \pg{str}

\e $a$-strongly regular, \pg{astrreg} 

\e topological Jordan decomposition, \pg{tjd}

\e topologically unipotent, \pg{topun}

\e quasi-isogeny, \pg{qi}

\e quasi-logarithm, \pg{qlog}

\e $\C{E}$-unstable, \pg{eunstable}

\e $a:T\hra G$, \pg{atg} 

\e $\ov{a}:\ov{T}\hra L_a$, \pg{ovatg} 

\e $\Ad G$, $\Ad g$, \pg{adg} 

\e $\wh{[a]}$, \pg{wha}

\e $a_x$, \pg{ax}

\e $\ov{a}$, \pg{ovg} 

\e  $\C{B}(G)$, \pg{Bg}

\e $[b]_G$, \pg{bg}

\e $b_{t,\dt}$, \pg{btdt}

\e  $C_c^{\infty}(X(E))$ \pg{ccinfty}

\e $\C{D}(X(E))$, \pg{dxe}

\e $\C{D}_G$, \pg{dg}

\e $\C{D}^0(X(E))$, \pg{d0xe}

\e $dg=|\om_G|$, \pg{omg} 

\e $dx=|\om_X|$, \pg{omx} 

\e $\ov{E}$, \pg{ove} 

\e $E^{\sep}$, \pg{esep} 

\e $E^{\nr}$, \pg{enr} 

\e $e(G)$,  $e'(G)$, \pg{eg}

\e $\C{E}=(H,[\eta],\ov{s})$, \pg{ehetas}

\e $\C{E}_{([a],\ka)}$, \pg{eaka}

\e $\C{E}_{t}=(H_t,[\eta_t],\ov{s}_t)$, \pg{et}

\e $F|_U$, \pg{fu}

\e $F(x,\ov{\xi})$, \pg{fxxi}

\e $F_{a,\theta}$, \pg{fatheta}

\e $F_{a_0,\ka,\theta}$, \pg{fa0katheta}

\e $F_{\tu}$, \pg{ftu}

\e $F_x$, $F_{x^*}$, \pg{fx}

\e $\C{F}$, \pg{f}

\e $G^0$, \pg{g0}

\e $G^{\ad}$, $G^{\der}$, \pg{gad} 

\e $G^{\ssc}$, \pg{gss}

\e $\C{G}$, \pg{lieg}
  
\e $G(E)_{\tu}$, $\C{G}(E)_{\tn}$, \pg{getu}

\e $\fq$, \pg{fq} 

\e  ${}^0G$, \pg{0g}

\e $G^*$, \pg{g*}

\e $\wh{G}$, \pg{dualg}  

\e $G_{a}$, $\C{G}_a$, \pg{ga}

\e  $G_{a^+}$, $\C{G}_{a^+}$, \pg{ga} 

\e $\wt{G_a}$, \pg{wtga}

\e $G_{x}$, $\C{G}_x$, \pg{stab}, \pg{gx} 

\e $G_{x^+}$, $\C{G}_{x^+}$, \pg{gx+}

\e $G_{x,\tu}$, $\C{G}_{x,\tn}$, \pg{gxtu}

\e $\un{G}_x$, $\ov{G}_x$, \pg{ovgx}

\e $\ov{g}$, \pg{ovg} 

\e $\C{H}$, \pg{lieh} 

\e $\C{I}^+$, \pg{I+}

\e $\Inn G$, $\Inn g$, \pg{intg} 

\e $\inv(G,G')$, \pg{invg'g}

\e $\inv(x,x')$, \pg{invxx'}

\e $\ov{\inv}(x,x')$, \pg{ovinvxx'}

\e $\ov{\inv}((x_1,x'_1);...;(x_k,x'_k))$, \pg{ovinvkxx'}

\e $\C{L}$, \pg{liel} 

\e $L^1_{loc}(X(E))$, \pg{l1loc}

\e  $L_{a}$, $\C{L}_a$, \pg{ga}

\e $L_x$, $\C{L}_x$, \pg{lx} 

\e $\frak{m}$, \pg{m} 

\e $\C{N}(\C{H})$, \pg{nh}

\e $\C{O}$, \pg{O} 

\e $O_{x}$, \pg{ox}

\e $O^{\ov{\xi}}_{x}$, \pg{oxix}

\e $O^{[y]}_x$, \pg{oyx} 

\e $p$, \pg{p} 

\e  $SO_{x}$,  $SO_{x^*}$, \pg{sox} 

\e $S_{([a],\ov{\ka})}$, \pg{saka} 

\e $S_{[b]}$, \pg{sb}

\e $\C{S}(X(E))$, \pg{sxe}

\e $\C{T}$, \pg{liet}  

\e  $\un{Tor}$, \pg{tor}

\e $T_{\C{O}}$, $\ov{T}$, \pg{to}

\e $t_{a,\theta}$, \pg{tath}

\e $\C{U}(H)$, \pg{uh} 

\e $W(G)$, \pg{wg}

\e $X^{\sr}$, \pg{sr}

\e $X^*$, \pg{g*} 

\e $[x]$, \pg{[x]} 

\e $[y]_{G}$,  \pg{yg} 

\e $Z(G)$, \pg{zg}
 
\e $Z(\C{E})$, \pg{ze} 

\e $Z(\C{E},[a],\ov{\ka})$, \pg{zeaka}

\e $Z_{[\wh{a}]}$, \pg{zwha}

\e $Z_{[\ov{\eta}]}$, \pg{zbareta}

\e $\Gm$, \pg{Gm} 

\e $\Dt_{a,t}$, \pg{Dtat}

\e ${\Dt}_{a_0,\ka,t}$, \pg{Dta0kat}

\e $\dt_1\sim_{(\gm_0,\ov{\xi})}\dt_2$, \pg{dt1simdt2}, \pg{dt'simdt}

\e $\dt_{a,t}$, $\ov{\dt}_{a,t}$, \pg{dtat}

\e $\theta:T(E)\to\B{C}\m$, \pg{th} 

\e $\ov{\theta}:\ov{T}(\fq)\to\B{C}\m$, \pg{ovth}

\e $\iota([b_1],[b_2])$, \pg{ib1b2}

\e $\ka([b_1],[b_2])$, \pg{kab1b2}

\e $\ka\left(\frac{[b_1]}{[b_2]}\right)$, \pg{ka'b1b2}

\e $\ov{\ka}_{[b]}$, \pg{ovkab}

\e $\ov{\ka}_{[y]}$, \pg{kay} 

\e $\La(\C{E})$, \pg{lae}

\e $\mu_H, \nu_H$, \pg{muh}

\e $\Pi_{\C{E}}$, \pg{Pie}

\e $\pi_{\C{E}}$, \pg{pie}

\e $\pi_{a,\theta}$, \pg{piath}

\e $[\wh{\pi}]$, \pg{whpi}

\e $\rho_{a,\theta}$, \pg{rhoath}

\e $\rho_{\ov{a},\ov{\theta}}$, \pg{rhoatheta}

\e $\rho_G:\Gm\to\Out(\wh{G})$, \pg{rhog}

\e $\Phi:G\to \C{G}$, \pg{phi}

\e $\Phi_x:\un{G}_x\to\C{G}_x$, \pg{phix}

\e $\ov{\Phi}_x:L_x\to\C{L}_x$, \pg{ovphix}

\e ${}^0\Phi$, \pg{0Phi}

\e  $\Phi_{\rho}$, \pg{phirho}

\e $\phi|_U$, \pg{phiu}

\e $\varphi:G\to G'$, \pg{it}

\e $\varphi_X:X\to X'$, \pg{varphiX}

\e $\varphi_{x,x'}$, \pg{varphixx'}

\e $\chi(\pi_{a,\theta})$, \pg{chipiath}

\e $\chi_{a_0,\ka,\theta}$, \pg{chia0kath}

\e $\Om_{a,t}$, $\ov{\Om}_{a,t}$, \pg{omat}

\e $\om_G$, $\om_{\C{G}}$, \pg{omg} 

\e $\om_X$,  \pg{omx} 

\e $\om_{X'}$, \pg{omx'}

\e $\lan\frac{a_1,a'_1;[b_1]}{a_2,a'_2;[b_2]}\ran$, \pg{brac}

\e $\lan\frac{a_1,a'_1;[b_1]}{a_2,a'_2;[b_2]}\ran_{\C{E}}$, \pg{brac}

\e $\lan\frac{a_1,a'_1;\ov{\ka}}{a_2,a'_2;[b_2]}\ran$, \pg{brac2}

\e $\lan\frac{a_1,a'_1;[b_1]}{a_2,a'_2;\ov{\ka}}\ran$, \pg{brac2}

\e $\lan\frac{x,x';\ov{\ka}}{a,a';[b]}\ran$, \pg{brac3}

\e $\lan\frac{x,x';[y]}{a,a';[b]}\ran$, \pg{brac4}

\e $\lan\cdot,\cdot\ran_x$, \pg{lanranx}

\e $\lan\cdot,\cdot\ran_{\rho}$, \pg{lanranrho}

\e

\e

\e

\e

\e

\e

\e

\e

\e

\e

\e

\end{multicols}

\end{document}